\documentclass[twoside,openright]{report}
\usepackage[a4paper,hmargin=2.2in,vmargin=2.5in]{geometry} %1.2, 1.1
\usepackage{definitions}
\usepackage{./subfiles/subfiles}
\usepackage{xmpincl}

\usepackage{hyperref}
\usepackage{cleveref}

\usepackage{thmtools}
\usepackage{thm-restate}

\usepackage{etoolbox}
\usepackage{cite}
\usepackage{mathrsfs}
\usepackage{epigraph}
\usepackage[small,nohug,heads=vee]{diagrams}
\usepackage{afterpage}

\makeatletter
\newcommand{\labitem}[2]{%
\def\itemlabel{\textbf{#1}}
\item
\def\currentlabel{#1}\label{#2}}
\makeatother

\newcommand{\xx}{\overline{x}}
\newcommand{\oa}{\overline{a}}

\newcommand{\QQ}{$Q$} %text command
\newcommand{\qg}{\QQ-group }   %text command
\newcommand{\qsub}[0]{\QQ-subgroup }

\newcommand{\pr}{$\pi$}  %text command
\newcommand{\gad}{geometric abelian decomposition }   %text command
\newcommand{\gat}{geometric abelian tree }  %text command
\newcommand{\gats}{geometric abelian trees }  %text command

\newcommand{\gr}{graded }   %text command
\newcommand{\gd}{graded }   %text command
\newcommand{\rs}[0]{restricted }   %text command
\newcommand{\res}{resolution } %text command
\newcommand{\ress}[0]{resolutions } %text command

\newcommand{\lvs}[1]{lv(#1)}   %% leaves of a frt
\newcommand{\trm}[0]{concretion }        	%text command
\newcommand{\cautious}[0]{cautious } %text command
\newcommand{\dt}[1]{\mathcal{D}(#1)}
\newcommand{\cfam}[0]{essential system of curves } %text command
\newcommand{\cfams}[0]{essential systems of curves } %text command
\newcommand{\inn}[1]{inn_{#1}}

\newcommand{\ws}[0]{well-separated } %text command
\newcommand{\fcr}[0]{factorization} %text command
\newcommand{\twist}[0]{twist } %text command

\newcommand{\fcrs}[0]{factorizations } %text command
\newcommand{\dsd}[0]{decomposed } %text command
\newcommand{\fcrm}[3]{(#2_{\lambda},#3_{\mu})_{\lambda\in#1,\mu\in\hat{#1}}}
\newcommand{\lft}[1]{\tilde{#1}}  	 % lift from the quotient graph to some transversal
\newcommand{\afld}[1]{#1^{f}} 		%% the operation folding edges adyacent to abelian subgroups of a group
\newcommand{\BS}[3]{(#1,#2,(#3_{e})_{e\in E#2})}
\newcommand{\BSP}[0]{Bass-Serre presentation } %tc
\newcommand{\vga}[4]{\bar{#1}_{(#2_{e})_{\substack{\alpha(e)=#3,\\e\in #4}}}}
\newcommand{\mvga}[3]{\bar{#1}_{(#2_{e})_{\alpha(e)=#3}}}
\newcommand{\svga}[2]{\bar{#1}_{(#2_{e})}}

\newcommand{\Act}[1]{\mathcal{T}(#1)}

\newcommand{\PAct}[1]{\mathcal{PA}(#1)}

\newcommand{\edin}[1]{\overline{#1}}

\newcommand{\LQ}[0]{\mathcal{L}^{Q}} % the language $L_q$

\newcommand{\lgp}[0]{level } %tc
\newcommand{\Ess}[1]{Ess_{Q}(#1)}   % the family of all systems of essential s.c.c. on a surface

\newcommand{\SEC}[0]{system of essential simple closed curves } %text command
\newcommand{\lvlgp}[0]{node group }  %text command
\newcommand{\TS}[2]{#1_{#2}} % the taut structure at the vertex group associated to vertex $v$
\newcommand{\treeact}[3]{(#1,#2_{(-)},p_{(-)})}   % notation for the tree of actions given in Rips theorem; arguments in this order; first the simplicial tree, then the group, then the

\newcommand{\maxim}[1]{\dsp{\max_{#1}} \,}
\newcommand{\pl}[3]{l_{#1}^{#3}(#2)}
\newcommand{\tl}[2]{tl^{#2}(#1)}
\newcommand{\ml}[3]{ml^{#3}_{#1}(#2)}  %maximum displacement of $x$ of a tuple of generators
\newcommand{\md}[2]{md(#1)^{#2}}

\newcommand{\egt}[0]{equivariant Gromov topology } %tc
\newcommand{\pegt}[0]{pointed equivariant Gromov topology } %tc
\newcommand{\looprep}[3]{#1_{0}t_{#2_{0}}\cdots t_{#2_{#3-1}}#1_{#3}} %% loop representation of an element
\newcommand{\vma}[1]{Aut_{\partial}(#1)}  			%vertex modular automorphisms

\newcommand{\atjsj}[0]{abelian* JSJ decomposition }  %tc
\newcommand{\adjsj}[0]{abelian* JSJ tree }  %tc

\newcommand{\ceq}[3]{#1^{#3}(#2)=1}     %constrained system of equations
\newcommand{\cneq}[3]{#1^{#3}(#2)\neq 1}
\newcommand{\syeq}[2]{#1(#2)=1}       %unconstrained system of equations
\newcommand{\syneq}[2]{#1(#2)\neq 1}

\newcommand{\NEQs}[0]{systems of inequations } %tc
\newcommand{\EQ}[0]{system of equations } %tc
\newcommand{\CEQ}[0]{\pr-system of equations } %tc
\newcommand{\CEQS}[0]{\pr-systems of equations } %tc
\newcommand{\CNEQS}[0]{\pr-systems of inequations } %tc
\renewcommand{\ml}[3]{ml_{#1}^{#3}(#2)}
\renewcommand{\md}[2]{md^{#2}(#1)}

\renewcommand{\sl}[3]{sl_{#1}^{#3}(#2)}
\renewcommand{\pl}[3]{l_{#1}^{#3}(#2)}
\renewcommand{\tl}[2]{tl^{#2}(#1)}

\newcommand{\fundamental}[0]{fundamental }   %tc
\newcommand{\disceq}[5]{\bvee{j=1}{(y\in #5 (j)\wedge\syeq{#2}{#1}\wedge \syneq{#3}{#1})}{#4}}    % disjunction

\newcommand{\pVAFor}[6]{\forall x\,((x\in q^{#4}\wedge\syeq{\Delta}{x,a})\sra \exists y\in #5
\wedge(\disceq{x,y,a}{#2}{#3}{#6}{q^{y}}))}
\newcommand{\pg}[1]{pgq(#1)}

\newcommand{\lo}[0]{\vartriangleleft}
\newcommand{\LL}[0]{\mathcal{L}}

\newcommand{\BTQ}[0]{T^{\pi}_{gp}}

\newcommand{\qG}[1]{\mathbf{\mathsf{#1}}}

\newcommand{\rqG}[2]{\qG{#1}_{#2}}
\newcommand{\gqG}[4]{\qG{#1}_{#2,(#3_{i})_{i\in #4}}}

\newcommand{\hDD}[1]{\hat{\Delta}_{\hat{#1}}}
\newcommand{\hDDa}[1]{\hat{\Delta}_{[\hat{#1}]}}
\newcommand{\ddn}[1]{\Delta_{#1}}

\newcommand{\DD}[1]{\hat{Delta}_{#1}}

\newcommand{\DDa}[1]{\hat{\Delta}_{[\hat{#1}]}}  % Those $2$ are the same.
\newcommand{\V}[0]{V}
\newcommand{\Va}[0]{V_{a}}

\newcommand{\hV}[0]{\hat{V}}

\newcommand{\E}[0]{E}

\newcommand{\abv}[1]{\hat{[#1]}}
\newcommand{\MT}[0]{Z}

\newcommand{\ctm}[2]{\iota_{#1}^{#2}}
\newcommand{\Sb}[2]{Stab_{#2}(#1)}

\newcommand{\RRA}[0]{\mathcal{R}_{A}}
\newcommand{\RSA}[0]{\mathcal{S}_{A}}
\newcommand{\RTA}[0]{\mathcal{T}_{A}}

\newcommand{\QR}[1]{\mathcal{#1}}
\newcommand{\RAl}[0]{\rqG{A}{R^{\lambda}}}
\newcommand{\TAl}[0]{\rqG{A}{S^{\lambda}}}

\newcommand{\TAm}[0]{\rqG{A}{S^{\mu}}}

\newcommand{\GRla}[0]{(\Delta\mathcal{R})^{\lambda}}

\newcommand{\GSla}[0]{(\Delta\mathcal{S})^{\lambda}}

\newcommand{\GTla}[0]{(\Delta\mathcal{T})^{\lambda}}

\newcommand{\Rl}[0]{\qG{R}^{\lambda}}
\newcommand{\Tl}[0]{\qG{T}^{\lambda}}

\newcommand{\Rm}[0]{\qG{R}^{\mu}}
\newcommand{\Tm}[0]{\qG{T}^{\mu}}
\newcommand{\Sm}[0]{\qG{S}^{\mu}}

\newcommand{\RLR}[3]{\qG{#1}_{#2}^{#3}}

\newcommand{\LRL}[3]{\QR{#1}_{#2}^{#3}}   %resolution
\newcommand{\LRLi}[4]{(\QR{#1}^{#2})_{#3}^{#4}}    %resolution
\newcommand{\LRLii}[3]{\QR{#1}^{#2}_{#3}}

\newcommand{\Mo}[2]{Mod(\qG{#1}^{#2},\Delta\QR{#1}^{#2})}
\newcommand{\RMo}[3]{Mod(\QR{#1}_{#2}^{#3})}
\newcommand{\GQR}[5]{\mathcal{#1}_{#2;(#4_{i})_{i\in #5}}^{#3}}

\newcommand{\DR}[2]{(\Delta\QR{#1})^{#2}}
\newcommand{\Title}[0]{On expansions of non-abelian free groups\\ by cosets of a finite index subgroup}
\newcommand{\Name}[0]{Javier de la Nuez Gonz\'alez}

\begin{document}
%\begin{comment}
\title{\Title}
\author{Javier de la Nuez Gonz\'alez\\Westf\"alische Wilhelms Universit\"at M\"unster}
\date{September, 2016}
%\end{comment}

%\maketitle

\titlepage{
\vspace*{5.5cm}
\centering
%\Large
\Name\\
\vspace*{1cm}
%{\LARGE
%A Duality Transform \\for Realizing Convex Polytopes \\with Small Integer Coordinates\\
%}
%\Large
\Title\\
\vspace*{1cm}
%\Large
2016

%\vspace*{7cm}
%\centering
%Alexander Igamberdiev\\
%A Duality Transform for Realizing Convex Polytopes with Small Integer Coordinates\\
%2015

}

\begin{comment}
The interest in the theory of non-abelian free groups goes back to Tarski, who first asked the question of whether all non-abelian free groups share the same order theory. While one needs to look no further than abelianity to find a first order property separating the free group in one generator from the non-abelian ones the so called "Tarski's problem" (question?), evaded the efforts of mathematicians for many decades, until it was solved affirmatively by Zlil Sela in $2001$ (and later independently by O.Kharlampovich and A. Myasnikkov). Sela's solution spans a series of $6$ papers, where he needs to develop an entirely new host of techniques.
\end{comment}

\pagestyle{empty}
\pagenumbering{gobble}
\afterpage{\null\newpage}

\begin{titlepage}

\vspace{1cm}

\centering

\center{
Mathematik
}

\vspace{4cm}

%\center{Dissertationsthema}

{\huge{
\Title
}}

\vspace{6cm}

\centering{
Inaugural-Dissertation\\
zur Erlangung des Doktorgrades\\
der Naturwissenschaften im Fachbereich\\
Mathematik und Informatik\\
der Mathematisch-Naturwissenschaftlichen Fakult\"at\\
der Westf\"alischen Wilhelms-Universit\"at M\"unster\\
}

\vspace{5cm}

\centering{
vorgelegt von\\
Javier de la Nuez Gonz\'alez\\
aus Madrid\\
2016
}

\end{titlepage}

\begin{titlepage}

\vspace*{18cm}

\begin{tabular}{p{.5\textwidth} p{.5\textwidth}}
Dekan:  &  Prof. Dr. Martin Stein \vspace*{0.2cm}\\
Erster Gutachter: & Prof. Dr. Dr. Katrin Tent \vspace*{0.2cm}\\
Zweiter Gutachter: & Prof. Dr. Zlil Sela \vspace*{0.2cm}\\
Tag der Disputation: & 24/11/2016 \vspace*{0.2cm}\\
Tag der Promotion:  
\end{tabular}
\end{titlepage}

%\setcounter{page}{1}
%\pagenumbering{none}
%\setcounter{page}{1}
%\pagenumbering{Roman}

\begin{abstract}
	Let $\F$ be a finitely generated non-abelian free group and $Q$ a finite quotient. Denote by $\LQ$ the language obtained by adding unary predicates $\{P_{q}\}_{q\in Q}$ to the language of groups.
	
	By generalizing some of the techniques involved in Zlil Sela's solution to Tarski's question on the elementary theory of non-abelian free groups, we provide a few basic results on the validity of first order sentences in the $\LQ$-expansion of $\F$ in which every $P_{q}$ is interpreted as the preimage of $q$ in $\F$.
	
	In particular, we prove an analogous result to Sela's generalization of Merzlyakov's theorem for $\forall\exists$-sentence and show that the positive theory of such a structure depends only on $Q$ and not on the rank of $\F$ nor the particular quotient map.
\end{abstract}

\newpage\null\thispagestyle{empty}\newpage

\chapter*{Acknowledgements}

%\addtocontents{toc}{\protect\thispagestyle{empty}}
First and foremost, I would never have come this far without the unconditional love and support of my family, who were always there in the highs and lows and had to sacrifice more than I did during these years. This is partly theirs as well.

It would also have been impossible without the support of my advisor Prof. Katrin Tent, the many opportunities for exchange with which she provided me.

I am very grateful for Prof. Zlil Sela for inviting me to visit to Jerusalem and his openness for discussion, often around what for him must have been fairly trite and boring technical details of his wonderful work. His patient and unassuming attitude is something to appreciate.

It is impossible for me to forget the truly exceptional hospitality of Montserrat Casals and Ilya Kazachkov, who were always ready to lend me an ear and helped me to gain focus and confidence I need of it.
%It is due to my experience in Bilbao that I believe utopia can exist within the mathematical community.

I must also thank Chloe Perin and Rizos Sklinos, who allowed this nobody to collaborate with them and  taught some invaluable lessons on how to put the book aside and start doing math.

I would also like to mention Prof. Gilbert Levitt, who I had the chance to visit for a week in Caen and with whom I came to experience the beauty of geometric group theory for the first time. Also Martin Bays, who helped me chose this dissertation's title and Franziska Janke and Daniel Palac\'in, who shared with me very valuable insights on the process of writing a thesis.

Finally, life in the math building would not have been the same without the warmth of my fellow graduate students in M\"unster:  Alex (who saved my life several times), Isabel and Zaniar.

\newpage\null\thispagestyle{empty}\newpage
%\addtocontents{toc}{\protect\thispagestyle{empty}}
\tableofcontents
\clearpage

\pagestyle{headings}
\setcounter{page}{1}
\pagenumbering{arabic}

  %  \subsection{Structure of Sela's proof}
    
    % \ltba{A}{group}
  %  \subsection{Content of this work}

\chapter*{Introduction}
\addcontentsline{toc}{chapter}{Introduction}
\markboth{INTRODUCTION}{INTRODUCTION}
\label{chap:intro}
      \paragraph{Some history}
        Model theory concerns itself with the interaction between syntactic objects, first order formulas and consistent collections thereof, and the mathematical structures which interpret them. A fundamental notion is that of the (first order) theory $Th(\mathcal{S})$ of a given structure $\mathcal{S}$.
        %In general, there is no to expect the existence of any simple description of said collection.
        %Since G\"odels result on the incompleteness of arithmetic we know this class might fail to be decidable.
        %and even if decidible, it might admit no simple 'description' in the human sense
        %the product operation being the operation on said quotient induced by 	cancellation of consecutive pairs of dual letters. Each element can be represented by a unique such word which is reduced, that is, one where such operation cannot be performed.
        
        % $z\bar{z}$ and $\bar{z}z$ .
        
        Free groups are a very natural and simple to describe mathematical object that plays a role in several areas outside of algebra, such as topological group theory. Interest in the model theory of free groups can be traced back to the forties by one of the fathers of modern logic, Alfred Tarski, who proposed the following question:
        \begin{question}
        	Do $\F_{m}$ and $\F_{n}$ have the same order theory for $2\leq m<n$?
        \end{question}
        
        The free group of rank $m$, $\F_{m}$ can be conceived as the collection of all finite words in $m$-letters and their iverses, up to the cancellation of consecutive pairs of dual letters, so that each element can be represented by a unique reduced such word. This point of view, within the larger context of what is known as "combinatorial group theory", is the one privileged by classical reference works like \cite{magnus2007combinatorial} and \cite{schupp1977combinatorial}  dominated mathematicians initial understanding of free groups \cite{schupp1977combinatorial}.
        
        Some early successes Tarski's problem was achieved using such methods by Merzlyakov \cite{merzlyakov1966positive}, who proved that if $a$ generates a free group $\F$, then the validity of any sentence of the form $\forall x\exists y\,\,\phi(x,y,a)$ true in $\F$, where $\phi(x,y,a)$ is a system of equations over $a$, is witnessed by what is known as a formal solution: some algebraic expression of the $y$'s in terms of $x$ and $a$, rendering the equations in the system trivial.  Generalizing this result, he was able to prove that all non-abelian free groups satisfy the same positive sentences.
        % ( those that can be formed without using the negation symbol).
        % is trivial for any interpretation of $x$ in a group containing $A$.
        
        Another early milestone was the work of Makanin (see \cite{makanin1983equations}),
        who provided an algorithm for deciding whether system of equations with parameters in the free group admits a solution. Later on by refining Makanin's technique Razborov (see \cite{razborov1985systems}) managed provide a description of the set of solutions to such systems.
        %The set of solution can be covered
        %by finitely many subsets, each one of which is parametrized by finitely
        
        Meanwhile P. Serre \cite{serre1977arbres}, established a duality between algebraic decompositions of a group as a graph of groups and particular actions on such groups on simplicial trees, providing a new approach to free groups and other notions in combinatorial group theory.
        %This approach would become extremely fruitful
        Later on there was an increased interest in the possibility of generalizing Serre's work to so called $\R$-trees, metric trees where the set of branching points is allows to be non-discrete. An important problem in its day was that of characterizing those finitely generated groups admitting a free action on an $\R$-tree (for simplicial trees the answer is precisely free groups).
        It was E. Rips' outstanding contribution to the subject was to realized that the combinatorial iterative procedure behind Makanin and Razborov's work could be understood in in a geometrical sense and hence generalized to such actions, providing a solution of the problem above. More general developments of the idea appeared later in works as \cite{bestvina1995stable} or \cite{guirardel2008actions} and its predecessors, which show how certain classes of actions always decompose into simpler pieces of a certain type.
        
        Another important idea, going back to the work Gromov, is that there is endow a class of metric spaces with a topology, that some form of compactness holds under certain circumstances and that one certain large scale properties of a metric space can be distilled from the limit of a sequence of rescalings of the space by smaller and smaller positive constants. Later Paulin \cite{paulin1988topologie} and Bestvina generalized this idea to isometric actions of a group $G$ on certain metric spaces.
        
        These are some of the main ingredients in Z. Sela remarkable positive answer to Tarski's question, which came at the turn of the century in the shape of a series of six papers (culminating in \cite{Sela6}). Tarski's problem is a consequence of the stronger result that the embedding of free group as a free factor of another one is an elementary embedding. A solution to Tarski's problem was reached independently by Kharlampovich and Myasnikov \cite{kharlampovich2006elementary} using combinatorial techniques.
        
      \paragraph{This thesis}
        
        This work provides a small illustration of the versatility of Sela's geometric approach. Fix a finite group $Q$ and consider the language $\LQ$ obtained by adding unary predicates $\{P_{q}\}_{q\in Q}$ to the language of groups.
        Given a group $G$, together with an epimorphism $\fun{\pi}{G}{Q}$. This is what we call a \pr-group. By interpreting $P_{q}$ as the preimage of $q$ in $G$ by $\pi$ we can promote $G$ to an $\LQ$ structure.
        
        This sort of expansion is fairly common move in model theory. Sometimes a mathematical object can be seen
        both as an $\mathcal{L}_{1}$-structure and as an $\mathcal{L}_{2}$-structure, where the theory is 'well-behaved' in both cases in a certain sense, while the $\mathcal{L}_{1}\cup\mathcal{L}_{2}$ structure
        %is very rich in definable sets and
        exhibits an extremely wild behaviour. In our case, the evidence collected so far points in the opposite direction.
        %For example, the poorness in definable sets of a particular structure might make it (more properly speaking, its theory) fall under a certain dividing line. Adding some natural (from the mathematicians point of view) predicate or function to the sructure might make it automatically fall into the other side of
        
        We are interested in the case in which $\F$ is a finitely generated free group and $\pi$ an epimorphism. What is the general picture? We conjecture the following generalizations of Sela's result:
        \begin{conjecture}
        	Given a finite group $Q$, non-abelian free groups
        	$F_{1},F_{2}$ and epimorphisms $\pi_{i}$ from $F_{i}$ to $Q$ for $i=1,2$
        	the pairs $(F_{1},\pi_{1})$ and $(F_{2},\pi_{2})$ are elementary equivalent as $\LQ$-structures.
        \end{conjecture}
        Unfortunately due to the complexity of the techniques involved, are not able to offer an aswer to the question in this work. We dapt several of his tools to the study of the language $\LQ$.
        
        A solution of a system of equations (for simplicity without constants) $\Sigma(x)=1$ in a given free group $\F$ can be corresponds to a homomorphism from the group $G_{\Sigma}$ with presentation $\subg{x|\Sigma(x)=1}$.
        Makanin showed the existence of finitely many finite chains of proper epimorphisms, terminating in a free group:
        \begin{align*}
        	G_{\Sigma}\xrightarrow{\eta_{0}} L_{1}\xrightarrow{\eta_{1}} L_{2}\xrightarrow{\eta_{2}} \cdots L_{m}=F_{k}
        \end{align*}
        Such that any $f\in Hom(G_{\Sigma},\F)$ can be written as $\eta_{0}\circ\tau_{1}\circ\eta_{1}\circ\tau_{2}\cdots\eta_{m-1}\circ h$
        for automorphisms $\tau_{j}$ of $L_{j}$ and a homomorphism $h$ from $F_{k}$ to $\F$. The groups $L_{i}$ happen to be so called limit groups: finitely generated models of the universal theory of a non-abelian free group. This is also one of the initial steps in Sela's work (see \cite{Sela1}), whose proof is based on the Bestvina-Paulin limiting method and Rips' analysis of axions of groups on real trees, which makes it easily generalizable to hyperbolic groups (see \cite{reinfeldtmakanin}).
        
        In the language $\LQ$ the basic quantifier free formula that needs to be considered is the conjunction of a system of equations together with conditions of the form $P_{q}(x_{j})$ (this is a special case of the notion of an equation with rational constraints, explored in works such as \cite{diekert2005existential} under a different angle, more akin to computer science).
        
        We work in the category finitely generated \pr-groups as a category,
        a morphism from $(G,\pi)$ to $(G',\pi')$ being any homorphism $\fun{f}{G}{G'}$ for which
        $\pi'\circ f=\pi$. Solutions to any condition as above can be identified with morphisms from a finitely generated \pr-group $(G,\pi_{G})$ to $(\F,\pi)$.
        % Say something about the first chapters.
        
        In chapter \ref{MR chapter} we provide a description of the family of such morphisms by generalizing Sela's version of the Makanin-Razborov procedure and one of its refinements (namely that of a taut-Makanin Razborov diagram) to our new category.
        
        The above mentioned morphism-classifying sequences of epimorphisms are what one usually refers to as a resolution. Our approach to resolutions is formally slightly different from the one found in Sela in the sense that we use a rooted tree to index what would be the free factors of a resolution in a classical sense. For the scope of this work this probably makes little difference, but we think it might be useful in later stages of the analysis of $\LQ$-sentences.
        
        The concept analogous to limit groups is that of what we call
        \pr-limit groups: finite models of the universal theory of our reference structure $(\F,\pi)$ or, equivalently, those finitely generated \pr-groups which are discriminated by morphisms to $(\F,\pi)$. As in the standard case, this notion can be eventually be shown to be independent from the particular reference structure. Not every homomorphism $\fun{\pi}{L}{Q}$ from a limit group to $Q$ makes it into a \pr-limit group, but we will formulate a sufficient (and necessary) condition for this to be the case, akin to that of a constructible limit group (see \cite{bestvina2009notes}).
        % One constraint is that any abelian subgroup of $L$ is mapped to a cyclic group
        
        In chapter \ref{Towers chapter} we present Sela's notions of a tower (a group together with a certain algebraic structure) and a test sequence (a sequence of homomorphisms from the tower to the free group), suitably adapted to our setting. Test sequences are a generalization of small cancellation sequences. They play a crucial role in Sela's generalization of Merzlyakov's theorem for $\forall\exists$-sentence, where universal variables to be constrained by certain algebraic relations.
        
        Our exposition differs from Sela's in that we isolate those properties of an individual sequence must satisfy from those that a collection of sequences has to satisfy in order for the theorem to hold. In the final section of the chapter we give a proof of the existence of such families.
        
        %        Test sequences play a crucial role in Sela's generalization of Merzlyakov result for $\forall\exists$ sentences. Sela's version allows for universal variables to be constrained by certain algebraic relation.
        In chapter \ref{Merzlyakov chapter} we adapt the aforementioned theorem to our setting, together with the following generalization of the other result by Merzlyakov mentioned before:
        \begingroup
        \def\thetheorem{\ref{equivalent}}
        \begin{theorem}
        	Let $A$ be a free factor of non-abelian free groups $F_{1}$ and $F_{2}$ and for $i=1,2$ let $\pi_{i}$ be a surjective homomorphism from $F_{i}$ to the finite group $Q$.
        	
        	Then $Th^{+}_{A}(F_{1},\pi_{1})=Th^{+}_{A}(F_{2},\pi_{2})$.
        \end{theorem}
        \addtocounter{theorem}{-1}
        \endgroup
        Here $Th_{A}^{+}(M)$ stands for the positive $A$-theory of $M$, that is, the collection of all sentences built without use of the negation symbol with constants in $A\subset M$.
        
        The results above seem to suggest there is little interaction between the newly introduced predicates and the old family of definable sets. We hope to be able to explore this line of thought in future work.

%CONSIDER GROUPS as trivially restricted graded

\chapter{Preliminaries}

%\section{The free group}
\section{Reminder: the free group}
  
  We start with a quick reminder of what free groups and free products are; see \cite{lyndon2015combinatorial} for a thorough introduction. Given a set of elements $X$ a group $F$ containing $X$ is free over $X$ if any map from $X$ to a group $G$ extends uniquely to some homomorphism from $F$ to $G$. This object is unique up to isomorphism relative to $X$ and we denote it by $\F(X)$. In this situation any element of $G$ can be uniquely written as a word in $\{x,\bar{x}\}_{x\in X}$, where $\bar{x}$ is a new letter (formal inverse) dual to $x$. The group $\F(X)$ can be constructed by taking the quotient of the set of such terms (where inverses are a priori formal) by the relationship generated
  by equivalences of the form $uxx^{-1}v\sim uv$, where $x\in X^{-1}$. Any element can be represented by a reduced word, i.e., one in which a consecutive pair $xx^{-1}$ does not appear. Such a word is cyclically reduced if in addition the first and last elements are not mutually inverse. A group $F$ which is free over a certain set $X$ is called free and any such $X$ is called a base of $F$. All bases of a given free group have the same cardinality which, if finite, is denoted as the rank of the free group in question. The free group of rank $k$ is often denoted as $\F_{k}$.
  
  Given a collection $\{G_{i}\}_{i\in I}$ of groups, its free product, denoted by an expression such as $\bfrp{G_{i}}{i\in I}{}$ or $G_{1}\frp G_{2}\cdots G_{k}$
  is a group containing each $G_{i}$ and satisfying a similar universal property. Namely, given any group $H$ any collection of homomorphisms $\{\fun{f_{i}}{G_{i}}{H}\}_{i\in I}$ extends to a unique $f$, which we will denote as $\bdcup{f_{i}}{i\in I}{}$. Any non-trivial element of $G$ can be written in a unique way as a product $a_{1}a_{2}\cdots a_{m}$, where $a_{l}$ and $a_{l+1}$ belong to $A_{i_{l}}$ and $A_{i_{l+1}}$ for any $1\leq l\leq m$.
  
  A free decomposition is an isomorphism (often thought as an identity) between a group and some non-trivial free product.
  
\section{A few words about first order logic}
  
  \newcommand{\LS}[0]{$\LL$-structure }  %tc
  \newcommand{\LSs}[0]{$\LL$-structures } %tc
  \newcommand{\UU}[0]{\mathfrak{U}}
  Let us introduce here some basic notions of first order logic with which to formulate some of the results in this work. We will gloss over many formalities.
  A complete yet swift introduction to the matter can be found in \cite{tent2012course}; for a more fleshed-out the reader can take for example \cite{shoenfieldbook}.
  
  A language $\mathcal{L}$ is a collection of symbols including constants, functions and relation symbols. Additionally, each function symbol and relation symbols are associated with a natural number
  $n\geq 1$, referred to as the arity of the symbol. When naming the constituents of $\mathcal{L}$ we will often include the $n$-arity of the symbols in a subscript between parenthesis.
  For most purposes, constants can be regarded as functions of arity $0$.
  In order to do anything meaningful with $\mathcal{L}$ one needs some additional logical symbols, namely:
  \elenco{
  	\item A set of infinitely many variables: $\{x_{n}\}_{n\in\N}$
  	\item The equality symbol $=$,
  	\item Negation $\neg$, conjunction $\wedge$ and disjunction $\vee$ symbols
  	\item An existential quantifier $\exists$ and a universal quantifier $\forall$
  }
  In addition to this we will also use commas and parenthesis as auxiliary symbols (in a somewhat loose way in the case of parenthesis), although they are not strictly necessary.
  %  for defining unambiguously interpretable terms and formulas (this can be achieved by means of the so called Polish notation).
  %  When writing down formulas (or rather, terms of the meta-language to be interpreted as formulas) we will drop
  Likely, one loses no expressive power by foregoing one of the two quantifier symbols, or one among $\{\wedge,\vee\}$. Expressions involving them can be safely regarded as abbreviations of expressions without them. In our case we will restrict that treatment to the logical connectors $\sla$,$\sra$ and $\slra$.
  %  Using these, we first construct what we call an $\mathcal{L}$-term is a string of symbols in $\mathcal{L}$ and
  %  An $\mathcal{L}$-term is a $s$ finite string of constants, variables and function symbols constructed inductively in the common reference
  $\mathcal{L}$-terms and $\mathcal{L}$-formulas are sets of finite strings of symbols of $\mathcal{L}$ and logical symbols \LS which are 'meaningful' or 'sound' in an intuitive way.
  An $\mathcal{L}$-terms is any such string that can be generated in finitely many steps using the following rules:
  \elenco{
  	\item The string containing a single variable or constant is an $\mathcal{L}$-term
  	\item Given a function symbol $f^{(m)}\in\mathcal{L}$, and terms $t_{1},t_{2},\cdots t_{m}$, the string $f(t_{1},t_{2},\cdots t_{m})$ is an $\mathcal{L}$-term.
  }
  %  Strictly speaking, terms as $a+b$ should be considered as abbreviations of 'properly formed' $\mathcal{L}$-terms. The latter can be thought as abbreviations of $\mathcal{L}$-terms as defined above: i.e. we can think $a+b$ as standing for $+(a,b)$ and so on. As
  for $\mathcal{L}$-formulas, the generating rules are:
  \elenco{
  	\item If $t_{1}$ and $t_{2}$ are terms, then $t_{1}=t_{2}$ is a formula.
  	\item If $R^{(m)}\in\LL$ is a relation symbol, and $t_{1},t_{2},\cdots t_{m}$ are terms, then $R(t_{1},t_{2},\cdots t_{m})$ is a formula.
  	\item If $\psi$ is an $\mathcal{L}$-formula, then $\neg\psi$ is one as well.
  	\item If $\psi_{1}$ and $\psi_{2}$ are $\mathcal{L}$-formulas, then $(\psi_{1}\vee\psi_{2})$ and $(\psi_{1}\wedge\psi_{2})$ are one as well.
  	\item If $\psi$ is a formula then $\exists x\,\,\psi$ and $\forall x\,\,\psi$ are one as well.
  }
  Those formulas as in the first two points are called \emph{atomic}.
  %  It can be shown that any $\mathcal{L}$-term or an $\mathcal{L}$-formula admits a unique construction tree.
  For any occurrence in a formula of a quantifier $Qx$ arising from an application of the third role, ones refers to the corresponding formula $\psi$ as the scope of the quantifier and to any occurrence of $x$ in it as a bounded occurrence of $x$.  A variable with at least one non-bounded occurrence in a certain formula is called a free variable. A formula without free variables is called a sentence.
  
  A formula with no occurrence of the negation symbol $\neg$ is called a \emph{positive formula} and one with no instances of quantifiers is called a \emph{quantifier free}.
  %  and one in which only the $\forall$ quantifier appears is said universal.
  Typically when introducing a term $t$ one writes $t(x_{1},x_{2},\cdots x_{k})$, where $x_{1},x_{2},\cdots x_{k}$ are  (distinct) variables containing all those appearing in $t$. Likewise $\phi(x_{1},x_{2}\cdots x_{k})$ will stand for a formula with free variables among $x_{1},x_{2}\cdots x_{k}$.
  The notation $\pi(x_{1},x_{2}\cdots x_{k})$ might also indicate a set of formulas whose free variables are among $x_{1},x_{2}\cdots x_{k}$.
  
  For $\mathcal{L}$ as above, a \LS $\UU$ is given by a set $U$, referred to as the universe of the structure, together with a tuple an interpretation $Z^{\UU}$ for each of the symbols $Z\in\mathcal{L}$, which is:
  \elenco{
  	\item an element of $U$ in case $Z$ is a constant
  	\item a function from $A^{m}$ to $A$, if $Z$ is a function symbol of arity $m$
  	\item some subset of $A^{m}$ if $Z$ is a function symbol of arity $m$
  }
  In the future we might incur in a certain lack of precision and refer to $\UU$ simply as $U$.
  Given an \LS $\UU$, one can interpret any term $t(x_{1},x_{2},\cdots x_{k})$ as a function $\fun{t^{\UU}}{U^{k}}{U}$. Likewise, a formula $\phi(x_{1},x_{2},\cdots x_{k})$ for $k\geq 1$ can be interpreted as a predicate in $\phi(\UU)\subset U^{k}$. The formal definition merely reflects the intuitive way such expressions are used in common mathematical practice; the interested reader can take a look at the references mentioned above.
  
  In the case of a $\mathcal{L}$-sentence $\phi$ one can ask whether $\phi$ holds (is true) in $\mathcal{M}$, denoted by the notation $\mathcal{M}\models\phi$. Again, the concept the formal definition merely reflects the informal commonsense notion of satisfaction.
  
  Given a term $t(x_{1},x_{2},\cdots x_{k})$ and any tuple $s_{1},s_{2},\cdots s_{k}$ of terms
  (possibly other variables) the expression $t(s_{1},\cdots s_{k})$ denotes the term obtained by replacing each variable $x_{i}$ by $s_{i}$. One can also 'substitute' $t_{1},t_{2},\cdots t_{k}$ in place of free variables of a formula, but this is slightly more subtle.
  %   It only works in the same way if we can assume no quantified variable in $\phi$ appears in one of the $t_{j}$, but one can always replace $\phi(x_{1},x_{2},\cdots x_{k})$ with an equivalent formula so that the latter holds.
  
  %   Given terms
  %  $t_{1},t_{2},\cdots t_{k}$ for any quantifier $Qx\phi$ in $\phi$ one needs to chose a variable $u$ not appearing in either any of the $t_{j}$ and replace $x$ by $u$ all over the scope of the quantifier, together with the quantifying expression itself.
  
  As it is common practice in model theory, we might use a single letter, such as $x$ to denote a whole tuple of $k=|x|$ many variables. In this case, given $Q\in\{\forall,\exists\}$, an expression like $Q\in\{\forall,\exists\}$ then $Qx$ will abbreviate the string $Qx_{1}\,Qx_{2}\,\cdots Qx_{k}$.
  A formula of the form $\forall x\,\phi(x)$, where $\phi$ is quantifier free is called a universal formula. One of the form
  $\forall x\,\exists y\,\,\phi(x,y)$ with $\phi(x,y)$ quantifier free is called an $\forall\exists$-formula.
  
  Given an \LS $\UU$ and a subset $A\subset U$, it is formally useful to consider the language $\mathcal{L}_{A}$ obtained by adding a constant $c_{a}$ to the language
  for any element $A$; in practice we will use the same letter $a,b\cdots$ denote both the element in $A$ and the associated constant. There is an obvious way one can extend $\UU$ to an $\mathcal{L}_{A}$-structure $\UU_{A}$.
  
  %  An $\LL$-theory is simply a collection of $\mathcal{L}$-sentences; an $\LL$-structure is said  model of $T$ if and only if $\mathfrak{M}\modelof\phi$ for any $\phi\in T$.
  The (positive) elementary theory of an \LS $\UU$, denoted by $Th(\UU)$ ($Th^{+}(\UU)$), is the collection of all (positive) $\mathcal{L}$-sentences valid in $\UU$.
  We say that two \LSs $\UU$ and $\UU'$ are (positively) elementarily equivalent, denoted by $\UU\equiv\UU'$ ($\UU\equiv^{+}\UU'$) if $Th(\UU)=Th(\UU')$ ($Th^{+} (\UU)=Th^{+}(\UU')$).
  
  Given \LSs $\UU$ and $\UU'$  such that their respective universes contain some common subset $A$, we say that they are elementarily equivalent over $A$,
  denoted by $\UU'\equiv_{A}\UU$, if $\UU_{A}'\equiv\UU_{A}$. The positive theory of $\mathfrak{U}$, denoted by $Th^{+}(\mathcal{U})$ is the collection of all positive sentences valid in $\UU$. We denote the corresponding weak elementarily equivalent relationships by $\equiv^{+}$ and $\equiv^{+}_{A}$. Considering only universal sentences instead of positive ones we obtain the notion of the positive theory and positive elementary equivalence respectively.
  
  %  Given two \LSs $\mathfrak{M}$ and $\mathfrak{N}$, where $M\subset N$, we say that $M$ is an elementary submodel of $\mathfrak{N}$ if $M\equiv_{M}N$.
  An $\LL$-theory $T$ implies an $\LL$-sentence $\phi$, something we denote as $T \vdash \phi$, if $\mathcal{M}\modelof \phi$ for any model $\mathcal{M}$ of $T$. We say that two $\LL$-formulas $\phi(x)$ and $\psi(x)$ are
  $T$-equivalent if $T\modelof\forall x\,(\phi(x)\lra\psi(x))$.  Any formula $\phi(x)$ is $\nil$-equivalent to one of the form:
  \begin{align*}
  	\forall y_{1}\,\,\exists y_{2}\,\,\cdots\forall y_{2m-1}\exists y_{2m}(\bvee{j=1}{(\bwedge{i=1}{\psi_{i}^{j}}{r_{i}})}{m})
  \end{align*}
  where each $\psi_{i}^{j}$ is either an atomic formula or the negation of one; we will refer to this as a normal form (for $\phi(x)$). If the formula we started with is positive, then the $\psi_{i}^{j}$ can be all taken to be atomic.

\section{Introducing our setting and some basic definitions}
  
  Groups can be seen as $\mathcal{L}_{gp}$ structures, where $\mathcal{L}_{gp}=\{1,\cdot,(-)^{-1}\}$, containing symbols for the identity, multiplication and inverse respectively.
  
  Fix a finite group $Q$ and let $\LQ$ be the language resulting from adding an unary predicate $P_{q}$ to $\mathcal{L}_{gp}$ for any $q\in Q$.
  
  By a \pr-\emph{group} we intend a pair $(G,\pi)$, where $\pi:G\to Q$ is a homomorphism. Such an object can be seen as an $\LQ$-structure in a natural way by interpreting $P_{q}$ as the preimage of $q$ by $\pi$. The epression as $T_{\pi}=Th(G,\pi)$ denotes the corresponding $\LQ$-theory.
  
  We are interested in the particular case of a \pr-group $(\F,\pi)$, fixed for most of this work, where $\F$ is a non-abelian free group and $\pi$ a surjection. This structure will be usually abbreviated as $\qG{F}$.
  
  A \emph{morphism}  between  \pr-groups $(G,\pi)$ and $(G',\pi')$ is a homomorphism $\phi:G\to G'$ satisfying $\pi_{G'}\circ \phi=\pi_{G}$. To help lighten the notation we will (for now) adopt the notation $\qG{G}$ to refer to a \pr-group $(G,\pi_{G})$, $\mathcal{H}$ one of the form $(H,\pi_{H})$ and so forth. Likewise, $\qG{F}$ will stand for our fixed pair $(\F,\pi)$.
  
  A fundamental notion in the study of the theory $T_{\F}$ is that of a limit group. Many basic questions about the theory $T_{\F}$ can be restated in terms of limit groups and the homomorphisms between them.  This notion and all of the main related results can be easily adapted to the theory $T_{\pi}$.
  
  \begin{comment}
  To help lighten the notation we will (for now) adopt the following convention.
  From now on, $\qG{G}$ will refer to a \pr-group $\qG{G}$, $\p{H}$ to one of the form $(H,\pi_{H})$ and so forth.
  \end{comment}
  \begin{comment}
  Likewise, from now now, $\qG{F}$ will stand for a fixed pair $(\F,\pi_{\F})$ where $\F$ is a finitely generated free group, $\pi_{F}$ a map (epimorphism) from $\F$ to some finite group $Q$.
  \end{comment}
  %	Although the definition leaves the possibility that two \qgs have
  %	the same underlying group open, we can in fact work within an isomorphic category
  %	in which this is never the case. Of course, we might still be interested in the case when two 	underlying groups might be wisomorphic. (IS THIS RELEVANT?)
  
  %  \begin{definition}
  	A sequence $(f_{n})_{n}$ of homomorphisms from a countable group $G$ to a group $H$ is said 	\emph{convergent} if and only if for any $g\in H$ eventually either $f_{n}(g)=1$ or $f_{n}(g)\neq 1$. 	By the limit kernel of such a sequence $(f_{n})_{n}$ we intend $\limker{f_{n}}{n}=\{w\in G\,|\,\#\setof{n}{f_{n}(w)=1}=\infty\}$.
  %  \end{definition}
  It is a trivial matter to check that the limit kernel is a subgroup of $G$. By the \emph{limit 	quotient} of $(f_{n})_{n}$ we intend the quotient map from $G$ to $G'=G/\limker{f_{n}}{n}$. Often we will use the term to refer to the quotient itself.
  If the $f_{n}$ are morphisms from some \pr-group $\qG{G}$ to $\qG{H}$, then we will be referring to the \pr-group $\qG{G'}$ obtained by pushing forward the \pr-structure of $G$ onto $G'$.
  %
  %  If all the $f_{n}$ are morphisms from $\pG{G}$ to $\pG{H}$ then the \pr-structure on $G$ can be pushed to a \pr-structure on $\pG{G}$ in an obvious way
  
  Given a countable set $X$, we can identify $\pow{X}$ with $2^{X}$, which is a compact separable metric space for the product topology. Convergence in the sense above is then equivalent to the  	topological convergence of $Ker(f_{n})$ in $2^{G}$ . In particular:
  \begin{observation}
  	Any sequence of homomorphisms between two groups contains a convergent subsequence.
  \end{observation}
  
  \begin{definition}
  	By a \pr-limit group we intend a \pr-group $\qG{L}$ which is the limit quotient of a finitely generated \pr-group $\qG{G}$ by a converging sequence of morphisms from $\qG{G}$ to $\qG{F}$.
  \end{definition}
  
  Letting $Q=\tg$ one recovers the standard notion of a limit group.
  We would like to remark that the notion of a \pr-limit group is strictly stronger than that of a limit group equipped with a morphism to $Q$.
  An obvious restriction is the fact that given a \pr-limit group $(L,\pi_{L})$ and any $x\in L$ the image by $\pi$ of the centralizer $Z_{L}(x)$ is necessarily cyclic.
  
  \begin{comment}
  \df{
  	By a \rs \pr-group we intend a pair $(\pmb{G},\iota)$, where $ \pmb{G}$ is a \qg and $\iota$ an isomorphism between a finitely generated \qsub $\pmb{A}\leq \F$ and a subgroup of $\pmb{G}$.
  }
  \end{comment}
  %  \begin{definition}
  	A \pr-subgroup of a \pr-group $\qG{G}$ is simply a subgroup of $\qG{G}$ endowed with a \pr-structure making the inclusion map a \pr-morphism.
  %  \end{definition}
  
  %  \begin{definition}
  	By a ($A$-)\emph{\rs }\pr-group we intend a pair $(\qG{G},\qG{A},\iota)$, where $\qG{G}$ is a \pr-group and $\iota$ an injective morphism from a finitely generated \pr-subgroup $\qG{A}$ of $\qG{F}$ to $\qG{G}$.
  	%  	\begin{comment}
  	%  	We say that the group is restricted over A.
  	%  	\end{comment}
  %  \end{definition}
  %(What about if the two gradings are different?)
  \begin{comment}
  Equivalently, a tuple $a$ of generators of $A$ is provided, in which case we can think of
  \end{comment}
  
  Since we will never be interested in comparing \rs groups differing only in the particular isomorphism $\iota$ and not in the underlying groups, we will often call such $(\qG{G},A,\iota)$ by the name $\rqG{G}{A}$
  %  	and refer to $\iota$ as
  %  	$\iota_{A}^{G}$,
  and possibly even identify $A$ with its image. %% (make it $\iota_{A}^{G}$)
  
  \df{
  	By a morphism from an $A$-\rs \pr-group $(\qG{G},A,\iota)$ to another one
  	$(\qG{G'},A,\iota')$ we intend a morphism $\phi:\qG{G}\to\qG{G'}$
  	such that $\iota'\circ \phi=\iota$.
  }
  
  %  In case of ambiguity, we will include the adjective \rs in front of the word morphism as well.
  %  We will in general avoid using the term \rs in front of the word 'morphism' but less so in front of names indicating particular types of such morphisms, such as retraction
  
  Of course, given a \rs \pr-group $\rqG{G}{A}$ and a morphism $\phi:G\to H$ injective on
  $A$, there is a unique way to push the $A$-\rs structure of $G$ to one on $H$, with respect to
  which $\phi$ is a \rs \pr-morphism.
  %  For any $A\leq F$ there is a natural \rs \pr-group structure of $\qG{F}$ over $\qG{A}$ and this, of course will be the meaning of the notation$\rrqG{F}{A}$.
  
  \df{
  	A \rs \pr-limit group is a \rs \qg $\rqG{L}{A}$ where $L$ is the quotient of a convergent sequence of morphisms $\phi_{n}$ from a finitely generated \rs \pr-group $\rqG{G}{A}$ to $\rqG{F}{A}$ and $\phi_{n}\circ\ctm{L}{A}=\ctm{F}{A}$ for all $m$.
  }

  \begin{definition}
  	A  \gr \pr-group is given by a tuple $(\rqG{G}{A},(P_{j})_{i\in I})$,
  	%  	 where $0\in i$
  	comprising:
  	\enum{i)}{
  		\item A \rs \pr-group $\rqG{G}{A}$
  		\item A family $\{P_{j}\}_{j\in I}$ of finitely generated subgroups of $G$, which we call parameter subgroups.
  	}
  	%  	In case $A$ is non-trivial we require the index set $I$ to contains $0$ and $A\leq P_{0}$.
  	%  	The $P_{i}$ will be called \emph{parameter subgroups} of $G$.
  	We will often refer to such a structure by an expression of the form $\gqG{G}{A}{P}{I}$.
  	%  	The images of the groups $P_{i}$ will be called \emph{parameter subgroups} of $G$.
  	A \gr \pr-limit group is just a \gr \pr-group for which the underlying \rs \pr-group is a \rs limit group.
  \end{definition}
  \begin{comment}
  Need a commen on the fact that $p$ is a valid name for a variable.
  \end{comment}
  
  By a marking of a group $G$ we intend a map $\fun{\iota}{x}{G}$
  where $x=(x_{i})_{=1}^{k}$ is a finite set of syntactic variables
  % form a finite set of syntactic variables to a generating set of the group
  . A marked group is a pair $(G,\iota)$ where $G$ is a group and $\fun{\iota}{x}{G}$ a marking of $G$ for which $\subg{\iota(x)}=G$.
  One should think of $\iota$ as an interpretation of the tuple $x$ in $G$. In accordance, we will always refer to $\iota(x)$, or its image in a quotient of
  $G$ as simply $x$. If instead of a group we have a \pr-group, \rs \pr-group we will talk of marked \pr-groups, marked \rs \pr-groups and so forth.
  %  	with the extra condition that $\subg{\iota(x)}=G$.  .
  %  	link surjectivity condition to some adjective
  %  	In most cases, such an object will be referred to simply as 'a group', but denoted by an expression such as $G$.
  The expression $x(G)$ will stand for $\iota(x)$ and $\iota(x_{i})$ for $x_{i}(G)$.
  A marked \rs group $G$ is a pair $(G,\iota)$ where $\iota$ is a marking such that $\subg{A,\iota(x)}=G$.
  A marked \gr group $A$ consists of a tuple $(G,(\iota_{i})_{i\in I\setminus\{0\}},\iota_{0})$, together with markings
  $\fun{\iota_{i}}{p_{i}}{P_{i}}$ and $\fun{\iota}{x}{G}$
  such that $\pi(\subg{A,\iota(x),\iota_{i}(p_{i})})=G$.
  \begin{comment}
  A marked \rs group $G$ is a pair $(G,\iota)$ where $\iota$ is a marking such that $\subg{A,\iota(x)}=G$.
  A marked \gr group $A$ consists of a pair $(G,(\iota_{i})_{i\in I})$, together with markings
  $\fun{\iota_{i}}{p_{i}}{P_{i}}$ and $\fun{\iota}{x}{G}$
  such that $\subg{A,\iota(x),\iota_{i}(p_{i})}$.
  \end{comment}
  A marked \rs group $G$ is a pair $(G,\iota)$ where $\iota$ an $x$-marking such that $\subg{A,\iota(x)}=G$.
  A marked \gr group $A$ consists of a tuple $(G,(\iota_{i})_{i\in I\setminus\{0\}},\iota_{0})$, together with markings
  $\fun{\iota_{i}}{p_{i}}{P_{i}}$ and $\fun{\iota}{x}{G}$
  such that $\subg{A,\iota(x),\iota_{i}(p_{i})}$.
  %  	We might refer to it in an abbreviated form as $\qG{G}_{a}[p_{i}]$.
  The notation might replace $A$ by some finite tuple $a$ of generators.
  Notice that given an epimorphism of (\rs, \gr) groups, a marking of the domain pushes forward to a marking of the target univocally.
  \begin{definition}
  	A group $G$ is called $CSA$ if and only if every maximal abelian subgroup $A$ oef $G$ is malnormal, i.e. $g\in A$ whenever $A\cap A^{g}\neq\tg$. It is called commutative transitive if and only if the relationship $R(x,y)\cong[x,y]=1$ is transitive when restricted to $G\setminus\{\tg\}$.
  \end{definition}
  
  The following is well-known and easy to prove:
  \begin{observation}
  	In a commutative transitive group any non-trivial element is contained in a unique maximal abelian subgroup. Any $CSA$ group is commutative transitive.
  \end{observation}
  For the following theorem see for example \cite{champetier2005limit}
  \begin{theorem}
  	Limit groups are $CSA$.
  \end{theorem}
  
  \begin{remark}
  	Given \pr-groups $\qG{G}$ and $\qG{H}$ we denote the unique extension $\pi_{K}$ of $\pi_{G}$ and $\pi_{H}$ to $K=G\frp H$ as $\qG{G}\frp \qG{H}$. If $\qG{G}$ can be promoted to an $A$-\rs \pr-group, the same injection map makes $\qG{K}$ into one.
  \end{remark}

\section{Preorders}
  
  \label{poset section} We will sometimes prove properties of a certain class of strucutres $\mathcal{C}$ by induction on some well-founded partial order of $\mathcal{C}$.
  Given a partially ordered set (poset) $(J,\leq)$, we denote by $\equiv$, as usual $\leq\cap\geq$. A strict poset will be one for which $\leq$ is irreflexive.
  %(do we need this? give a name to the case when equivalence classes are trivial)
  Recall a poset $(J,\leq)$ is said well-founded if there are no infinite chains $\lambda_{1}>\lambda_{2}>\cdots\lambda_{n}\cdots$ of elements of $P$. For $\lambda\in J$ we let $J\rst_{\lambda}=\setof{\mu\in P}{\mu\geq\lambda}$.
  Any poset determines a successor relationship, $\succ$, where $\lambda\prec\mu$ if and only if  $\lambda<\mu$ and there is no element between $\lambda$ and $\mu$.
  %    By a discrete poset we intend one for which any element which is non-maximal has a predecessor and any which is not minimal has a successor.
  %The poset $(J,\leq)$ will be said to have finite width if it is lower-discrete and every element has only finitely many successors.
  %Given posets $(I,\leq)$ and $(J,\leq)$ a map $\fun{f}{I}{J}$ is said monotonous
  By a rooted simplicial tree we will intend a poset $(J,\geq,r)$ together with a distinguished element $r$, the root such that:
  \enum{i)}{
  	%    	\item $r\geq j$ for any $j\in J$
  	\item If $i,j\in J$ have a common lower bound, they are comparable.
  	\item Any $i,j\in J$ have a least common upper bound.
  	\item For any $\lambda\in J$ there is a unique sequence $r=\lambda_{0}\succ \lambda_{1}\succ\cdots \lambda_{n}=\lambda$ for some $n\in\N$.
  }
  By a branch of it we mean a sequence $\lambda_{0}\succ \lambda_{1}\succ\cdots \lambda_{n}$ or $\lambda_{0}\succ \lambda_{1}\succ\cdots \lambda_{n}\succ\cdots$.
  %% Denote it by $(J,r)$
  %    	In this context, we will use the terms 'child' and 'parent' instead of 'successor and 'predecessor'.
  Whenever $\lambda\prec\mu$ in this context we will refer to $\lambda$ as a child of $\mu$ and to $\mu$ as the (unique) parent of $\lambda$, or $p(\lambda)$. Let
  $Ch(\lambda)$ stand for the set of all the children of $\lambda$. Given $\lambda\geq\mu$, we will refer to $\lambda$ as an ancestor of $\mu$ and to $\mu$ as a descendant of $\lambda$.
  We call \emph{leaves} those elements of $J$ which are minimal and denote the set of all of them by $\lvs{J}$, while $\hat{J}$ will stand for the complement of the latter in $J$.
  %    A branch ending in a leaf will be called complete.
  %    Given $(J,r)$, the expression $(\hat{J})$ will stand for the restriction of $(J,r)$ to the set of all the non-leaves.
  %Whenever $\lambda\succ\mu$ we will refer to $\lambda,\mu$ as standing in a parent-child relationship and whenever 	$\lambda\geq\mu$, as them standing in the ancestor-descendant relationship.
  % Some comment on restrictions.
  % We say that the the tree is finite branching if any $\lambda\in J$ has only finitely many children.
  \begin{comment}
  Given $\lambda\in (J,r)$,
  %% by $J\restriction_{\lambda}$ we intend the set $\setof{\mu}{\lambda\geq\mu}$
  by $(J\restriction_{\lambda})$ or $(J)$ we will intend the substructure generated by the set $\setof{\mu}{\lambda\geq\mu}$.
  \end{comment}
  The following is a consequence of K\"onig's lemma.
  \begin{lemma}
  	\label{Koenig's lemma}
  	Every finitely branching well-founded rooted simplicial tree is finite.
  \end{lemma}
  
  Given preorders $\leq_{1}$ and $\leq_{2}$ on a set $P$, we denote by $\leq_{0}=\leq_{1}\times\leq_{2}$ the relationship defined by $p\leq_{0}q\lra\,p\leq_{1}q\vee(p\equiv_{1}q\wedge p\leq_{2} q)$. It
  is well-known and easy to check this is a preorder and
  %with $\equiv_{0}=\equiv_{1}\cap\equiv_{2}$
  and that it is well-founded whenever both the $\leq_{i}$ are
  % and strict whenever of the $\leq_{i}$ are.
  \begin{notation}
  	Given a class $\mathcal{C}$ of structures and a function $\fun{f}{\mathcal{C}}{\N}$ we let $\leq_{f}$ be the preorder given by $\leq_{f}=f^{-1}(\leq_{\N})$.
  	%    	 Given functions $f_{1},f_{2},\cdots f_{k}$ from $\mathcal{C}$ to $\N$, we let $\leq_{(f_{1},f_{2},\cdots f_{k})}=\leq_{f_{1}}\cdot \leq_{f_{2}}\cdot\cdots \leq_{f_{k}}$.
  \end{notation}
  
  We say that a rooted tree $(S,r)$ extends some subtree $(T,r)$ in case each $v\in S\setminus T$ is a descendant of some $w\in\lvs{T}$. The $l$-th level of a rooted tree $(T,r)$, denoted by $[T]_{l}$ consists of those nodes of $T$ at distance $l$ from $r$. Let also $[T]_{\leq l}=\bunion{k\leq l}{[T]_{k}}{}$.

  Given a well founded preorder $(P,\leq)$ let $Tr(P,\leq)$ be the set of pairs $((T,\leq,r),p)$ where $(T,\leq,r)$ is a finite rooted tree and $p$ a monotonous map from $T$ to $P$
  taking root to root and whose restriction to $T\setminus\{r\}$ is strictly monotonous.
  
  We define a partial order $\leq_{Tr}$ on $Tr(P,\leq)$ as follows: we say that
  \begin{align*}
  	((T,\leq,r),p)<^{0}_{Tr}((T',\leq,r'),p')
  \end{align*}
  if and only if there is a subtree $ S\subset T$ and an isomorphism $f$ between $S$ and some subtree $S'$ of $T'$ such that $f(r)=r$ and:
  %% this you have to modify
  \enum{i)}{
  	\item  \label{preorder extension} $T$ extends $S$ and $T'$ extends $S'$
  	\item  \label{str existence} Either $S=T$ or $p'(f(u))<p(u)$ for some $u\in S$
  	\item  \label{str necessity}$p'(f(u))<p(u)$ if $u\in\lvs{S}\setminus\lvs{T}$
  	\item  \label{domination} Each $v\in T'\setminus f(S)$ lies below $f(u)$ for some $u\in \lvs{S}$ for which $ii)$ holds (in particular, $T'$ extends $f(S)$).
  }
  It is easy to check that $<_{Tr}$ is anti-reflexive and transitive.
  \begin{lemma}
  	\label{tree preorder}
  	The preorder $\leq_{Tr}$ is well-founded.
  \end{lemma}
  \begin{proof}
  	Suppose there is some infinite descending chain
  	${(((T_{k},\leq,r),p_{k}))_{k\in\N}}$ in $Tr(P,\leq)$. For each $k\in\N$ take $S_{k}\subset T_{k}$ and a map $f_{k}$ witnessing $((T_{k},\leq,r),p_{k})>((T_{k+1},\leq,r),p_{k+1})$.
  	For each $m\leq k$, let $g_{k}^{m}=f_{k}\circ f_{k-1}\circ\cdots\circ f_{m}$  (its domain of definition might be smaller than that of the $f_{m}$).
  	%Let also $U_{k}=\lvs{S_{k}}\setminus\lvs{T_{k}}\subset S_{k}$.
  	Let $U_{k}=\setof{u\in S_{k}}{p_{k+1}(f_{k}(u))<p_{k}(u)}$, as in $ii)$.
  	\begin{claim}
  		For each $L\in\N$ there is $N_{L}\in\in\N$ such that for all $k\geq N_{L}$
  		$S_{k}\supset[T_{k}]_{L}$, $f_{k}$ restricts to an isomorphism between $[T_{k}]_{\leq L}$ and $[T_{k}]_{\leq L}$ and $U_{k}\cap [T_{k}]_{\leq L-1}=\nil$.
  	\end{claim}
  	\begin{proof}
  		We will prove the existence of $N_{L}$ by induction on $L$.
  		Take $L=0$ as a base case. Then $N_{-1}=0$ does the job. The first property is satisfied by all $f_{k}$, as they take roots to roots. The second one follows from the well-foundedness of $\leq$.
  		
  		We now deal with the induction step. Suppose $N_{L}$ is given such that for any $k\geq N_{L}$ the map $f_{k}$ maps $[T_{k}]_{\leq L}$ isomorphically onto $[T_{k+1}]_{\leq L}$ and $\lvs{U_{k}}\cap[T_{k}]_{\leq L}=\nil$. Now, for any $\lambda\in[T_{N_{L}}]_{L}$ either its image $v_{k}$  in $T_{k}$ belongs to $\lvs{T_{k}}$ for any
  		$k\geq N_{L}$ or eventually this image has children.
  		
  		Suppose this occurs for $T_{M}$. Then property (\ref{str necessity}) implies that $Ch(v_{k})\subset S_{k}$ and therefore
  		$f_{k}$ restricts to a bijection between $Ch(v_{k})$ and $Ch(v_{k+1})$ from that point on (it is surjective because $T_{k}$ extends $S_{k}$ as opposed to merely containing it as a subtree).
  		We conclude that $f_{k}$ restricts to an isomorphism between $[T_{k}]_{\leq N+1}$ and $[T_{k+1}]_{\leq N+1}$ for all $k$ greater or equal than some $M_{L}$.
  		The rest of the induction hypothesis for $L+1$ clearly follows from the well-foundedness of the base partial order.
  	\end{proof}
  	This, together with property (\ref{str existence}) of the definition of the order $<$ implies that $[T_{k}]_{\geq L+1}\neq\nil$ for $k\geq N_{L}$.
  	Now consider the direct limit $T_{\infty}$ of the restriction of the maps $g^{N_{L}}_{N_{L+1}}$ to $[T_{N_{L}}]_{L}$ with $[T_{N_{L+1}}]_{L}$.
  	Let $p_{\infty}$ be the map from $T_{\infty}$ to $(P,\leq)$ which restricts to $p_{N_{L}}$ on the image of $[T_{N_{L}}]_{L}$.
  	It follows from the properties of the $N_{L}$, that $p_{\infty}$ is strictly monotonous. Since the other hand $T_{\infty}$ is clearly finite-branching, so by K\"onig's lemma it has to contain some infinite branch. This contradicts well-foundedness of $(P,\leq)$.
  \end{proof}
  
  Strictly speaking, this result needs only be applied to linear trees in $Tr(P,\leq)$ (which is simpler), but we believe it can be of use in proving the iterative procedure for analyzing general $\forall\exists$ formulas, which unfortunately goes beyond the scope of this work.

\chapter{Actions on trees}
\label{Trees chapter}
\section{$\R$-trees}
  \subsection{Basic definitions}
    \newcommand{\Mod}[1]{Mod(#1)}
    So called $\R$-trees are  class of metric spaces which generalize the notion of a simplicial tree which play a crucial role in the study of the first order theory of free groups.  We will give a quick survey of the main results
    that will be needed later and point the reader to the main source for the material of this chapter, \cite{chiswell2001introduction}, for a more in detail and general account (see also \cite{bestvina2002r}).
    
    Given a metric space $(X,d)$, by an oriented segment in $X$ we denote an isometric embedding $i$ from a bounded interval $[a,b]\subset\R$, into $X$; we say that such $i$ is a segment from
    $i(a)$ to $i(b)$. If we are only interested in $i$ up to precomposition with an isometry of $[a,b]$, then we talk about an unoriented segment. In this case we will often blur the distinction between $\iota$ and its image, as this hardly generates any ambiguity. For example, we might say that the union of two segments is a segment, when what we mean is that formally the union of their images is the image of another segment.
    If instead of a bounded interval we take the entire real line $\R$ (and we allow ourselves to precompose by translations of $\R$) we obtained a geodesic line; if we take a semi-infinite interval $[a,\infty)$, a geodesic ray. A simple path in $X$ is a continuous embedding of an interval $[a,b]$ into $X$, defined up to precomposition by a homeomorphism of $[a,b]$ relative to $a,b$.
    \begin{definition}
    	A metric space is said geodesic if for any $x,y\in X$ there is a segment $\fun{\iota}{[a,b]}{X}$, with $\iota(a)=x$ and $\iota(b)=y$.
    \end{definition}
    \begin{definition}[Gromov]
    	Given a metric space $(X,d)$, and $x\in X$, for $y,z\in X$ let $(y,z)_{x}=\nicefrac{1}{2}(d(x,y)+d(x,z)-d(z,y))$. Given a positive constant $\delta$, a metric space is called $\delta$-hyperbolic if for some (equivalently for each) $x\in X$, and for all $y,z,w\in X$:
    	\begin{align*}
    		(y,w)_{x}\geq\,\min\{(y,z)_{x},(w,z)_{x}\}-\delta
    	\end{align*}
    \end{definition}
    We collect some background results.
    %  We collect some results from chapter $2$ in \cite{chiswell2001intoduction}. Note that some of them are phrased there with a greater degree of generality.
    \begin{lemma}[{\cite[p.43]{chiswell2001introduction}}]
    	Let $(X,d)$ be a metric space. The following are equivalent:
    	\enum{a)}{
    		\item Any two distinct points of $X$ are joined by unique simple path in $X$, represented by a segment.
    		\item $(X,d)$ is geodesic and contains no subspace homeomorphic to the circle.
    		\enum{i)}{
    			\item $(X,d)$ is geodesic
    			\item If a point is a common endpoint of two points and their only intersection point, then their union is a segment.
    			\item The intersection of two segments with a common endpoint is again a segment.
    		}
    	}
    \end{lemma}
    We refer to any space satisfying the equivalent conditions above as an $\R$-tree.  We will often use the expression $[x,y]$ to denote the unique oriented segment from $x$ to $y$.
    \begin{comment}
    \begin{lemma}
    	Given any metric space $(X,d)$ satisfying $a)$ and $b)$, for any two points $x,y\in X$ there is a unique path from $x$ to $y$ and it
    \end{lemma}
    \end{comment}
    \begin{lemma}
    	Let $(X,d)$ be a tree and $x,y,z\in X$. Then either
    	\enum{i)}{
    		\item   Exactly one of $\{x,y,z\}$ belongs to the segment between the other two.
    		\item   There is a unique point $w$ in the interior of each of the segments $[x,y]$, $[x,z]$, $[y,z]$, so that $[x,z] \cap [y,z] = [w,z]$, $[y,x]\cap[z,x]=[w,x]$, $[z,y]\cap[x,y]=[w,y]$.
    	}
    	In both cases $[x,y]\cap[x,z]\cap[z,y]$ contains a unique point, which we denoted by $Y(x,y,z)$.
    \end{lemma}
    Call a subset of an $\R$-tree
    convex if it contains the segment between any of its points. The following is well-known:
    \begin{lemma}
    	Given a closed convex $C\subset X$ and a point $x\in X\setminus C$ there is a unique point $x'\in C$ such that $d(x,C)=d(x,x')$.
    \end{lemma}
    \begin{proof}
    	Given any segment $[x,y]$, where $y\in C$, there is clearly a point $x_{0}\in[x,y]\cap C$ closest to $x$. All we need to show is that $x_{0}$ is contained in $[x,z]$ for any $z\in C$. So consider any $z\in C\setminus\{y\}$ and let $u=Y(x,y,z)$. Clearly $u\in[y,z]\subset C$, so $c\nin[y,u]$. Therefore $c\in[x,u]\subset[x,z]$.
    \end{proof}
    For the following result, see the discussion in chapter $3$ of \cite{chiswell2001introduction}, starting on p.79.
    \begin{lemma}
    	\label{dichotomy} Let $g$ be an isometry of a $\R$-tree $X$. Then either:
    	\elenco{
    		\item $g$ fixes a point.
    		\item There is a unique line invariant under $g$, on which $g$ acts by proper translations. .
    	}
    \end{lemma}
    In the first case we will say that $g$ acts \emph{elliptically} on $X$. In the second case we say that $g$ acts hyperbolically on $X$ and we refer to the aforementioned line as the \emph{axis} of $g$, or $Ax(g)$ and to through constant distance $d(x,g\cdot x)$ for $x\in Ax(g)$ as the
    translation length of $g$, or $tl(g)$.
    \begin{lemma}
    	Given an isometry $g$ of an $\R$-tree $Y$ and $x\in Y$ we have $\pl{x}{g}{}=\tl{g}{}+2d(x,Ax(g))$.
    \end{lemma}
    We denote by $Fix(g)$ the set of points fixed by $g$. Note that this is a convex set.
    \begin{comment}
    %  Given an isometry $g$ of a real tree $Y$ and $x\in Y$, we let $\pl{x}{g}{}=d(x,g(x))$
    %  by the translation length of $g$ or $\tl{x}{g}{}$ we mean $\inf_{x\in Y}\pl{x}{g}{}$.
    %  The previous lemma implies the infimum above is always attained.
    %  We define:
    \begin{comment}
    \begin{align}
    	Ax(g)=\setof{x\in Y}{d(x,g(x))=tr(g)}
    \end{align}
    \end{comment}
    %  It can be proven that if $g$ does not act elliptically on $Y$, then $Ax(Y)$ is a linear subtree of $Y$ on which $g$ acts as a translation.
    
  \subsection{Group actions on $\R$-trees by isometries }
    %    \df{
    An action $\lambda$ of a group $G$ on a real tree $X$ is called \emph{minimal} if $X$ contains no proper non-empty $G$-invariant sub-trees.
    It is called faithful if no $g\in G\setminus \{1\}$ fixes $X$ point-wise. It is called \emph{trivial} if $G$ fixes some $x\in X$. We say that $X$ is \emph{minimal} if it contains no proper invariant subtree.
    A subgroup $H\geq G$ is called \emph{elliptic} with respect to $\lambda$ if the restriction of $\lambda$ to $H$ is trivial.
    In the absence of ambiguity as to what the intended action is we will tend to refer to the tree itself as minimal. Likewise, we might refer to elements of $G$ itself as hyperbolic or elliptic.
    By the \emph{kernel} of $\lambda$ we intend the subgroup of all those elements fixing the whole space $X$.
    
    An action is relative to a family $\mathcal{A}$ of subgroups of $G$ if and only if every $A\in \mathcal{A}$ is elliptic.

    We denote by $\fun{tl^{\lambda}}{G}{\R}$ the map which to each $g\in G$ associates the translation length of the isometry $\lambda(g)$.
    Likewise, given $x\in X$ we can consider the function $\fun{l_{x}}{G}{\R_{\geq 0}}$ defined by $l_{x}^{\lambda}(g)=d(\lambda(g)\cdot x,x)$.
    %     We refer to the function from $G$ to $\R_{\geq 0}$ thus defined as the length function with respect to $x$.
    In both cases $\lambda$ might be replaced by $X$ or even fall from the notation altogether when the action is unambiguous.
    Recall that in general, given actions of a group $G$ on sets $S$ and $S'$ a map $\fun{S}{S'}$ is equivariant if $f(g\cdot s)=g\cdot s$ for any $s\in S$ and $g\in G$.
    Given isometric actions $\lambda$ and $\lambda'$ of a group $G$ on real trees $X$ and $X'$ we will say that $X$ and $X'$ are said \emph{equivariantly isometric} if some equivariant isometry between the two exists.
    
    \begin{lemma}
    	\label{elliptic products}Suppose we are given a finite set $S$ of generators of a group $G$ which acts by isometries on a $\R$-tree $T$. If $x$ and $x\cdot y$ are elliptic for each $x,y\in S$, then the action is trivial.
    \end{lemma}
    \begin{proof}
    	This follows from the fact that $Fix(gh)=Fix(g)\cap Fix(h)$ for any two elliptic elements $g$ and $h$ together with the well-known fact that the intersection of a finite set of convex subsets of an $\R$-tree with pair-wise non-empty
    	intersection is itself non-empty (Helly's thoerem).
    \end{proof}
    As a consequence, a subgroup $H\leq G$ is elliptic with respect to $\lambda$ if and only if all of its elements are. We now summarize the fundamental classifications of isometric actions $\lambda$ of a group $G$ on a real tree $X$.
    
    To begin with, an action $\lambda$ as above is called \emph{abelian} if the translation length map $\fun{tr^{\lambda}}{G}{\R}$ is a homomorphism.
    \begin{lemma}
    	The action $\lambda$ is abelian if and only if one of the following possibilities occurs:
    	\enum{i)}{
    		\item It is trivial.
    		\item It restricts to an action by translations on some invariant axis.
    		\item There is some infinite ray $\rho$ of $X$ such that $\rho\setminus Fix(g)$ is a bounded segment for any $g\in G$, in which case $G$ is not finitely generated.\footnote{One usually describes this by saying that $G$ is not elliptic in $X$ but fixes and end of $X$.}
    		%    		$G$ fixes some end of $T$ but no point of $T$. In this case $G$ cannot be finitely generated.
    	}
    \end{lemma}
    \begin{corollary}
    	\label{abelians acting} An isometric action of a finitely generated abelian group $A$ on a tree $X$ is either trivial, or restricts to an action by translations on a geodesic line of $X$.
    \end{corollary}
    The following lemma characterizes all remaining actions.
    \begin{proposition}[{\cite{chiswell2001introduction}[Proposition 3.7]}]
    	\label{irreducibility} Let $G$ be a group acting by isometries on an $\R$-tree $X$. The following are equivalent:
    	\enum{i)}{
    		\item There are hyperbolic elements $g,h\in G$ such that $tr(gh)>tr(g)+tr(h)$
    		\item There are hyperbolic elements $g,h\in G$ such that $tl([g,h])\neq 0$.
    		\item There are hyperbolic elements $g,h\in G$ such that $Ax(g)\cap Ax(h)$ is a segment of length strictly less than $tr(g)+tr(h)$.
    		\item $G$ contains a free subgroup of rank $2$ which acts freely, without inversions and properly discontinously on $X$.
    	}
    \end{proposition}
    One refers to any action satisfying the properties above \emph{irreducible}.
    \begin{lemma}[see {\cite[Theorem 4.1]{chiswell2001introduction}}]
    	\label{minimality and translation length} Suppose we are given an isometric action $\lambda$ of a group $G$ on a $\R$-tree $X$ which is either irreducible or restricts to a non-trivial action by translations on a line of $X$. Then the union of the axis of hyperbolic elements of $G$ is a minimal $G$-invariant subtree $X_{min}\subset X$. Given any other action $\lambda'$ of $G$ on an $\R$-tree $X'$, if the functions $tl^{\lambda}$ and $tl^{\lambda'}$ coincide, then the trees $X_{min}$ and $X'_{min}$ are equivariantly isometric.
    \end{lemma}
    In the context of an action of a supergroup $H$ of $G$ on $X$ we will usually denote the minimal tree of $G$ in $X$ as $X_{G}$.
    
    The following is well-known.
    \begin{lemma}
    	Suppose that a group $G$ acts by isometries on an $\R$-tree in such a way that some subgroup $H$ of $G$ of finite index fixes a point $y$ of $Y$. Then the action is trivial.
    \end{lemma}
    \begin{proof}[sketch]
    	\label{ellipticity and finite index} Let $y_{1}$ and $y_{2}$ be two points in the orbit of $y$ by $G$ at maximum distance from each other and $m$ the midpoint of the segment between them. Given $g\in G$ of course $g\cdot y_{j}\in G\cdot y$. From the maximality of $d(y_{1},y_{2})$ one can easily see that $[y_{1},y_{2}]\cap[g\cdot y_{1},g\cdot y_{2}]\neq\nil$ and in fact $g\cdot m=m$.
    \end{proof}
    We say that a real tree on which $G$ acts by isometries is $K$-acylindrical, for $K>0$, if the point-wise stabilizer of any segment of length greater than $K$ is trivial.
    
\section{Simplicial trees: Bass-Serre theory}
  
  Bass-Serre theory provides a correspondence between each from a very general class of actions of a group $G$ on a simplicial trees and certain types of presentations of $G$ in terms of small building blocks (finitely many in case $G$ is finitely generated). It is the natural setting for many essential tools in the study of the first order theory of the free group. We will roughly follow the approach found in \cite[Chapter 2]{dicks1989groups}.
  %  Define the star of a vertex
  By a \emph{graph} $Y$ we intend a structure consisting of two disjoint sets. One $EY=E$ of (oriented)
  edges and another one $VY=V$ of vertices, together with a function $\alpha:E\to V$ and an involution
  $\overline{\cdot}:E\to E$ without fix points. We will use $\omega(e)$ as an abbreviation of $\alpha(\edin{e})$. If $(\alpha(e),\omega(e))=(u,v)$ we will
  say that $e$ originates at $u$ and terminates at $v$, or that it is an edge from $u$ to $v$ and we will refer to $u$ and $v$ as adjacent.
  
  A \emph{path} in $Y$ from a vertex $u$ to a vertex $v$ is a finite sequence $u=u_{0},e_{0},u_{1}\cdots u_{m}=v$ ($m\geq 0$), where
  $\alpha(e_{i})=u_{i}$ and $\omega(e_{i})=u_{i+1}$ for $0\leq i\leq m-1$.
  The path is closed if $u=v$ and simple if $u_{i}=u_{j}$ and $i<j$ implies $i=0,j=m$.
  A graph is called connected when for any two distinct vertices $u,w$ there is a path between them. A \emph{simplicial tree} (or simply a tree) is a connected graph without simple closed paths.
  In particular, given two adjacent vertices $u,v$ of a simplicial tree, there is a unique edge between them, which we will often denote as $(u,v)$.
  For edges $e,f$ in a tree we say that $e\succ f$ if and only if $\alpha(f)\in[\alpha(e),\omega(f)]$.
  
  Given graphs $X$ and $Y$, a cellular map $f:X\to Y$ is a map $VX\cup EY\to VY\cup EY$
  preserving $VX$ and $VY$ and commuting with $\alpha$ and $\edin{\cdot}$. In case $f$ is an inclusion, we refer to $X$ as a subgraph of $Y$.
  %  The sugraph $X$ is full if two vertices of $X$ are joined.
  Given a graph $(V,E)$ and $v\in V$, by the star around $v$ we intend the subgraph consisting of $v$, all edges at distance $1$ from $v$, together with all the edges between any two of two said vertices.
  
  By an action without inversions of a group $G$ on a graph $X=(V,E)$ we intend an action of $G$ on $E\cup V$
  respecting the sets $V$ and $E$, commuting with $\alpha$ and $\edin{\cdot}$ and such that
  $g\cdot e\neq \edin{e}$ for no $g\in G$ and $e\in E$.
  Given such an action, $\alpha$ and $\edin{-}$ push to the origin and inverse function of a graph structure with edge st $G\backslash E$ and vertex group $G\backslash V$, the quotient 	graph, or $G\backslash Y$.
  
  We will refer to a simplicial tree $T$ endowed with an action without inversions of a group $G$ as a $G$-tree. This is called \emph{trivial} if there is some vertex of $T$ fixed by the entire $G$. It is called \emph{reduced} if the stabilizer of a vertex is never contained in that of an adjacent edge. It is called minimal if $T$ does not contain a proper $G$-invariant subtree.
  
  \begin{fact}
  	If $G$ is finitely generated and $T$ minimal, then $G\backslash T$ is finite for any minimal $G$-tree $T$.
  \end{fact}
  
  A Bass-Serre presentation compresses in a certain sense
  all the information contained in a $G$-tree into a smaller object, finite whenever $G$ is finitely generated.
  %Given a group $ G$, By a $G$-graph $X$ we will intend a graph $X$ together with an
  %action of $G$ on it.
  %  	with respect to which the quotient map is a graph map.
  %  	 The quotient map on $E\cup V$ is then a graph map $p_{G}:X\to G\backslash X$.
  %  We will often refer to the image of a subgraph $Z$ by the quotient map $p_{G}$ as $\widecheck{Z}$, and to the preimage by $p_{G}$ of a subgraph $Y$ as $\lft{Y}$.
  
  %  Change this notation HAT is used for copletions
  \begin{definition}
  	A \emph{Bass-Serre presentation} of a $G$-tree $T$ consists of a triple $(Y_{0},Y,(g_{e})_{e\in Y})$, where $Y_{0}\subset Y\subset T$ and $Y_{0}$ is a subtree of $T$
  	%  	$Y=Y_{0}\cup F$ for some set $F$ of edges
  	satisfying the following conditions:
  	%  	whose entries are a sub sets $Y_{0}$ of $T$, a set of edges $Y$, $EY\subset S\subset ET$
  	%  	and a collection $\mathcal{G}=\{g_{e}\}_{e\in EY\setminus EY_{0}}$ of elements of $G$ satisfying the following conditions.
  	\begin{itemize}
  		\item $p_{G}$ restricts to a bijection between $EY$ and $G\backslash ET$ and to a graph isomorphism between $Y$ and a maximal subtree of $G\backslash T$
  		\item $\alpha(e)\in VY_{0}$ for every edge $e\in Y$
  		\item The edge $g_{e}\cdot \edin{e}$ belongs to $Y$ and $g_{f}=g_{e}^{-1}$.
  		\item $g_{e}=1$ for $e\in Y_{0}$
  	\end{itemize}
  	%  	One can extend  the notation by setting $g_{e}=1$, for $e\in EY$. This preserves the validity of the two points above.
  \end{definition}
  The elements $g_{e}$ in the definition above are usually referred to as Bass-Serre elements.
  Given $x\in G\backslash T$ denote by $\lft{x}$ be the unique preimage of $x$ in $Y$.
  
  \begin{definition}
  	A \emph{graph of groups} consists of a connected graph $(E,V)$,
  	%  	\abs{\Gamma}
  	together with the following assignments:
  	\begin{itemize}
  		\item A group $\Gamma_{v}$ for each $v\in V$,
  		\item For all $e\in E$ a subgroup $\Gamma_{e}\leq \Gamma_{\alpha(v)}$
  		and an isomorphism $i_{e}:\Gamma_{e}\cong\Gamma_{\edin{e}}$, such that $i_{\edin{e}}=(i_{e})^{-1}$.
  	\end{itemize}
  \end{definition}
  
  \begin{definition}
  	Let $\Gamma$ be a graph of groups and $Z$ a maximal subtree of $\abs{\Gamma}$.
  	By the fundamental group of $\Gamma$ with respect to $Z$, denoted by $\pi(\Gamma,Z)$ we intend
  	the group defined by the following relative presentation:
  	\begin{equation*}
  		\pi(\Gamma,Z)=\langle\bigcup_{v\in V} \Gamma_{v}\cup\{t_{e}\}_{e\in E}|
  		\{t_{e}=1\}_{e\in Z}\cup\{ i_{e}(g)=g^{t_{e}}
  		,\;t_{e}=t_{\edin{e}}^{-1}\}_{e\in E, g\in \Gamma_{e}}\rangle
  	\end{equation*}
  	\begin{comment}
  	\begin{equation*}
  		\pi(\Gamma,Z)=\langle\bigcup_{v\in V} \Gamma_{v}\cup\{t_{e}\}_{e\in E}|
  		\{t_{e}=1\}_{e\in Z}\cup\{ i_{e}(g)^{-1}g^{t_{e}}=1
  		,\;t_{e}t_{\edin{e}}=1\}_{e\in E, g\in \Gamma_{e}}\rangle
  	\end{equation*}
  	\end{comment}
  \end{definition}
  It can be shown that if $G$ is finitely generated, then $\Delta_{v}$ is generated by finitely many elements together with $\Delta_{e}$ for all $e$ with $\alpha(e)=v$.
  %  \begin{remark}
  	The natural homomorphism from each $G_{v}$ to $\pi(\Gamma,Z)$ can be shown to be
  	injective and we will think of it from now on as an inclusion.
  	%  	For $e\in Z$ one can identfy $G_{e}$ with $G_{\edin{e}}$ and identiy it with $G_{\alpha(e)}=G_{\alpha(e)}\cap G_{\omega(e)}$.
  	One can prove that given two maximal trees $Z$, $Z'$ there is an isomorphism
  	$\theta:\pi(\Gamma,Z)\equiv \pi(\Gamma,Z')$ compatible with the embeddings of
  	$G_{v}$ into $\pi(\Gamma,Z)$ and $\pi(\Gamma, Z')$, up to inner automorphism of one of the fundamental groups.
  %  \end{remark}
  %%  \begin{definition}
  	Let $P=\BS{Y_{0}}{Y}{t}$ a Bass-Serre presentation the
  	action $\cdot$ of a group $G$ on a tree $T=(V,E)$.
  	To this presentation we can associate a graph of groups $\Gamma(P)$ with underlying graph $X=G\backslash T$ and assigments:
  	\begin{itemize}
  		\item  $\Gamma_{v}=Stab_{G}(s(v))$ for $v\in V$
  		\item  $\Gamma_{e}=Stab_{G}(\lft{e})$ for $e\in G\backslash E$.
  		\item	 For $e \in S$ the map $i_{e}$ is the restriction of the inner automorphism $\inn{g_{e}}$ to the subgroup $G_{e}$ (in particular, the identity when $e\in EY$).
  	\end{itemize}
  %%  \end{definition}
  Let $Z_{0}=p_{G}(Y_{0})$. This is a maximal subtree of $G\backslash T$. Consider the map from $\{t_{e}\}_{e\in EX\setminus Z_{0}}\cup(\bunion{v\in\in VX}{\Gamma_{x}}{})$ which sends each
  $t_{e}\in\pi(\Gamma,p_{G}(Y_{0}))$ to $g_{e}\in G$ and restricts to the inclusion on each $\Gamma_{v}$. One can easily check that all relations in the presentation of
  $\pi(\Gamma,p_{G}(Y_{0}))$ hold for their images in $G$, so it extends to a homomorphism $\fun{\phi_{P}}{\pi(\Gamma,Z_{0})}{G}$.
  %  \df{
  By a \emph{graph of groups decomposition} of a group $G$ we intend a triple $(\phi,\Delta,Z)$ where $\phi$ is an isomorphism between $G$ and $\pi_{1}(\Delta,Z)$.
  %  a graph of groups $\Delta$ where each $\Delta_{v}\subset G$ and the inclusion maps of $\Delta_{v}$ into $G$ induced an isomorphism from
  %  In the future we will tend to denote such an object simply by $\Delta$,
  \begin{comment}
  Since in the future we will not be really interested in comparing different of isomorphisms as above with a common domain and target,
  we will refer to the graph of groups decomposition simply as $\Delta$ or $\Delta$.
  \end{comment}
  
  The fundamental theorem of Bass-Serre theory establishes a strong correspondence between graph of groups decompositions of a group and Bass-Serre presentations of actions of the group on a tree.
  \begin{theorem}
  	For any pair $(\Gamma,Z)$ where $\Gamma$ is a graph of groups and $Z$ a maximal subtree of its underlying graph $X$ we can construct an action
  	$\lambda(\Gamma,Z)$ of $\pi(\Gamma,Z)$ on a tree $\lft{(\Gamma,Z)}$ and a Bass-Serre presentations $P(\Gamma,Z)$ of it so that the following properties are satisfied:
  	\begin{itemize}
  		\item  %Modulo identifying vertex groups of $\Gamma$ with their image,
  		The graph of groups associated to the presentation $P(\Gamma,Z)$ in the way described coincides with $\Gamma$.\footnote{In other words, $\pi(\Gamma,Z)\backslash\lft{(\Gamma,Z)}$ is identical to $X$ and $\Gamma_{v}$ coincides with the stabilizer of the lift of $x$ in the maximal subtree of $P(\Gamma,Z)$.}
  		%  		 In particular, the map taking a Bass-Serre element
  		%  		Let $\lft{\Gamma}=\Gamma((\lambda,P)(\Gamma,Z))$. Notice that the partial maps, subgroups and $t$ elements determining the graph of groups structure on both $\Gamma$ and $\lft{\Gamma}$ all live inside $\pi(\Gamma,Z)$.
  		%  		There is a graph isomorphism $\theta$ between $X\Gamma$ and $X\lft{\Gamma}$ seen through which the two graph of groups structures $\Gamma$ and $\lft{\Gamma}$ are 	identical, in the sense that the $G_{v}$, $G_{e}$ and the maps $i_{e}$ are
  		\item Suppose one is given an action $\lambda$ of a group $G$ on a tree $T$ and a Bass-Serre presentation
  		$P=\BS{Y_{0}}{Y}{t}$ of it and let $Z=p_{G}(Y_{0})$. Then the map $\phi_{P}:\pi(\Gamma(\lambda,P),Z)\to G$ is an isomorphism and
  		there is a graph isomorphism $f:T(\Gamma(\lambda,P),Z)\to T$ such that $\phi(h)\cdot v=f(h\cdot v)$ and
  		pair $(f,\phi)$ sends $P(\Gamma,Z)$ to $P$.
  	\end{itemize}
  \end{theorem}
  \begin{comment}
  \item Suppose one is given an action $\lambda$ of a group $G$ on a tree $T$ and a Bass-Serre presentation
  $P=\BS{Y_{0}}{Y}{t}$ of it and let $Z=p_{G}(Y_{0})$. Then the map $\phi_{P}:G\to \pi(\Gamma(\lambda,P),Z)$ compatible with the embeddings of the $G_{v}$ on each side and such that $\phi(g_{e})=t_{p_{G}(e)}$; also a unique $\phi$-equivariant graph isomorphism $f:T\to T(\Gamma(\lambda,P),Z)$ such that $f(S)=S(\Gamma(P),Z)$ and $f(Y)=Y(\Gamma,P,Z)$.
  \end{comment}
  
  %  Given an action $\lambda$ of a group on a simplicial tree $T$, the target of the isomorphism $\phi,\phi'$ between $G$ and the fundamental group of a graph of groups given by two different presentations $P=(Y,S,\{g_{e}\}_{e\in S\setminus EY})$, $P'=(Y',S',\{g'_{e}\}_{e\in S\setminus EY'})$ will be even formally the same as long as the subtrees $Y$ and $Y'$ of the two presentations are translates of each other.
  %  		 but the particular isomorphism might differ in each case.
  %  One can prove that if a group $G$ is finitely generated, then $G\backslash T$ is finite for any $G$-tree. This will be the assumption for the rest of the section.
  %  A $G$-tree (or an associated graph of groups decomposition) is called reduced if the stabilizer of an edge is not contained in the family.
  
  The two most basics examples of fundamental groups of graph of groups are so called amalgamated products and HNN extensions.
  Given two groups $A,B$ and $i$ an embedding of some subgroup $C$ of $A$ into $B$ (usually one tends to think of $C$ as a common subgroup of $A$ and $B$),
  by the \emph{amalgamated product} of $A$ and $B$ over $C$, denoted by $A\frp_{B}C$ we intend
  the quotient of $A\frp B$ by the normal subgroup generated by all the elements of the form $i(c)c^{-1}$ for $c\in C$. This is merely the fundamental group of a graph of
  groups with a single pair of mutually iverse edges and vertex groups $A$ and $B$. The fundamental group of a graph of groups with a single pair of mutually inverse edges and one single vertex, with associated group $A$, is called an \emph{HNN extension} of $A$. It is the quotient of $A\frp\subg{t}$ by the normal subgroup generated by all the elements of the form
  $i(c)^{-1}c^{t}$ for some $C\leq A$ and some embedding $i$ of $C$ into $A$.
  
  The defining property of the free group $\F(x)$ can be restated as saying that the $Cayl(\F,x)$ is a tree. In fact any group which acts freely on a tree is a free and the collection of Bass-Serre elements in any presentation of this action a basis for it. More generally, any $G$-simplicial tree with trivial edge stabilizers is associated to a decomposition of $G$ as a free product
  \begin{align*}
  	(\bfrp{G_{i}}{i=1}{m})\frp F
  \end{align*}
  where each $G_{i}$ is the stabilizer of some vertex $v_{i}$ and $F$ is a group acting freely on $T$.
  We say that a \gr \pr-group $\gqG{G}{A}{P}{I}$ is freely indecomposable if $G$ admits no non-trivial action without inversions on a simplicial tree with trivial edge stabilizers relative to $\{A,P_{i}\}_{i\in I}$.
  
  A metric realization of a connected graph $X=(V,E)$ is the metric space obtained by gluing together
  a segment of a certain positive length for each pair of mutually inverse edges according to the maps $\alpha$ and $\omega$. The distance between two points is given by the shortest total length of a path between those points.
  
  In case of a $G$-tree $T$, we require the assigned length to be the same for edges in the same $G$-orbit. The resulting space $|T|$ is easily seen to be an $\R$-tree on which $G$ acts by isometries.
  
  Each of the notions used to describe the action of $G$ or particular elements of $G$ on $T$ defined in the previous section  have discrete counterparts which hold if and only if the original holds for the geometric realization.
  
  Let us only remark that the fact that the action of $G$ on $T$ is without inversions guarantees that an element of subgroup of $G$ is elliptic in $T$ if and only if it fixes a vertex of $T$.
  In the standard geometric realization every edge is assigned length $1$.
  
      \paragraph{Normal forms}
        
        \label{loop representation}The following appears within the proof of the fundamental theorem of Bass-Serre theory.
        Let $T$ be a $G$-tree, $\BS{Z_{0}}{Z}{t}$ a Bass-Serre presentation for it and $v_{0}$ a vertex in $Z_{0}$.
        Let $g$ be any element of $G$ and consider the path from $u_{0}$ to $g\cdot u_{0}$:
        \begin{align*}
        	u_{0},f_{0},u_{1}\cdots u_{m}=g\cdot u_{0}
        \end{align*}
        then we can write $g$ as
        %  there is a path $p: v_{0},e_{0},v_{1}\cdots v_{m}=v_{0}$ and
        %  for each $0\leq i\leq m-1$ some element $g_{i}\in\ddn{v_{i}}$ such that
        \begin{align*}
        	g=g_{0}t_{e_{0}}g_{1}\cdots t_{e_{m-1}}g_{m}
        \end{align*}
        where $g_{i}\in Stab(v_{i})$ for the unique translate of $u_{i}$ in $v_{i}$ and $e_{i}$ is a translate of $f_{i}$
        (in the language of graph of groups the same can be expressed as saying that $g_{i}\in \Gamma_{v_{i}}$
        for some closed path $u_{0},f_{0},u_{1}\cdots u_{m}$ in the underlying graph).
        This word is reduced in the sense that for $1\leq i\leq m-1$ we have $g_{i}\nin Stab(e_{i})$ whenever $e_{i+1}\neq t_{e_{i}}\cdot\edin{e_{i}}$
        (equivalently $g_{i}\nin\Delta_{\edin{e_{i}}}$ whenever $\edin{e_{i}}=e_{i+1}$  )
        %        	The path $p$ is uniquely determined by $g$. It is precisely the projection to $G\backslash T$ of the unique simple path between $v_{0}$ and $g\cdot v_{0}$ in $T$.
        
      \paragraph{Expansions of groups acting on trees by elliptic elements}
        
        The following lemma should be well-known, but we have failed to find an adequate reference.
        \begin{lemma}
        	\label{elliptic expansion}
        	Suppose we are given groups $G\leq H$, as well as a simplicial $G$-tree $T$ and a subtree $S\leq T$ invariant under the action of $G$,
        	whose translates by elements of $H$ span $T$. Assume moreover that $H$ is generated by $G$ and a family
        	$\mathcal{E}$ of elliptic elements of $G$, invariant under conjugation by elements of $G$. Then the following holds:
        	\enum{i)}{
        		\item $\bigcup_{h\in H}h\cdot S=T$ or, equivalently, the quotient $\fun{q}{G\backslash S}{H\backslash T}$ induced by inclusion is onto. \label{el_one}
        		\item	For any vertex $v\in S$ denote by $K_{v}$ the subgroup generated by $\Sb{v}{G}$ and $\mathcal{E}\cap\Sb{v}{H}$, let $\,\bar{\mathcal{E}}=\abunion{v\in VS}{K_{v}}{}$ and assume that for any $v\in VS$ and $e,f\in ES$ originating at $v$ one has $e\in K_{v}\cdot f$ if and only if $e\in \Sb{v}{G}\cdot f$. Then:
        		\enum{a)}{
        			\item  $\Sb{v}{H}$ is generated by $\Sb{v}{G}\cup(\mathcal{E}\cap \Sb{v}{H})$ For any $v\in VS$ . \label{el-two}
        			\item  $q$ is an isomorphism. \label{el_three}
        		}
        	}
        \end{lemma}
        \begin{proof}
        	Clearly $\bunion{h\in H}{h\cdot S}{}$ is connected. Given a finite product $g_{1}g_{2}\cdots g_{m}$, where each $g_{j}\in\mathcal{E}$, let $\bar{g}_{j}=g_{1}g_{2}\cdots g_{j}$ for $1\leq j\leq m$, clearly 	$\bar{g}_{j}\cdot S\cap\bar{g}_{j+1}\cdot S\neq\nil$ for $1\leq j\leq m-1$. We have shown (\ref{el_one}), since we assume $T$ is minimal.
        	We now assume that two edges adjacent to $v\in VS$ are in the same $\Sb{v}{G}$-orbit if and only if they are in the same $K_{v}$-orbit.
        	\begin{lemma}
        		Suppose that $h\cdot S\cap S\neq\nil$ for some $h\in H$. Then $h=kg$, for some $k\in\bar{\mathcal{E}}$ and $g\in G$.
        	\end{lemma}
        	\begin{proof}
        		Clearly, since $\mathcal{E}$ is invariant under the action of $G$ by conjugation, the subgroup $N$ generated by $\bar{\mathcal{E}}$ is normal in $H$, so one can write $h=ng$ for $n\in N$ and some $g\in G$. 	All we have to show is that $n\in \bar{\mathcal{E}}$. So write it as $a_{1}a_{2}\cdots a_{r}$, where $a_{j}\in\bar{\mathcal{E}}$; we can assume that $r$ is minimal among all such possible representations of $n$.
        		We deal with case $k=2$ first. For the sake of contradiction we assume that $a_{1}a_{2}\nin\bar{\mathcal{E}}G$.
        		Let $F_{j}$ be the fixed point set of $a_{j}$ for $j=1,2$.
        		%        					Then $F_{i}\cap F_{j}=\nil$ and we know that
        		Take $u_{1},u_{2}\in S$ be such that $a_{1}^{-1}\cdot u_{1}=a_{2}\cdot u_{2}:=v$
        		and $L(a_{1},a_{2})$ the minimal value of $d(u_{1},F_{1})+d(u_{2},F_{2})$ for any such $u_{1},u_{2}$.
        		We can assume that $L(a_{1},a_{2})$ is minimal among all $a_{1},a_{2}$ such that $a_{1}a_{2}\in h\cdot G$.
        		Denote by $w$ be the vertex of $S$ closest to $v$.
        		For each $j=1,2$ one of the following two possibilities occurs:
        		\elenco{
        			\item $w\in F_{j}$
        			\item $w\nin F_{j}$ and there are two consecutive edges $e,f\in S$ in the path
        			from $u_{j}$ to $v$ such that $a_{j}^{-1^{j}}\cdot e=\edin{f}$. In particular, $a_{j}$ fixes $\omega(u_{1})=u_{2}$.
        		}
        		The first possibility cannot be the case for both values of $j$, since then $F_{1}\cap F_{2}\neq\nil$, which implies that $a_{1}a_{2}\in\bar{\mathcal{E}}$.
        		By symmetry we can assume the second holds for $j=1$ (otherwise work on $a_{2}^{-1}a_{1}^{-1}$ and then use the $G$-invariance of $\bar{\mathcal{E}}$).
        		Our condition on $K_{v}$ implies the existence of some $g_{0}\in Stab(u)$ sending $\edin{f}$ to $e$, which in turn implies the inequality $d(g_{0}a_{1}^{-1}\cdot u_{1},u_{1})=d(g_{0}\cdot 	v,u_{1})<d(v,u_{1})$.
        		
        		On the other hand
        		$d(a_{2}^{g_{0}^{-1}}\cdot(g_{0}\cdot u_{2}))=d(g_{0}a_{2}^{-1}\cdot u_{2},g_{0}\cdot u_{2})=d(v,u_{2})$.
        		This means that $L(a_{1}g_{0}^{-1},a_{2}^{g_{0}^{-1}})<L(a_{1},a_{2})$, contradicting the minimality assumption, since $(a_{1}g_{0}^{-1})a_{2}^{g^{-1}}=a_{1}a_{2}g^{-1}$.
        		
        		Suppose now that $k>2$ and let $\bar{a}_{j}=a_{1}\cdots a_{j}$ (so $\bar{a}_{0}=1$), $S_{j}:=\bar{a}_{j}\cdot S$ and $d_{l}=d(S_{l},S)$ for $0\leq j\leq m$.
        		Pick any $1<j_{0}<m$ maximizing $d_{j_{0}}$ and $D=d_{j_{0}}$. In view of the previous case we can assume that $D>0$.
        		We claim that $S_{j_{0}-1}\cap S_{j_{0}+1}\neq\nil$. Indeed, both intersect $S_{j_{0}}$ and therefore they must contain the unique ($D>0$) point of $S_{j_{0}}$
        		at minimal distance from $S$.
        	\end{proof}
        	Using this the following additional claim, from which properties (\ref{el_one}) and (\ref{el_three}) are easy to prove.
        	\begin{lemma}
        		Let $h\in H$ and suppose we are given vertices $u,v\in S$ such that $h\cdot u=v$. Then $h=kg$ for some $k\in K_{v}$ and $g\in G$ such that $g\cdot u=v$.
        	\end{lemma}
        	\begin{proof}
        		Given any such $h$ consider a pair $(k_{0},g_{0})\in \bar{\mathcal{E}}\times G$ such that $k_{0}g_{0}=h$ and $d=d(g_{0}\cdot u,v)$ is minimal for this property.
        		If $d=0$ we are done, by the previous case, so assume this is not the case. Let $v=g_{0}\cdot u=n_{0}\cdot u$. Just as in the proof of the previous lemma, the segment from $g\cdot u$ to $v$ must contain two consecutive edges $e,f$ such that $\omega(e)=\alpha(e):=w$ and $n_{0}\cdot e=\edin{f}$ and some $g_{1}\in \Sb{G}{w}$ such that $g_{1}\cdot e=\edin{f}$ as well. It is clear then that $d(g_{0}g\cdot u,v)<d(g\cdot u,v)$. Since $k$ fixes $w$ as well, $kg_{0}^{-1}\in\bar{\mathcal{E}}$ which contradicts the minimality of $d$, since of course $g_{0}g\in G$.
        	\end{proof}
        \end{proof}
        
\section{Operations on trees}
  
  \subsection{Collapses and blow-ups}
    \begin{comment}
    \df{
    	Given graphs $X,Y$, by a \emph{cellular map} from $X$ to $Y$ we intend   % use cellular map instead
    	a map $f: X\to Y$ such that $f(VX)\subset VY$,  and for all $e\in EY$ one has
    	$f(e)\in VX \slra f(\edin{e})\in VX$ and $f(\alpha(e))$ is either:
    	\elenco{
    		\item $\alpha(f(e)) $ in case $f(e)\in EX$
    		\item $f(e)$ in case $f(e)\in VX$
    	}
    }
    \end{comment}
    %    \df{
    Given a graph $X$ and $F=\edin{F}\subset EX$
    %    	one can define a graph $C_{F}(X)$ and
    let $p_{F}: X \to C_{F}(X)$ be the cellular map obtained by collapsing the edges in $F$.
    More precisely, let $X_{F}=F\cup\alpha(F) $ and $\mathcal{C}$ the
    collection of its connected components.
    We let
    \begin{align*}
    	V(C_{F}(X))=&V(X)/\{u\sim v|u,v \in V(C),\,C\in \mathcal{C}\} \\
    	E(C_{F}(X))=&E(X) \setminus \bunion{C\in\mathcal{C}}{EC}{}
    \end{align*}
    where the new incidence function is obtained by post-composing the previous one with the quotient map and the inverse operation the adequate restriction of the old one.
    The cellular map $p_{F}$ is the obvious quotient mapping each component $C$ to a point.
    
    For simplicity, in this situation above let $Y^{F}=C^{F}(Y)$ for $F\subset EX$ and $Y$ a subgraph of $X$. Assume now that $X$ is a $G$-tree and $F$ a $G$-invariant family.
    \newcommand{\gqt}[1]{(G\backslash #1)}
    Clearly $X^{F}$ inherits a $G$-graph structure from $X$ and there is an isomorphism $\gqt{X}^{\gqt{F}}\equiv\gqt{X^{F}}$ compatible with $p^{F}$ and $p_{\gqt{F}}$.
    Suppose we are given a $G$-tree $T$, as well as some $G$-invariant $F\subset ET$ take any Bass-Serre presentation $\BS{Z_{0}}{Z}{t}$ of it
    such that whenever $e\in\gqt{F}$ for some edge $e\in Z\setminus Z_{0}$ we have $f\in\gqt{F}$ for each $f$ in the unique path in $Z$ between the endpoints of $e$.
    Then $\gqt{Y^{F}}$ is a maximal subtree of $\gqt{X^{F}}$, and $P'=(Z_{0}^{F},Z^{F},\{g_{e}\}_{e\in EZ^{F}})$ a Bass-Serre presentation of the $G$-tree $T^{F}$.
    It is easy to see that at the level of graph sof groups we can describe the operation as replacing each subgraph in the image of a connected component of $F$
    by the group generated by the corresponding generators.
    \begin{comment}
    Each connected component $\tah{C}$ of $\tah{F}$ is in bijective correspondence to an equivariant disjoint family of connected components of $F$. In virtue of the assumption, there is one, $\lft{C}$, intersecting the transversal $Y$.
    It is possible to show that the set-wise stabilizer of $C$ and, hence, the point-wise stabilizer $Stab(x_{C})$ of the corresponding point $x_{C}$ in the collapsed tree is precisely the subgroup of $G$
    
    %   (under the identification with $\pi_{1}(\Delta,Y)$ given by the Bass-Serre presentation)
    generated by the stabilizers of the vertices of $Y\cap C$ and the Bass-Serre elements $g_{e}$ for $e\in EY\cap E$.
    This implies the vertex group associated to $x_{C}$ in the graph of groups decomposition $\Delta'$ induced by $P'$, is isomorphic to $\pi_{1}(\tah{C},\tah{Y})$.
    \end{comment}
    
    The inverse operation is what one calls a \emph{blow-up}, or refinement.
    The pertinent data in this case is some $v_{0}\in VT$ and some $Stab(v_{0})$-tree $S$ in which $Stab(e)$ is elliptic for each edge $e$ incident to $v_{0}$.
    Using this one can construct $G$-tree	in the following way. First, for each $[g]\in G/Stab(v_{0})$ take a copy $S^{[g]}$ of $S$. We can assume that $S^{1}=S$.
    Pick any set of representatives $\mathcal{R}$ of $G/Stab(v_{0})$ and for each $r\in\mathcal{R}$ an isomorphism $\fun{\theta^{r}}{S^{[1]}}{S^{[r]}}$.
    There is a unique way of extending this assignment to a collection of graph isomorphisms
    $\theta^{g}_{[h]}$ between $S^{[h]}$ and $S^{[gh]}$ for any $[g]\in G/Stab(v_{0})$ and any $g\in G$ in such a way
    that $r^{g}_{[1]}=\lambda(g)$ for any $g\in Stab(v_{0})$ and $r^{gh}_{[1]}=r^{g}_{[h]}\circ r^{h}_{[1]}$.
    %    Up to here the picture does not depend on the choice of $\mathcal{R}$, only on the action on $S$.
    Take any system of representatives $\mathcal{ER}$ of orbits by $G$ of edges of $T$ in the family $\mathcal{E}_{v_{0}}$ of all the edges originating at a
    at a translate of $v_{0}$. Any assignment of a vertex $p_{e}\in S^{\alpha(e)}$ to each $e\in\mathcal{ER}$ can be extended uniquely to the whole
    $\mathcal{E}_{v_{0}}$ in a way that the equality $p_{h\cdot e}=\theta^{h}_{[g]}$ holds for any $e\in\mathcal{E}$ with $\alpha(e)=g\cdot v$ and $h\in G$.
    The tree $T_{S,(p_{e})_{e\in\mathcal{E}}}$ (often simply $T_{S}$, when context allows) is defined by:
    \begin{align*}
    	VT_{S,(p_{e})_{e\in\mathcal{E}}}=(VT\setminus (G\cdot\{v_{0}\}))\cup\bunion{g\in G}{VS^{[g]}}{}
    	\\ ET_{S,(p_{e})_{e\in\mathcal{E}}}=\setof{\bar{e}}{e\in ET}\cup(\bunion{g\in G}{ES^{[g]}}{})
    \end{align*}
    where for each $e\in ET$ we have $\edin{\bar{e}}=\bar{\edin{e}}$ and the copy $\bar{e}$ of $e$ originates either at $v=\alpha(e)$ in case $\alpha(e)\neq G\cdot v_{0}$ or else
    at $p_{e}$. Note that although there is no canonical bijection between the vertex and edge sets of $S$ and $S^{[h]}$, the maps $\theta$ do induce a canonical isomorphism between  $Stab(v_{0})\backslash S$ and $Stab(v_{0})\backslash S^{[h]}$.
    % replace $[g]$ by $g\cdot v$
    
  \subsection{Lifting decompositions through blow-up and collapse}
    
    Clearly, if we collapse all edges inherited from $S$ in the blow-up $T_{S}$ described in the previous subsections we recover the original tree $S$. One can also collapse those edges inherited from $T$ instead of from $S$. This operation has very nice properties under certain circumstances, a fact that we will be useful later.
    \begin{lemma}
    	\label{blow and collapse} Suppose we are given a non-trivial simplicial $H$-tree $T$, $G\leq H$, a $G$-invariant subtree $U$ of $T$ and for some $v_{0}\in V$ a simplicial $\Sb{v_{0}}{H}$-tree $S$.
    	Assume furthermore that:
    	\enum{a)}{
    		\item  Two edges of $U$ are in the same $H$-orbit only if they are already in the same $G$-orbit\label{hip1}
    		%    		The natural graph map form $G\backslash T_{G}$ to $H\backslash T$ is injective on edges
    		%    		Translates of the minimal tree $T_{G}$ cover $T$ and
    		\item $\Sb{e}{H}$ is elliptic in $S$ for any edge $e$ adjacent to $v_{0}$. \label{hip2}
    		\item For any $h\in H$ such that $h\cdot v_{0}\in T_{G}$ (equivalently, for all $h$ in some system of representatives of $G\backslash H$) \label{hip3}
    		if we let $\mathcal{E}_{v}$ be the set of edges of $T_{G}$ originating at $v$ then the subgroup generated by
    		\begin{align*}
    			\Sb{h\cdot v_{0}}{G}^{h}\cup\{\Sb{e}{H}^{h}\}_{e\in\mathcal{E}_{h\cdot v_{0}}}
    		\end{align*}
    		%    			, where  (equivalently a system of representatives of their orbits by the action of $Sb{v_{0}}{G}$)
    		fixes a vertex in some $Stab(v_{0})$-invariant family $\mathcal{O}\subset VS$.
    	}
    	Then a simplicial $H$-tree $S'$ exists together with an equivariant embedding $\psi$ of $S$ into $T$ with the following properties:
    	%    	Then there is a way way of blowing up each point in the orbit of $v$ to a copy of $S$ such that that if we let $T_{S}$ be the resulting simplicial $H$-tree
    	%    	and $T^{c}_{S}$ the one obtained collapsing the images in $T_{S}$ of the edges of $T$, the group $G$ is elliptic in $T^{c}_{S}$.
    	\enum{i)}{
    		\item  \label{prop1}The union of translates of $\psi(S)$ cover $S'$.
    		\item  For each $v\in VT\setminus H\cdot\{v\}$ the group $Stab(v)$ fixes some vertex in $H\cdot\psi(\mathcal{O})$.
    		\item  \label{prop3} $Stab(\psi(e))=Stab(e)$ for any $e\in ES$
    		\item  \label{prop5} $Stab(v)=Stab(\psi(v))$ for any $v\in VS$
    		\item  \label{prop2} $G$ fixes some vertex $v_{G}$ in $T\setminus H\cdot\psi(\mathcal{O})$. Let $H_{0}=Sb(v_{G})$
    		\item  \label{prop4} If  $e,e'\in ET$ are in distinct $H_{0}$-orbits, then $\psi(e),\psi(e')$ are in distinct $H$-orbits.
    		\item \label{prop6} Any $v,v'\in VT\setminus\mathcal{O}$ in distinct $H_{0}$-orbits are also in distinct $H$-orbits.
    	}
    	Furthermore, for some $v_{G}\in\psi(\mathcal{O})$ fixed by $G$ if we let $H_{0}$ be the stabilizer of $v_{G}$ and $U_{0}$ be the convex closure in $T$ of the union of translates of $U$ by elements of $H_{0}$. The intersection with $U_{0}$ of either an $H$-orbit of edges or an $H$-orbit of vertices not containing $v_{0}$ is either empty or an $H_{0}$-orbit.
    	%    	, in which case the stabilizer in $H$ of any of its members is also contained in $G$. The group $H_{0}$ contains the stabilizer in $G$
    \end{lemma}
    \begin{proof}
    	We first construct an $H$-tree $T_{S}:=T_{S,(p_{e})_{e}}$ by equivariantly blowing up all vertices in $H\cdot v_{0}$ to a copy of $S$. Condition (\ref{hip2}) guarantees this can be done (the result is generally non-unique). Afterwards, we collapse all edges of $T_{S}$ inherited from $T$. As $\psi$ is simply the isomorphism between $S$ and $S^{[1]}\subset T_{S}$.
    	The resulting $H$-tree $T^{c}_{S}$ together with the obvious map from $S$ into $T^{c}_{S}$ clearly satisfies properties (\ref{prop1}),(\ref{prop2}), (\ref{prop3}), (\ref{prop5}) and the last part of (\ref{prop6}) by construction.
    	%    	, as it depends on the choice of a system of attaching points on the new copies of $S$ for the family of edges originating at vertices of $H\cdot v_{0}$.
    	%    	set of representatives of the family of all edges of $U$ originating at $r\cdot v_{0}$.
    	Let $\mathcal{V}^{G}$ be a set of representatives of $G\backslash (H\cdot v_{0})$. For any $u\in\mathcal{V}^{G}$ chose
    	a set of representatives $\mathcal{ER}^{G}_{r}$ of the orbits by $Stab(v_{0})$ of the family of edges of $T$ originating at $u$.
    	Our first hypothesis implies that for any $u\in\mathcal{V}^{G}$ and any
    	distinct $u,u'\in\mathcal{V}^{G}$ and $h\in H$ the sets $h\cdot\mathcal{ER}^{G}_{u}$ and $\mathcal{ER}^{G}_{u'}$ are disjoint.
    	
    	In particular, we can assume that $\bunion{u\in\mathcal{U}^{G}}{\mathcal{ER}^{G}_{u}}{}$ is a subset of the transversal $\mathcal{ER}^{H}_{v_{0}}$ of the set of edges originating in $H\cdot v_{0}$
    	used in the blow-up construction.
    	In virtue of the third hypothesis, for any $r\in\mathcal{R}$ we can arrange so that all $e\in\mathcal{ER}^{G}_{u}$ have the same attaching point $p_{r}=p_{e}\in\mathcal{O}_{u}$,
    	where $\mathcal{O}_{u}$ is the image of $\mathcal{O}$ in the copy of $S$ at $u$. This $p_{r}$ is fixed by $\Sb{r\cdot v_{0}}{G}$.
    	
    	Clearly, every pair of edges $(e,f)$ of $U$ originating at a common translate of $v_{0}$ is an $H$-translate of a pair of the form
    	$(e_{0},k\cdot f_{0})$, where $k\in\Sb{u}{G}$, $e,f\in\mathcal{ER}^{G}_{u}$ and $u\in\mathcal{V}^{G}$, so $p_{e}=p_{f}$ in that case as well.
    	This implies that subtree spanned by the images of the edges of $U$ in $T_{S,(p_{e})_{e}}$ is connected, and hence maps to a unique point $v_{G}\in T^{c}_{S}$,
    	fixed by $G$. The preimage $C$ of $v_{G}$ in $T_{S}$ has to contain an attaching point (since it does not fill the whole $T_{S}$), so $v_{G}\in H\cdot\psi(\mathcal{O})$. The same argument
    	applies to $Stab(v)$ for $v\in VT\setminus H\cdot v_{0}$.
    	On the other hand $C$ is $H_{0}$-equivariantly isomorphic to soem subtree of $T$ containing $U$, since the collapse map to $T$ does not touch any of its edges.
    	%  This implies that $U$ collapses to a point in $T^{c}_{S}$, clearly fixed by $G$.
    	%    	(a set which is invariant by the action of $H$).
    	% It is easy to see that the minimal tree $T_{H_{0}}$ is $H_{0}$-equivariantly isomorphic to $C$ using the fact they are both locally isometric.
    	If two vertices $v,v'\in C$ inherited from $T$ or two edges $e,e'\in C$ differ by translation by some $h\in H$ then $h\cdot C\cap C\neq\nil$, which implies $C=h\cdot C$,
    	since $C$ is a connected component of the complement of the union of the collection of all copies of $S$ in $T_{S}$ we conclude that $h\in Stab(v_{G})$, which proves the last claim.
    	%    	Since $U_{0}$ is isomorphic to an $H_{0}$ invariant subtree of $C$.
    \end{proof}

  \subsection{Foldings and slides}
    
    Given a graph $X=(V,E)$, and $e,f\in E$ such that $\alpha(a)=\alpha(f)$ by the \emph{fold} of $e$ and $f$ we intend the quotient map $\fun{F_{e,f}}{X}{X'}$, where $X'=(V',E')$, $V'$ obtained from $V$ by identifying $\omega(e)$ and $\omega(f)$ and $E'$ is obtained from $E$ by identifying $e$ and $f$. Given a family $\mathcal{P}$ of such pairs we can consider the (inverse) limit graph
    obtained by quotiening each of those pairs one by one (in successive steps we use the image of the initial pairs in the current quotient tree). It can be shown this does not depend on the order in which the different foldings are performed. There is a natural cellular map from $X$ onto the resulting graph preserving the set of edges.
    
    If $X$ is a tree, the result of folding $\mathcal{P}$ is a tree as well. If we are dealing with a simplicial $G$-tree the operation is to be interpreted equivariantly. Namely,
    we not only fold $(e,f)$ but also any image of the pair by the action of $G$. The result is again a $\mathcal{G}$ tree. See \cite{stallings1991foldings} or \cite{dunwoody1998folding} for more details.
    The effect of an equivariant folding at the level of graph of groups has a varied typology. We will restrict ourselves to the case relevant to us.
    %    It can be easily checked that letting $\alpha'([e])=[\alpha(e)]$ is a well-defined graph structure on $X'$ with respect to which the map $F_{e,f}$ is a cellular map which does not collapse  edges.
    \begin{lemma}
    	Let $T$ be a simplicial $G$-tree and suppose we are given $v\in VT$ and $e_{1},\cdots e_{m}$ origiating at $v$ and in distinct orbits by the
    	action of $G$ and suppose that for $1\leq j\leq m$ we are given $Stab(e_{j})< H_{j}\geq Stab(\omega(e_{j}))$. Assume that for any $1\leq j\leq m$ the vertex $\omega(e_{j})$ is not in the orbit of $v$.
    	Let $\fun{f_{\bar{H}}}{T}{T^{f}_{\bar{H}}}$ the equivariant fold of all the pairs
    	$(\edin{e_{j}},h\cdot\edin{e_{j}})$, where $h\in H_{j}$. Then $Stab(f_{\bar{H}}(v))$ is equal to the amalgamated product of $Stab(v)$ with $H_{j}$ over $Stab(e_{j})$ and the stabilizer
    	of $f_{\bar{H}}(e_{j})$ is equal to $H_{j}$.
    \end{lemma}
    \begin{proof}
    	That the stabilizer of the image of $e_{i}$ is equal to $H_{i}$ is immediate from the fact that no fold is performed at $\epsilon(e_{i})$.
    	
    	The stabilizer $H$ of $f_{\bar{H}}(v)$ coincides with the set-wise stabilizer of its preimage $\mathcal{P}$ in $T$, which is precisely the set of all vertices $v'$  in the orbit of $v$ that can be joined to $v$ by some path of the form
    	$v_{0}=v,f_{0},u_{0},f'_{0},v_{1},e_{1}\cdots v_{m}=v'$, where for each $0\leq j\leq m-1$ edges $f_{j}$ and $\edin{f'_{j}}$ are folded together (we are using the assumption that none of the ends of an edge which is folded is in the orbit of $v$). The union of all those paths is a subtree $S$ whose intersection with $G\cdot v$ is precisely $\mathcal{P}$, so that the set-wise stabilizer of $S$ is equal to $H$. We claim that the action of $H$ on $S$ is dual to the graph of groups decomposition described above. It is easy to check that $\mathcal{P}$ is a single orbit by the action of $H$. This implies that the orbits by the action of $H$ of edges in $S$ are in bijective correspondence with $e_{1},\cdots e_{m}$.
    	Since the star around $\omega(e_{j})$ in $S$ contains only edges in the orbit of $e_{j}$ for $1\leq i<j\leq m$ vertices $\omega(e_{j})$ and $\omega(e_{i})$ must be in distinct orbits by the action of $H$. This implies the action of $H$ on $G$ is of the required form.
    \end{proof}
    Note that the situation above $f_{\bar{H}}$ induces an isomorphism from $G\backslash T$ to $G\backslash T^{f}_{\bar{H}}$.
    \begin{comment}
    \begin{lemma}
    	\label{orbit folding}
    	Let $G$ be a group and $T$ a $G$-tree. Let $v\in VT$ and $\mathcal{W}$ a set of representatives of orbits for the action of $Stab(v)$ on the set of edges originating at $v$. Suppose that for each $w\in\mathcal{W}$ we are given some subgroup $K_{w}$ such that $Stab(e)\subset K_{w}\subset Stab(w)$. Let $\fun{q}{T}{T'}$ the quotient map obtained from folding any pair of edges $(e,f)$ in the orbit by $G$ of a pair of the form $(e,h\cdot e)$ where $\alpha(e)=v$, $\omega(e)\in\mathcal{W}$ and $h\in K_{w}$. Then, if we let $v'=q(v)$, $Stab(v')$ is the subgroup $H$ of $G$ generated by $Stab(v)$ and $Stab(w)$ for $w\in\mathcal{W}$. Moreover, $v\cup\{\mathcal{W}\}$ spans a fundamental domain of the minimal tree $T_{v}$ of $Stab(v')$ in $T$. In particular, the stabilizer of an edge of $T_{v}$ of the form $h\cdot e$ for an edge $e$ with $\alpha(e)=v$ and $\omega(e)\in\mathcal{W}$ and $h\in H$ is precisely $Stab(e)^{h^{-1}}$.
    \end{lemma}
    \end{comment}
    
    Another equivariant move that can be performed on a tree is a \emph{slide}. It moves the origin of of some edge $e$ with $\alpha(e)=v$  with $v'$ (as well as that of all other edges in the orbit of $e$, in an equivariant way) on the condition
    that there is some edge $f=(v,v')$ such that $Stab(e)\leq Stab(f)$.

\section{Trees and surfaces}
  \label{surfaces section}  Let $\Sigma$ be a non-simply connected compact surface with boundary. We quickly recall some notions on the topology of compact surfaces. For more details the reader is referred to \cite{farb2012primer}.
  \begin{comment}
  We start by recalling the very well known classification theorem:
  It is a well-known fact that for any compact surface with boundary $\Sigma$ is homeomorphic to either $S^{2}_{b}$ or $S^{2}_{b}\#\Sigma'$, where $\Sigma'$ is either isomorphic to either:
  \elenco{
  	\item  The connected sum of $m$ tori for $m\geq 1$ in case $\Sigma$ is orientable.
  	\item  The connected sum of $m$ projective planes in case $\Sigma$ is non-orientable.
  }
  where $S^{2}_{b}$ for sphere minus $b$ disjoint open disks.
  In the first case the fundamental group admits a presentation of the following type:
  \begin{align*}
  	\subg{x_{1},y_{1}\cdots x_{m},y_{m},s_{1},\cdots s_{k}\,|\,[x_{1},y_{1}][x_{2},y_{2}]\cdots[x_{m},y_{m}]=s_{1}\cdots s_{k}}
  \end{align*}
  While in the second it is:
  \begin{align*}
  	\subg{x_{1}\cdots x_{m},s_{1},\cdots s_{k}\,|\,x_{1}^{2}x_{2}^{2}\cdots x_{m}^{2}=s_{1}\cdots s_{k}}
  \end{align*}
  In both cases the $s_{i}$ are classes of boundary loops of $\Sigma$. The Euler characteristic $\xi(\Sigma)$ is equal to $2-2m-b$ in the first case and $2-m-b$ in the second.
  \end{comment}
  %  	A surface $\Pi$ of finite homology dimension can be obtained by removing some boundary components from a uniquely determined compact surface with boundary, which we by denote $\Pi^{*}$.
  We will refer to those compact surfaces with (possibly empty) boundary, $\Sigma$, for which either $\xi(\Sigma)\leq-2$ or $\Sigma$ is a torus minus an open disk as \emph{big}. We say that a surface is \emph{closed} if contains no boundary.
  By a simple closed curve (abbreviated as s.c.c.) in a surface $\Sigma$ we intend an homeomorphic embedding from a
  circle into $\Sigma$. We will usually be interested in the class $\alpha=[a]$ of a curve $a$ modulo homotopy and reparametrization (so not an oriented curve), rather than $a$ itself.
  
  There are two mutually excluding possibilities. Either $\alpha$ has a neighbourhood homeomorphic to an annulus or one homeomorphic to a M\"obius band; in the first case we refer to the curve as two-sided.
  % Possibly talk about separating curves.
  A s.c.c. will be called \emph{essential} if it is two-sided and is not homotopic to a boundary component of $\Sigma$. A subsurface $\Sigma'\subset\Sigma$ will be called essential if none of its boundary components are nullhomotopic in $\Sigma$. By an \cfam of $\Sigma$, $\delta$, we intend a collection of distinct homotopy classes of essential s.c.c. in $\Sigma$. Sometimes it is convenient to think of it as the union of disjoint representatives of such classes.
  %  We might sometimes regard $\delta$ as the finite set of all its components. Denote $\pi_{1}(\Sigma)$ by $S$.
  %   though we might tend to blur the distinction between the two.
  %%Recall that there is a bijection between homotopy classes of curves in $\Sigma$, $[\alpha]$ and conjugacy classes $C(\alpha)$ of $\pi_{1}(\Sigma)$.
  %  Properly speaking $S$ has to refer to the fundamental group with respect to some particular basepoint $*$.
  
  Recall that for any path connected subspace $X\subset\Sigma'$ of $\Sigma$, $*'\in X$, any homotopy class of continuous paths $[p]$ from $*$ to $*'$
  defines a homomorphism $\iota_[p]$ from $\pi_{1}(X,*')$ into $\pi_{1}(\Sigma,*)$. The effect of chosing a different $p$ is that of postcomposing $\iota_{[p]}$ by some
  inner automorphism of $\pi_{1}(\Sigma,*)S$ and unless $*=*'$, there is no canonical one, but in most of the applications we are interested in we will be usually able to ignore this issue.
  The homomorphism $\iota_{*}$ is injective in case of an incompressible subsurface (one none of whose boundary connected compoenents is null-homotopic in $S$).
  We refer to the cyclic group in one of the conjugacy classes associated with the boundary of $\Sigma$ as a boundary subgroup and to any of its members as a boundary element.
  The collection $Homeo_{\partial}(\Sigma)$ of all homeomorphisms of $\Sigma$ fixing its boundary components, quotiened by the equivalence relation any two such homeomorphisms
  which are related by an isotopy preserving the boundary components set-wise. We will refer to the resulting group as the modular group of $\Sigma$, or $Mod(\Sigma)$.\footnote{In the literature this is usually defined as the modular group of the corresponding punctured surface (in which a point instead of a closed disk is removed from the closed surface).
  Likewise, when referred to an orientable surface, in absence of punctures the modular group stands for the set of classes of orientation preserving homeomorphism exclusively. If we include non-orientable homeomorphisms as well (in the empty boundary is empty) we are left with what is usually called the extended mapping class group.}
  It follows from by the Dehn-Nielsen Baer theorem that the group, which we denote by $Mod(\pi_{1}(\Sigma))$ of all automorphisms of $\pi_{1}(\Sigma)$ induced by homeomorphisms of the surface relative to the boundary has finite index in the group of all the automorphisms of $\pi_{1}(\Sigma)$ which restrict to inner automorphisms on peripheral subgroups. In truth, the classic result (see \cite[Theorem 8.1]{farb2012primer}) is only valid for orientable surfaces; an appropriate generalization to non-orientable surfaces can be found in \cite[p.278]{fujiwara2002}.
  
  %  Given a class $\eta\in Mod(\Sigma)$ and an essential simple closed curve $\alpha$ in $\Sigma$
  %  The action of $Homeo_{\partial}(\Sigma)$ on the
  The action of $Homeo_{\partial}(\Sigma)$ on $\Sigma$ descends to a well-defined action of $Mod(\Sigma)$ on the set of homotopy classes of essential simple closed curves in $\Sigma$ and on the set
  of \cfams on $\Sigma$. It is possible to show the latter has finitely many orbits. Given an \cfam $\delta$, we denote by $[\Sigma]_{\delta}$ the collection of pieces in which $\delta$ cuts $\Sigma$,
  %  minimal subsurfaces of $\Sigma$ bounded by curves in $\delta$,
  i.e., the (closures of) the connected components of $\Sigma\setminus\delta$. Of course, these are only well defined up to isotopy of the surface.
  The system $\delta$ dertermines an action $\lambda_{\delta}$ without inversions of $S$ on a simplicial tree $T_{\delta}$
  whose associated graph of groups $\Delta$ has the following properties:
  \enum{i)}{
  	\item There is a bijective correspondence by which each $\alpha\in\delta$ corresponds to some orbit $S\cdot e$ of edges of $T_{\delta}$, so that the family $\{Stab(s\cdot e)\}_{s\in S}$ coincides with that of cyclic groups associated to $\alpha$.
  	\item There is a bijective correspondence associating every $\Sigma_{0}\in[\Sigma]_{\delta}$ with an orbit $S\cdot v$ of vertices of $T_{\delta}$, so that the family $\{Stab(s\cdot v)\}_{s\in S}$
  	coincides with all the mutually conjugate images of $\pi_{1}(\Sigma_{0})$ in $S$ by the embeddings associated to the inclusion map.
  	%  		\item The cyclic subgroups associated with boundary components of $\Sigma$ are elliptic in $T_{\delta}$.\label{boundaries are elliptic}
  }
  %  Accordingly, the corresponding graph of groups has a vertex with associated group $\pi_{1}(\Sigma_{0})$ (one of its conjugate images in $\pi_{1}(\Sigma)$) for each subsurface $\Sigma_{0}$ of $[\Sigma]_{\delta}$ and an edge $e$ between $\Sigma_{0}$ and $\Sigma_{1}$ for each $\alpha\in\delta$ contained in $\Sigma_{0}\cap\Sigma_{1}$, the edge group being one of the mutually conjugate groups corresponding to the curve $\alpha$.
  %also referred to as the tree dual to $\delta$.
  Alternatively, consider the preimage $\mathcal{D}$ of $\delta$ by the universal covering map $p:\widetilde{\Sigma}\to\Sigma$. Necessarily $\widetilde{\Sigma}$ is homeomorphic to an infinite plane.
  %  Here $\widetilde{\Sigma}$ is homeomorphic to an open disk $D$.
  The family $\mathcal{D}$ is invariant under the action of $S$ and consists of infinitely many disjoint bi-infinite curves. Each one of them, say $\tilde{\alpha}$, projecting to $\alpha$, separates $D_{2}$ into two connected components and is stabilized by some $Z\in C_{\alpha}^{\Sigma}$. The tree $T_{\delta}$ can be obtained by taking a vertex $V_{U}$ for each connected component $U$ of $D\setminus(\cup\mathcal{D})$ and a pair of edges $e_{\widetilde{\alpha}},\edin{e_{\widetilde{\alpha}}}$ for each $\tilde{\alpha}\in\mathcal{D}$ joining the two connected components meeting at $\tilde{\alpha}$.
  
  Clearly $S$ acts isomorphically on $T$, with edge stabilizers in $C_{\alpha}^{\Sigma}$. That $T$ is a tree comes from the fact that the $\tilde{\alpha}$ separate $D$. The hypothesis that $\alpha$ is one-sided guarantees that the action is without inversions.
  
  Alternatively, the dual graph of groups decomposition can be obtained using Van Kampen's theorem, which provides the required presentation for $\pi_{1}(\Sigma)$.
  
  Given two homotopy classes of essential s.c.c., $\alpha,\beta$, the intersection number of $\alpha$ and $\beta$, denoted by $i(\alpha,\beta)$ is the minimum $n\in\N$ such that there are representatives $a\in\alpha$ and $b\in\beta$ which intersect in exactly $n$ distinct points.
  \begin{fact}
  	%  	Given $\alpha,\beta$ homotopy classes of essential s.c.c., $i(\alpha,\beta)\neq 0$ if and only if the elements of $\pi_{1}(\Sigma)$ associated to $\alpha$ are hyperbolic in the tree
  	%  	$T_{\beta}$.
  	The intersection number $i(\alpha,\beta)$ coincides with the translation length of a generator of a subgroup associated to $\alpha$ in $T_{\beta}$.
  \end{fact}
  We will refer to those compact surfaces with (possibly empty) boundary, $\Sigma$, for which either $\xi(\Sigma)\leq-2$ or $\Sigma$ is a torus minus an open disk as big.
  %  These are precisely those admitting two essential simple closed curves $\alpha$ and $\beta$ with $i(\alpha,\beta)\neq0$.
  
  Part of the important role played by surface groups in the theory of simplicial actions of trees has to do with the following property (see \cite[Theorem III.2.6]{morgan1984valuations} or the proof of \cite[Lemma 9.4]{guirardel2010jsj}):
  \begin{lemma}
  	\label{geometric splittings}Let $S$ be the fundamental group of a compact surface with boundary $\Sigma$ supporting a hyperbolic metric (i.e. neither a torus nor the connected sum of one or two projective planes) and $T$ a minimal non-trivial $\pi_{1}(\Sigma,*)$-tree in which all peripheral subgroups of $S$ are elliptic.
  	Then there is an equivariant map $f$ from $T$ to $T_{\delta}$ for some \cfam $\delta$ on $\Sigma$. If all stabilizers of edge groups in $T$ are cyclic then
  	$f$ can be taken to be an isomorphism.
  \end{lemma}
  
  \subsection{Geometric abelian decompositions}
    We will borrow the notion of a geometric abelian decomposition from \cite{bestvina2009notes}, with only slight modifications.
    Let $G$ be a finitely generated group and $\mathcal{R}$ a family of subgroups of $G$.
    By a geometric-abelian $G$-tree  we intend a simplicial $G$-tree $T$ with abelian edge stabilizers together with a partition of the set of vertices $V$ of $T$, into $G$-invariant sets $V_{a}$,$V_{s}$ and $V_{r}$,
    the set of abelian, surface and rigid type vertices respectively,
    such that the following properties are satisfied:
    %    	, such that:
    \enum{i)}{
    	\item For $v\in V_{a}$ the group $\Sb{v}{G}$ is abelian. Additionally, we distinguished what we will call the generalized peripheral subgroup of $Stab(v)$, or $Per^{*}(v)$. This is
    	a free summand of $Stab(v)$ containing the peripheral subgroup $Per(v)$, which is the smallest summand containing the stabilizer of each of the edges incident at $E$.
    	%  			 	is referred to as an abelian type vertex.
    	%the intersection $ker(f)$ for all homomorphisms $\fun{f}{Sb{v}{G}}{\Z}$. of $\Delta_{v}$
    	\item Each $v\in V_{s}$ is endowed with an isomorphism $\theta_{v}$ between $Stab(v)$ and the fundamental group of a compact surface with boundary (possibly empty) $\Sigma_{v}$, which admits essential simple closed curves \footnote{In other words, } . Each edge stabilizer of an edge incident to $v$ is mapped by $\theta_{v}$ to some finite index subgroup of the peripheral subgroup
    	of $\pi_{1}(\Sigma)$ and the map thus induced from the set of orbits of edges incident to $v$ to the family of boundary components of $\Sigma$ is injective.
    }
    We say that such $T$ is relative to $\mathcal{R}$ in case for any $R\in\mathcal{R}$:
    \elenco{
    	\item $R$ is elliptic in $T$.
    	\item For any $v\in V_{s}$ the intersection of $R$ with $Stab(v)$ is contained in a boundary subgroup .
    }
    In this context, if the preimage of each boundary subgroup of $\pi_{1}(\Sigma_{v})$ contains a finite index subgroup which is either conjugate into $\mathcal{R}$ or the stabilizer of an edge incident to $v$ we say that all boundaries are used. A sufficient condition for all boundaries to be used is for $G$ to be freely indecomposable relative to $\mathcal{R}$ \footnote{The fundamental group $S$ of a surface with boundary is a free group and given any boundary component any collection $\mathcal{B}$ of generators of each of the remaining boundary components can be extended to a base of $S$. In particular $S$ admits a free decomposition relative to $\mathcal{B}$. }) .
    \begin{comment}
    \enum{i)}{
    	\item Is the stabilizer of an edge incident to $v$. \label{edge condition}
    	\item Can be conjugated into a subgroup from $\mathcal{R}$ \label{relative condition}
    }
    \end{comment}
    %  	Dually, we will refer to any graph of groups decomposition of $G$, $\Delta$, together with the associated partition of the edges of its underlying graph as a \gad.
    A graph of groups decomposition $\Delta$ dual to $T$, together with the induced partition of the vertices of the underlying graph and peripheral subgroups of the abelian type vertex groups is called a \gad of $G$. In this context we will often talk about rigid/abelian/surface vertex stabilizers or vertex groups.
    
    By a trivial \gat we intend one in which $G$ stabilizes a vertex of rigid type. The surface $\Sigma_{v}$ can have empty boundary, in which case the definition implies the tree $T$ consists of a single point and $\mathcal{R}$ can only contain the trivial group.
    %% possibly cite Bestvina here
    %  	Given, $v\in V_{a}$, by the peripheral subgroup of $Stab(v)$ or $Per(v)$  we intend the smallest free summand of $Stab(v)$ containing the stabilizer of any edge incident to $v$
    %  	and the intersection with $Stab(v)$ of any conjugate of a group in $\mathcal{R}$.\footnote{In oder words, the root closure in $Stab(v)$ of the group generated by all of those groups.}
    %  	The last property above implies that the intersection of any $R\in\mathcal{R}$ with $Stab(v)$ is contained in
    %  	$Per(v)$ for $v\in V_{a}$ and in a boundary subgroup of $Stab(v)$ for $v\in V_{s}$.
    %  		A \gad will be called commutative-transitive if for each $e\in E\Delta$ for some $f\in\{e,\edin{e}\}$ the image $i_{e}(\Delta_{e})$ is maximal abelian in $\Delta_{\alpha(e)}$.
    
    We will say that a \gat (or the corresponding \gad) is commutative transitive if every stabilizer is its own centralizer in the stabilizer of at least one of the two endpoints.
    
  \subsection{Pinching}
    
    \label{pinching subsection}Let $\Delta$ be a \gad of a group finitely generated group $G$.
    By an \cfam in $\Delta$ we will intend simply the union of an \cfam $\Sigma_{v}$ for each surface type vertex $v$ of the underlying graph
    (for that purpose, we regard surfaces corresponding to distinct vertices as distinct and disjoint).
    Denote by $\Delta_{\delta}$ the graph of groups decomposition obtained by blowing each surface type vertex $v$ in the tree dual to $\Delta$
    to a copy of $T_{\Sigma_{[v]}\cap\delta}$ (since edge groups incident to $v$ fix a unique point in $T_{\Sigma_{[v]}\cap\delta}$, the result is unique).
    This can be seen as a \gad as well, whose surfaces are the pieces in which $\delta$ cuts each of the surfaces in $\Delta$.
    Let $\fun{pinch_{\delta}}{G}{G^{p}_{\delta}}$ be the quotient by the normal subgroup generated by the cyclic groups corresponding to simple closed curves in $\delta$.
    This amounts to collapsing to a point each curve in some disjoint set of representatives of the classes in $\delta$
    %   within each $\Sigma_{0}\in [\Sigma_{v},\Theta_{\lambda}]$
    making the corresponding edge groups.
    %   $\Sigma'_{0}$ the surface resulting in each case and
    Denote by $\Delta^{p}_{\delta}$ the resulting graph of groups decomposition of $G^{p}_{\delta}$ and by $Pinch(\Delta,\delta)$ the union of all surfaces resulting from those
    of $Cut(\Sigma_{v},\delta\cap\Sigma_{v})$ by collapsing the curves in $\delta$ to a point.
    %  We use the notation and $Pinch(\Delta,\delta)$ and $p^{\delta}$ to refer to the this kind of operation applied to a single surface $\Sigma$ with boundary and its fundamental group.
    We distinguish between the set $Pinch^{in}(\Sigma,\Delta)$ of all of those which contain some boundary component of the
    original $\Sigma_{v}$ and that $Pinch^{ex}(\Sigma,\Delta)$ of those which do not.
    Let $\Delta_{\delta}^{c}$ be the free decomposition obtained by collapsing the non-trivially stabilized edges of $\Delta_{\delta}$ and $X^{c}_{\delta}$ its underlying graph. This is of the form:
    \begin{equation}
    	(\bfrp{A_{i}}{i=1}{m})\frp(\bfrp{\pi_{1}(\Pi)}{\Pi\in Pinch^{ex}(\Delta,\delta)}{})\frp\F_{b_{1}(X^{c}_{\delta})}
    \end{equation}
    where the last factor is the free group generated by the Bass-Serre elements and $b_{1}(X)$ stands for the first betti number of the graph $X$, while $A_{1},A_{2}\cdots A_{m}$ are
    what we call the internal pinching factors, that is,
    the images of the fundamental groups of the finitely many graph of groups into which $\Delta^{p}_{\delta}$ decomposes after removing the edges associated with the curves in $\delta$. By the exterior rank of $\Delta^{p}_{\delta}$ we intend the sum of the maximum rank of a free images of the remaining factors, what we call external factors.
    
    %  The following can be seen as a consequence of \ref{geometric splittings}.
    %    A \cfam on an compact surface is maximal if
    \begin{lemma}
    	\label{pinching to a free group}
    	\leavevmode
    	Suppose we are given a closed compact surface $\Sigma$ which is not a sphere. Then any homomorphism $f$ from $\pi_{1}(\Sigma)$ to a free group $\F$ pinches some maximal (in the set theoretical sense) \cfam
    	in $\Sigma$.
    \end{lemma}
    \begin{proof}[sketch]
    	In case $\xi(\Sigma)=0$, i.e., when $\Sigma$ is either a klein bottle or a torus, it follows immediately from the presentation of $\pi_{1}(\Sigma)$
    	and the fact that free groups are $CSA$ with cyclic maximal abelian groups that any such morphism must have a cyclic image. The sought \cfam consists of a single curve that be then be found by hand.
    	
    	As for the other cases, a straightforward induction argument implies that it is enough to show the existence of some essential simple closed curve in $\Sigma$ which is pinched by $f$. This can be deduced from  \ref{geometric splittings}, applied to the action of $\pi_{1}(\Sigma)$ on the Cayley graph of $\F$ induced by $f$ (see also \cite[Theorem 5.2]{reznikov1997quadratic}).
    \end{proof}
    %  This means $f$ can be written as $\bar{f}\circ q\circ p^{\delta}$, where $q$ is the quotient map resulting from killing the fundamental group of any projective plane appearing in $Pinch(\Sigma,\delta)$. The fact that $\delta$ is maximal implies that
    %  the image of $q\circ p^{\delta}$ is a free group.
    
    %  The modular group of $\Sigma$, $Mod(\Sigma)$ is the group of all homeomorphisms of $\Sigma$ which preserve the boundary componenents, up to isotopy (equivalently, homotopy).
    
  \subsection{Seifert type actions of surface groups on real trees}
    %  punctured surface
    The notion of a singular measured foliation was introduce by Thurston in his groundbreaking work on the classification of homeomorphisms of surfaces (see \cite{thurston1988geometry} and \cite{fathi2012thurston}). It generalizes that of an \cfam on a surface.
    %For hyperbolic surfaces, a equivalent notion commonly found in the literature is that of a measured lamination.
    Discrete versions of that notion are often useful: this is the approach found in \cite{levitt1997geometric}, which uses simplicial complexes, or \cite{bestvina1995stable}, which uses band complexes (a finite graph to which foliated squares are attached along the sides perpendicular to the foliation, plus some additional) are used instead.
    
    %  We briefly introduce the former
    %  \begin{comment}
    What follows is only an intuitive sketch. Roughly speaking, given a surface with boundary $\Sigma$, a singular foliation on $\Sigma$
    %  is given by finitely many points $x_{1},x_{2}\cdots x_{m}$ of $\Sigma$ and a differentiable choice, for each $x\in\Sigma\setminus\{x_{1},x_{2},\cdots x_{m}\}$ of an element in the projective tangent bundle
    is a partition $\mathcal{F}=\{F_{\lambda}\}_{\lambda\in\Lambda}$ of $\Sigma$ into path-connected sets called leafs, admitting a differential atlas $\mathcal{A}$ of $\Sigma$, such that for some finite set $x_{1},x_{2}\cdots x_{m}$ of points of the surface (possibly on its boundary) the following properties are satisfied. First of all, any boundary component is contained in a leaf.
    %  For the sake of simplicity we consider only the case in which $\partial\Sigma=\nilthe$.
    For any chart $\fun{\phi}{U}{\phi(U)\subset\R^{2}}$ in $\mathcal{A}$ such that $U\cap\{x_{i}\}_{i=1}^{m}=\nil$ the intersection of each
    $F_{\lambda}$ with $U$ is equal to $\phi^{-1}(\bunion{i\in I}{\setof{(x,y)}{y=c_{i}}}{})$ for some countable collection $\{c_{i}\}_{i\in I}$ of constants. For any other chart 	$\fun{\phi}{U}{\phi(U)\subset\R^{2}}$ of $\mathcal{A}$ the set $U$ contains a unique point $x_{i_{0}}$, sent to $(0,0)$ by $\phi$  and any $F_{\lambda}$ is the pullback of a similar set by the map
    $g_{k}\circ\phi$, where $\fun{g_{k}}{\R^{2}}{\R^{2}}$, is some branch of the complex exponentiation $z\mapsto z^{\frac{k}{2}}$, for $k=1,3,4\cdots$ up to identifying $\R^{2}$ with $\mathcal{C}$.
    %	  In particular, assuming $\phi(U)$ is convex,
    %    This mean sone can assume $F_{\lambda_{s}}\cap U$ to be homeomorphic to the union of $k$ semi-open segments with a common end in $\phi(\mathring{U})$, while for any other $\lambda$ the intersection
    %	  $U_{\lambda}\cap F_{\lambda}$ is an open segment.
    A leaf that contains one of the $x_{i}$ is called singular.
    An arc $\alpha$ in $\Sigma$ is said \emph{transverse} to $\mathcal{F}$ if it can be covered by a union of domains of charts $\fun{\phi}{U}{\phi(U)}$ in $\mathcal{A}$ with the property that for any of them
    $\phi(\alpha\cap U)$ does not intersect any of the level sets $\{y=c\}$ twice (resp. the preimage of such a set by $f_{k}$).
    A measured foliation on $\Sigma$ is given by a pair $(\mathcal{F},m)$, together with a non-negative real valued function $m$ on the set of transverse paths between points of $\Sigma$ such that $m(\alpha)=m(\alpha')$ for
    paths $\alpha$ and $\alpha'$ whenever there is a homotopy between $\alpha$ and $\alpha'$ keeping any point in the same leaf.
    
    Any measured foliation lifts to one, $(\tilde{\mathcal{F}},\tilde{m})$ on the fundamental cover $\tilde{\Sigma}$ of $\Sigma$ invariant under the action of $\pi_{1}(\Sigma)$ on $\tilde{\Sigma}$.
    %  It can be proven that given any two points $x,y\in\tilde{\Sigma}$ transverse paths between $x$ and $y$, that they intersect any leaf of $\tilde{F}$ at most once and
    %  that their measure depends only on the leafs containing $x$ and $y$,
    %  $x$ and $y$.
    One can define a pseudo-metric on the set $\tilde{\mathcal{F}}$ by letting the distance between to leaves be equal to the infimum of the measure of the paths between them and that the corresponding quotient metric space $X=\tilde{\Sigma}/\mathcal{F}$ is an $\mathcal{R}$-tree.
    %	  One says that $\mathcal{F}$ is arational if there is no essential simple closed curve contained in a leaf of $\tilde{F}$. This implies that any element fixing a leaf of
    %
    It is clear that if no leaf of $\mathcal{F}$ contains a path homotopic to an essential simple closed curves the only elements of $\pi_{1}(\Sigma)$ fixing a $\tilde{\mathcal{F}}$ are the boundary subgroups. Under certain additional assumptions the quotient becomes nice enough, so that an element $\pi_{1}(\Sigma)$ fixes a point of $Y_{\mathcal{F}}$ if and only if it is contained in a boundary subgroup. One can also assume each boundary component for contains a singularity.
    
    The resulting action of $\pi_{1}(\Sigma)$ on $X$ is what is known as a Seifert type action. Its main properties are the fact that
    %     an element $g\in\pi(\Sigma)$ is elliptic if and only if it belongs to a boundary subgroup of $\pi(\Sigma)$, that
    no non-trivial element fixes a non-degenerate segment and that distinct boundary subgroups fix distinct points of $X$.

\section{Simplicial trees and group automorphisms}
  \label{trees and automorphisms section}
  %  We will temporarily drop the composition symbol for the scope of this section.
  \begin{comment}
  \begin{definition}
  	Given a connected graph $X$, and edges $e,f$ of $T$, we say that $e\succeq f$ if and only if there is a simple path in $T$ such whose first edge is $e$ and whose last edge is $f$ (possibly coinicdent).
  \end{definition}
  It is easy to check that both $e\succ f$ and $f\succ e$ cannot be the case in case $X$ is a tree.
  \end{comment}
  Let $\BS{Y_{0}}{Y}{t}$ be a Bass-Serre presentation for an action without inversion of a finitely generated group $G$ on a simplicial tree $T$. Fix $e\in Y$ and $c\in Z_{Stab(e)}(Stab(\omega(e)))$ and
  for each $v\in VY$ let $c_{v}$ be equal to $1$ if $e$ is in the unique path from $\alpha(e)$ to $v$ in $Y_{0}$ and $c$ otherwise. Likewise, let $c_{f}$ be equal to $c$ if $e\succeq f$ and $1$ otherwise. Consider the map which sends:
  \elenco{
  	\item  $g\in Stab(v)$ to $g^{c_{v}}$ for $v\in VY$.
  	\item  $t_{f}$ to $c^{-1}_{\edin{f}}t_{f}c_{f}$ for $f\in E\setminus Z$
  }
  It is easy to check that this induces a homomorphism from $G$ to itself, $\tau_{e,c}$, which we refer to as the Dehn twist over $e$ by $c$. It is also easy to check that $\tau_{e,c}\circ\tau_{e,c^{-1}}=\tau_{e,c^{-1}}\circ\tau_{e,c}=Id$, so that $\tau_{e,c}$ is in fact an automorphism of $G$. We will refer to $e$ as the base of $\tau_{e,c}$, denoted by $B(\tau_{e,c})$.
  
  Now suppose we are given $v_{0}\in VY$ and an automorphism $\sigma$ of $Stab(v_{0})$
  compatible with the incidence structure, by which we mean
  that $\sigma\rst_{Stab(e)}$ restricts to conjugation by some element $c_{e}\in Stab(v_{0})$ for each $e\in E$ with $\alpha(e)=v_{0}$.
  
  Then there is a natural way of extending $\sigma$ to an endomorphism $\bar{\sigma}=\bar{\sigma}_{(c_{e})_{e}}$ of $G$ such that:
  %  	\alpha(e)=v,\,e\in EY
  \enum{i)}{
  	\item $\bar{\sigma}$ restricts to $\sigma$ on $Stab(v_{0})$
  	\item $\bar{\sigma}$ restricts to $\inn{c_{v}}$ on $Stab(v)$, for $v\in VY$, where $c_{v}=c_{in(e)}$ for $in(e)$ the first edge in $[v_{0},v]$
  	\item $\bar{\sigma}(g_{e})=c_{e}^{-1}t_{e}c_{\edin{e}}$ for $e\in E$ where $c_{e}=c(u)$ for the unique vertex $u\in Z_{0}\cap G\cdot\omega(e)$.
  }
  One can check, just as before, that this induces an isomorphism of $G$ with inverse $\bar{\sigma^{-1}}_{(c_{e}^{-1})_{e}}$, to which we will refer to as a natural extension of a vertex group automorphism. In this case we let $\{v\}$ be the base $B(\sigma)$ of $\bar{\sigma}$. We will refer to the elements $c_{e}$ above as twisting elements.
  
  We will refer to Dehn twists and extensions of vertex group automorphisms as elementary automorphisms of $G$.
  
  \begin{comment}
  \begin{fact}
  	If $G$ admits a commutative transitive \gad such that each rigid type vertex group is $CSA$ then $G$ is $CSA$.
  \end{fact}
  \end{comment}
  %  \begin{definition} \label{mod generators}
  	Suppose we are given a \gd \pr-group $\gqG{G}{A}{P}{I}$. By \gat (resp. \gad ) of it we intend a \gat (resp. \gad) $T$ of $G$ relative to $\{A\}\cup\{P_{i}\}_{i\in I}$.
  	By a \gat of $\rqG{G}{A}$ we intend one relative to $\{A\}$.
  	%  		, by the \pr-modular group of $\qG{G}$ with respect to
  	%  		$T$.
  	Fix a \BSP $\mathcal{P}=\BS{Y}{S}{t}$ for the action of $T$. By $Mod(\rqG{G}{A},T)$ or the \emph{modular group} of $\qG{G}_{A}$ with respect to $T$ we intend the subgroup of $Aut(G)$ generated by all those elementary automorphisms $\sigma$ of the form:
  	\enum{i)}{
  		\item A Dehn \twist over an edge of $Y$.
  		\item A vertex group automorphism $\vga{\sigma}{c}{v}{S}$ of some $v\in VY$, such that:
  		\enum{i)}{
  			\item  $\sigma$ is equal to $Id_{Stab(v)}$ in case $v\in V_{r}$.
  			\item  $\sigma$ is induced by a homeomorphism of $\Sigma_{v}$ in case $v\in V_{s}$
  			\item	 $\sigma$ fixes $Per^{*}(v)$ for any $v\in V_{a}$
  		}
  		\item An inner automorphism $\inn{c}$ of $G$ (conjugation by $c$ ).
  	}
  	and satisfying properties:
  	\enum{i)}{
  		\item  $\pi\rst_{Stab(v)}\circ\sigma=\pi\rst_{Stab(v)}$ and $\pi(c_{e}),\pi(c)\in Z(Q)$
  		\item $\sigma\rst_{A}=Id\rst_{A}$
  	}
  	\newcommand{\PMGen}[0]{MGen^{\pi}}
  	Denote by $MGen(\rqG{G}{A},\mathcal{P})$ the set of generators listed above; we call them elementary automorphisms. It can be shown that the subgroup they generate does not depend on the presentation $\mathcal{P}$, justifying the use of the alternative notation $\PMGen(\rqG{G}{A},T)$ or $\PMGen(\gqG{G}{A},\Delta)$ for any associated \gad. We denote the subgroup generated by
  	$\PMGen(\gqG{G}{A},\Delta)$ as $Mod(\qG{G}_{A},\Delta)$ (or dually $Mod(\qG{G}_{A},T)$, where.  Alternatively, we might use the notation $Mod^{\pi}_{A}(\Delta)$.
  	
  	Ignoring the \pr-structure yields the standard modular group $Mod_{A}(G,T)$. The definition extends to \rs \pr-groups and \pr-groups by regarding them as \gr \pr-groups in the obvious way.
  	Note that in particular, $\sigma$ will fix $Per^{*}(v)$ for any $v\in V_{a}$. Given $v$, we will let $Mod^{\pi}(v)$ stand for the subgroup of all possible $\sigma\in Aut(Stab(v))$ satisfying the restrictions above.
  	%  	 $Mod(v)$ in case we ignore $\pi$.
  	%  		Alternatively we will denote such group by $PMod(\Delta;A)$, where $\Delta$ is the graph of groups decomposition of $G$ associated with $P$.
  	%  		Removing the condition on the invariance of $\pi$ leaves us with the standard modular group $\Mod{G,A,T}$.
  	
  	If the given $T$ is merely a pure $G$-tree over abelian edge stabilizers and not a \gat then in an expression such as $Mod(G,T)$ we will be implicitly promoting $T$ to a \gat by declaring all its vertices as rigid (so that only Dehn twists are taken into account).
  	%  	Given a vertex $v$ of surface type, denote by $PMod(Stab(v))$ the subgroup consisting of all the automorphisms
  %  \end{definition}
  %  	Notice that $\sigma$ fixes $Per(v)$ in case $v$ is of abelian type. The following fact is not essential for our arguments, but is worth pointing out.
  \begin{observation}
  	$PM:=Mod(\rqG{G}{A},T)$ is a subgroup of finite index in $M:=Mod_{A}(G,T)$.
  \end{observation}
  \begin{proof}
  	The group $M$
  	%  	the group $Inn_{G}(Z_{G}(A))$
  	\newcommand{\MA}[0]{Mod_{A}(G,T)}
  	\newcommand{\M}[0]{Mod(G,T)}
  	contains the group $I:=Inn(Z_{G}(A))$ of inner automorphisms of $G$ by elements of $Z_{A}(G)$ as a normal subgroup. Clearly $PM$ intersects
  	$I$ in a subgroup of finite index. On the other hand, given a transversal $\mathcal{R}$ of $G\backslash VT$, for any $v\in\mathcal{V}$ and
  	$\tau\in M$ there is $c$ such that $\inn{c}\circ\tau$ preserves $Stab(v)$ and it is possible to check that its projection in
  	$Mod(G_{v})/Inn(G_{v})$ is well defined. This determines a map $\fun{\eta_{v}}{M/I}{Mod(G_{v})}$. It is well known that the collection of all of them induces an isomorphism
  	between $M/I$ and $\bcsum{Mod(Mod(G_{v})/Inn(G_{v})}{v\in\mathcal{R}}{}$.
  	Indeed, this is well-known to be true in case $A=1$ (see \cite{levitt2005automorphisms}), and it is easy to see that the inclusion map induces an isomorphism between $M/I$ and $Mod(G,T)/Inn(G)$.
  	Clearly the image of $Mod{\rqG{G}{A},T}$ intersects each of the direct summands in a finite index subgroup.
  	Of course, given any group $G$, $N\unlhd G$ and $H\leq G$, if both $[N,N\cap H]$ and $[G/K,NH/K]$ are finite, then $[G,H]$ is finite as well.
  \end{proof}
  
  %  	Observe that in the definition we might have replaced the first class of generators with that of Dehn twists over edges of $E$.
  The following can be easily checked by inspection.
  %  \begin{fact}
  	\begin{observation}
  		\label{elementary automorphisms}
  		If $e$ is a vertex from $u$ to $v$, and $Stab(v)$ is abelian, then $\tau_{e,c}=\tau_{e_{1},c}\cdots\tau_{e_{k},c}$, where $e_{i}$, $i=1,2\cdots k$ enumerates the edges $f\neq\edin{e}$ with $\alpha(f)=v$.
  		Notice also that $\inn{c^{-1}}\tau_{e,c}=\tau_{\edin{e},c^{-1}}$ and for every $v\in VT$ and vertex group automorphism $\bar{\sigma}_{(c_{e})_{\alpha(e)=v}}$ and $g\in Stab(v)$ one has $\inn{g}\bar{\sigma}_{(c_{e})_{\alpha(e)=v}}=\bar{\inn{g}\sigma}_{(c_{e}g)_{\alpha(e)=v}}$. In particular, for any $e_{0}$ with $\alpha(e_{0})=v$ we have
  		$e_{0}\nin Supp(\inn{c_{e_{0}}^{-1}}\circ\bar{\sigma}_{(c_{e})_{\alpha(e)=v}})$.
  	\end{observation}
  %  \end{fact}
  
  \begin{definition}
  	The \pr-\emph{modular group} of a freely indecomposable $\gqG{G}{A}{P}{I}$, denoted by $Mod(\gqG{G}{A}{P}{I})$ is the subgroup of $Aut(G)$ generated by $Mod(\rqG{G}{A},\Delta)$ where $\Delta$ ranges among all \gad of $\gqG{G}{A}{P}{I}$.
  	The modular group of $\rqG{G}{A}$ is simply the modular group of the associated \gr \pr-group (with trivial grading).
  	%  		Given a freely indecomposable \rs \pr-limit group $\rqG{G}{A}$, $PMod(G;A)$ is the subgroup generated $PMod(\Delta)$, where $\Delta$ is any \rs \gad of $\rqG{L}{A}$.
  	%  		In the case of a freely indecomposable \gd \pr-limit group $\gqG{G}{A}{P}{I}$, we let $PMod()$ be the subgroup generated by all $\tau\in Mod(\rqG{G}{A}}{}\Delta)$, where $\Delta$ is any \gd \gad of $PMod()$.
  \end{definition}
  \begin{comment}
  The fact that modular groups of surfaces are generated by Dehn twists for the one edge splittings coming from a s.c.c. and that transvections generate the point-wise stabilizer of a direct summand of a free abelian group in its automorphism group yields the following fact (see \cite{wilton2009solutions} for a detailed proof):
  \begin{fact}
  	The \pr-modular group of a freely indecomposable \pr-group $\mathcal{G}$, \rs group $\rqG{G}{A}$ or \gr group $\gqG{G}{A}{P}{I}$ is generated by Dehn twists over one edge splittings of $G$, which are \rs or \gr in the second and third case respectively.
  \end{fact}
  We now recall a very well known fact. Abelian decompositions can be made relative to non-cyclic subgroups without any substantial loss of information.
  \end{comment}
  \newcommand{\NZA}[0]{\mathcal{A}^{nz}}
  \begin{lemma}
  	\label{non cyclic abelians} Suppose we are given a \gad $\Delta$ of a $CSA$ \gr \pr-group $\gqG{G}{A}{P}{I}$. Assume furthermore that any maximal abelian subgroup of $\qG{G}$ is finitely generated. Then any $\sigma\in Mod(\rqG{G}{A},\Delta)$ is contained in $\sigma\in Mod(\rqG{G}{A},\Delta')$, for a \gad $\Delta'$ of $\gqG{G}{A}{P}{I}$ whose underlying tree is relative to the family of all abelian non-cyclic subgroups of $G$.
  \end{lemma}
  \begin{proof}
  	Take any non-cyclic maximal abelian subrgroup $M$ of $G$ which is not elliptic. Since $M$ is finitely generated, in virtue of \ref{abelians acting}, it acts by translations on some line $L$ of
  	$T=\tilde{(\Delta,Z)}$.
  	
  	Since $M$ is non-cyclic, the group $K=Fix(L)\cap M$ is non-trivial. Let $T'$ be the result of collapsing all edges in $L$. If $\sigma$ is based on a vertex outside $L$
  	(in particular, if this vertex is of surface type), then it is clear it belongs to $Mod(\rqG{G}{A},T')$ as well. If not, observe that by \ref{factoring lines} we can decompose
  	$G$ as an amalgamated product of $M$ over the edge group $K$. By inspection one can easily check that when $\sigma$ is a Dehn twist over an edge in $L$ or the extension of an automorphism of the stabilizer of an abelian type subgroup of $L$ the latter decomposition does the job. %
  	%  		It is easy to check any
  	%  		elementary $\sigma$ as above belongs to $Mod(\gqG{G}{A}{P}{I},T')$ as well (we regard the vertices resulting from the collapse as abelian or not depending on whether the group is actually abelian).
  	%  		 (as, for example no surface type vertex can belong to $L$), with the only exception of
  	%  		There are two possibilities. If $\Delta$ is an amalgamated product,
  	%  		Suppose there is $e\in ET$ with $\alpha(e)=v$ and $g\in Stab(v)\setminus\{e\}$ such that $g^{-1}\cdot e,e\in L$. It follows that $K\subset Stab(e),Stab(e)^{g}$. Sine
  	%  		$G$ is $CSA$, this implies that $g\in M$. This is impossible, sice all elliptic elements of $M$ fix $L$.
  	%  	Consider the image of $L$ in the quotient graph $G\backslash T$.
  \end{proof}
  Let $G$ be a finitely generated $CSA$ group and let $T$ be a \gat for $G$. We will say that $T$ is \emph{normalized} if:
  \enum{i)}{
  	\item It is relative to the family of all non-cyclic abelian subgroups of $G$.
  	\item Exactly one of the endpoints of each edge $e\in T$ is contained in an abelian type vertex, the other endpoint being non-abelian.
  	\item The intersection of $Stab(e)$ and $Stab(e')$ for each rigid type vertex $v\in T$ and any two distinct edges $e,e'$ originating at $v$.
  }
  % 	 Given a \gad, it is generally possible to take it into a normalized form without any essential loss of information.
  %  In the literature this is often done by hand; an alternative is the tree of cylinders construction introduced in \cite{guirardel2011trees}.
  \begin{observation}
  	A normalized \gat in which the stabilizer of every rigid type vertex is $CSA$ is $2$-acylindrical.
  \end{observation}
  %  [Given a relative decomposition, how to get the \gad  structure; we need to exclude vertices which of surface type which are not hyperbolic-hyperbolic, etc.]
\section{Grushko and JSJ decompositions}
  
  In this section we give a quick overview of JSJ theory, as laid out in \cite{guirardel2009jsj} and \cite{guirardel2010jsj}.
  Or main goal is to show the following:
  \begin{proposition}
  	\label{modular group universality} Any \gr \pr-limit group $\gqG{L}{A}{P}{I}$ admits some \gad $\Delta_{JSJ}$ such that any \pr-modular automorphism of $\gqG{L}{A}{P}{I}$ is contained $Mod(\rqG{L}{A},\Delta_{JSJ})$.
  \end{proposition}
  
  Let $G$ be a fixed finitely generated group and $\mathcal{R}$ some family of subgroups of $G$.
  \begin{comment}
  We will assume that $\mathcal{R}$ is invariant under conjugation. This hypothesis is not essential. The definitions presented here make sense without it and any non-invariant family $\mathcal{R}$ could be replaced by that of all subgroups conjugate to one in $\mathcal{R}$; the statement of some of the conclusions, however might need a slight (obvious modification).
  \end{comment}
  %  All notions introduced in this section can be defined for a non-invariant family in exactly the same way. The only thing th
  The following result was first proven by Grushko in case $\mathcal{R}=\nil$ (a nice topological proof due to Stallings can be found in \cite{stallings1965topological}).
  \begin{theorem}
  	There is a free decomposition relative to $\mathcal{R}$ of the form:
  	\begin{align*}
  		G=(\bfrp{G_{i}}{i=1}{m})*F
  	\end{align*}
  	Where $l,m\geq 0$, $F$ is a (possibly trivial) free group and each $G_{i}$ is freely indecomposable relative to $\mathcal{R}\cap G_{i}$. The decomposition is finest in the sense that given any other decomposition of the same form $(\bfrp{G'_{i}}{i=1}{m'})*F'$ then every $G'_{i}$ is the free product of a free factor of a conjugate of $F$ and some conjugates of groups $G_{j}$, while $F'$ is a free factor of a conjugate of $F$.
  \end{theorem}
  \newcommand{\ER}[0]{$(\mathcal{E},\mathcal{R})$}    %text command
  As a matter of fact Bass-Serre presentations are pretty much irrelevant to the discussion, which will thus take place in the language of trees rather than in that of graph of groups decompositions. Given a conjugacy invariant family $\mathcal{E}$ of subgroups an \ER-tree is a simplicial $G$-tree relative to $\mathcal{R}$ with edge stabilizers in $\mathcal{E}$. The notion of an \ER-JSJ generalizes the notion of a Grushko decomposition to the class of \ER-trees.
  \begin{comment}
  %  Proofs of the existence of JSJ decompositions for particular cases started to appear before any consensus on which of their properties would enter the definition of any general definition (see \cite{FuPaJSJ},\cite{RipsSelaJSJ}). We survey here the approach found in \cite{GuiLevJSJ}, whose definition is more general than will be needed here (it leaves the characterization of flexible vertices to the particular case), stroke us as being particularly clear and amenable to application (see for example \cite{PerinRizos forking} for an excellent example). Further work of the same authors on the construction of an $Aut_{\mathcal{R}}(G)$ invariant JSJ decomposition will be also useful (see also \cite{GuiLevJSJ2}).
  %  \cite{SelaJSJ}
  %  From now on in this section $\mathcal{E}$ will stand for some conjugacy invariant family of subgroups of $G$. %All the notions to be defined depend on $\mathcal{E}$, though we might
  \end{comment}
  A subgroup $H\leq G$ will be called \ER-universally elliptic in case it is elliptic in any \ER-tree. A simple example of an universally elliptic subgroup is that of some subgroup $H\leq G$ containing some $R\in\mathcal{R}$ as a finite index subgroup.
  %  \footnote{Let $v$ be the vertex fixed by the corresponding $R\in\mathcal{R}$. The orbit of $v$ by $F$ is finite. One can easily check  } .
  Given $G$-trees $S$ and $T$, one says that $S$ dominates $T$ if and only if $Ell(S)\subset Ell(T)$, where $Ell(T)$ stands for the familiy of all elliptic subgroups of $G$ elliptic in $T$.  It is easy to check this is equivalent to the existence of an equivariant map from $VS$ to $VT$.
  
  A deformation space is the class of all \ER-trees $T$ for which $Ell(T)$ has a particular value.
  %  for which $Ell(T)$ has a certain common value.
  %  is a certain fixed family of subgroups of $G$.
  A \ER-universally elliptic tree is an \ER-tree all of whose edge stabilizers are \ER-universally elliptic. We say that an universally elliptic tree $T$ is an \ER-JSJ tree if $T$ is \ER-universally elliptic and $Ell(T)$ is minimal among universally elliptic trees.
  %  The deformation space of $T$ is the class of all \ER-trees $S$ such that $Ell(S)=Ell(T)$.
  The following follows from \cite[3.2,3.8]{guirardel2009jsj} (the second of those results
  is formulated there only in case $\mathcal{R}=\nil$ but the difference is inessential).
  \begin{lemma}
  	\label{UE property}Let $T$ be an universally elliptic \ER-tree and $S$ an \ER-tree. Then there is a refinement %still to define
  	$\hat{T}$ of $T$ such that $Ell(\hat{T})=Ell(T)\cap Ell(S)$. If $S$ is a universal \ER tree then $\hat{T}$ is an universal \ER tree as well.
  \end{lemma}
  As a consequence, if the class of \ER-JSJ trees is non-empty, it is a deformation space, referred to as the \ER-JSJ deformation space, and every \ER-universally elliptic $T$ tree has a \ER-JSJ refinement.
  %  	Given any universally elliptic \ER-tree $T$ and some \ER-tree $S$, therer is some refinement of $T$ which dominates $S$.
  Using this it is possible to give an outline of the basic strategy used in proving existence. The following is lemma 3.9 in \cite{guirardel2009jsj}
  \begin{lemma}
  	Let $\{T_{i}\}_{i\in I}$ be any family of $G$-trees, where $G$ is finitely generated. There exists a countable subset $J\subset I$ such that if $T$ is elliptic with respect to every $T_{i}$ for $i\in I$ and dominates every	$T_{j}$ for $j\in J$ then $T$ dominates $T_{i}$ for all $i\in I$.
  \end{lemma}
  Let $\{T_{n}\}_{n\in \N}$ be the countable subfamily extracted from the family of all \ER-universally elliptic trees using the previous lemma. Iterative application of \ref{UE property} yields a sequence of $(T'_{n})_{n}$ of \ER-universally elliptic trees, where $T'_{n}$ is a collapse of $T'_{n+1}$ and dominates every $T_{m}$ for $m\leq n$.
  All it is then left to show is that the deformation space of any chain of refinements of \ER-trees eventually stabilizes. This is what is commonly referred to as an accessibility result.
  
  The following is a relative version of \cite[Proposition 4.4]{guirardel2009jsj}, which was originally proved in \cite{dunwoody1985accessibility}
  \begin{lemma}
  	Suppose $G$ is finitely presented and let $\ssq{T}{n}$ be a sequence of
  	$(\mathcal{E},\mathcal{R})$-trees, where $T_{n+1}$ is a refinement of $T_{n}$ for any $n\in\N$. Then there is some $G$-tree $S$ dominating each of the $T_{n}$.
  \end{lemma}
  If one takes for granted, as we do, that limit groups are finitely presented, this is good enough for our purposes. Without this assumption one needs to apply acylindrical accesibility.
  See \cite{Sela1997acylindrical} and \cite{weidmann2002nielsen} for the basic combinatorial result and \cite{reinfeldtmakanin} for its application to limit groups.
  
  Given an \ER-JSJ tree $T_{JSJ}$ one distinguishes two types of vertices (vertex stabilizers)
  of $T$. We say that $v$ (resp. $Stab(v)$) is \emph{rigid} if $Stab(v)$ is \ER-universally elliptic. Otherwise we will refer to it as \emph{flexible}
  %  It follows from the properties of the blow-up construction that:
  The following lemma is easy to prove using the blow-up construction.
  %  Otherwise, we say it is flexible. For the following
  %  see \cite[section 2]{guirardel2009jsj}.
  \begin{lemma}
  	If $Stab(v)$ is rigid then there is no non-trivial $Stab(v)$-tree $S$ over edge stabilizers in $\mathcal{E}$ in which the stabilizer of any $e$ incident to $v$ in $T$ and the intersection with $Stab(v)$ of any conjugate of a subgroup from $\mathcal{R}$ are elliptic.
  \end{lemma}
  \begin{comment}
  \begin{proof}
  	Indeed, the direction from right to left is trivial. For the direction from left to right, blowing-up of $v$ to a copy of $S$ produces a \ER-tree which contradicting the fact that $Stab(v)$ is universally elliptic.
  \end{proof}
  \end{comment}
  Let us now restrict to the case in which $\mathcal{E}$ contains only abelian subgroups.
  Take any \ER-tree $T$ and suppose $v$ is of surface type, with $Stab(v)\cong\pi_{1}(\Sigma)$.
  %   Suppose that a \ER-universally elliptic tree $T$
  %  is given a \gat structure relative to $\mathcal{R}$
  We claim that for any minimal non-trivial $Stab(v)$-tree $S$ relative to the boundary subgroups over edge stabilizers in $\mathcal{E}$ and any edge $e\in S$ the stabilizer of $e$ is not universally elliptic. Indeed, (up to subdivision) we know $S$ to be dual to some \cfam $\delta$ on $\Sigma$. On the other hand there is some essential simple closed curve $\alpha$ such that $i(\alpha,\beta)\neq 0$ for any $\beta\in\delta$. The result of blowing up $v$ to $T_{\alpha}$ in $S$ witnesses the fact $Stab(e)$ is not \ER-universally elliptic.
  
  Under certain weak conditions $Stab(v)$ can be proved to be elliptic in any \ER-universally elliptic tree (see \cite[Proposition 7.13]{guirardel2009jsj}). Note that this is obvious if we assume $T$ to be \ER-universally elliptic to begin with, since in that case it can be refined to a \ER-JSJ tree.
  %   Using the fact that essential simple closed curves fill $\Sigma$, it can be proved that
  %  any universally-elliptic element of $Stab(v)$ is peripheral and,
  %  HYPO
  %  $Stab(v)$ is contained in some flexible vertex stabilizer of any \ER-JSJ tree (this is obvious in case $T$ is universally elliptic, since ).
  %  Now, we know that some \ER-JSJ tree $T_{JSJ}$ can be obtained by blowing-up some of the vertices of $S$. The previous argument shows $v$ is not one of them, and hence $Stab(v)$
  %  is a flexible vertex of $T_{JSJ}$.
  %   In fact $Stab(v)$ is elliptic in the JSJ-deformation space
  %  even if the
  %  It turns out
  
  The crux of JSJ theory is the fact that under certain assumptions
  all flexible vertex groups of a JSJ-tree are of this form (see \cite[theorem 7.7]{guirardel2009jsj}.
  With a different formulation this was proven for example in \cite{rips1997cyclic} using
  a finer version of the technique used to construct test sequences presented in this work.
  %  and the phenomenon above explains all \ER-splittings of $G$ over non \ER-universally elliptic edge groups.
  In truth, in the general case the relevant notion is not that of a surface type vertex group but the broader one of a quadratically hanging (QH) vertex group (see definition 8.17 in\cite{guirardel2009jsj}). Both are equivalent, however, in case $G$ is a torsion free $CSA$ (see \cite[8.12]{guirardel2009jsj}).
  %  \footnote{Slender and abelians subgroup coincide, since any non-abelian subgroup of a limit group contains a free group therefore infinitely generated subgroups. On the other hand $G$ has no nontrivial normal abelian subgroups. Any orbifold with torsion free fundamental group is a surface.}.
  %  introduce the word refine somewhere
  
  An \ER-tree $T$ is called universally compatible if $T$ and $T'$ admit a common refinement $\hat{T}$ for any \ER tree $T'$.
  %  Compare this with the case in which satisfied by an \ER-JSJ tree $T$,
  This is more restrictive than being an \ER-JSJ tree and more convenient (but not essential) when dealing with modular automorphisms. The existence of a universally compatible \ER-JSJ tree is proven in \cite[11.1]{guirardel2010jsj}, under certain hypothesis, in particular in case of a limit group which is freely indecomposable relative to $\mathcal{R}$ and $\mathcal{E}$ the family of all abelian subgroups of $G$.
  Let us collect this result and the above discussion.
  %  	and phrase it in the language of \gats.
  \begin{comment}
  An universally compatible \ER-JSJ tree with the property exists in the case of a torsion-free $CSA$ group $G$ which is freely indecomposable relative to $\mathcal{R}$, where $\mathcal{R}$ contains all non-cyclic abelian subgroups of $G$  (see \cite[11.1]{guirardel2010jsj}). Let
  \end{comment}
  \begin{proposition}
  	\label{JSJ summary} Let $G$ be a limit group which is freely indecomposable with respect to some family $\mathcal{R}$ of subgroups. Let also $\mathcal{A}^{nz}$ be the family of all the non-cyclic abelian subgroups of $G$.
  	%  	let $\mathcal{A}$ be the family of all abelian subgroups
  	%  	of $G$ and $\mathcal{NCA}$ that of all its non-cyclic abelian subgroups.
  	Then $G$ admits a (necessarily unique) universally compatible $(\mathcal{E},\mathcal{R}\cup\mathcal{A}^{nz})$-JSJ tree $T_{JSJ}$ with the following properties:
  	\elenco{
  		\item $T_{JSJ}$ is normalized
  		\item A vertex group of $T_{JSJ}$ is flexible if and only if it is of surface type (relative to $\mathcal{R}$).
  		\item Up to edge subdivision, any minimal reduced \gat $S$ of $G$ relative to $\mathcal{R}\cup\mathcal{A}^{nz}$  can be obtained from $T_{JSJ}$ by equivariantly blowing up some of the flexible vertices of $T_{JSJ}$ according to a \cfam $\delta$ and collapsing (equivariantly) some of the edges inherited from $T_{JSJ}$.
  		\item Given a tree $S$ as above the stabilizer of any surface type vertex of $S$ is the fundamental group of some of the pieces in which the  associated to vertices of $T$.
  	}
  \end{proposition}
  \begin{proof}
  	Flexible vertices of $T_{JSJ}$ are the only ones that can be affected by the blow-up (as we assume $S$ is reduced).
  	%  	\ref{geometric splittings} the action of to a \cfam in the corresponding surface
  	Regarding the last claim, notice that for any vertex stabilizer $Stab(v)$ of $S$ which is not of the form specified above either:
  	\elenco{
  		\item  $Stab(v)$ equals some rigid vertex stabilizer of $T_{JSJ}$
  		\item  $Stab(v)$ splits relative to stabilizers of incident edges over universally elliptic vertex stabilizers.
  	}
  	Both are impossible for surface type vertex stabilizers.
  \end{proof}
  We can regard $T_{JSJ}$ itself as a \gat relative to $\mathcal{R}$ by declaring all its flexible vertices to be of surface type and all those which are abelian and as being of abelian type, with
  extended peripheral subgroup the smallest summand containing the peripheral subgroup and all the conjugates of $\mathcal{R}$ in $Stab(v)$. All the remaining ones are considered to be of rigid type.
  Proposition \ref{modular group universality} can be proven by inspection using \ref{JSJ summary} (see \cite{perinforking}) .
  
  \begin{remark}
  	In the propostion above normality comes from the fact that $T_{JSJ}$ can be taken to be a 'tree of cylinders' in the sense of Guirardel and Levitt.
  \end{remark}
  %
  %  \begin{remark}
  	%  	In the proposition above, one can also assume $T_{JSJ}$ to be normalized in the sense defined above. One approach is to use the 'tree of cylinders' constructions introduce in .
  	%  	A more hands on approach takes any JSJ tree can be also taken to a normal form by performing a series of operations that preserve the modular group, as well as the family of elliptic subgroups. This is the approach
  %  \end{remark}
  %  This is not essential, but makes the construction of the completion of a resolutoi
  %  We can always make this into a nice \gat by attaching to
  Given a \gr \pr-group we will refer to the JSJ of $G$ relative to $\{A\}\cup\{P_{i}\}_{i\in I}$ as above as its \adjsj.
      %	\paragraph{Acylindricity and normalized \gad}
        
        %  [maybe add in which all the abelian subgroups are free of finite rank]
        %  	Fix $G$ a finitely generated $CSA$ group and let $\mathcal{NCA}$ the family of its abelian subgroups of rank $\geq 2$. Following \cite{GuiLev cylinders} we will describe a way of obtaining from each minimal abelian $G$-tree $T$ a tree $T_{c}$ preserving most of the information of $T$ but without any of its redundancies. We might regard it as a canonical form of $T$. A hands-on construction of an isomorphic object for this particular case can be found in \cite{Reinfeldt thesis}.
        
        %  $(\mathcal{AB},\mathcal{NCA})$-$JSJ$
        %   $T_{c}^{*}$ in the same deformation space as $T$
        %  	When the construction is applied to an $(\mathcal{Ab}^{*},\mathcal{NCA})$-$JSJ$ it outputs a canonical $(\mathcal{Ab},\mathcal{NCA})$-$JSJ$ tree.
        %  	This is done in \cite{GuiLev cylinders} with a wider degree of generality

\newcommand{\PP}[0]{\mathscr{P}}
\newcommand{\pe}[1]{p_{e_{#1}}}
\newcommand{\qe}[1]{p_{f_{#1}}}
\newcommand{\ce}[1]{c_{e_{#1}}}
\newcommand{\de}[1]{c_{f_{#1}}}
\newcommand{\pn}[1]{p_{e_{#1},n}}
\newcommand{\qn}[1]{p_{f_{#1},n}}
\chapter{Limits of trees and the shortening argument}

\section{Limits of actions on real trees: the Bestvina Paulin method}
  \label{limit actions section}
  \subsection{Definitions}
    % and \cite{Paulin Gromov topology 2},
    %    Credit vyto Reifieldt Bestvina
    In \cite{paulin1989gromov} Paulin adapted Gromov's notion of convergence among metric spaces to the context of isometric actions of a group $G$ on a hyperbolic space.
    In our case we are mostly concerned with actions on real trees that arise as limits of sequences of actions of a group on the Cayley graph of $\F$ induced by homomorphisms to $\F$.
    In this section we give an overview of the techniques and state the main compactness result needed later. Most of the definitions and results come from \cite{bestvina2002r}.
    \begin{definition}
    	An $\epsilon$-approximation between metric spaces $(X_{1},d_{1})$ and $(X_{2},d_{2})$ is a relation $R\subseteq X_{1}\times X_{2}$ such that:
    	\elenco{
    		\item For all $x_{1},y_{1}\in X_{1}$ and $x_{2},y_{2}\in X_{2}$ we have:
    		\begin{align*}
    			R(x_{1},x_{2})\cap R(y_{1},y_{2})\ra|d_{1}(x_{1},y_{1})-d_{2}(x_{2},y_{2})|<\epsilon
    		\end{align*}
    		
    		\item $R(X_{1})=X_{2}$ and $R^{-1}(X_{2})=X_{1}$
    	}
    	If $X_{1}$ and $X_{2}$ are equipped with an action of a group $G$ by isometries, for any $S\subset G$ we say that the approximation is $S$-equivariant if, in addition
    	for any $(x_{1},x_{2})$ and $s\in S$, if $s\cdot x_{i}\in X_{i}$ then $(s\cdot x_{1},s\cdot x_{2})\in R$.
    	%
    	%    	In case have pointed spaces $(X_{1},x_{1})$ and $(X_{2},x_{2})$ instead, we require the approximation to contain the pair $(x_{1},x_{2})$.
    \end{definition}
    %%In \cite{Paulin Gromov 1}, the author proves a criteria for a sequence of metric spaces to have a convergent subsequence, which is enough for what what we want to prove.
    %% The tools are coming from
    \begin{comment}
    (or $X_{1}$ and $X_{2}$, if no ambiguity regarding the actions considered ensues)
    \end{comment}
    \begin{comment}
    \begin{definition}
    	Suppose isometric actions $\cdot_{1}$ and $\cdot_{2}$ of a group $G$ on metric spaces $Y_{1}$ and $Y_{2}$, respectively, $X_{i}$ a subspace of $Y_{i}$ for $i=1,2$ and $S\subseteq G$ are given.
    	Let $\cdot=\cdot_{1}\times\cdot_{2}$ be the composite action on $Y_{1}\times Y_{2}$. For $g\in G,\,i=1,2$ let $Z_{i}(g)=X_{i}\cap g^{-1}\cdot_{i}X_{i}$ and $Z(s)=Z_{1}(g)\times Z_{2}(g)$.
    	An $S$-equivariant $\epsilon$-approximation between $(X_{1},Y_{1},\cdot_{1})$ and $(X_{2},Y_{2},\cdot_{2})$ is an $\epsilon$-approximation between $X_{1}$ and $X_{2}$ such that for all $s\in S$ one has $(s\cdot R)\rst_{Z(s)}=R\rst_{Z(s)}$.
    	%    	such that for each
    	If we consider pointed spaces
    \end{definition}
    \end{comment}
    Let $G$ be a group and suppose we are given $\mathcal{X}$ a set of tuples $(X,d,\lambda)$, where $(X,d)$ is a metric space, $\lambda$ an action by isometries of $G$ on $X$
    and $\mathcal{Y}$ a set of tuples of the form $(X,d,\lambda,\cdot)$ with $(X,d,\lambda)$ as before and
    $\bar{x}\in X^{m}$ a tuple of $m$ points of $X$ (a marking). In both cases we can assume that $\lambda$ an action by isometries of $G$ on $X$.

    %    Let $\bar{\mathcal{X}}$ ($\bar{\mathcal{Y}}$) the quotient of $\mathcal{X}$ ($\mathcal{Y}$) up to isometries which are equivariant and compatible with the markings.
    \begin{comment}
    Let $\bar{\mathcal{X}}=\mathcal{X}/\sim$ be the equivalence relation on $\mathcal{X}$ relating
    $(X,\cdot_{1})$ and $(X',\cdot_{2})$ if and only if the two spaces are equivariantly isometric, i.e.,
    if there is an isometry $\fun{\phi}{X_{1}}{X_{2}}$ which is equivariant with respect to $\cdot_{1}$ and $\cdot_{2}$. The goal is to topologize the space $\mathcal{X}/\sim$.
    \end{comment}
    
    %    Fix $\pmb{X}=[(X_{0},d_{X_{0}},\lambda_{0})]\in\mathcal{X}$ and $\pmb{Y}=[((Y_{0},d{Y_{0}}_{0}\bar{y_{0}}),\rho_{0})]\in\mathcal{Y}$ .
    
    For any $S\in\fpow{G}$ and $K_{0}\subset S$ a finite subspace of $Y_{0}$, we define $N_{K,S,\epsilon}(X)$ as the set of all $(Y,\bar{y},\cdot)\mathcal{Y}$ such that there exists a
    finite subspace $K\subseteq Y$ and a $S$-equivariant $\epsilon$-approximation $R$ between $K_{0}\cup\bar{y}$ and $K\cup\bar{y}$ such that $(y_{0,j},y_{j})\in R$ for all $1\leq j\leq m$.
    Let $N_{K,S,\epsilon}(X)$ be the analogous for $\mathcal{X}$, without the condition on the marking. Both sets are well defined in the sense that they only depend on the equivalence class up to
    equivariant isometry.
    
    %    Obviously $N_{K,S,\epsilon}(\pmb{X})$ is well defined.
    %     and that the family of all $N_{K,S,\epsilon}(\pmb{X})$ for $\epsilon>0$, $S\in\fpow{G}$ and $K\subseteq X$
    The (marked) equivariant Gromov topology on $\mathcal{X}$ ($\mathcal{Y}$)
    %    	, denoted $\mathfrak{T}(\mathcal{X})$,
    is the topology defined by taking as a neighbourhood system of $\mathcal{X}$ ($\mathcal{Y}$)
    the family of all $N_{K,S,\epsilon}(X)$, where $S$ ranges over all $S\fpow{G}$ and $K$ over all finite $K\subseteq X$. When the marking consists of a single point we use the adjective \emph{pointed} instead.
    
    Observe that in case the metric spaces involved are $\R$-trees any $\epsilon$-approximation between $K\subset X\in\mathcal{X}$ and $K'\subset X'\in\mathcal{X}$  can be extended to an $\epsilon$-approximation between the convex envelopes $conv(K)$ and $conv(K')$, which are compact. So in the definition of the Gromov Hausdorff topology we might as well have used 'compact' instead of 'finite'.
    %  This is very dependent v
    %    Make a comment to omit the word equivariant.
    Let $\Act{G}$ be the space of non-trivial minimal actions of $G$ on $\R$-trees which are either linear or irreducible, up to equivariant isometry, equipped with the equivariant Gromov topology.
    Let $\PAct{G}$ be the space of non-trivial minimal actions of $G$ on pointed real trees which are spanned by the basepoint, equipped with the \pegt.
    %    	 and $\Irr{G}$ that of those which are irreducible, both equipped with the Gromov topology.
    \begin{definition}
    	A function $\fun{l}{G}{\R_{\leq 0}}$ is called a Lyndon length function if it satisfies the following conditions:
    	\enum{i)}{
    		\item $l(1)=0$
    		\item $\forall g\in G\,\,l(g)=l(g^{-1})$
    		\item For all $g,h,k\in G$, the inequality $c(g,h)\geq \min\{c(g,k),c(h,k)\}$ holds,
    		where $c(g,h)=\frac{1}{2}(l(g)+l(h)-l(g^{-1}h))$
    	}
    \end{definition}
    We say that a point $y\in Y$ \emph{spans} $Y$ in case $Y$ is the convex closure of the orbit of $y$ by the action of $G$.
    \begin{theorem}[\cite{chiswell2001introduction}, 4.6 Chapter 2]
    	\label{from length to action}Given any length function $l$, there is an action $\lambda$ of $G$ by isometries on a $\R$-tree $Y$ and $y\in Y$  such that $l$ is equal to the length function $l^{\lambda}_{y}$. The tuple $(Y,d,\lambda,y)$ is unique up to basepoint preserving equivariant isometry.
    \end{theorem}
    The proof of \ref{from length to action} yields also the stronger result:
    \begin{proposition}
    	\label{length topology is fine}The map $\fun{\pmb{l}}{\PAct{G}}{\R_{\geq 0}^{G}}$ taking any $(X,d,\cdot,x)$ to $(l_{x}^{\cdot}(g))_{g\in G}$ is a homeomorphism between $\PAct{G}$ and its image.
    \end{proposition}
    Assume that $G$ is countable. Then in particular the pointed equivariant Gromov topology is separable. In fact a similar non-pointed result is true:
    \begin{proposition}
    	\label{tl topology is fine}The map $\fun{\pmb{tl}}{\PAct{G}}{\R_{\geq 0}^{G}}$ taking any $(X,d,\cdot)$ to $(tl^{\cdot}(g))_{g\in G}$ is a homeomorphism between the subset of $\Act{G}$ of all those actions which are either linear or irreducible and its image.
    \end{proposition}
    
    Suppose we are given a sequence $((Y_{n},d_{n},y_{n},\rho_{n}))_{n}\subset \PAct{G}$
    %    of pointed actions converging in the \pegt
    converging to some $(Y,d,y,\rho)\neq\PAct{G}$ in the equivariant Gromov topology. We say that a sequence of points $\ssq{y}{n}\in Y_{n}$ converges to a point $y\in Y$ if the sequence $(Y_{n},d_{n},\rho_{n},(x_{n},y_{n}))$ converges to $(Y,d,\rho,(x,y))$ in the marked equivariant Gromov topology. This behaves as expected:
    \begin{lemma}[see \cite{reinfeldtmakanin}]
    	In the situation above, any point of $y$ is the limit of some convergent sequence of points. If the sequence of points $(y_{n})_{n}$ is convergent,
    	then for any other sequence $\ssq{y'}{n}$ of points if $(d(y_{n},y'_{n}))_{n}$ converges to $0$ then $\ssq{y'}{n}$ is convergent and has the same limit as $\ssq{y}{n}$.
    \end{lemma}

  \subsection{Compactness, rescaling and apriori basepoints}
    
    \label{basepoints section} The last result above is enough to prove compactness of the projectivization of $\Act{G}$ in the cases we are interested in. Otherwise formulated,
    any sequence $\ssq{\rho}{n}$ of minimal nontrivial actions of our finitely generated $G$ on real trees contains a subsequence which, after rescaling converges to another non-trivial action.
    
    For the shortening argument one needs however to be able to make an apriori choice of a point $*^{\rho}\in Y$ for every $(Y,d,\rho)\in\PAct{G}$ in such a way that
    when extracting a subsequence we can guarantee the points $*^{\rho_{n}}$ converge as well.
    
    \begin{definition}
    	\label{useful length functions}Given an action $\lambda$ by isometries of a group $G$ on an $\R$-tree $Y$ and $U\subset G$ a finite set
    	we define the following functions:
    	\begin{align*}
    		\ml{y}{U}{\lambda}=\maxim{u\in U}\;l^{\lambda}_{y}(u) \\  	\sl{y}{U}{\lambda}=\serprop{u\in U}{\pl{y}{U}{}} \\ \md{U}{\lambda}=\min_{y\in Y}\ml{y}{U}{} \\ .
    	\end{align*}
    \end{definition}
    %    or a subsequence we have convergence of the corresponding pointed spaces as wellbe obtained in $\PAct{G}$ when adding the $*^{\rho_{n}}$ as basepoints.
    %        for any sequence $((Y_{n},d_{n},\rho_{n}))_{n}\subset\PAct{G}$ some subsequence $((Y_{n_{k}},d_{n_{k}},y^{\rho_{n_{k}}},\rho_{n_{k}}))_{n_{k}}$ is convergent in the projectivization of $\PAct{G}$.
    %    A more careful choice of $*^{\rho}$ results in $((Y_{n},d_{n},*^{\rho_{n}},\rho_{n}))_{n}$ projectively converging every time $(Y_{n},d_{n},*^{\rho_{n}},\rho_{n})$ does. But the previous
    %    weaker statement, which we shall now formally state, is enough for our purposes.
    %    is implied by convergence of the functor taking an action $\rho$ to $(\tl{g}{\rho})_{g\in G}$
    \begin{comment}  % 2016/04/17/16:14:13_:
    \begin{definition}
    	Let $s$ be a finite tuple of generators of a group $G$ equipped with an action $\rho$ by isometries on a real tree $Y$. We say that a point $y\in Y$ is minimally displaced by $s$ if and only if
    	$\ml{y}{s}{\rho}=\md{n}{\rho}$. From now on, given such $G,s$, we work with some choice for any such $\rho$ of some $*_{min}^{\rho,s}$.
    \end{definition}
    \end{comment}
    %    	 We can assume a fixed choice of some $y_{min}^{\rho,s}$ minimally displaced by $s$ is given for any $(G,s)$ and any isometric action $(Y,d,\rho)$ of $G$ on
    %    	a real tree. In case the action is simplicial, we can assume further more that $y_{min}^{\rho}$ comes from a simplicial vertex.
    \begin{lemma}
    	\label{compactness of actions} Let $s$ be a finite tuple of generators of a group $G$ and suppose that for each $n\in\N$ we are given a minimal non-trivial action $(Y_{n},d_{n},\rho_{n})$ of $G$ on a real tree and some $y_{n}$ for which $\ml{y_{n}}{s}{\rho_{n}}=\md{s}{\rho_{n}}=:\mu_{n}$.
    	Then some subsequence of $((Y_{n},\frac{1}{\mu_{n}}d_{n},y_{n},\rho_{n}))_{n}$ converges in the equivariant pointed Gromov topology to a minimal non-trivial pointed action $(Y,d,y,\rho)$ of $G$ on a real tree $Y$.
    \end{lemma}
    \begin{proof}
    	First of all, for any $g\in G$ the length functions $\pl{y_{n}}{g}{\nicefrac{1}{\mu_{n}}Y_{n}}$ are uniformly bounded, since in general $\pl{x}{gh}{}\leq\pl{x}{g}{}+\pl{x}{h}{}$.
    	This implies some subsequence of $(l_{y_{n_{k}}}^{\nicefrac{1}{\mu_{n_{k}}}Y_{k}})_{k\in\N}$ converges to some non-constant function $l$, which is necessarily the length function associated to some
    	%    	which
    	%    	in virtue of \ref{from length to action} and \ref{length topology is fine} the sequence of pointed actions
    	tuple $(Y,d,y,\rho)$ to which
    	$((Y_{n_{k}},\frac{1}{\mu_{n_{k}}}d_{n_{k}},y_{n_{k}},\rho_{n_{k}}))_{n_{k}}$ converges in the \pegt.
    	We claim this action is non-trivial. Indeed, suppose that $z\in Y$ is fixed by $G$ and let $(z_{n})$ a sequence of points in $Y_{n}$ converging to $z$.
    	Since $\lim_{k}\pl{z_{n_{k}}}{g}{(\nicefrac{1}{\mu_{n_{k}}}d_{n_{k}},\rho_{n_{k}})}=\pl{z}{g}{d,\rho}$ for any $g\in G$, eventually ${\pl{z_{n_{k}}}{s}{(d_{n_{k}},\rho_{n_{k}})}<\mu_{n_{k}}}$, contradicting the choice of $z_{n_{k}}$. We claim that the action is minimal. Indeed, in case it is not, let $Z$ be a proper $G$-invariant subtree of
    	$Y$. Since $y$ spans $Y$, clearly $y\nin Z$. But then, given any $g\in s$ which moves $y$, it is clear that there must be some $z\in Z$ such that $\pl{z}{g}{}<\pl{y}{g}{}$.
    	Any sequence of approximating $z$ will eventually contradict the choice of $y_{n}$.
    \end{proof}
    \begin{definition}
    	From now on, given a finitely generated $G$ we fix some finite tuple $s$ of generators. For any $\rho$ as above let $\mu^{\rho}=\md{s}{\rho}$ and fix some
    	$*_{min}^{\rho}$ for which $\ml{*_{min}^{\rho}}{s}{\rho}=\mu^{\rho}$. If $Y$ is the geometric realization of a simplicial tree, one can assume that $*_{min}^{\rho}$ corresponds to
    	a simplicial vertex.
    \end{definition}
    
    In certain situations, in which $G$ is equipped with a distinguished subgroup $H$ the point $*_{min}^{\rho}$ is not good enough.
    In order to be able to apply the shortening argument we will need to be able to make an a priori choice of points that will converge (perhaps after subsequence extraction)
    to a point in the minimal tree or in the fixed point set of $H$ in the limit tree. A small complication arises from the fact that we cannot tell in advance whether $H$ will be elliptic or not in the limiting action. In the final writing of this work it turns out the case in which $H$ is not elliptic in the limiting tree can be ignored, but we will present the other case in case it might be useful in other circumstances.
    
    Given non-empty convex subsets $X,Z$ of a real tree $Y$ and $y\in Y$, let $pr_{X}(y)$ be the point in $X$ closest to $Y$, and $Lpr_{X,Z}(y)$ the points
    in $\{pr_{X}(y),pr_{Z}(y)\}$ whose distance to $y$ is minimal. Let $Lpr_{X,Z}(y)$ be either the unique point in $Lpr_{X,Z}(y)$ or $pr_{X}(y)$ in case of a draw.
    Given $h\in G$, we let $Ax^{\rho}(h)$ stand for the set of points minimally displaced by $h$ with respect to $\rho$.
    \begin{observation}
    	\label{limiting axis projection}
    	Suppoe a group $G$ endowed with a finite set of generators $s$ and a sequence
    	$((Y_{n},d_{n},y_{n},\rho_{n}))_{n}$ of minimal pointed isometric actions of $G$ on real trees,
    	converging to some $(Y,d,y,\rho)$ in the \pegt is given. Fix some $h \in G$. For each $n$ let $A_{n}=Ax^{\rho_{n}}(h)$, and $A=Ax^{\rho}(h)$. Then $pr_{A_{n}}(y_{n})$ converges to $pr_{A}(y)$.
    	
    	If we are also given $h'\in G$ and we let $A'_{n}=Ax^{\rho_{n}}(h')$ and $A'=Ax^{\rho}(h')$, then some subsequence of $Lpr_{A_{n},A'_{n}}(y_{n})$ converges to a point in $Lpr_{A}(y)$.
    \end{observation}
    \begin{proof}
    	Let $z_{n}=pr_{A_{n}}(y)$ and $z=pr_{A}(y)$. Let $\ssq{w}{n}$ be a sequence of points converging to $w$.
    	Now, since both $\pl{w_{n}}{h}{\rho_{n}}$ and $\tl{h}{\rho_{n}}$ tend to $\tl{h}{\rho}$, the sequence of distances $d(w_{n},A_{n})=\nicefrac{\pl{w_{n}}{h}{\rho_{n}}-\tl{h}{\rho_{n}}}{2}$ converges to $0$, and thus we can assume $w_{n}\in A_{n}$.
    	
    	Now, on one hand $d(y_{n},w_{n})$ converges to $d(y,z)$ by definition. On the other $d(y_{n},z_{n})=d(y_{n},A_{n})=\nicefrac{\pl{y_{n}}{h}{\rho_{n}}-\tl{h}{\rho_{n}}}{2}$ must also converge to  $d(y,z)=d(y,A)=\nicefrac{\pl{y}{h}{\rho}-\tl{h}{\rho}}{2}$. But then $|d(y_{n},w_{n})-d(y_{n},z_{n})|$ tends to $0$ (regardless of whether $y\in A$), so $\ssq{z}{n}$ converges to $z$.
    	
    	The second statement follows almost immediately from this.
    \end{proof}
    
    \begin{lemma}
    	\label{kitten killer} Let $\rho$ be an isometric action of a group $G$ on a real tree $Y$ be given, together with $h_{1},h_{2}\in H\leq G$ and $y\in Y$.
    	Let $A_{i}=Ax^{\rho}(h_{i})$. Assume that if both $h_{1}$ and $h_{2}$ are elliptic then $A_{1}\cap A_{2}$ has at most one point. Then $Lpr_{A_{1},A_{2}}(y)$ belongs to the
    	fix point set of $H$ if $H$ acts elliptically, or to its minimal tree otherwise.
    \end{lemma}
    \begin{proof}
    	Let $w:=lpr_{A_{1},A_{2}}(y)=pr_{A_{i_{0}}}(y)$. Of course, each $A_{i}$ contains some point $a_{i}$ in $Y_{H}$.
    	If either $h_{1}$ or $h_{2}$ is hyperbolic, then $w\in Y_{H}$, since then $w$ either belongs to the axis of a hyperbolic element in $H$
    	or belongs to the path between some $a_{i}\in A_{i}$ and some axis contained in $Y_{H}$.
    	
    	If both $h_{1}$ and $h_{2}$ are elliptic, then $w$ is closest to $A_{3-i_{0}}$ among the points in $A_{i_{0}}$ and hence must belong to $[a_{1},a_{2}]$, since $A_{1}$ and $A_{2}$ intersect in at most
    	one point in that case.
    \end{proof}
    Observe that if the $A_{n}^{1}$ and $A_{n}^{2}$ in the lemma intersect in a single point $w_{n}$, this coincides with $z_{n}^{i}$.
    \begin{definition}
    	\label{basepoints} %   We will say that the sequence of actions above is unbounded if the sequence $\ssq{\mu}{n}$ contains some unbounded subsequence and that it is bounded otherwise.
    	Now, suppose we are given some finitely generated limit group $G$ and $1\neq H\lneq G$.
    	For each action $\rho$ of $G$ on some real tree $X$ we assign a point $*_{H}^{\rho}$ as follows:
    	\enum{i)}{
    		\item If $H$ is trivial, take as $*=*^{\rho}_{min}$
    		\item If $H$ is abelian, pick some $h\in H\setminus\{1\}$. For any $\rho$ as above we let $*_{H}^{\rho}=pr_{Ax^{\rho}(h)}(*_{min}^{\rho})$.
    		\item If $H$ is non-abelian, chose $h_{1},h_{2}\in H$ such that $[h_{1},h_{2}]\neq 1$. For any $\rho$ as above we let
    		$*_{H}^{\rho,s}=lpr_{Ax^{\rho}(h_{1}),Ax^{\rho}(h_{2})}(*_{min}^{\rho,s})$.
    	}
    	As before, when the tree $X$ is a geometric realization of some simplicial tree, we take the basepoints to correspond to simplicial vertices.
    	%    	Sometimes we will omit explicit mention to $s$ or $\rho$ in the notation when allowed by the context.
    \end{definition}
    %    The previous results guarantee that:
    %    \begin{comment}
    \begin{comment}
    \begin{lemma}
    	Let $(Y_{n},d_{n},y_{n},\rho_{n})$ converge in $\PAct{G}$ to some  $(Y,d,\rho,y)\in\PAct{G}$ and $\tg\neq H\leq G$. Then the sequence of points $(*_{H}^{\rho_{n}})_{n}$
    	converges to either a point in $Fix_{Y}(H)$, if $H$ is elliptic in $Y$, or a point in the minimal tree of $H$ in $Y$ otherwise.
    \end{lemma}
    \end{comment}
    %    \end{comment}
    
  \subsection{Limit actions of increasingly acylindrical sequences of actions}
    \begin{remark}
    	\label{homotety} Suppose we are given a sequence $((Y_{n},d_{n},\rho_{n}))_{n}$ of minimal isometric actions of a finitely generated group $G$ on real trees and two different sequences of rescalings $((Y_{n},\lambda_{n}d_{n},\rho_{n}))_{n}$ and $((Y_{n},\lambda'_{n},d_{n},\rho_{n}))_{n}$ ($\lambda_{n},\lambda'_{n}\in\R$) converging to minimal actions $(Y,d,\rho)$ and $(Y',d',\rho')$ respectively in the \egt. If both $Y$ and $Y'$ are non-trivial, then there is an equivariant homotety between $(Y,d)$ and $(Y',d')$.
    \end{remark}
    One says that an action by isometries of a group $G$ on a tree $Y$ is \emph{superstable} if and only if for any non-degenerate segments $I\subsetneq J$ either $Fix(I)=\tg$ or $Fix(I)=Fix(J)$.
    %    For a proof of the following theorem see for example \cite{Reinfeldt}.
    \begin{comment}
    \begin{observation}
    	\label{axis preservation} Suppose that some sequence $(Y_{n},d_{n},\rho_{n})$ of non-trivial minimal actions a group $G$ on real trees converges to some non-trivial minimal
    	$(Y,d,\rho)$ in the equivariant Gromov topology and that an element $g$ acts by translations on a line $L$ of $Y$ and lines $\epsilon_{n}$ of $Y_{n}$ for all $n$. Then any element
    	$h$ preserving $\epsilon_{n}$ for each $n$ must preserve $L$ in the limit.
    \end{observation}
    \end{comment}
    The following is a well now result (see \cite{rips1994structure} or \cite{wilton2009solutions}) that can be generalized to sequences of action on hyperbolic spaces (see \cite{perinElementary} or \cite{reinfeldtmakanin}).
    \begin{lemma} \label{properties limit}
    	Let $G$ be a torsion free finitely generated group and $((Y_{n},d_{n},\rho_{n}))_{n}$ a sequence of minimal actions by isometries of $G$ on real trees,
    	converging to some action $(Y,d,\rho)$ in the equivariant Gromov topology.
    	%    	where for some sequence $\ssq{L}{n}$ converging to $0$ the action $\rho_{n}$ is $\epsilon_{n}$-acylindrical for all $n$.
    	Assume that the stabilizer of any segment of length greater than $\geq \epsilon_{n}$ is contained in $K_{n}:=\ker \rho_{n}$, for some sequence $(\epsilon_{n})_{n}$ of positive constants converging to $0$,
    	and that for any $g\in G\setminus\{1\}$ we have $n\nin K_{n}$ for $n$ big enough. Then:
    	%   where for a certain $L>0$ the action of $G$ on $T_{n}$ is a $L$-acylindrical.
    	\enum{i)}{
    		\item  Point-wise stabilizers of non-degenerate segments of $Y$ are free abelian and each of their non-trivial elements is hyperbolic in $T_{n}$ for $n$ big enough.
    		\item  Fix point stabilizers of non-degenerate tripods are trivial.
    		\item  The action of $G$ on $Y$ is super-stable.
    		\item  Either $Y_{n}$ is linear for $n$ big enough (so that $Y$ is as well) or the kernel of $\rho$ is trivial.
    	}
    	If we assume no subgroup of $G$ acts on a line of $Y_{n}$ dihedrally, then the same is true for the action of $G$ on $Y$.
    \end{lemma}

  \subsection{Limit actions induced by a sequence of homomorphisms}
    
    \label{homomorphisms limit tree} Suppose finitely generated groups $G$, and $K$ are given, as well as an acylindrical action
    $\rho$ of $K$
    %    $K$ on a simplicial tree, which we regard as a
    on an $\R$-tree $(T,d)$.
    Given a homomorphism $\fun{f}{G}{K}$, let $\rho^{f}$ denote the action of $G$ by the pullback of $\rho$ by $f$ on its minimal tree $Y^{f}=Y_{f(G)}$ in $T$.
    %        on which $K$ acts without inversion
    % 		on whihch $K$ acts by isometries
    
    Suppose that a sequence $(f_{n})_{n}$ of homomorphisms from $G$ to $K$ with trivial limit kernel is given. For each $n\in\N$ let $Y_{n}:=Y^{f_{n}}$, $\rho_{n}=\rho^{f_{n}}\rst_{Y_{n}}$, $*_{n}=*_{min}^{\rho_{n}}\in Y$ the minimum displacement point and $\mu_{n}=\mu^{\rho_{n}}$ the corresponding maximal displacement by some fixed generating set, as in subsection \ref{basepoints section}.
    Then we know that some subsequence of rescaled pointed actions $((Y_{n_{k}},\frac{1}{\mu_{n_{k}}}d\rst_{Y_{n_{k}}}),*_{n},\rho_{n_{k}})_{k}$ will converge to an action of $G$ on a pointed real tree $(Y,d_{\infty}\,*)$.
    
    If the sequence $(\mu_{n_{k}})_{n_{k}}$ tends to infinity, we say that $(f_{n_{k}})_{n_{k}}$ is unbounded with respect to the action of $K$ on $Y$. Otherwise we say that it is \emph{bounded}.
    %    In the first case some subsequence $(f_{n_{k}})_{n_{k}}$ will be what we will refer to as geometrically convergent with respect to $\rho$.
    We refer to any (pointed) tree obtained from a subsequence of $\ssq{f}{n}$ in this way as the one as a \emph{limit tree} for $(f_{n})_{n}$ and the action $\rho$.
    %    	 When the target of $f_{n}$ is the reference model, in absence to any explicit mention to $\rho$, the action of $\F$ on its Cayley graph will be intended.
    %        limit
    In the first case, the sequence of acylindricity constants of the sequence of rescaled actions tends to $0$, so that the hypothesis of theorem \ref{properties limit} hold.
    In particular, point-wise stabilizers of non-degenerate segments are abelian and hence, by \ref{kitten killer}, given $H\leq G$,
    the sequence consisting of each of the points $*_{H}^{\rho_{n_{k}}}$ chosen as specified in \ref{basepoints} converges to a point of $Y$ which is either in the minimal tree of $H$,
    in case $H$ does not fix any point of $Y$, or it is otherwise fixed by $H$.
    For simplicity, in case the action $\rho$ to which we are referring is clear we will  $*_{H}^{\rho_{f}}$ simply by $*_{H}^{f}$.
    In the second case, the limit action is a simplicial tree.
    
    We will be mainly concerned with the particular case in which $H$ is our reference free group $\F$, and $\rho$ the action of $\F$ on its Cayley graph with respect to a fixed base.
    In this situation we might refer to $Y$ simply as a limit tree for the sequence $\ssq{f}{n}$.
    
    For $h,k\in F$ we have that $\pl{k}{h}{}=\pl{1}{h^{k}}{}$. Hence in this case $\mu_{n}$ can be expressed as $\min_{h}\max_{g\in s}|\inn{k}\circ f_{n}(g)|$, for $k\in H$.
    Since for any positive constant $N$ there are finitely many homomorphisms from $G$ to $H$ for which $\max_{g\in s}|f(g)|<N$ is bounded, any sequence $\ssq{f}{n}\subset Mor(G_{A},\F_{A})$
    contains an unbounded subsequence if and only if it contains infinitely many distinct members up to postcomposition by an inner automorphism of $\F$.
    
    \begin{definition}
    	Let $G_{1}$ and $G_{2}$ be finitely generated groups and consider sequences $\ssq{f^{i}}{n}\subset Hom(G^{i},\F)$ for $i=1,2$. We say that
    	$\ssq{f^{2}}{n}$ grows faster than $\ssq{f^{1}}{n}$ if for any finitely generated $K$ containing $G_{1},G_{2}$ and any sequence $\ssq{h}{n}\subset Hom(K,\F)$
    	such that $h_{n}\rst_{G_{i}}=\inn{c_{n}}\circ f^{i}_{n}$ for some $c_{n}\in\F$ the group $G_{1}$ fixes a point in any limit action associated to $\ssq{h}{n}$.
    	Alternatively, given a single sequence $\ssq{f}{n}\subset Hom(K,\F)$, where $G_{1},G_{2}\leq K$, we say that $\ssq{f}{n}$ makes $G_{2}$ grow faster than $G_{1}$
    	if  $(f_{n}\rst_{G_{1}})_{n}$ grows faster than $(f_{n}\rst_{G_{2}})_{n}$.
    \end{definition}
    It is perhaps not a bad idea to formally state the following (well-known) lemma.
    \begin{observation}
    	\label{growth and ellipticity} Let $G_{1}$ and $G_{2}$ be finitely generated groups and consider sequences $\ssq{f^{i}}{n}\subset Hom(G^{i},\F)$ for $i=1,2$.
    	The following are equivalent:
    	\elenco{
    		% 	 		\item For any finitely generated $K$ containing $G_{1},G_{2}$ and any sequence $\ssq{h}{n}\subset Hom(K,\F)$  \label{one}
    		%  			such that $h_{n}\rst_{G_{i}}=\inn{c_{n}}\circ f^{i}_{n}$ for some $c_{n}\in\F$ the group $G_{1}$ fixes a point in any limiting action associated to $\ssq{h}{n}$.
    		\item $\ssq{f^{2}}{n}$ grows faster than $\ssq{f^{1}}{n}$ \label{onee}
    		%  			\item Given any finite sets of generators $X_{i}$ of $G_{i}$ the sequence
    		\item There is some $g\in G_{2}$ and some tuple  $x$ of generators of $G_{1}$ such that the sequence 	$\tl{f_{n}(x_{i}x_{j})}{\lambda}$   \label{two}
    		tends to infinity with $n$ for any $x_{i}\neq x_{j}$, where $\lambda$ is the action of $\F$ on its Cayley graph.
    		\item Some $g_{2}\in G_{2}$ is not killed by $f^{2}_{n}$ for $n$ big enough and for every $g_{1}\in G_{1}$ \label{three}
    		the quotient $\frac{\tl{f^{1}_{n}(g_{1})}{\lambda}}{\tl{f^{2}_{n}(g_{2})}{\lambda}}$ tends to $0$ with $n$.
    	}
    \end{observation}
    \begin{proof}
    	%  	Clearly \ref{three} implies \ref{two}.
    	Suppose we are given a group $G$ which acts by isometries on a real tree $Y$ and a set $\mathcal{X}\subset G$, the group generated by $\mathcal{X}$ fixes a point of $Y$ if and only if
    	for each $x,y\in\mathcal{X}$ both elements $x,xy$ fix a point of $Y$. The equivalence between \ref{onee} and \ref{two} follows easily from continuity of translation length
    	with respect to the equivariant Gromov topology and its invariance by conjugation. Using continuity it is also easy to show that \ref{onee} implies \ref{three}, while the direction from
    	\ref{three} to \ref{two} is trivial.
    \end{proof}
    
  \subsection{Trees of actions}
    
    We borrow some tools from \cite{guirardel2008actions}, the term "graph of actions" instead of "tree of actions":
    \begin{definition}
    	A \emph{tree of actions} can be given as a tuple $(T^{G},Y_{(-)},p_{(-)},\cdot)$, comprising:
    	\enum{i)}{
    		\item  A simplicial $G$-tree $T^{G}$, which we call the skeleton of the tree of actions.
    		\item  A map assigning to each $v\in VT$ a real tree $Y_{v}$.
    		\item  For each $e\in ET$ a point $p_{e}\in Y_{\alpha(e)}$.
    		\item  An action of $G$ on $\bdcup{Y_{v}^{}}{v\in VT}{}$. So that for each $g\in G$ and $v\in VT$ the action of $g$ restricts to an isometry between $Y_{v}$ and $Y_{g\cdot v}$.
    	}
    	In particular the data makes each $Y_{v}$ into a $Stab(v)$-tree.
    \end{definition}
    
    Given a tree of actions as above, consider the quotient $Y$ of $\bdcup{Y_{v}^{}}{v\in VT}{}$ by the equivalence relation generated by $\{(p_{\alpha(e)},p_{\alpha(\edin(e))})\}_{e\in ET}$.
    It can be shown there is a metric structure on $Y$ making it an $\R$-tree with respect to which	the inclusion map of $Y_{v}$ into $Y$ is an isometric embeddings for every $v$. We denote the resulting real tree
    by $\dt{(S^{G},Y_{(-)},p_{(-)},\cdot)}$.
    
    \begin{definition}
    	Let $Y$ be an $\R$-tree on which a group $G$ acts by isometries. A decomposition of the action as a tree of actions comprises:
    	\elenco{
    		\item A tree of actions $(S^{G},Y_{(-)},p_{(-)},\cdot)$
    		\item An equivariant isomorphism $\phi$ between
    		$\dt{(S^{G},Y_{(-)},p_{(-)},\cdot)}$ and $Y$.
    	}
    \end{definition}
    It is easy to see that $S$ is minimal whenever $Y$ is. We will usually identify the vertex trees $Y_{v}$ with their preimage in $Y$ and by an abuse of notation refer to the decomposition as $(T^{G},Y_{(-)},p_{(-)},\cdot)$.
    \newcommand{\candtoa}[0]{normalized } %text command
    We say that a family $\{Z_{j}\}_{j\in J}$ of subtrees of $Y$ is \fundamental if $Z_{i}\cap h\cdot Z_{j}$ is non-degenerate only in case $i=j$. We say that a subtree $Z$ is \fundamental if $\{Z\}$ is.
    \enum{i)}{
    	\item If $Y,Z\in\mathcal{Y}$ intersect in more than two points, then $Y=Z$.
    	\item Every arc of $X$ is covered by finitely many members of $\mathcal{Y}$.
    }
    \begin{observation}
    	There is a transverse covering $\mathcal{Z}$ of $Y$ containing $\{Z_{j}\}_{j\in J}$, with $Z_{i}$ and $Z_{j}$ in different orbits for $i\neq j$ if and only if $\{Z_{j}\}_{j\in J}$ is fundamental. \label{coverB}
    \end{observation}
    
    Given a decomposition of an action of $G$ on $Y$ as a tree of actions $\treeact{S}{Y}{}$, the family of images in $Y$ of non-degenerate vertex trees $Y_{v}$ constitutes a transverse covering of $Y$. Viceversa: given a transverse covering $\mathcal{Y}$ of a $Y$ on which a group $G$ acts by isometries, let $V=\mathcal{Y}\dcup\mathcal{I}$, where $\mathcal{I}$ is the family of all points of $Y$ belonging to more than one member of $Y$; notice how $G$ acts on both sets. Put an edge $(p,Z)$ between $p\in\mathcal{I}$ and $Z\in\mathcal{Y}$ whenever $p\in Z$. Clearly $G$ acts without inversion in the resulting graph $S$. For each $Z\in\mathcal{Y}$ there is an associated action of the set-wise stabilizer of $\mathcal{Y}$ on $Y$. For each $p$ we can consider the trivial action of $Stab(p)$ on $\{p\}$. It is easy to see that the dual of the corresponding tree of actions is isomorphic to $Y$. We call decompositions as trees of actions of this form normalize.
    %    We will denote by $\mathcal{V}^{deg}$ the set of vertices of the skeleton associated
    \begin{comment}
    Given a decomposition of an action of $G$ on $Y$ as a tree of actions $\treeact{S}{Y}{}$, the family of images in $Y$ of non-degenerate vertex trees $Y_{v}$ constitutes a transverse covering of $Y$ and viceversa: given a transverse covering $\mathcal{Y}$ of a $Y$ on which a group $G$ acts by isometries, let $V=\mathcal{V}_{c}\dcup\mathcal{V}_{i}$, where $\mathcal{V}_{i}$
    contains a vertex $v_{p}$ for any point $p$ belonging to more than one member of $Y$ and $\mathcal{V}_{comp}$ contains a point $p_{Y}$ for each $Z\in\mathcal{Y}$. Put an edge between $v_{p}$ and $v_{Z}$ whenever $p\in Z$. The resulting graph is a tree $S$ on which onviously $G$ acts without inversion and is obviously the skeleton of a tree of actions,
    whose vertex actions are those of $\SwS{Z}=Stab(v_{Z})$ on $Z\in\mathcal{Y}$ and that of $Stab(p)$ on $\{p\}$. It is easy to see that the dual of the corresponding tree of actions is isomorphic to $Y$. We call decompositions as trees of actions of this form \candtoa.
    \end{comment}
    \begin{observation}
    	Let $\treeact{S}{Y}{}$ a normalized decomposition as a tree of actions of the action of $G$ on $Y$, as above. Then any $H\leq G$ fixing a point of $Y$ fixes a vertex of $S$.
    \end{observation}
    
    \begin{observation}
    	\label{fundamental subtrees}Suppose we are given a decomposition as a tree of actions $\treeact{S}{Y}{}$ of an action of a group $G$ on a $\R$-tree $Y$. Let $T$ a subtree of $S$ which is invariant under the action of some $H\leq G$, and such that for any $v\in VT$ such that $Y_{v}$ is not degenerate condition
    	$G\cdot v\cap T=H\cdot v$ holds and $Stab(v)\leq H$. Then $Y_{T}:=\bunion{v\in VT}{Y_{v}}{}\subset Y$ is fundamental and $Stab(T)=H$, ie., $g\cdot Y_{T}\cap Y_{T}$ is non-degenerate for $g\in G$ only if $g\in H$.
    \end{observation}
    \begin{proof}
    	Given $g\in G$ such that $g\cdot Y_{T}\cap Y_{T}$ is non-degenerate, our assumption on the action of $H$ on $S$ implies $g\cdot Y_{v}\cap h\cdot Y_{v}$ is non-degenerate for some $v\in VS$. Since $Y_{v}$ is itself fundamental, $g^{-1}h\in Stab(v)\subset Stab(Y_{v})\subset H$. But $h\in H$, so this proves the result.
    \end{proof}
    
    We note in passing the following:
    \begin{observation}
    	\label{free skeleton}Suppose that we are given a normalized decomposition as a tree of actions of an isometric action of a group $G$ on a real tree $Y$, $\treeact{S}{Y}{}$, where the
    	action on $S$ is non-trivial. If two distinct degenerate components $Y_{u},Y_{v}$ intersect in a point $p$ with
    	$Fix_{Stab(v)}(p)=\tg$ then $S$ contains trivially stabilized edges and, therefore $G$ is freely decomposable relative to $\setof{Stab(v)}{v\in VS}$.
    \end{observation}
    \begin{corollary}
    	\label{free skeleton 2} Given $G$, $Y$ and $\treeact{S}{Y}{}$ as above, suppose that $G$ is freely indecomposable relative to $\setof{Stab(v)}{v\in VS}$. Suppose that some non-degenerate component $Y_{v}$, with stabilizer $G_{v}$ and $p\in Y_{v}$ such that $\Sb{p}{G_{v}}=\tg$. Then $\Sb{p}{G}=\tg$ and $G\cdot p\cap Y_{v}=G_{v}\cdot p$.
    \end{corollary}
    %    Let us briefly discuss the variant of the construction employed in the simplicial setting, where one is mostly interested in the complementary
    %    point of view. That is, rather than decomposing an action into smaller pieces more often than not we wish to combine global and local actions.
    %     (dually, for some vertex group $\Delta_{v}$ a graph of groups decomposition of $\Delta_{v}$ relative to
    %    $\Delta_{e}$ for each $e\in \abs{\Delta}$ originating at $v$)
    
    One property of decompositions of the type above worth mentioning is the following:
    \begin{lemma}
    	\label{factoring lines} Let $G$ be a $CSA$ group equipped with an action by isometries on a real tree $Y$ in which segment stabilizers are abelian.
    	%    	admitting a decompositon as a tree of actions $\treeact{S}{Y}{}$, where each component is either	of simplicial, axial or Seifert type.
    	Let $Z$ be a line of $Y$ and suppose that $M\leq G$ consist of those elements preserving $Z$ set-wise, with its orientation.
    	Let also $K\leq M$ be the kernel of the action of $M$ on $Z$. If $M\neq K$, then  $M$ is a maximal abelian subgroup of $G$. If moreover $G$ is torsion-free, then $M$ coincides with the set-wise stabilizer of $Z$. If $K\neq\tg$ as well, then $Z$ is fundamental and $G$ decomposes as an amalgamated product of the form $G'\frp_{K}M$,
    	where every subgroup of $G$ elliptic in $Y$ can be conjugated into $G'$.
    	%    	, where $K$ is the kernel of the action of $M$ on $Z$.
    \end{lemma}
    \begin{proof}
    	First of all we claim that $M$ is abelian. Indeed, first of all clearly $M/K$ is free abelian. On the other hand $K$ is abelian by our assumption on the action on $Y$.
    	If $K\neq\tg$, since also $K\normal{}M$, the $CSA$ property implies that $M$ is abelian.
    	
    	Now, $M$ is maximal as well. Indeed, pick any $g\in G$ commuting with $M$, let $Z'=a\cdot Z$ and suppose that $Z'\neq Z$.
    	Any $m\in M\setminus K$ has to act as a proper translation on both $Z$ and $Z'$, which is impossible in any action by isometries on a real tree. So $a\cdot Z=Z$.
    	It is clear as well that $g$ has to preserve the orientation of $Z$, since otherwise $m^{g}=m^{-1}$ so $g\in M$.
    	%    	\footnote{The same argument can be applied to the case in which instead of $K\neq M$ we assume that $K\neq\tg$ and that the stabilizers non-degenerate tripods of $Y$ are trivial.}
    	Finally, assume that $G$ is torsion free and $s\in G$ stabilize $Z$ set-wise. Then $s^{2}\in M$. Since $s^{2}\neq 1$, commutative transitivity implies that
    	$s\in M$.
    	
    	We now add the assumption that $K\neq\tg$. We claim that $Z$ is fundamental in $Y$.
    	Indeed, pick any $k\in K\setminus\{1\}$. Given any $g\in G$ for which $J:=g^{-1}\cdot Z\cap Z$,
    	$k^{g}$ fixes $J$ which implies that $k^{g}\in M$, since segment stabilizers are abelian. The $CSA$ property implies that $g\in M$ as well, so that $g\cdot Z=Z$.
    	This implies that there is a decomposition of the action on $Y$ as a normalized tree of actions $\treeact{S}{Y}{}$, where $Z=Y_{v}$ for some $v\in VS$. And observe that the edge group
    	associated to any edge $e$ incident on $v$ is precisely $K$. This implies that $G$ admits the required decomposition.
    	%    	Collapse all the edges of $S$ not incident on a translate of $v$. And then perform a fold
    	Any subgroup of $G$ elliptic in $Y$ is elliptic in $Y$. If it is contained in a conjugate of $M$, then it must be contained in a conjugate of $K$. This proves the last claim.
    \end{proof}
  \subsection{Rips' decomposition}
    
    The combinatorial core of the following result goes back to the work of Makanin and Razborov on solutions to equations in the free group. Rips' contribution was the realization that their ideas could be read geometrically and hold in a much wider context. The initial goal was to study free actions on $\R$ trees but it turned out it could be applied to a larger class of actions (see \cite{bestvina1995stable}).  We now state the version found as theorem 4.1 of \cite{guirardel2008actions}, with some extra detail borrowed from its proof. See theorem 5.1 in that same paper for the most general result of this kind presented there.
    \newcommand{\ZZ}[0]{Z}
    \newcommand{\Gv}[0]{Stab(Z)}
    \newcommand{\Nv}[0]{K}
    \newcommand{\SwS}[1]{Stab(#1)}  %% set-wise stabilizer
    \newcommand{\PwS}[1]{Stab(#1)}  %% set-wise stabilizer
    \begin{theorem}\label{Rips decomposition}
    	Any minimal superstable action by isometries of a group $G$ on an $\R$-tree $Y$, where $G$ is finitely generated relative to the family of its elliptic subgroups
    	admits a decomposition as a tree of actions $\treeact{S}{Y}{}$ where either:
    	\enum{a)}{
    		\item There is $v\in V$ such that $Y_{v}$ contains an infinite tripod and some $e\in S$ incident to $v$ such that $Stab(e)$ coincides with the kernel $\Nv$ of the vertex action on $Y_{v}$.
    		In particular there is a splitting of $G$ over $\Nv$ relative to the stabilizer of any point in $Y$.\label{exotic}
    		\item Each vertex action is either:
    		\enum{i)}{
    			\item Simplicial, i.e., isomorphic to the geometric realization of a simplicial tree.
    			%    		\item A segment whose interior does not intersect any other vertex action. $(Y\in\mathcal{Y})$
    			\item \label{Seifert} Of Seifert type. If we denote by $\Nv$ the kernel of the action of $\Gv$ on $Y_{v}$, then $\Gv/\Nv$ is dual to the fundamental group of an orbifold with boundary and the action on $Y_{v}$ dual to an arational measured foliation on it.
    			\item Of axial type: $\ZZ$ is a line and the image of $\Gv$ in $Isom(\ZZ)$ is a non-discrete group of isometries of $\ZZ$.  \label{Axial}
    		}
    	}
    \end{theorem}
    
    One can assume that the decomposition as a tree of actions above is normalized. ALso, after further decomposing any simplicial action, that in the second alternative of the theorem all simplicial vertex trees are either a point or a non-trivially stabilized segment with its endpoints are the only attaching points $p_{e}$. The first alternative in the theorem can be understood in terms of the existence of exotic components. See for example \cite{bestvina1995stable}. If $G$ is torsion free, then in case (\ref{Seifert}) the orbifold is actually a surface and the action of Seifert type as defined earlier.
    %    Of course, if $Y$ is not linear and the point-wise stabilizers of non-degenerate tripods are trivial, then $\Nv$ above is trivial.
    \newcommand{\absorbing}[0]{absorbing } %text command
    %    	Given an action of a group $G$ on a $\R$-tree $Y$ by isometries, we say that a finite family $\{Z_{j}\}_{j\in J}$ of subtrees of $Y$ is \absorbing if and only if
    %    The following lemma will be useful.
    \begin{comment}
    Let $S^{*}$ be the tree obtained by collapsing all the edges of $S$ corresponding to vertices in $\mathcal{O}$. Clearly $H$ fixes a vertex of $S^{*}$ and in virtue of the above all
    its edges are trivially stabilized. This gives right away a decomposition as above, where $F$ is generated by Bass-Serre elements of a presentation of $S^{*}$.
    Of course, if $Y_{v}$ contains degenerate components outside of $\mathcal{O}$, then there are edges which survive the collapse, so the decomposition above is not trivial.
    Any
    \end{comment}
    \begin{comment}
    Any pair of distinct points $p',p\in Y_{v}$ such that $Stab(v)\cdot p\neq Stab(v)\cdot p'$ but $G\cdot p=G\cdot p'$ and $\{p,p'\}\nsubset\mathcal{O}$
    is associated to a couple of edges $e,e'$ originating at $v$, such that $\omega(e)$ and $\omega(e')$ are mutual translates.
    but which themselves lie in different orbits by the action of $G$. This corresponds to a loop (of length 2) in $G\backslash S$. Since at least one of $e,e'$ survives the collapse
    the loops survives in $G\backslash S^{*}$, implying that $F\neq\tg$.
    \end{comment}
    
    The result above can be slightly refined:
    \begin{corollary}
    	\label{Rips refined} Consider an isometric action of a finitely generated group $G$ on a $\R$-tree $Y$, which is superstable and has trivial tripod stabilizers
    	%    	Let $Z_{1},Z_{2},\cdots Z_{m}$ be sub-trees of $Y$ with the property that for any $g\in G$ and $1\leq i,j\leq m$ if $Z_{i}\cap g\cdot Z_{j}$ is a non-degenerate segment then $i=j$ and $g\cdot Z_{i}=Z_{i}$. Let $H_{i}=\SwS{Z_{i}}$. Then either:
    	and let $\{Z_{j}\}_{j\in J}$ be a \fundamental family of subtrees of $Y$. Then either:
    	\enum{i)}{
    		\item $G$ is freely decomposable relative to $\{Stab(Z_{1}),Stab(Z_{2}),\cdots Stab(Z_{m})\}\cup\mathcal{E}$ where $\mathcal{E}$ is the family of all the subgroups of $G$ elliptic in $Y$.
    		\item The action decomposes as a tree of actions $\treeact{S}{Y}{}$, where $Y_{v_{j}}=Z_{j}$ for some $v_{j}\in VS$ and
    		for any $w\in VS\setminus (G\cdot\setof{v_{j}}{j\in J})$ the vertex action on $Y_{v}$ is either simplicial with non-trivially stabilized edges, of axial type or of Seifert type.
    		Given a non-degenerate component $Y_{v}$, all its points of intersection with other non-degenerate components have non-trivial stabilizer in $Stab(v)$.
    	}
    \end{corollary}
    \begin{proof}
    	%    	Let $\mathcal{X}$ consist of the closures of all connected components of $Y\backslash\bunion{g\cdot Z_{j}}{g\in G\\1\leq j\leq m}{}$. Let $\mathcal{Y}=\mathcal{X}\cup\{g\cdot Z\}_{g\in G}$. Our hypothesis on the family $\{Z_{j}\}_{j=1}^{m}$ implies that this is a transverse covering of $Y$.
    	Let $\treeact{S}{Y}{}$ a graph of actions, where each $Z_{j}=Y_{v_{j}}$ for some $v_{j}$ in a different orbit.
    	For all non-trivial $v\in VS\setminus G\cdot\{Z_{j}\}_{j\in J}$ we can apply \ref{Rips decomposition} to $Y_{v}$, which yields a decomposition of the action of $Stab(v)$ as a normalized tree of actions $\treeact{S^{v}}{Y^{v}}{}$, with associated transverse covering $\mathcal{Y}_{v}$.
    	Denote by $Y'_{v}$ the minimal tree of $Y_{v}$.	If one of $Y'_{v}$ falls into first alternative of theorem  \ref{Rips decomposition}, then
    	$Stab(v)$ splits over a tripod stabilizer relative to its elliptic subgroups; in particular, relative to the stabilizer of edges of $S$ incident to $v$. Hence this free splitting extends to a free splitting of $G$ relative to $\{H_{j}\}_{j=1}^{m}$.
    	If the second alternative, observe that all the components of $Y_{v}\setminus Y'_{v}$ are segments. Indeed, a segment with and endpoint in $Y'_{v}$ and another in $\partial (Y\setminus Y_{v})$ cannot be trivially stabilized, since otherwise $S$ would contain trivially stabilized edges. Superstability tells us that its point-wise stabilizer coincides with that any of its non-degenerate initial subsegments. The fact non-degenerate tripods have trivial stabilizer then implies any two such segments can only intersect in their initial points.
    	
    	We conclude that $Y_{v}$ itself admits a transverse covering $\mathcal{S}_{v}$ as in the previous lemma.
    	The covering $\{Y_{g\cdot v_{j}}\}_{\substack{g\in G\\j\in J}}\cup\bunion{\mathcal{Y}_{v}}{v\in VS\setminus (G\cdot\{Z_{j}\}_{j\in J}}{})$ of $Y$ satisfies the required conditions. The final claim follows from observation \ref{free skeleton}.
    \end{proof}
    %    \end{comment}

%( chosing the basepoint (conjugating element) in the CG of \F  also in H seems to be completely innecessary; is it?)

\section{The shortening argument}
  
  The following theorem contains the geometric engine of the shortening argument. Except for the constraint on the modular group, the slightly looser hypothesis and the additional clause regarding the limit of the shortening argument. This is usually stated in the more general context of actions on hyperbolic spaces (see \cite{rips1994structure} or \cite{perin2008elementary}). For the sake of completeness we provide an almost self contained exposition of the more restrictive version used in this work. The original result could have been enough for most of the applications in this thesis, had we used a different approach but in view of future work, we have chosen to state the result involving a finite index subgroup of the modular group.
  
  \begin{proposition}
  	\label{general shortening}
  	Let $G$ be a group and $(Y_{n},d_{n},y_{n},\lambda_{n})$ a sequence of minimal pointed actions by isometries of $L$ on real trees, converging to some faithful pointed action $(Y,d,y,\lambda)$ in the \pegt, such that any segment stabilizer is abelian.
  	Let $U$ be a finite tuple of elements from $G$.
  	Suppose that $\lambda$ admits a decomposition as a graph of actions $\treeact{S}{Y}{}$
  	and $\mathcal{O}\subset VS$ an equivariant set such that for each $v\in\mathcal{O}$ the vertex action of $Stab(v)$ on $Y_{v}$ is either simplicial with non-trivially stabilized edges or of Seifert or axial type.
  	%  	The following trivial general remark won't is not strictly needed, but possibly worth doing.
  	
  	Let  $S'$ be the tree obtained by collapsing all those edges of $S$ none of whose ends belong to $\mathcal{O}$ and blowing up into $Y_{v}$ every $v\in VS$ for which $Y_{v}$ is simplicial.
  	
  	Regard $S$ as a \gat by declaring of surface type all vertices associated to actions of Seifert type, of abelian type those associated with mixing components and regarding all the rest as rigid.
  	Suppose that a normal finite index subgroup $\PP\trianglelefteq Mod(\Delta)$ is given.
  	
  	Then there is $\tau\in\PP$ such that for infinitely many $n\in\N$:
  	\enum{i)}{
  		\item $\tau(g)=g$ for any $g\in G$ such that $[y,s\cdot y]$ does not intersect
  		$Y_{v}$ for $v\in \mathcal{O}$ in a non-degenerate segment.
  		\item $\pl{y_{n}}{g}{\lambda_{n}\circ\tau}<\pl{y_{n}}{g}{\lambda_{n}}$ for any $g\in U$ for which $[y,g\cdot y]$ intersects $Y_{v}$ in a non-degenerate segment for some $v\in \mathcal{O}$.
  	}
  	%	In all the vertices in $\mathcal{O}$ are associated with actions of simplicial type, then the second point holds for $g$
  	%	Moreover, for some $\tau\in\PP$
  	%	and any function $\pmb{w}$ from the family of simplicial edges to $(0,1]$ which is equivariant
  	%	the $\sigma$ can be chosen in such a way that the sequence of actions $\lambda_{n}\circ\sigma$ converges to
  	%	the action $\lambda\circ\tau$ for some $\tau\in\PP$.
  	%	the action obtained from $\lambda\circ\tau$ by replacing the length $l$ assigned to any simplicial edge $e$ by $\pmb{w}$.
  	%  	have to make the comment in the end that the presentation you chose is irrelevant
  \end{proposition}
  \newcommand{\fai}[1]{#1^{-}}   %%% false inverse (in the transversal)
  \begin{proof}
  	We can assume that the decomposition is normalized, so that $VS=V_{0}\dcup V_{1}$, with $Y_{v}$ a single point in case $v\in V_{0}$ and a non-degenerate tree in case $v\in V_{1}$, as well that
  	any attaching points on an action of simplicial type correspond to simplicial vertices.
  	%	$x\in Y_{v_{0}}$ for some $v_{0}\in S\setminus\mathcal{O}$
  	For convenience, we will assume that $x\in Y_{v_{0}}$ for some $v_{0}\in S\setminus\mathcal{O}$ (the other case does not represent any fundamental complication).
  	We will start with some preliminary considerations, valid for any decomposition of $Y$ as a graph of actions.
  	Fix some presentation $P=\BS{W_{0}}{W}{t}$ of the action of $G$ on $S$ for which $v_{0}\in Z$.
  	Let $g$ be any element of $g$; it has a loop normal form $\looprep{g}{e}{m}$ with respect to this presentation, as in \ref{loop representation}, where $g_{i}\in Stab(v_{i})$ and
  	$v_{0},e_{0},v_{1}\cdots e_{m-1},v_{m}\subset W_{0}\cup W$ projects to a closed loop in $G\backslash S$.
  	For each $0\leq j\leq m-1$ we let $f_{j}$ be the unique edge in $W$ whose inverse is in the orbit of $e_{j}$.
  	Suppose we are given $v\in V_{1}\cap W_{0}$ and  a natural extension $\bar{\sigma}=\mvga{\sigma}{c}{v}$ of some $\sigma\in Aut(Stab(v))$ to $G$. with $c_{v_{0}}=1$ (recall that $\bar{\sigma}$ restricts to conjugation by $c_{u}$ on $Stab(v)$ for $u\in W_{0}\setminus \{v\}$.
  	Then the element $\bar{\sigma}(g)$ has a normal form $\looprep{g^{\bar{\sigma}}}{e}{m}$, where
  	\elenco{
  		\item $g^{\bar{\sigma}}_{i}=g_{i}$ in case $v_{i}\neq v$
  		\item $g^{\bar{\sigma}}_{i}=c_{e_{i-1}}\sigma(g_{i})c_{e_{i}}^{-1}$ in case $v_{i}=v$
  	}
  	Recall that a vertex $v\in\mathcal{O}\cap W_{0}$ appears in the loop representation of $g$ if and only if $[g\cdot v_{0},v_{0}]$ intersects the orbit of $w$.
  	Since we assume our decomposition $\treeact{S}{Y}{}$ to be normalized, this happens precisely when $[y,g\cdot y]$ intersects some translate of $Y_{v_{0}}$ in a non-degenerate segment.
  	Let $N(v)$ be the subgroup of all those $g\in G$ for which that is not the case. Clearly $\bar{\sigma}$ fixes $N(v)$.
  	Likewise, observe that the axis of $g$ intersects the orbit of $v$  in $S$ if and only if
  	its axis in $Y$ intersects some translate of $Y_{v}$ in a non-degenerate segment.
  	The automorphism $\bar{\sigma}$ restricts to an inner automorphism of $G$
  	on any subgroup $H\leq G$ which is either elliptic in $S$ or whose minimal tree does not contain the vertex $v$. By the previous assertion, the latter is equivalent to the fact that
  	none of the axis of hyperbolic elements in $H$ intersects $Y_{v_{0}}$ in a non-degenerate segment.
  	
  	We define:
  	\begin{align*}
  		\pl{y}{g}{v}=|[y,g\cdot y]\cap\bunion{u\in[v]}{Y_{u}}{}|=\serprop{0\leq j\leq m\\v_{j}=v}{|J_{j}|}
  	\end{align*}
  	\begin{comment}
  	\begin{align*}
  		\pl{y}{g}{\mathcal{Q}}=|[y,g\cdot y]\cap\bunion{v\in\mathcal{Q}}{Y_{v}}{}|=\serprop{0\leq j\leq m\\v_{j}\in \mathcal{Q}}{|J_{j}|}
  	\end{align*}
  	In general, given $\mathcal{Q}\subset VS$, let $N(\mathcal{Q})$ be the set of all those $g\in G$ for which $[y,g\cdot y]$ intersects some $Y_{v}$, for $v\in\mathcal{Q}$ in a non-degenerate subsegment and
  	\begin{align*}
  		\pl{y}{g}{\mathcal{Q}}=|[y,g\cdot y]\cap\bunion{v\in\mathcal{Q}}{Y_{v}}{}|=\serprop{0\leq j\leq m\\v_{j}\in \mathcal{Q}}{|J_{j}|}
  	\end{align*}
  	\end{comment}
  	More precisely, the segment $[y,g\cdot y]$ can be decomposed into mutually non-overlapping sub-segments of the form	$I_{j}=h_{j}\cdot J_{j}$ for $0\leq j\leq m$, where $J_{j}\subset Y_{v_{j}}$ and $h_{j}$ are given by:
  	\elenco{
  		\item  $h_{0}$, $\,\,J_{0}=[y,p_{e_{0}}]$
  		\item  $h_{m}=g$, $\,\,J_{m}=[p_{f_{m-1}},x]$
  		\item $h_{j}=g_{0}t_{e_{0}}g_{1}\cdots t_{e_{j-1}}$ and $J_{j}=[p_{f_{j-1}},g_{j}\cdot p_{e_{j}}]$ for any $1\leq j\leq m-1$.
  	}
  	\begin{comment}
  	\begin{align*}
  		h_{0}&=1\text{, $h_{m}=g$,\,\, $J_{0}=[y,p_{e_{0}}]$,\,\, } J_{m}=[p_{f_{m-1}},x] \\
  		h_{j}&=g_{0}t_{e_{0}}g_{1}\cdots t_{e_{j-1}}\text{ and }J_{j}=[p_{f_{j-1}},g_{j}\cdot p_{e_{j}}] \text{ for any $1\leq j\leq m-1$.}
  	\end{align*}
  	\end{comment}
  	Let $\bar{\sigma}=\mvga{\sigma}{c}{v}$ be a natural extension of some $\sigma\in Aut(Stab(v))$ to $G$ with $c_{v_{0}}=1$, as before.
  	Let $h^{\bar{\sigma}}_{j}$ be defined in the same way as $h_{j}$, but with 	$g^{\bar{\sigma}}_{j}$ in place of $g_{j}$. Exactly as before, $[y,\bar{\sigma}(g)\cdot y]$ is
  	the union of the non-overlapping segments $h^{\bar{\sigma}}_{j}\cdot J^{\bar{\sigma}}_{j}$, where $J^{\bar{\sigma}}_{j}=[\qe{j-1},g^{\bar{\sigma}}_{j}\cdot \pe{j}]$ for $1\leq j\leq m-1$.
  	In particular, for $\bar{\sigma}$ as above  $\pl{y}{\bar{\sigma}(g)}{w}=\pl{y}{\bar{\sigma}(g)}{w}$ for $v\neq w\in W_{0}$.
  	
  	The proof of our proposition reduces to the following two lemmas:
  	\begin{lemma}[Mixing case]
  		\label{mixing shortening} Let $v\in W_{0}$ be associated to an action of Seifert or axial type, $U$ a finite subset of $G$ and $\epsilon>0$. Then there is a natural extension $\bar{\sigma}=\svga{\sigma}{c}\in\PP$ of an automorphism $\sigma$ of $Stab(v)$ such that for all $g\in U\setminus N(v)$ we have:
  		\begin{align*}
  			\pl{y}{\bar{\sigma}(g)}{[v]}<\epsilon
  		\end{align*}
  	\end{lemma}
  	%	Up to further splitting simplicial actions, one can assume each of them is a non-trivially stabilized segment $I=[p,q]$, whose endpoints are the only attaching points (if one of them were not
  	%	an attaching point then of course $Y$ would not be minimal).
  	As we know, we can assume any simplicial component is an edge containing no attaching points in its interior.
  	Given any such $I$, one consider the simplicial tree $T_{I}$ which has a vertex $v_{C}$ in correspondence to each connected component of 	$Y\setminus\bunion{g\in G}{g\cdot I}{}$ and a simplicial edge between $v_{C_{1}}$ and $v_{C_{2}}$ whenever $C_{1}$ and $C_{2}$  intersect some common translate of $I$.
  	%	and that in correspondence to each translate $J$ of $I$ has a simplicial edge between the two
  	\begin{lemma}[Simplicial case]
  		\label{edge shortening} Let $v\in W_{0}$ and suppose $Y_{v}$ is a non-trivially stabilized edge $[p,q]$ and $e$ the corresponding simplicial edge in $T^{[p,q]}$. Then there is a Dehn twist $\tau_{e,c}\in\PP$  such that for some diverging sequence $\ssq{k}{n}\subset\N$, and all $g\in G\cap L(v)$ we have:
  		\begin{align*}
  			\pl{y_{n_{k}}}{\tau(g)}{}<\pl{y_{n_{k}}}{g}{}
  		\end{align*}
  		\begin{comment}
  		\begin{align*}
  			\pl{y_{n_{k}}}{\tau(g)}{}<\pl{y_{n_{k}}}{g}{}
  		\end{align*}
  		\end{comment}
  		for $k$ big enough (the particular bound dependent on $g$ )
  	\end{lemma}
  	\begin{comment}  % 2016/05/10/23:56:41_: edge shortening
  	\begin{lemma}[Simplicial case]
  		\label{edge shortening} Let $v\in W_{0}$ be associated to an action of simplicial type. Then there is a natural extension $\svga{\sigma}{c}\in\PP$ of
  		an automorphism $\sigma$ of $Stab(v)$ such that some diverging sequence $\ssq{k}{n}\subset\N$, and all $g\in G\cap L(v)$ we have:
  		\begin{align*}
  			\pl{y_{n_{k}}}{\bar{\sigma}(g)}{}<\pl{y_{n_{k}}}{g}{}
  		\end{align*}
  		for $k$ big enough (the particular bound depending on $g$ )
  	\end{lemma}
  	\end{comment}
  	Let us see how to prove the proposition using the lemmas above. Denote by $M$ be the set of all those $v\in W_{0}$ associated to an axion of Seifert or axial type. For each $v\in M$ let $\bar{\sigma}^{v}$ be the automorphism obtained by applying the first lemma to $\epsilon_{v}=\mn{g\in U}\,\pl{y}{g}{v}$. We know that $\bar{\sigma}^{v}(g)=g$ for $g\in N(v)$, while in case of $g\in U\setminus N(v)$ we have
  	\begin{align*}
  		\pl{y}{\bar{\sigma}^{v}(g)}{}=\pl{y}{\bar{\sigma}^{v}(g)}{[v]}+\pl{y}{\bar{\sigma}^{v}(g)}{[v]^{c}}=\\
  		=\pl{y}{\bar{\sigma}^{v}(g)}{[v]}+\pl{y}{g}{[v]^{c}}<\pl{y}{g}{[v]}+\pl{y}{g}{[v]^{c}}=\pl{y}{g}{}
  	\end{align*}
  	If we let $\theta$ be the product of all the $\bar{\sigma}^{v}$ as we let $v$ range among all vertices in $W_{0}$,
  	(the order is irrelevant)
  	then for any $v\in U\setminus\bunion{v\in M}{N(v)}{}$ we have  $\pl{y}{g}{}<\pl{y}{\theta(g)}{}$.
  	Which implies that $\pl{y_{n}}{\tau(g)}{\lambda_{n}}<\pl{y_{n}}{\tau(g)}{\lambda_{n}}$ for $n$ large enough, since of course $(Y_{n},d_{n},y_{n},\lambda_{n})$ has to converge to $\lambda\circ\sigma$ in the \pegt. The proof can be completed by applying the second lemma to each of the orbits of simplicial edges in $Y$, starting with the sequence of actions $\lambda_{n}\circ\sigma$. What remains is dedicated to the proof of the aforementioned lemmas. For simplicity, denote $N_{v}\cap U$ by $U_{v}$.
  	
    	\subsubsection{Mixing cases}
      	Let now $v\in W_{0}$ be a vertex associated with an action of Seifert or axial type. In view of the discussion above, it will be enough for us to find, given $\epsilon>0$ and any finite $H\subset Stab(v)$, some extension $\svga{\sigma}{c}\in\PP$ such that for all $h\in H$:
      	\begin{equation*}
      		d(c_{e}^{-1}\cdot p_{e},\sigma(h) c_{f}^{-1}\cdot p_{f})< \epsilon
      	\end{equation*}
      	%      	For all $e,f\in W$ originating at $v$ and $h\in V$. If we let $\sigma_{n}=\sigma$ and 	$c_{e,n}=c_{e}$ for all $n$ inequality \ref{atomic shortening strong} then follows by convergence.
      	In general, given an action by isometries of a group $H$ on a $\R$-tree, and $A\subset Aut(H)$, let us say that $A$ has the shortening property if for any $V\subset H$, any $z\in G$ and any $\epsilon>0$ some $\sigma\in\mathcal{A}$ exists for which $\pl{y}{\sigma(g)}{}<0$ for all $g\in H$.
      	Ww claim that if $H$ has the shortening property and $K\leq H$ has finite index in $H$ then $K$ has the shortening property as well. Let $\mathcal{R}$ be a  finite set of representative of $H/K$. Given $\epsilon$, $V$ and $z$ as above, chose $\theta\in Mod(\Sigma)$ such that $\pl{y}{\theta(g)}{}<\epsilon$ for all $h\in\bunion{\tau\in\mathcal{R}}{\tau^{-1}(v)}{}$.
      	Then for some $\tau\in\mathcal{R}$ the composition $\sigma=\theta\circ\tau^{-1}\in B$ witnesses shortening for $(\epsilon,V)$. We set $z=p_{e^{*}}$.
      	Denote by $Mod(Stab(v))$ the group of automorphisms of $Stab(v)$  which coincide with an inner automorphism on each of the stabilizers of an edge incident to $e$.
      	
      	\paragraph{Seifert case}
        	Here $Stab(v)\cong\pi_{1}(\Sigma)$, where $\Sigma$ is a compact surface with boundary. For a proof of the following result see
        	\cite[Lemma 5.13]{perin2008elementary} or \cite[p.347-348]{rips1994structure}
        	\begin{proposition}
        		In the situation above, given $v\in W_{0}$ associated to an action of surface type $Mod(Stab(v))\leq Aut(Stab(v))$  has the shortening property.
        		\begin{comment}
        		\enum{i)}{
        			\item The only points of $Y_{v}$ with non-trivial stabilizer are precisely $\{p_{e}\}_{\alpha(e)=v}$ (in which case the stabilizer corresponds to some boundary component of $\Sigma$).\label{non trivially stabilized}  				%	        		\item Orbits of points of $Y_{v}$ are dense. \label{dense orbits}
        			\item   \label{surface shortening}
        		}
        		\end{comment}
        	\end{proposition}
        	Let $M_{\PP}$ be the subgroup consisting of all those elements of $Mod(\Sigma)$ extending to $\PP$; this clearly has finite index in $Mod(\Sigma)$, hence has the shortening property as well.
        	Let $\mathcal{E}=\setof{e\in W}{e\neq E^{*},\alpha(e)=v}$ and for each $e\in\mathcal{E}$ let $s_{e}$ generate $Stab(e)$.
        	
        	Let $V=\{s_{e}\}_{\substack{\alpha(e)=v\\e\in W}}\cup U_{v}$.
        	Pick $\sigma\in M_{H}$ such that $\pl{x}{h}{}<\frac{\epsilon}{6}$ for all $h\in V$.
        	If it restricts to conjugation by some $c_{e}$ on each $\subg{s_{e}}$, then $\sigma(\subg{s_{e}})$ stabilizes $c_{e}^{-1}p_{e}$ and no other point.
        	Hence $d(z,c_{e}^{-1}\cdot p_{e})=\pl{z}{\sigma(s_{e})}{}<\frac{\epsilon}{5}$. Given any other $f\in\mathcal{E}$ we get:
        	\begin{align*}
        		d(h c_{f}^{-1}\cdot p_{f},z)\leq
        		d(hc_{f}^{-1}\cdot p_{f},c_{f}^{-1}\cdot p_{f})+d(c_{f}^{-1}\cdot p_{f},z)\leq\\
        		\leq (\pl{z}{h}{}+2d(c_{f}^{-1} p_{f},z))+\pl{z}{\sigma(s_{f})}{}\leq\frac{4\epsilon}{6}% d(h\cdot z,z)
        	\end{align*}
        	It follows that $d(c_{e}^{-1}\cdot p_{e},hc_{f}^{-1}\cdot p_{f})\leq\frac{5\epsilon}{6}<\epsilon$ and we are done.
        	
      	\paragraph{Axial case}
        	\begin{fact}
        		The action of  $Stab(v)$ on $Y_{v}$ satisfies the following properties:
        		\enum{i)}{
        			\item Orbits are dense. \label{dense orbits}
        			\item  $Mod(Stab(v))$ has the shortening property. \label{axial shortening}
        		}
        	\end{fact}
        	First of all, we claim it was possible for us to chose our transversal $W$, containing $v_{0}$, for the action of $G$ on $S$ in such a way that $d(p_{e},z)<\frac{\epsilon}{3}$ for any $e\in W$ originating at $v$.
        	Indeed, if $W$ already satisfies the property for all $v'$ in the path from $v_{0}$ to $v$ associated to an axial type action, let $W^{1}, W^{2},\cdots W^{m}$ consist of the connected components of $W\setminus\{v\}$ which do not contain $v_{0}$. Using (\ref{dense orbits}), we can guarantee that the property holds for $v$, by replacing each of the $W^{i}$ by its translate by an appropriate element of $Stab(v)$, an operation which preserves the condition on those vertex on which it has already been arranged for.
        	%        	In the section about trees: use $(Y_{0},Y)$ notation and define $Mod(Stab(w))$
        	Let $M_{\PP}$ be the subgroup of $\sigma\in\vma{Stab(v)}$ such $\bar{\sigma}_{\bar{1}}\in \PP$. Being of finite index in $Stab(v)$, $M_{\PP}$ has the shortening property, so $\sigma\in M_{\PP}$ exists such that $\pl{z}{\sigma(h)}{}<\frac{\epsilon}{3}$ for any $h\in U_{v}$. Now, for any $e,f\in\mathcal{E}$ we have $c_{e}=c_{f}=1$ while
        	\begin{align*}
        		d(p_{e},\sigma(h)p_{f})\leq d(p_{e},z)+d(p_{f},z)+\tl{\sigma(h)}{}<\epsilon
        	\end{align*}
        	
        	%        IN THE\,\SECTION\,ON\,TREeS\,DEFINE\,CONN COMPONENTS\,FOR\,THINGS\,THAT\,ARE\,NOT\,SUBTReES
    	\subsubsection{Simplicial case}
      	
      	Let $[p,q]$ be our simplicial edge and $\ssq{p}{n}$ and $\ssq{q}{n}$ approximating sequences for $p$ and $q$ respectively. We can assume that $p$ is closest to $y$.
      	Set $l_{n}=d(p_{n},q_{n})$ and $l=d(p,q)$, $t_{n}=\tl{c}{\lambda_{n}}$ and $\zeta_{n}=\max\{d(p_{n},h\cdot p_{n}),d(q_{n},h\cdot q_{n})\}$. By convergence both $\ssq{t}{n}$ and $\tl{h}{\lambda_{n}}\neq0$ tend to $0$ with $n$.
      	Choose $c\in Stab(v)=Fix([p,q])\setminus\{1\}$ hyperbolical for the action $\lambda_{n}$. The acylindricity condition implies that for $n$ big enough $c$ has to act hyperbolically for $\lambda_{n}$.
      	Indeed there are points $p'_{n}$ and $q'_{n}$ in $Fix_{\lambda_{n}(c)}$ at distance strictly less than $\frac{\zeta_{n}}{2}$ from $p_{n}$ and $q_{n}$ respectively, so that the diameter of $Fix(c)$ is greater or equal than $l_{n}-\zeta_{n}$. Since this value tends to $d(p,q)\neq 0$ for $n$ big enough this cannot be the case.
      	
      	Up to restricting to a subsequence, we can assume that $c$ always translates $Ax(c)$ in the direction from $p'_{n}$ to $q'_{n}$ and up to replacing $h$ with some power we can assume that $\tau=\tau_{e,c}\in\PP$.
      	
      	Fix now $g\in G$ and suppose that $[y,g\cdot y]$ intersects some translate of $[p,q]$ non-trivially. Let $\looprep{g}{e}{m}$ be a loop representation of $g$ with respect to
      	some presentation of the tree $T_{[p,q]}$ described above, associated to some sequence $u_{0},e_{0}\cdots u_{m}\subset T_{[p,q]}$.
      	Notice that in this case $e_{j}$ is either the edge $e\in T_{[p,q]}$ associated to $[p,q]$ or $t_{e}^{-1}\cdot\edin{e}$ for all $i$ (possibly $t_{e}=1$).
      	For each $0\leq j\leq m-1$ let $k_{j}=g_{0}t_{e_{0}}g_{e_{1}}\cdots t_{e_{i-1}} g_{j}$.
      	For any $0\leq j\leq m$ let $p^{j}=p$, $q^{j}=q$ if $e^{j}=e$ and $p^{j}=q$, $q^{j}=p$ otherwise.
      	
      	Observe that we encounter the sequence of points $p^{0},q_{0},k_{1}\cdot p^{1},k_{1}\cdot q^{1}\cdots k_{m}\cdot q^{m}$ in that same order along the path between $y$ and $g\cdot y$.
      	\begin{comment}
      	\footnote{Indeed by assumption no translate of $I$ lies in the path between $y$ and $p$ or $t_{e}\cdot q$;
      	since $t_{e_{i}}^{-1}q^{i}\in C$, neither it does between $k_{i}q^{i}$ and $k_{i+1}p^{i+1}=k_{i}t^{e_{i}}p^{i}$. On the other hand, for $0\leq i\leq m-1$
      	if the two translates $h_{i}\cdot I$ and $h_{i+1}\cdot I$ coincide, then $[e_{i}]$ and $[e_{i+1}]$ are inverses in $G\backslash T_{I}$ and $g_{i}\in Stab(e)$, contradicting the properties of $\looprep{g}{e}{m}$.}
      	\end{comment}
      	By convergence in the equivariant Gromov topology,
      	\begin{align*}
      		\lmt{n}{\,d(p_{n},q_{n})+\serrang{0}{L}{m-1}d(t_{e_{j}^{-1}\cdot q^{j}_{n}},g_{j+1}p^{j+1}_{n})+d(y_{n},p^{0}_{n})+d(g_{m}\cdot y_{n},t_{e_{m}}^{-1}q^{0,m}_{n})}=\pl{y}{g}{}
      	\end{align*}
      	Now, let $k_{0}=1$ and for $0\leq j\leq m-1$ let $k_{j+1}=k_{}$. $k_{m}=\tau(g)$.
      	This implies the existence of sequences of points $\ssq{a}{n},\ssq{b}{n}\subset[p'_{n},q'_{n}]$ with $\ssq{a}{n}$ converging to $p$ and $\ssq{b}{n}$ converging to $q$ respectively such that if we let :
      	$a_{n}^{j}=a_{n}$ and $b_{n}^{j}=b_{n}$  in case $e_{j}=e$ and $a_{n}^{j}=b_{n}$ and $b_{n}^{j}=a_{n}$ otherwise, then
      	the sequence of points $w_{n}^{0},z_{n}^{0},w_{n}^{1}\cdots z_{n}^{m}$ defined by $w_{n}^{j}=k_{j}a_{n}^{j}$ and $z_{n}^{j}=k_{j}b_{n}^{j}$ appears in that order along the segment $[y_{n},g\cdot y_{n}]$ for $n$ big enough.
      	
      	\newcommand{\wjn}[0]{\hat{w}_{n}^{j}}
      	\newcommand{\zjn}[0]{\hat{z}_{n}^{j}}
      	\newcommand{\wjjn}[0]{\hat{w}_{n}^{j+1}}
      	\newcommand{\zjjn}[0]{\hat{z}_{n}^{j+1}}
      	Let the pair $(\wjn,\zjn)$ be equal to $(k')_{n}^{j}\cdot (a_{n},c^{-1}\cdot b_{n})$ if $e_{j}=e$ and $(b_{n},c\cdot a_{n})$ otherwise.
      	Since clearly $d(\wjn,\zjn)=d(w_{n}^{j},z_{n}^{j})-t_{n}$, in order to show that $\pl{y_{n}}{\tau(g)}{}=\pl{y_{n}}{\tau(g)}{}-t_{n}m$ for $n$ big enough it is enough to prove:
      	\elenco{
      		\item The segments $[z_{n}^{j},w_{n}^{j}]$ and $[\hat{z}_{n}^{j},\hat{w}_{n}^{j}]$ are translates of each other for $1\leq j\leq m-1$
      		\item $\hat{w}_{n}^{0}\in[y_{n},\hat{z}_{n}^{0}]$
      		\item $\hat{z}_{n}^{m}\in[\hat{w}_{n}^{m}\tau(g)\cdot y_{n}]$
      		\item $\zjn\in[\wjn,\wjjn]$, $\wjjn\in[\zjn,\zjjn]$ for each $1\leq j\leq m-1$
      	}
      	First of all, for each $1\leq j\leq m$ we have $\zjn=k'_{j}c_{e_{j}}^{-1}\cdot b_{n}^{j}$, while
      	$\wjjn=k'_{j}\tau(t_{e_{j}})g_{j+1}a_{n}^{j+1}=k'_{j}c_{e_{j}}^{-1}t_{e_{j}}g_{j+1}a_{n}^{j+1}$. It is clear from this that $(\zjn,\wjjn)=\mu_{j}\cdot(z_{n}^{j}w_{n}^{j+1})$, where $\mu_{j}=(h'_{j})^{-1}c_{e_{j}}^{-1}$.
      	
      	Likewise,	the triple $(\wjn,\zjn,\wjjn)$ is a translate of $(a_{n}^{j},c_{e_{j}}^{-1}\cdot{b}^{j}_{n},g_{j+1}\cdot a_{n}^{j+1})$. Now, $b_{n}^{j}\in [a_{n}^{j},g_{j+1}\cdot a_{n}^{j+1}]$, since $z_{n}^{j}\in[w_{n}^{j},a_{n}^{j+1}]$. On the other hand $[a_{n}^{j},b_{n}^{j}]\subset Ax(c)$, so as soon as the translation length of $c$ with respect to $\lambda_{n}$ is smaller than $d(a_{n}^{j},b_{n}^{j})$, we have $\zjn\in[\wjn,\wjjn]$. The argument showing $\wjjn\in[\zjjn,\zjn]$ is entirely symmetrical.
      	%      	(using the fact that oriented segment $[y\cdot \tau(g)\cdot y]$, together with the sequence of points in its interior is a translate of the inverse of $[y\cdot\tau(g)^{-1}\cdot y]$ ).
      \end{proof}
  \subsection{Proper quotients through shortening}
    
    \label{shortening morphism section} It is time to draw some consequences of a more algebraic nature from the previous theorem.
    \begin{comment}
    Let $\qG{L}_{A}[P_{i}](x)$ be a marked \gr \pr-limit group, and $\Delta$ an \gd \gad of $\rqG{G}{A}$. If $A\neq\tg$, fix some $h\in A\setminus\{1\}$.
    Let $\rqG{H}{A}$ be \pr-group, endowed with an action without inversion $\rho$ on a $K$-acylindrical $\R$-tree $Y$ and $f$ a morphism from $\rqG{L}{A}$ to $\rqG{H}{A}$.
    We say that $f$ is $\Delta$-short (or simply 'short' in case $\Delta$ is a \rs \atjsj of $\gqG{L}{A}{P}{I}$)
    with respect to the given marking if and only if
    $\sl{*^{\lambda_{f}}_{A}}{x}{\lambda\circ\tau}\geq\sl{*}{x}{T}$, where
    $*^{\lambda_{f}}_{A}$ is the point chosen in \ref{basepoints}.
    \end{comment}
    
    Let $\gqG{L}{A}{P}{I}(x)$ be a marked non-cyclic \gr \pr-limit group, and $\Delta$ an \gd \gad of $\gqG{G}{A}{P}{I}$. If $A\neq\tg$, fix some $h\in A\setminus\{1\}$.
    Let $\rqG{H}{A}$ be \pr-group, endowed with an isometric action without inversion $\lambda$ on a $K$-acylindrical $\R$-tree $(Y,d)$ and $f$ a morphism from $\rqG{L}{A}$ to $\rqG{H}{A}$.
    
    We say that $f$ is $\Delta$-\emph{short}
    with respect to $x$-marking if and only if
    $\sl{*^{\lambda_{f}}_{A}}{x}{\lambda\circ\tau}\geq\sl{*^{\lambda_{f}}_{A}}{x}{T}$, for any $\tau\in Mod(\rqG{G}{A},\Delta)$, where
    $*^{\lambda_{f}}_{A}$ is the point as defined in \ref{basepoints section}.
    %  	maybe make this more concrete later
    
    Given a property $\mathcal{P}$ (and an $L$-invariant family $\mathcal{A}$ of subgroups of $L$) we say that $f$ is $\Delta$-short (relative to $\mathcal{A}$) with respect to the $x$-marking among those satisfying property $\mathcal{P}$ if it satisfies $\mathcal{P}$ and $\sl{*^{\lambda_{f}}_{A}}{x}{\lambda\circ\tau}\geq\sl{*}{x}{T}$ for any $\tau\in Mod(\gqG{L}{A}{P}{I},\Delta)$ for which $f\circ\tau$ also satisfies $\mathcal{P}$.
    %    	[In case no explicit reference to $\Delta$ is given, we take by default the \gr \atjsj of $\gqG{L}{A}{P}{I}$.]
    %        By a shortening \pr-quotient of a $\qG{G}_{A}[P_{i}](x)$ we will intend a limit
    %         Let $undefined$ a sequence of morphisms from $\rqG{G}{A}$ to $\rqG{H}{A}$. We say that
    \begin{corollary}
    	\label{morphism shortening}
    	Let $\gqG{G}{A}{P}{I}(x)$ be a marked \gr \pr-limit group and $\rqG{H}{A}$ a \rs limit group which is equipped with an acylindrical isometric action $\lambda$ on a real tree.
    	Suppose we are given a system of inequalities $\syneq{\Psi}{x,a}$ and $\ssq{f}{n}$ an unbounded sequence of morphisms from $\qG{G}_{A}(x)$ to $\rqG{H}{A}$ which are short among those preserving $\syneq{\Psi}{x,a}$.
    	Let $Y$ a limiting tree for the sequence $\ssq{f}{n}$ and $\lambda$, as described in \ref{properties limit}.
    	If each $\gqG{G}{A}{P}{I}$ is freely indecomposable and each of the $P_{i}$ is elliptic in $Y$, then $\limker{f}{n}\neq\tg$.
    	In particular the latter is the case whenever the domain is a freely indecomposable \rs group (trivial grading case).
    	%        	In particular, any sequence of morphisms from a \rs \pr-limit group $\rqG{G}{A}$ to a free group $\rqG{F}{A}\frp\qG{H}$ which are short in the sense of the Cayley graph of the latter has non-trivial limit kernel.
    \end{corollary}
    \begin{proof}
    	Suppose for the sake of contradiction that $\limker{f_{n}}{n}=\tg$. Up to replacing $f_{n}$ with a subsequence we can assume that the sequence of rescalings of the actions induced by the $f_{n}$ converges to a limiting action $\rho$ on a real tree $Y$. In virtue of \ref{properties limit}, the action of $G$ on $Y$ is superstable with trivial tripod stabilizers. Notice that since the sequence is unbounded, the subgroup $A$ is elliptic in $Y$. In particular it fixes the limit $*$ of the sequence $(*_{n})_{n\in\N}$ given by $*_{n}=*_{A}^{\lambda^{f_{n}}}$.
    	%  	[maybe not necessary if we assume that some parameter subgroup contains $A$ ]
    	
    	Notice $[h\cdot *,*]=\{*\}$ for any $h\in A$. At the same time, $[*,y_{j}\cdot *]$ is necessarily non-degenerate for some $y_{j}\in x$. Taking $*$ as a basepoint in \ref{general shortening}, we conclude that an automorphism $\sigma\in Mod(\rqG{G}{A},S)$ exists such that $\sl{*_{n}}{f_{n}\circ\sigma(x)}{}<\sl{*_{n}}{f_{n}(x)}{}$ for $n$ big enough. Since $\limker{f_{n}\circ\sigma}{n}=\tg$, for $n$ big enough
    	the map $f_{n}\circ\sigma$ has to preserve the system $\syneq{\Psi}{x,a}$, contradicting the assumption on $f_{n}$.
    \end{proof}
    
    Let $\gqG{L}{A}{P}{I}(x)$ be a marked \gd \pr-limit group and $\Delta$ a graded \gad of it.
    %     of it and $\rqG{K}{A}=\rqG{F}{A}\frp H$ a free \rs \pr-limit group.
    We say that $(\gqG{L}{A}{P}{I}(x),\Delta)$ is a solid pair if and only if a sequence of $\Delta$-short restricted morphisms $\fun{f_{n}}{\rqG{G}{A}}{\rqG{K}{A}}$ with trivial limit kernel exists.
    We say that a freely indecomposable \gr \pr-limit group $\gqG{L}{A}{P}{I}(x)$ is \emph{solid} if the pair $(\gqG{L}{A}{P},\Delta_{JSJ})$ is solid, where $\Delta_{JSJ}$ is	the \adjsj of $\gqG{L}{A}{P}{I}$.

    Suppose that for each $i\in I$ some finite tuple of generators $p^{i}$ of $P^{i}$ is given. A \emph{flexible sequence} of morphisms from $\gqG{L}{A}{P}{I}(x)$ to $\rqG{K}{A}$ is
    a sequence $(\fun{f_{n}}{\rqG{L}{A}}{\rqG{K}{A}})_{n\in\N}$ of short morphisms satisfying the condition:
    \begin{align*}
    	\tl{f_{n}(p^{i}_{k}p^{i}_{l})}{}\leq n\cdot \max_{1\leq k,l\leq |x|}\tl{f_{n}(x_{l_{0}}x_{k_{0}})}{}
    \end{align*}
    For any $i\in I$ and $1\leq k,l\leq|p^{i}|$. By a \emph{flexible \pr-quotient} of $\gqG{G}{A}{P}{I}$ we will intend the limit \pr-quotient of some flexible sequence.
    The fact that parameter subgroups must be elliptic in any limiting tree associated to such a sequence implies, by the results above, that its limit kernel must be non-trivial. Notice that the collection of all flexible sequences is closed under subsequences and diagonal sequences.
    
    A morphism $f$ from $\qG{L}_{A}(x)$ to $\rqG{K}{A}$ is said $\Delta$-\emph{flexible} if it factors through some flexible \pr-quotient after precomposing by some $\tau\in Mod^{\pi}_{A}(\Delta)$.
    Otherwise we call it a $\Delta$-solid morphism. By a $\Delta$-solid \pr-quotient of $\gqG{L}{A}{P}(x)$ we intend the limit \pr-quotient of a sequence of solid morphisms.
    In absence of an explicit reference to $\Delta$ we will tacitly assume it coincides with \adjsj of $\gqG{G}{A}{P}{I}$.
    %    	If we are given a family $\mathcal{G}$ of morphisms from $\rqG{G}{A}$ to $\rqG{K}{A}$
    %    	We say $\gqG{G}{A}{P}{I}$ is flexible if it admits a sequence of short morphisms with trivial limit kernel.
    %    	We can likewise define a $\Delta$-flexible sequence relative to $\mathcal{F}$ as a $\Delta$-flexible sequence whose members are of that same form, flexible quotients
    %    	localized to $\mathcal{F}$ as limit quotients of such sequences and classify morphisms of $\mathcal{F}$ in solid or flexible in $\mathcal{F}$ according
    %    	to whether up to \pr modular automorphism they factor through one of these quotients or not.
    %    	    	we will intend $(\qG{L}_{A}(x),\Delta)$
    %  \end{definition}
    It can be proven that the distinction between solid and flexible morphisms is independent from the choice of $x$ and $p_{i}$, but we will not prove it here.
    In the degenerate case of a \rs \pr-group one obtains the definition above has a simple interpretation.
    %    together with a particular case of a well known result (see \ref{Perin}).
    \begin{observation}\leavevmode
    \elenco{
    	\item $\rqG{L}{A}$ is solid if and only if it injects into $\rqG{F}{A}$
    	\item A morphisms from a freely indecomposable \rs \pr-limit group $\rqG{L}{A}$ to some free $\rqG{K}{A}$ is solid if and only if it is injective. \label{solid restricted}
    	% \item Any non-injective quotient is flexible
    	\item If $A\neq\tg$ or $\rqG{L}{A}$ is not cyclic, then up to precomposition with an element of $Mod(\rqG{L}{A})$, there are only finitely many distinct injective morphisms from a solid $\rqG{L}{A}$ into $\rqG{F}{A}$.
    	%  		If $A=\tg$ and $G$ cyclic, then the same result holds if we allow ourselves to postcompose by conjugation by elements of $\F$.
    }
    %    	\item A purely \rs (trivially graded) freely indecomposable \rs \pr-limit group $\rqG{L}{A}$ is solid if and only if an injective morphism from $\rqG{L}{A}$ to $\rqG{F}{A}$ exists. \label{solid restricted}
    \end{observation}
    \begin{proof}
    For the first point, suppose a converging sequence of short morphisms $\ssq{f}{n}$ from $\rqG{L}{A}$ to $\rqG{K}{A}$ exists such that $\limker{f_{n}}{n}=\tg$. By \ref{morphism shortening} this sequence has to be bounded, which implies that all morphisms in some subsequence $(f_{n_{m}})_{m}$ can be obtained from each other by
    postcomposing with an inner automorphism of $\F$, implying they all have trivial kernel.
    But then $\ker{f_{n_{m}}}=\limker{f_{n_{m}}}{m}=\tg$.	Any $f\in Mor(\rqG{G}{A},\rqG{F}{A})$ which is non-injective is flexible; here the condition on the growth of parameter groups is
    empty, so if $f\circ\sigma$ is short for some $\sigma\in Mod(\rqG{G}{A})$ the constant sequence $f_{n}=f\circ\sigma$.
    \end{proof}
    
    In case $A$ is a free factor of $\F$, the situation above can only be the case if $L=A$, since otherwise the pullback of a free splitting of $\F$ containing $A$ as a factor is a free splitting
    of $L$ relative to $A$. This becomes crucial later.

\chapter{Makanin-Razborov diagrams  }
\label{MR chapter}
\section{Basic notions}
  %  A fundamental step in the study of the theory $T$ is the analysis of varieties, i.e., solutions sets of finite systems of equations, as these definable sets generates the boolean algebra quantifier free definable sets (with parameters).
  By a system of equations in the tuple of variables $x$ and the tuple of parameters $a$ we mean a condition of the form:
  \begin{align*}
  	\bwedge{i\in I}{w_{i}(x,a)=1}{ }
  \end{align*}
  where for $w_{i}(x,a)$ is a term in the language of groups. This is mostly abbreviated as $\syeq{\Sigma}{x,a}$. By a system of inequations we mean a condition of the form:
  \begin{align*}
  	\bwedge{i\in I}{w_{i}(x,a)\neq 1}{ }
  \end{align*}
  mostly abbreviate by an expression such as $\syneq{\Sigma}{x,a}$. We can think of $\Sigma$ as a subset of elements of any particular $\subg{a}$-\rs group. Any quantifier-free definable formula is equivalent to one of the form $\bvee{i=1}{\ceq{\Sigma}{(x,a)}{}\wedge\cneq{\Pi}{x,a}{}}{m}$.
  When passing on to the language $\LQ$ the same property holds for the intersection of the two types of formulas above with those of the form $\binters{i=1}{P_{q_{i}}(x_{i})}{k}$, for $k=|x|$ and $q=(q_{i})_{i=1}^{k}\in Q^{k}$, which we will often abbreviate the latter as $x\in q$.
  This will be be referred to as \CEQS and \CNEQS respectively
  and denoted by an expression of the form $\ceq{\Sigma}{x,a}{q}$  (resp. $\cneq{\Sigma}{x,a}{q}$).
  %  Any quantifier free formula in the language $\LQ$ is equivalent to one of the form $\bvee{i=1}{\ceq{\Sigma}{(x,a)}{q(i)}\wedge\cneq{\Pi}{x,a}{r(i)}}{m}$.
  %($\Sigma_{\qG{G},x,a}$)
  \begin{definition}
  	A group $G$ is said to be equationally noetherian if any system of equations with parameters in $G$ is equivalent to some finite subsystem, i.e., they both have the same set of solutions in $G$.
  \end{definition}
  The property holds for free groups, see \cite{guba1986equivalence}.
  %  Suppose for a second all our notions are defined in terms of an arbitrary \pr-group $\qG{K}$ in place of $\qG{F}$
  Given a group $K$, $A\leq K$ and a marked \rs group $G_{A}(x)$ the variety $Var(G,A,x)[K]$ is defined as the set of tuples $x\in K^{\pmb{x}}$ which are the image of $x(G)$ by some morphism from $G_{A}$ to $\F_{A}$.
  Let $Var(G_{A},x)$ stand for $Var(G,A,x)[\F]$.
  
  For any  surjective \res homomorphism $\fun{\phi}{G_{A}}{H_{A}}$ and compatible $x$-markings on both $G_{A}$ and $H_{A}$, clearly
  $Var(H_{A},x)\subset Var(G_{A},x)$ and it can be easily seen this inclusion is proper in case that $G$ is residually free, that is, whenever for any
  $g\in G$  there is some homomorphism (morphism) from $G$ ($\qG{G}$) to $A$.
  Observe that equational noetherianity implies that for every $A\leq K$ and any chain of epimorphisms between finitely generated marked groups:
  \begin{align*}
  	G_{A}^{1}(x)\to G_{A}^{2}(x)\to\cdots G_{A}^{m}(x)\to\cdots
  \end{align*}
  the sequence
  \begin{align*}
  	Var(G^{1}_{A};x)[K]\supset Var(G^{2}_{A};x)[K]\supset\cdots Var(G^{m}_{A};x)[K]\cdots
  \end{align*}
  contains finitely many proper inclusions. Another consequence is that all the homomorphisms in some convergent sequence with $\F$ as target eventually factor through the limit quotient of the sequence.
  %  We shall now state the important consequences of the equational noetherianity of the free group.
  %  Clearly, given any relative presentation $G=\subg{x,A|w_{i}(x,A)}_{i\in I}$ of the finitely generated* \rs group $G$
  \begin{corollary}
  	Every chain of epimorphisms between residually free groups (hence, in particular, between limit groups) contains only finitely many which are proper.
  \end{corollary}
  
  If we are given a marked \rs \pr-group, $\qG{G}_{A}(x)$ instead, we can define a \emph{\pr-variety} $Var(\rqG{G}{A},x)$ as the set of tuples in $ \F^{x}$ which are the image of $x(G)$ by some morphism from $\rqG{G}{A}$ to $\rqG{K}{A}$.
  We might refer to the points in the \pr-variety as solutions of $\qG{G}_{A}(x)$.
  %  For an unrestricted \pr-group $\qG{G}=\subg{s}$, let $PVar(G;s)=PVar(\rpG{G,1})$.
  %  If $x=(x_{i})_{i=1}^{k}$, and $\pi(x_{i})=q_{i}$, $PVar(G;A;s)=Var(G;A;s)\cap\binters{i=1}{\{x_{i}=q_{i}\}}{k}$.
  
  Notice that the correspondence $f\mapsto f(x)$ is a bijection between $Hom_{A}(G,\F)$ ($Mor(\rqG{G}{A},\rqG{F}{A})$) and $Var(G_{A},x)$ ($Var(\rqG{G}{A},x)$).
  Given a system of equations $\syeq{\Sigma}{x,a}$, let $G^{\Sigma}_{A}$ be the $A$-\rs marked group with presentation $\subg{x,A|\syeq{\Sigma}{x,a}}$ (relative to $A$), where $A=\subg{a}$. The set of morphisms from it to $\rqG{F}{A}$
  determines in turn the same variety defined by $\Sigma$.
  %  If we also have conditions specifying the predicates to which each of the $x_{i}$ belong, yields a \CEQ $\ceq{\Sigma}{x,a}{p}$ instead.
  Given a \CEQ we can expand the previous group to a \pr-group by sending each generator in the marking to its specified value in $Q$. We denote the resulting \rs \pr-group as $\qG{G^\Sigma}_{A}(x)$, where $A=\subg{a}$. The same comment applies.
  
  %  Note that if $G$ ($\qG{G}_{A}(x)$) is isomorphic to $G_{\Sigma}$ ($\qG{G_{\Sigma}}_{a}(x)$), where $a$ is a tuple of generators of $A$ then $Var(G;a;x)[K]$ ($PVar(G;a;x)[K]$) is the set of solutions in $K$ of $\ceq{\Sigma}{x,a}{}$ ($\ceq{\Sigma}{x,a}{}$).
  
  %  The above remark has an obvious generalization to our context. Given a surjective \res homomorphism $\fun{\phi}{\rqG{G}{A}}{\rqG{H}{A}}$,
  %  $PVar(H;A;x)\subset PVar(H;A;s,x)$. Viceversa, any \pr-subvariety $V\subset PVar(H;A;s,x)$ is associated with some \pr-quotient $\fun{\phi}{\rqG{G}{A}}{\rqG{H}{A}}$.
  %  The following remark should be taken into account as well:
  \begin{comment}
  \begin{remark}
  	Let $\pmb{x}$ and $\pmb{z}$ two different markings of a \rs \pr-group $\rqG{G}{A}$. Then the varietis $PVar(G;A;x)$ and $PVar(G;A;z)$ are interdefinable.
  	They are in fact interalgebraic: given a generating tuple $a$ of $A$, there is an $|x|$-tuple $u(z,a)$ of words in $z$ and and a $|z|$-tuple $v(x,a)$ of words in $x$ and $a$ such that
  	\begin{align*}
  		\forall x\in\F^{|x|}\,\,x\in PVar(G;A;x)\slra v(x,a)\in PVar(G;A;z)
  		\\ \forall z\in\F^{|z|}\,\,z\in PVar(G;A;z)\slra u(z,a)\in PVar(G;A;x)
  		\\ u(v(x,a),a)=x\,\,\,v(u(z,a),a)=z
  	\end{align*}
  	This implies, in particular, that any subvariety of $PVar(G;A;x)$ pulls back to a subvariety of $PVar(G;A;z)$ and viceversa.
  	%  the two word maps determine a correspondence between definahle subsets of both $PVar(G;A;x)$ and $PVar(G;A;z)$.
  \end{remark}
  \end{comment}
  Given a family $\mathcal{C}$ of \pr-limit groups or one $\mathcal{C}_{A}$ of $\qG{A}$-\rs \pr-limit groups, define the partial order $\leq_{Z}$ on $\mathcal{C}$ by saying that $\qG{G}\leq\qG{H}$ ($\rqG{G}{A}\leq\rqG{H}{A}$) if and only if there is a surjective morphism from $\qG{H}$ to $\qG{G}$ ($\rqG{G}{A}$ to $\rqG{H}{A}$). The discussion above implies $\leq_{Z}$ is is well-founded.
  
  Equational noetherianity implies the following characterization of \pr-limit groups, entirely analogous to that of limit groups.
  \begin{lemma}
  	\label{pilimit characterization}Let $\rqG{G}{A}$ be any finitely generated \rs \pr-group. The following are equivalent:
  	\enum{i)}{
  		\item  $\rqG{G}{A}$ is a \rs \pr-limit group. \label{islimit}
  		\item  $\rqG{G}{A}$ is \pr-discriminated by $\rqG{F}{A}$, that is: for any finite set $S\subset G$ there is a morphism $\fun{f}{\rqG{L}{A}}{\rqG{F}{A}}$ which is injective on $S$. \label{isfully}
  		\item $\rqG{G}{A}$ is a model of the universal part of $Th(\qG{F}_{A})$.  \label{ismodeluniversal}
  	}
  \end{lemma}
  %  The implication from (\ref{ismodeluniversal}) to (\ref{islimit}) is the only one requiring equational noetherianity.
  %  \begin{comment}
  \begin{proof}
  	The implication from (\ref{islimit}) to (\ref{isfully}) follows from equational noetherianity. The homomorphisms in a ny convergent sequence targetting must factor through the limit quotient after a certain point.
  	
  	%  (Observe that the equivalence between \ref{ismodeluniversal} and \ref{isfully} and the fact that they are implied by \ref{})
  	
  	Let us look now at the implication from (\ref{isfully}) to (\ref{ismodeluniversal}). Take a \rs \pr-group $\rqG{G}{A}$ which is \pr-discriminated by $\rqG{F}{A}$. All one has to show is that for any tuple $b$ of elements from $G$ and any atomic formula $\phi(x,a)$ ($|x|=|b|$) with constants $a$ in $A$ such that $\phi(b,a)$ holds in $\rqG{L}{A}$ some tuple $c$ in $F$ exists such that $\phi(c,a)$ holds in $\rqG{F}{A}$. Observe that it is enough to verify this for a formula of the form:
  	\begin{align*}
  		(\bwedge{i=1}{P_{q_{i}}(x_{i})}{|x|})\wedge\syeq{\Sigma}{x,a}\wedge\syneq{\Delta}{x,a}
  	\end{align*}
  	since any atomic formula with free variables in $x$ is the disjunction of finitely many formulas of that same form. Property (\ref{isfully}) implies the existence of a morphism $f$ from $\qG{G}_{A}$ to $\rqG{F}{A}$ which does not kill the words in $\syneq{\Delta}{x,a}$. Of course, since this is a homorphism fixing $A$ any element which can be written as $u(b,a)$ for some word $u(x,y)$ is sent by such $f$ to $u(f(b),a)$, so that $\syneq{\Delta}{f(b),a}$ holds as well. That $\phi(f(b),a)$ holds is now obvious, the fact that a morphism preserves the positive part of such a formula being equally easy to prove.
  	
  	The implication from (\ref{ismodeluniversal}) to (\ref{islimit})
  	follows the inverse path. Let $s$ be a tuple of generators of a given finitely generated model $\rqG{G}{A}$ of the universal theory of $\rqG{F}{A}$ and $u$ a finite tuple of elements of $G$. Let $\rqG{H}{A}$ the \pr-group with underlying group $A\frp\F(s)$, together with the obvious homomorphism to $Q$. There is a natural morphism $p$ from $\rqG{H}{A}$ to $\rqG{G}{A}$. Let $\syeq{\Sigma_{1}}{x,a},\syeq{\Sigma_{2}}{x,a}\cdots$
  	an increasing chain of systems of equations whose union contains the collection of all the equations with parameters in $a$ satisfied by the tuple $(s,a)$. Let $\syneq{\Delta_{1}}{x,a},\syneq{\Delta_{2}}{x,a}\cdots$ a similar exhaustion, but with inequations in place of equations.
  	
  	Denote by $\phi_{n}(x,a)$ be the conjunction of $\syeq{\Sigma_{n}}{x,a}$, $\syneq{\Delta_{n}}{x,a}$ and the condition $P_{\pi(s_{i})}(x_{i})$. Since $s$ satisfies $\phi$ in $\rqG{G}{A}$ and $\rqG{G}{A}$ is a model of the universal theory of $\rqG{F}{A}$,
  	some tupe $s'$ must also satisfy $\phi$ in $\rqG{F}{A}$. The map sending $s_{i}$ to $s'_{i}$ determines a morphism $f_{n}$ from $\rqG{K}{A}$ to $\rqG{F}{A}$. The choice of $f_{n}$ implies that the limit kernel of the sequence $\ssq{f}{n}$ is the same as the kernel of the projection $p$.

  	%  	[quote Abderezak notes on limit groups of EN groups]
  \end{proof}
  %  \end{comment}
  \begin{comment}
  Larsen Louder has shown in \cite{Krull dimension}, that in fact the length of a descending $\leq_{Z}$-chain of limit groups admits an affine bound in terms of the rank of $G_{1}$. This allows to define a finite-valued dimension $\fun{dim}{\mathcal{C}_{A}\times\mathcal{C}_{A}}{}$ as the sharpest such bound. In our case, well-foundedness of$\leq_{Z}$ will be enough for the kind of inductions arguments we are going
  %  can be in principle phrased in terms of $dim$, although this is just a formality.
  \end{comment}
  
  While equational noetherianity guarantees the termination of the Makanin-Razborov procedure, the next result takes care of the finite branching of the process. We provide a topological interpretation of the diagonalization argument found in \cite{Sela1}.
  %  references
  %  to draw from the topological approach in \cite{Guirardel compactifying}.
  
  \subsection{Finite width}
    
    Let $\mathcal{F}$ be the family of homomorphisms from a group $G$ to an equationally noetherian group $H$ and
    let $\mathcal{K}$ a subset of the closure of $\{ker(f)\,:\,f\in\mathcal{F}\}$ in $2^{G}$. Recall that the product topology on $2^{G}$ is compact and metrizable, hence sequentially compact.
    One can easily check that in this case $\mathcal{K}$ consists of limit kernels of convergent sequences of homomorphisms in $\mathcal{F}$.
    Equationally noetherianity conditions the way $\mathcal{K}$ sits in $2^{G}$, namely:
    \begin{equation*}
    	\text{For each $K_{0}\in 2^{G}$ there is a neighbourhood $V$ of $K_{0}$ such that $\forall K\in V\cap\mathcal{F}\,\,K_{0}\subseteq K$  }
    \end{equation*}
    Observe that $\subset$ is a closed relationship in $2^{G}$, so one can strengthen the previous statement to:
    \begin{equation}
    	\tag{*} \label{topological noetherianity}
    	\text{For each $K_{0}\in 2^{G}$ there is a neighbourhood $V$ of $K_{0}$ such that $\forall K\in V\cap\mathcal{K}\,\,K_{0}\subseteq K$  }
    \end{equation}
    The intersection of the elements of a chain $\mathfrak{C}=\{\mathcal{K}_{i}\}_{i\in I}$ of $(\overline{\mathcal{K}},\subseteq)$ belongs to $\overline{\mathfrak{C}}\subset\overline{\mathcal{K}}$, since the restriction of $\overline{\mathcal{K}}$ to the elements in any ball of the Cayley graph of $G$ eventually stabilizes. Hence $\binters{i\in I}{K_{i}}{}\in\overline{\mathcal{K}}$. The observation above, and Zorn's lemma imply that for any element $K\in\overline{\mathcal{K}}$ there is some $L\in\mathcal{L}$ such that $L\neq K$.
    \begin{claim}
    	If $\mathcal{K}$ is closed then $\mathcal{L}$ is finite.
    \end{claim}
    \begin{proof}
    	Since the space $2^{G}$ is sequentially compact, it is enough to show that $\mathcal{L}$ has no accumulation points. Suppose there was a sequence $(L_{n})_{n\in\N}\subset\mathcal{L}$ and $L\in\overline{\mathcal{K}}$, such that $L\neq L_{n}$ for all $n$ and $\ssq{L}{n}$ converges to $L$. Then (\ref{topological noetherianity}) implies that eventually $L\subsetneq L_{n}$, contradicting the fact that $L_{n}\in\mathcal{L}$.
    \end{proof}
    %        We are really just interested in the case in which $\mathcal{K}$ is the set of limit kernels of a certain family $\mathcal{S}$ of sequences of homomorphisms from
    %        $G$ to $H$.
    If $\mathcal{K}$ is the set of limit kernels of a family of sequences of homomorphisms from $G$ to $H$, there is a simple criterion on $\mathcal{S}$ ensuring that $\mathcal{K}$ is compact.
    %     Let $d$ be some metric compatible with the product topology on $2^{G}$.
    
    We say that a family $\mathcal{S}$ of sequences is closed under diagonal sequences if for any sequence $(\ssq{f^{k}}{n})_{k}\subset\mathcal{F}$
    %    	some subsequence of
    the sequence $(f^{n}_{n})_{n}$ belongs again to $\mathcal{S}$.
    \begin{lemma}
    	\label{diagonals}
    	If $\mathcal{S}$ is closed under subsequences and diagonal sequences, then $\mathcal{K}$ is closed in $2^{G}$.
    	\begin{comment}
    	%    	Give the first condition a name.
    	%    	 suppose that for any sequence $(\ssq{f^{k}}{n})_{k}$ of sequences from $\mathcal{S}$
    	%    	each sequence has fast converging kernel
    	Suppose that for each sequence $\ssq{f}{n}\in\mathcal{S}$ is of rapidly converging kernel.
    	By this we mean that for some tuple $s$ of generators of $G$, $\limker{f_{n}}{n}\cap B_{n}=\ker{f_{m}}\cap B_{n}$ for any $m\geq n\geq 0$,
    	where $B_{n}$ is the ball of radius $n$ around the identity in the Cayley graph of $G$ with respect to $s$.
    	Assume furthermore, that for any sequence $(\ssq{f^{k}}{n})_{k}\subset\mathcal{F}$
    	%    	some subsequence of
    	the sequence $(f^{n}_{n})_{n}$ belongs to $\mathcal{S}$.
    	\end{comment}
    	%    	there are subsequences
    \end{lemma}
    \begin{proof}
    	This follows easily from the fact that $2^{G}$ is metrizable. For any $K_{n}=lker_{n} f^{n}_{m}$ in a sequence $\ssq{K}{n}$   converging to $K$, take some $f^{n}_{k_{n}}\in\mathcal{F}$
    	in the corresponding sequence $(f^{n}_{m})_{n}$ such that $Ker\,f^{n}_{k_{n}}$ is as close to $K_{n}$ as $K_{n}$ is to $K$.
    \end{proof}
    We will refer to a family satisfying the property above as being closed under subsequences.
    Here are some first consequences of the above discussion of importance to us:
    \begin{lemma}
    	\label{factor sets}
    	Let $\rqG{G}{A}$ a finitely generated \rs \pr-group and $\mathcal{F}$ a family of morphisms from finitely generated \rs \pr-group $\rqG{G}{A}$ to $\rqG{F}{A}$. Then any $f\in\mathcal{F}$ factors through one of finitely many  \pr-quotients which are limits of sequences of elements of $\mathcal{F}$.
    	%    	The family $\mathcal{M}(\mathcal{F})$ of all the limit quotients of morphisms of $\mathcal{F}$ is finite.
    \end{lemma}
    In particular we have:
    \begin{corollary}
    	\label{absolute factor sets}
    	For any \rs group $\rqG{G}{A}$ there is a finite family $\mathcal{MQ}(\rqG{G}{A})$ of \rs \pr-quotient maps from $\rqG{L}{A}$ onto an $A$-\rs \pr-limit group such that any morphism from $\rqG{G}{A}$ to a given \pr-group $\rqG{H}{A}$ factors through some $q\in \mathcal{MQ}(\rqG{G}{A})$. %If $\rqG{G}{A}$ is not a limit group, then all such $q$ are proper.
    \end{corollary}

  \subsection{\pr-Resolutions}
    
    One of the most important notions in the work of Sela is that of a resolution. One of their uses is that of parametrizing the set of solutions of a  systems of equations. They generalize straightforwardly to our context.
    
    \begin{definition}
    	A \pr-\res  $\QR{R}=\LRL{R}{}{(J,r)}$ consists of a finite rooted tree $(J,r)$ and a series of assignments:   %use finite branching instead ?
    	\elenco{
    		\item To each $\lambda\in J$:
    		\enum{i)}{
    			\item A non-trivial \pr-group $\Rl$, which is a \pr-limit group in case $\lambda\neq r$ or $\GRla$ is non-trivial.
    			\item A \gad $\DR{\Delta}{\lambda}$ of $\Rl$ relative to the family of all non-cyclic abelian subgroups.
    		}
    		\item To each $\lambda\in \hat{J}$, a surjective \pr-morphism $\eta^{\lambda}_{R}$ from $\Rl$ to the free product:
    		\begin{align*}
    			\bfrp{\RLR{R}{}{\mu}}{\lambda\succ\mu}{}
    		\end{align*}
    	}
    	We will refer to $\RLR{R}{}{\lambda}$ as the top of the resolution $\QR{R}$ and to $\QR{R}$ as a resolution of $\qG{G}$.
    \end{definition}
    
    \begin{definition}
    	By a \rs \pr-resolution $\RRA=\LRL{R}{A}{(J,r)}$ we intend a resolution of the underlying \pr-group $\qG{G}$ in which for some $A\leq\F$ each of the node groups appearing along some branch $\mathcal{B}$ of $J$
    	%  	, the \prbr,
    	are all endowed with an additional \rs structure: $\RAl$. Moreover:
    	\elenco{
    		\item The group of constants $A$ is elliptic in the \gad $\DR{\Delta}{\lambda}$ for any $\lambda\in\mathcal{B}$.
    		\item For $\lambda\in\hat{\mathcal{B}}$ the map $\eta^{\lambda}$ is a morphism of \rs \pr-groups onto the free product $\bfrp{\RLR{R}{}{\lambda}}{\lambda\succ\mu}{}$.
    	}
    	Formally we will consider the \lvlgp at any node $\lambda$ as an $A_{\lambda}$-\rs \pr-limit groups, where $A_{\lambda}=A$ in case $\lambda\in\mathcal{B}$ and $\tg$ if not.
    	Again we will refer to $\RLR{R}{A}{\lambda}$ as the top of the \pr-\res and to it as a \rs \pr-\res of $\RLR{R}{A}{\lambda}$.
    \end{definition}
    
    Given a \pr-\res $\LRL{R}{}{(J,r)}$, we say that $\lambda$ is of free product type in case $\GRla$ is trivial and $\eta^{\lambda,R}_{\lambda,R}$ is injective.
    %  	quotient type if $\eta^{\lambda,R}_{\lambda,R}$ is non-injective. We say it is of free group type.
    A  (restricted or graded) \pr-\res will be called proper in case $\eta^{\lambda}_{\lambda}$ is non-injective 	%$J$ has at least one edge and
    whenever $\mu$ is the only child of $\lambda$.
    %It will be called \emph{non-mixed} if for each $\lambda\in\hat{J}$ either $\lambda$ has only one child or it is of free product type.
    %Of course, a resolution can be always taken into a non-mixed one by adding an intermediate single child when required.
    Given $\mathcal{R}^{J}$, and $\lambda\in J$, by $\mathcal{R}^{J}\restriction_{\lambda}$ we intend its restriction of $\mathcal{R}$ to $J\restriction_{\lambda}$.
    
    %Given a free product $\bfrp{A_{i}}{i\in I}{}$ and a family of homomorphisms $\{\fun{f_{i}}{A_{i}}{G}\}$, let
    %$\bdcup{f_{i}}{i\in I}{}$ stand for the unique homomorphism $\fun{g}{\bfrp{A_{i}}{i\in I}{}}{G}$ such that for each $\in I$,
    %$g\restriction_{A_{i}}=g_{i}$.
    
    Given a \pr-\res $\LRL{R}{}{(J,r)}$, a \trm of it
    to a \pr-group $\qG{H}$
    is a pair $(f_{\lambda},\tau_{\mu})_{\lambda\in J,\mu\in J}$,
    where the $\fun{f_{\lambda}}{\Rl}{\qG{H}}$ is a morphism  and
    $\tau_{\lambda}\in \Mo{R}{\lambda}$ are such that for any $\lambda\in (\hat{J})$ :
    \begin{align*}
    	f_{\lambda}=(\bdcup{f_{\mu}}{\lambda\succ\mu}{}) \circ \eta^{\lambda}_{\lambda}\circ\tau_{\lambda}
    \end{align*}
    An obvious induction argument shows that any choice of morphisms $\{f_{\lambda}\}_{\lambda\in \lvs{J}}$ and of automorphisms $\{\tau_{\lambda}\}_{\lambda\in (\hat{J})}$ determines a unique \trm of $\LRL{R}{}{(J,r)}$.
    
    The morphism $f_{r}$, will be said to factor through $\QR{R}$
    %  	We will say that a morphism $\fun{f}{\qG{G}}{\qG{F}}$ \emph{factors through} $\LRL{R}{}$ if there is a \trm as before such that $f=f_{r}$.
    and we will refer to $f_{\lambda}$ as a \fcr for $f$ at $\lambda$.
    Given a \pr-\res $\LRL{R}{}{(J,r)}$, we let $Fct(\mathcal{R})$ be the set of morphisms from $\qG{G}$ to $\qG{F}$ which factor through it.
    
    Given a \rs \pr-\res , $\LRL{R}{A}{(J,r)}$, a \trm of $\QR{R}$
    to a \rs \pr-group $\rqG{H}{A}$ is a \trm $(\{f_{\lambda}\}_{\lambda\in J},\{\tau_{\lambda}\}_{\lambda\in (\hat{J})})$,
    of the underlying \pr-\res,
    to $\qG{H}$ such that
    the $f_{\lambda}$ are morphisms from
    $\RAl$ to $\rqG{H}{A}$ and $\tau_{\lambda}\in \Mo{R}{\lambda}$.
    We will say that a morphism $\fun{f}{\rqG{G}{A}}{\rqG{F}{A}}$ factors through $\LRL{R}{A}{}$ if there is a \rs \trm $(f_{\lambda})_{\lambda\in J}$ such that $f=f_{r}$.
    Once more, we denote by $Fct(\QR{R}_{A})$ the set of morphisms from $\rqG{G}{A}$ to $\rqG{F}{A}$ factoring through $\RRA$.
    %If $\rqG{G}{A}$ is marked with a tuple of variables $x$, we might refer to the image of $x$ by a morphism which factors through $\RRA$ as a solution factoring through $\QR{R}$.
    
    \begin{definition}
    	A \gd \pr-\res $\QR{R}=\GQR{R}{A}{}{P}{I}$ of a \gr \pr-group $\gqG{G}{A}{P}{I}$ is a \rs \res of $\rqG{G}{A}$ together with an additional \gr \pr-structure $\gqG{R^{\lambda}}{A}{P}{I}$ on each of its node groups in such a way that the following conditions are satisfied:
    	\enum{i)}{
    		\item For each $\lambda\in\hat{J}$ the index set $I_{\lambda}$ is equal to $\bunion{\mu\prec\lambda}{I_{\lambda}}{}$.
    		\item For each $j\in I_{\lambda}$, the parameter group $P^{\lambda}_{j}\leq \Rl$ is mapped by $\eta^{\lambda}_{\Q R{R}}$ onto a conjugate in $\bfrp{\RLR{R}{}{\mu}}{\mu\prec\lambda}{}$ of $P^{\mu}_{j}$	for some $\mu\prec\lambda$ and $j\in I_{\mu}$. This defines a bijective map from $I_{\lambda}$ to $\abunion{\mu\prec\lambda}{I_{\mu}}{}$.
    		\item For each $\lambda\in\hat{J}$ the \gad $\GRla$ is relative to $\{P^{\lambda}_{j}\}_{j\in I_{\lambda}}$.
    	}
    	We say that the node group at $\lambda$ carries parameter groups in case $I_{\lambda}\neq\nil$.
    \end{definition}
    
    We use the notation $\QR{R}$ indistinctively for \gd,\rs or simple \pr-\res as long as the context makes clear the kind of structure we are referring to.
    To abbreviate, given a \gd  \pr-\res $\mathcal{R}_{A}$ as above and a node $\lambda\in J$ denote the group $Mod(\RAl,\GRla)$ by $\RMo{R}{A}{\lambda}$. Same comment applies.
    %    \pr-\res $\LRL{R}{A}{} $.
    \newcommand{\strm}[0]{solid \trm } %text command
    %\begin{definition}
    	We say that a \gr \pr-\res $\QR{R}$ is closed if for any $\lambda\in\lvs{J}$ either:
    	\enum{i)}{
    		\item  $R^{\lambda}$ carries the constants or some parameter subgroup, in which case $\gqG{R^{\lambda}}{A}{P}{I}$ is freely indecomposable and solid
    		\item  $R^{\lambda}$ is free and carries no parameters nor constants.
    	}
    %\end{definition}
    In the last case we refer to $\lambda$ as a solid leaf. Recall that
    in the case of a merely \rs \pr-\res $\rqG{G}{A}$ if the group $A$ of constants is a free factor of $\F$
    then $A$ occupies a whole group at one of the leafs.
    
    By a solid \trm of a closed \gr \pr-\res (in particular of a \rs \pr-\res) we mean a \trm $\fcrm{J}{f}{\tau}$ of it such that for any $\lambda\in\lvs{J}$ not carrying any constants for parameters the morphism $f_{\lambda}$ is solid. In the case of a \rs \pr-\ress this means the groups carrying the constants are injected.
    In this case we say that $f_{r}$  factors solidly through $\RRA$ and to $Fct^{sld}(\QR{R})$ as the set of morphism factoring solidly through $\RRA$.
    
    We will entertain for some time the the possibility that any \pr-group along a  (\gd or \rs) \pr-\res might not be a (\rs, \gd) \pr-limit group. We will refer to such an object as a 'weak' (\rs, \gd) \pr-\res.
    %  obtained by lifting the restriction that \lvlgps must (\gd or \rs) \pr-limit groups.
    It turns out that if a weak resolution is closed and strict in the sense below then it is also in fact a \pr-\res in the standard sense, as we will show later.

\section{The Makanin-Razborov procedure}
  
  %  \begin{definition}
  	We say that a \pr-\res $\LRL{R}{}{(J,r)}$ is \emph{\cautious} if $\DR{R}{r}$ is trivial. This will be the case by default for us whenever $\RLR{R}{}{\lambda}$ is not a limit group.
  As we saw, understanding the set of solutions of a \CEQ over some $A\subset\F$ amounts to that of the family of all morphisms from a
  certain $\qG{A}$-\rs \pr-group $\rqG{G}{A}$ to $\rqG{F}{A}$.
  %   Recall that here by $\qG{F}$ our reference model is intended.
  It will be in fact convenient to analyze the family of all morphisms from $\rqG{G}{A}$ to any \pr-group $\rqG{F'}{A}$ of the form $\rqG{F}{A}\frp\qG{H}$. For a fixed $\rqG{L}{A}$ we might as well restrict to a single $\qG{H}$.
  \newcommand{\K}[0]{\rqG{K}{A}}
  \begin{claim}
  	Given a finitely generated $\rqG{G}{A}$, there is some \pr-group $\qG{H_{0}}$ such that the image of
  	any morphism $\fun{f}{\rqG{G}{A}}{\rqG{F'}{A}}$ for $\rqG{F'}{A}$ as above factors through $\rqG{F}{A}\frp\qG{H_{0}}$.
  \end{claim}
  %  Given $H$, consider the free splitting of $\F\frp H$ which has $\F$ as single vertex group and $rk(H)$ edges with trivial edge groups.
  \begin{proof}
  	Given any free \pr-group $\qG{F'}=\rqG{F}{A}\frp \qG{H}$, consider the $\F'$-tree $T$ with trivial edges stabilizers and a single orbit of vertices, stabilized by conjugates of $\F$.
  	Now, suppose that $\fun{f}{\rqG{G}{A}}{\rqG{F}{A}\frp\qG{H}}$ is also given. The map $f$ induces an action of $Q=G/ker(f)$ on $T$ with trivial edge stabilizers. Consider the action of $G'$ on its minimal tree $S\subset T$.
  	%  a free decomposition $\Gamma$ of $G$ given by a presentation of
  	There cannot be more than $rk(G')\leq rk(G)$ distinct orbits of non-trivially stabilized vertices of $S$ (by the action of $Q$). Since they are all contained in the same $\F'$-orbit of $VT$, by adding at most $rk(G)-1$ elements of $\F'$ to $Q$, we can get $Q\leq Q'\leq\F'$ with $rk(Q')\leq 2rk(G)-1$ acting on its minimal tree $S'$ with a single non-trivially stabilized orbit.
  	Now, fix some maximal spanning tree of the quotient graph $Q'\backslash S'$ and collapse all the edges in its preimage in $S'$. The tree so obtained is dual to that of a free decomposition of $G'$ of the form $K\frp H'$ where $K\subset\F$ and $H'$ is free. Now, $rk(H')\leq rk(G')\leq 2rk(G)-1$. Let $\qG{H^{1}}\cdots\qG{H^{m}}$ be all the non isomorphic \pr-structures on a free group of rank $\leq 2rk(G)-1$. Hence  $\qG{H^{0}}=\bfrp{\qG{H^{i}}}{i=1}{m}$ satisfies the properties we need.
  	%  Use the relative rank of $G$?
  \end{proof}
  The idea behind this is that of preserving some more information on the set of solutions in $\qG{F}$, rather than studying the set of solutions in a larger model.
  %  As we shall see, the family of retractions from $\qG[\rqG{F}{A}_{i}]\frp\qG{H}$ to $\rqG{F}{A}$ is discriminating.
  We are now in a position to state a very basic Makanin-Razborov type result:
  
  \newcommand{\mrd}[0]{\mathcal{MR}}
  \begin{proposition}
  	\label{MR} For any finitely generated \gr \rs group $\gqG{G}{A}{P}{I}$ there is a finite family $\mrd(\mathcal{G})$ of closed \gr \pr-\ress such that
  	\begin{align*}
  		\mathcal{G}\subset\bunion{\QR{R}\in\mrd}{Fct(\QR{R})}{}
  	\end{align*}
  \end{proposition}
  \begin{proof}
  	The proof is by induction on the well-founded partial order $\leq^{\pmb{r}}\cong\,\leq_{rk(G)}\times\leq_{Z}$.
  	If $\gqG{G}{A}{P}{I}$ is freely decomposable, let
  	\begin{align*}
  		G=\left(\bfrp{G_{l}}{l=0}{k}\right)\frp\F_{l}
  	\end{align*}
  	be a Grushko decomposition of $G$ relative to $\{A,P_{i}\}_{i\in I}$, where $A\leq G_{0}$ and no $P_{i}$ is conjugate into the free factor $\F_{l}$ .
  	Let $A(l)$ be $A$ in case $l=0$ and $\tg$ otherwise and $I(l)$ the set of those $i$ for which $P_{i}$ is conjugate into $A(i)$.
  	We can assume that $P_{j}\leq G_{l}$ for any $j\in I(l)$.
  	For each $0\leq l\leq k$ there is by induction a finite family of \gd \pr-\res $\{(\mathcal{R}^{j}_{l})_{A(l)}^{(J^{j}_{l},r^{j}_{l})}\}_{j\in\Lambda_{l}}$ such that
  	\begin{align*}
  		Mor(\qG{G_{l}}_{A(l)},\qG{K}_{A(l)})=\bunion{j\in\Lambda_{i}}{Fct^{sld}((\mathcal{R}^{j}_{l})_{A(l)}}{})
  	\end{align*}
  	Let $\gqG{(G_{l}^{j})}{A(l)}{P^{l}}{I_{j}}$ be the top \gr \pr-limit group appearing in $\QR{R^{j}}{}{}$. \
  	
  	For each choice of a tuple $\bar{j}=(j_{l})_{i=0}^{k}\subset\Lambda_{0}\times\cdots\Lambda_{k}$ we can construct a
  	\gd \pr-\res $(\mathcal{R}^{\bar{j}})_{A}^{(J^{\bar{j}},r^{\bar{j}})}$ of $\gqG{G}{A}{P}{I}$, where $\Delta\mathcal{R}^{r^{\bar{j}}}$ is trivial,
  	the children of $r^{\bar{j}}$ contain the groups $\{G_{l}^{j_{i}},\lambda\}_{l=1}^{k}$ and $F_{l}$ with the appropriate additional structure, the map $\eta^{R^{\bar{j}}}_{R^{\bar{j}}},r^{\bar{j}}$ is
  	the obvious quotient, $\mathcal{R}^{\bar{j}}_{A}\rst_{r^{j_{l}}_{l}}=(\mathcal{R}^{j_{l}}_{l})_{A_{l}}^{(J_{l},r^{j_{l}}_{l})}$
  	and $F_{l}$ is terminal.
  	Clearly any morphism from $\rqG{G}{A}$ to $\rqG{K}{A}$ factors solidly through some of the \pr-\ress constructed this way.
  	
  	Suppose now $\gqG{G}{A}{P}{I}$ is freely indecomposable. If $\gqG{G}{A}{P}{I}$ is solid add the \gr \res $\QR{R}^{sol}$ whose single node contains the pair $(\gqG{G}{A}{P}{I},\Delta_{JSJ})$, where $\Delta_{JSJ}$ is an abelian* JSJ decomposition of $\gqG{G}{A}{P}{I}$.
  	%  	and at least some of the morphisms in $\mathcal{G}$ are solid,
  	Chose a marking $x$ of $\rqG{G}{A}$ and let $\mathcal{F}$ be the finite family of maximal flexible limit quotients of
  	sequences from $\gqG{G}{A}{P}{I}$ to $\rqG{K}{A}$ with respect to $x$ and the action of $K$ on its Cayley graph.
  	
  	Each $\fun{q}{\rqG{G}{A}}{\qG{L^{q}}_{A}}\in\mathcal{F}$ is proper, so we can apply the induction hypothesis to $\gqG{L^{q}}{A}{P}{I}$, which results in a
  	finite set $\mrd^{q}$ of resolutions of $\rqG{L}{A}$ such that $Mor(\rqG{L^{q}}{A},\rqG{K}{A})=\bunion{q\in\mathcal{M}}{Fc^{sol}(R^{q})}{}$. For each $\LRL{R}{A}{(J,r)}\in\mrd^{q}$ construct a \gd \pr-resolution $\LRL{R'}{A}{(J',r')}$ of $\qG{G^{q}}{A}{P}{I}$ by attaching a new root $r'$ associated to $\rqG{L}{A}$ and letting $\DR{R'}{r'}=\Delta_{JSJ}$ and taking $q$ as $\eta^{r'}_{\QR{R'}}$.
  	Given any flexible $f\in Mor(\gqG{G}{A}{P}{I},\rqG{K}{A})$, there is by definition some $\sigma_{r}\in Mod(\gqG{G}{A}{P}{I})$ such that $f\circ\sigma_{r}=f'\circ q$ for some $q\in\mathcal{M}$.
  	
  	It is obvious how to extend some \cfam of $f'$ through $\QR{R}\in\mrd^{q}$ to some \cfam of $f$ through the corresponding $\LRL{R'}{}{}$.
  	In case $\rqG{G}{A}$ is not a \rs \pr-limit group we know that there is some finite set $\mathcal{MQ}(\rqG{G}{A})$ of quotients of $\rqG{L}{A}$ to $A$-\rs limit groups such that every
  	morphism from $\rqG{G}{A}$ to $\rqG{K}{A}$ factors through one of them. This clearly reduces this case to the previous one.
  \end{proof}
  
\section{Strict \pr-\ress }
  %  \begin{comment}
  \begin{definition}
  	Suppose we are given a \rs \pr-limit groups $\rqG{G}{A}$ and $\qG{H}{A}$, as well as a \gad $\Delta$ of $\rqG{G}{A}$, we say that a morphism $\fun{g}{\qG{G}}{\qG{H}}$ is \emph{strict} with respect to $G,\Delta$ if
  	for any sequence $\ssq{h}{n}$ of morphisms from $\rqG{G}{A}$ to $\qG{F}_{A}$ there is a sequence of automorphisms $\sigma_{n}\in Mod(\qG{G}_{A})$ such that $(h_{n}\circ g\circ\sigma_{n})_{n}$ has trivial limit kernel as well.
  	%  	In the case of a morphism of \rs \pr-limit groups $\fun{g}{\rqG{G}{A}}{\rqG{H}{A}}$, we require $G,\Delta$ to be a restricted graph of groups decomposition and  $\sigma_{n}\in PMod(\Delta;G;A)$.
  	%  	In the case of a morphism of \gr \pr-limit groups $\fun{g}{\gqG{G}{A}{P}{I}}{\gqG{H}{A}{P}{I}}$, we require $G,\Delta$ to be a graded graph of groups decomposition and $\sigma_{n}\in PMod(\Delta;G;A)$.
  	%	In absence of any explicit mention to $G,\Delta$, the canonical JSJ decomposition of the (restricted, graded) group will be considered.
  \end{definition}
  %  It is clear that considering \rs \pr-groups does not yield a stronger notion. It is slightly less easy to see that the \pr-structure does not play any role either.
  %  \end{comment}
  \begin{remark}
  	The former is equivalent to the fact that for any (some) sequence $(\fun{h_{n}}{\qG{G}}{\qG{F}})_{n}$ with trivial limit kernel
  	there are $\sigma_{n}\in Mod(\rqG{G}{A},\Delta)$ such that $\limker{h_{n}\circ h\circ\sigma_{n}}{n}=\tg$.
  \end{remark}
  %  \begin{definition}
  	\label{envelope}
  	Given a \gat $T$ of $G$ we let $\afld{T}$ be the $G$-tree obtained by folding together for each abelian type vertex $v$ all the edges incident to $v$ which be long to the same $Per^{*}(v)$-orbit. Denote by $\afld{e}$ and $\afld{v}$ the image in $\afld{T}$ of an edge $e$ and a rigid type vertex $v$ respectively. All those edges and vertices of $T$ not involved in the folding can be thought of as identical to their image in $\afld{T}$.
  	%  	 The presentation for $T$ associated to $\Delta$ projects to one of $\afld{T}$. We call the resulting decomposition $\afld{\Delta}$;
  	%  	the underlying graph of $\afld{\Delta}$ is identical to that of $\Delta$; the operation only produces an enlargement of groups associated to rigid type groups and edges between them and abelian type groups.
  	
  	Given a rigid vertex $v$ of $T$, $Stab(\afld{v})$,
  	is the subgroup of $G$ generated by $Stab(v)$ and $Per(v)$ for all $w\in VT_{a}$ of abelian type adjacent to $v$, what one
  	usually refers to as the envelope of $v$.  %and which we denote by $@$.
  	%[ this only makes sense for CmT groups ]
  %  \end{definition}
  
  \begin{definition}
  	\label{weakly strict} We say that a homomorphism between groups $G$ and $H$ is formally strict with respect to a \gat $T$ of $G$ if:
  	%  	($\rqG{G}{A}$) with covering tree $T$ is formally strict if and only if:
  	\enum{i)}{
  		\item For any rigid type vertex $v$ of $\afld{T}$, $f$ is injective on $Stab(v)$.
  		\item For any $e\in E$, $f$ is injective on $Stab(e)$
  		\item For any surface type vertex $v$ of $\afld{T}$ the image $f(Stab(v))$ is non-abelian.
  		%  		\item For any abelian type vertex $v$, $f$ is injective on $Per^{*}(v)$. %(peripheral subgroup of $Stab(\Delta;v)$, yet to define)
  	}
  	%% quote Abderezak
  \end{definition}
  
  The following fact is well-known in the standard case (see for example \cite{bestvina2009notes} and \cite{wilton2009solutions})).
  \begin{lemma}
  	%  	A homomorphism $\fun{f}{\qG{G}}{\qG{H}}$ ($\fun{f}{\rqG{G}{A}}{\rqG{H}{A}}$) between (\rs) limit groups is formally strict with respect to a \gad $\Delta$ of $\qG{G}$ ($\rqG{L}{A}$) if and only if it is strict.
  	A morphism $f$ from a \rs \pr-group $\rqG{G}{A}$ and a \rs \pr-limit group $\rqG{L}{A}$ which is formally strict is also strict.
  	
  \end{lemma}
  \begin{proof}
  	The implication from right to left is immediate. Given any of the standard generators of $\tau\in Mod(\rqG{R}{A},\Delta)$ it is easy to check by inspection that $\tau$ restricts to an inner automorphism on:
  	\enum{i)}{
  		\item $Stab(v)$ for any rigid type vertex $v\in\afld{T}$.
  		\item $Stab(e)$ for any $e\in E$
  		%  		\item $Per(T;v)$ for any abelian type vertex $v$
  	}
  	Hence the kernel of the restriction of $f$ and $f\circ\tau$ to any of those groups coincide. Moreover, for any surface type vertex $v\in\afld{T}$,
  	$\tau$ takes the image of $Stab(v)$ to a conjugate in $G$, hence $f\circ\tau$ sends it to an abelian group if and only if $f$ does.
  	
  	The converse can be shown using a generalization of the method fleshed out in \cite{wilton2009solutions}. See section \ref{discrimination section} for more details.
  	%  	For the converse can be seen as a consequence of the existence of test sequences of \pr-towers and the embeddability of the top limit group of a strict resolution into its completion (see \label{strict is strict}),
  	%  	although the more direct method used in \cite{wilton2009solutions} easily adapts to the constrained case as well.
  	%  	For the general case the reader is referred to \cite{Bestvina Feighn notes} and \cite{solutions Bestvina Feighn}.
  \end{proof}
  In particular, if such a formally strict morphism exists $\rqG{G}{A}$ is also a \rs \pr-limit group.
  \begin{definition}
  	A  \rs \pr-\res $\LRL{R}{A}{(J,r)}$ will be called strict if and only if each map $\eta^{\lambda}_{\lambda}$ is strict with respect to $\GRla$. It will be called quasi-strict if the property holds for any $\lambda\neq r$ and the root $r$ has a single child and is associated to a trivial decomposition.
  \end{definition}
  %  It turns one can factor every morphism to the free group through finitely many \pr-\ress which only fail to be strict at the top:
  \begin{lemma}
  	\label{quasi strict MR}
  	For any \gd \pr-group $\gqG{G}{A}{P}{I}$ there is a finite family $\mathscr{S}$ of closed quasi-strict \gd \pr-\res of $\gqG{G}{A}{P}{I}$ such that any morphism from $\rqG{G}{A}$ to $\rqG{K}{A}$ factors solidly through one of them.
  \end{lemma}
  \begin{proof}
  	\newcommand{\M}[0]{}
  	Start with each $\LRL{R}{(J,r)}\in\mathcal{MR}(\gqG{G}{A}{P}{I},\mathcal{F})$ of $\gqG{L}{A}{P}{I}$ which is not strict.
  	For each $\lambda\in (J,r)$ let $\mathcal{H}_{\lambda}=Fct(\QR{R}\rst_{\lambda})\subset Mor(\RAl,\qG{F}_{A_{\lambda}})$. The fact that the resolution is not
  	strict implies that for some $\lambda\in J$ any limit quotient of morphisms in $\mathcal{H}_{\lambda}$ is proper;
  	let $\lambda_{0}$ be some maximal index in $J$ for which this is the case and $\mathcal{MQ}_{\lambda}$ the associated finite set of maximal quotients.
  	
  	For each $\fun{p}{\RAl}{\qG{L}_{A_{\lambda}}}\in\mathcal{MQ}_{\lambda}$ we consider all possible resolutions obtained by replacing %$\RAl$ with L$\lambda\in\Lambda\lambda_{0}$,
  	$\QR{R}\rst_{\lambda}$ with some \pr-\res $\LRL{R'}{}{(J',r')}\in\mathcal{MR}(\qG{L^{p}}_{A_{\lambda}})$ in the obvious sense.
  	%  	\footnote{In the following sense: if $\lambda$ has a parent $\mu$, then we need to postcompose the map
  	%  	at $\mu$ with the map $\eta\circ p\frp Id$, where $\eta=\eta_{R'}^{r'}$ if $\DR{R'}{r'}$ is trivial and the identity otherwise. In the first case the new child of $\lambda$ is not $r'$ but its unique child.}
  	\begin{comment}
  	%  	 and postcomposing $\eta^{\mu}&} with the appropriate quotient for any parent of a node in $J_{0}$.
  	%  	replacing $\eta^{\mu}&} for $\mu$ the father of some $\lambda\in\Lambda_{0}$ by postcomposing it with the quotient restricting to $Id_{Rm}$ on $Rm$ for $\lambda\nin\Lambda_{0}$ and $q_{\lambda}$ whenever $\lambda\in\Lambda_{0}$.
  	%the groups associated with those of its children not in $\Lambda_{0}$ and to the map $q_{\lambda}$ of its children in $\Lambda_{0}$.
  	%  	It is a very simple check that any morphism factoring through $\QR{R}$ factors through one of the resolutions $\LRL{R^{\bar{q}}}{}$ so obtained.
  	\end{comment}
  	Let $Rep_{\lambda}(\QR{R})$ be the finite family of all the \pr-\ress obtained as we let $p$ range within $\mathcal{MQ}_{\lambda}$.
  	Clearly any morphism factoring through $\QR{R}_{A}^{\lambda}$ factors through some $\LRL{R'}{A}{}\in Rep_{\lambda}(\QR{R})$.
  	
  	Iterate the operation in parallel on each $\QR{R}\in Rep_{\lambda}(\QR{R})$ as long as it is possible. We can represent the course of this substitution procedure using a rooted tree $(\mathfrak{T},\rho)$ each of whose nodes $\tau$ indexes some \pr-\res  $\LRL{R^{\tau}}{}{}$. Each pair  $\tau\prec\tau'$ corresponds to the derivation of
  	$\mathcal{R}^{\tau}$ from $\mathcal{R}^{\tau'}$ by performing a substitution as described above.
  	
  	If we manage to prove that $\mathfrak{T}$ is finite, then we are done. Any resolution appearing at a leaf of $\mathfrak{T}$ will be strict.
  	Moreover, in that case the union of the \pr-\ress which appear at the leaves of $\mathfrak{T}$ at
  	the starting \pr-\res ranges within $\mathcal{MR}(\gqG{G}{A}{P}{I})$ will meet our requirements.
  	%  	Since the tree is maximal, the resolutions appearing at any of its leaves are strict. Observe also that if $\mathfrak{T}$ is finite $\mathscr{F}(\LRL{R^{\mu}}{},\LRL{R^{\mu}}{})=\bunion{\mathscr{F}(\LRL{R^{\lambda}}{},\LRL{R^{\lambda}}{})}{\lambda\in\lvs{\mathfrak{T}}}{}$.
  	In virtue of K\"onig's lemma this amounts to showing that $\mathfrak{T}$ contains no infinite branches. A \pr-\res resolution can be seen as an element of the partial order $Tr:=Tr(PL,\leq_{rk}\times\leq_{Z})$ described in subsection \ref{poset section}, so the result follows from the somewhat more general fact that $(Tr,\leq_{Tr})$ is well-founded (the map $\tau\mapsto \LRL{R^{\tau}}{}{}\in Tr$ is strictly monotonous by construction).
  	%  	A resolution $\LRL{R'}{}$ obtained from another one $\QR{R}$ by the replacement described above is strictly smaller than $\QR{R}$ with respect to $\leq_{Tr}$, where the set of $\mu$ above some $\lambda\in\Lambda_{0}$ plays the role of $S$ in the definition of $\leq_{Tr}$.
  \end{proof}
  %  In the proof above the use of \ref{order on trees} on its full generality is superfluous, since one ca  always gather a branching resolution as defined here into a linear one.
  
  Although in general it is not possible to capture all morphisms from $\rqG{L}{A}$ to a free $\rqG{K}{A}$ using finitely many strict \rs \pr-\ress, there is always at least one.
  \begin{lemma}
  	Let $\rqG{L}{A}$ be a freely indecomposable \rs limit group and $\Delta_{JSJ}$ its JSJ decomposition. Then $\rqG{L}{A}$ admits a non-injective $\Delta_{JSJ}$-strict morphism $\fun{\phi}{\rqG{L}{A}}{\rqG{K}{A}}$. As a consequence, any such $\rqG{L}{A}$ admits a strict \pr-\res $\LRL{R}{A}{(J,r)}$.
  \end{lemma}
  \begin{proof}
  	Fix a tuple of generators of $L$. For each $n\in \N$ we know a morphism $\fun{f_{n}}{\rqG{L}{A}}{\rqG{F}{A}}$ which is inyective on the ball of radius $n$ in the Cayley graph. Now take
  	$\tau_{n}\in Mod(\rqG{L}{A},\Delta_{JSJ})$ such that $f_{n}\circ\tau_{n}$ is short in the sense of \ref{shortening morphism section}. Any convergent subsequence of $(f_{n})_{n}$ will have non-trivial limit kernel.
  	The second statement can be proven by induction using equational noetherianity, as usual.
  \end{proof}
  
\section{Well-separated \pr-\ress }
  \begin{comment}
  Have to dfine restrictions of graphs of groups and the embedding of their fundamental group.
  In the relevant section. Be careful to distinguish between surface groups of the graph and those of the tree
  \end{comment}
  %  \begin{definition}
  	A \ws \gr \pr-\res is a special type of quasi-strict \gr \pr-\res $\LRL{R}{}{(J,r)}$ to which is possible to assign an completion (see section \ref{completions subsection}).
  	%  	, together with some extra information.
  	%  	refines the notion of a \gr \pr-\res.
  	More precisely, for any $\lambda\in \hat{J}$ not of free product type there is a distinguished set
  	\begin{align*}
  		DFt(\lambda)\subset Ch(\lambda)\cap\lvs{J}
  	\end{align*}
  	so that the group at $\mu\in DFt(\lambda)$ is free and the following conditions hold:
  	%  	let $\GRla=(\Delta^{\lambda},Y^{\lambda},\phi^{\lambda})$.
  	%  	Then either $\eta^{\lambda}_{\QR{R}}$ is an isomorphism and $\GRla$ is trivial or:
  	\enum{i)}{
  		%      	\item  $\GRla$ is normalized in the sense of $\label{cylinders}$  (this assumption is not essential an essential one, but will make the following more convenient)
  		%      	\item For every $e\in E$ the image of $\GRla_{e}$ by $\eta^{\lambda}_{\lambda}$ is elliptic in $\bfrp{Rm}{\mu\prec\lambda}{}$.
  		\item For every vertex $v$ of non-surface type the image of $\GRla_{v}$ by $\eta^{\lambda}_{\lambda}$ is contained in $\Rm$ for some $\mu_{v}\prec\lambda$.\footnote{It is important that it is contained as opposed to merely conjugate into it. }
  		\item For $\mu\prec\lambda$, if $\Rm$ does not intersect the image of a non-surface type vertex $\mu\in\lvs{(J)}$. \label{dropped}
  		Moreover, $\Rm$ is a free factor of the image of the subgroup generated by some unique surface type vertex group $\GRla _{v}$ and the Bass-Serre elements associated to edges incident to $v$.  \label{situ one}
  		\item If $\GRla$ consists of a single vertex of surface type     %\,I M $P$ O $R$ T\,A\,N
  		or non-cyclic abelian type, then $Ch(\lambda)=DFt(\lambda)$. \label{situ two}.
  		\item If $\Rl$ is free and does not carry a constant or parameter group, then $\lambda\in\lvs{J}$.
  	}
  	We will refer to $\lambda\in DFt(\lambda)$ as a dropping node at $\lambda$ and to $\Rl$ as a factor dropping at $\lambda$.
  	Let also $nDFT(\lambda)=Ch(\lambda)\setminus DFt(\lambda)$.
  	
  	%  	If the resolution is endowed with additional structure we can talk of \rs or \gr \pr-\res.
  	Observe that the second condition implies	that no $\lambda\in DFt(\lambda)$ can carry any parameters or constants. We will let $DFt(T)=\abunion{\lambda\in \hat{J}}{DFt(\lambda)}{}$ and
  	$DFT=(J)\setminus DFt(T)$.
  	
  %  \end{definition}
  
  Note that leaves in $DFt(\lambda)$ are not necessarily the only ones for which $\Tm$ is free.
  Also that the fact that the resolution is strict implies $\eta^{\lambda}_{\QR{R}}$ is injective on $P=Per^{*}(v)$, whose image is necessarily elliptic in the free decomposition
  $\bfrp{\Rm}{\lambda\succ\mu}{}$, as otherwise the set of its elliptic elements is the kernel of a map onto $\Z$ containing the finite index subgroup generated by all incident edge groups (see \ref{abelians acting}). If $Per^{*}(v,\GRla)$ is cyclic, the fact that $\QR{R}$ is \ws implies that the image of some finite index subgroup of $P$ is elliptic, implying that the same holds for $P$ itself.
  This leads to the following observation:
  \begin{remark}
  	\label{connectedness remark}Given any path $[u_{1},u_{2},\cdots u_{m}]$ of the underlying graph $(V,E)$ of $\GRla$ containing only abelian an rigid subgroups, all $\Delta_{u_{j}}$ for $u_{j}\in V_{s}$ or $Per(u_{j})$ for $u_{j}\in V_{a}$ project into the same factor $\Rm$. Moreover, the Bass-Serre element associated to the edge $e_{i}$ from $u_{i}$ to $u_{i+1}$ maps into $\Rm$ as well.
  \end{remark}
  %      [move this to right after the definition of well structured; define what it means for a map to be compatible]
  
  Define $DFt(R,\lambda)$ as the set of $\mu$ such that $\Rm$ does not contain the image of a non-surface type vertex of $\GRla$, and $nDFT(\lambda)$ as $Ch(\lambda)\setminus DFt(R,\lambda)$.
  
  By the dropping component at $\lambda$ we intend $DCt(R,\lambda):=\bfrp{\Rm}{\mu\in DFt(R,\lambda)}{}$ and by the non-dropping component at $\lambda$ we intend $NDCt(R,\lambda):=\bfrp{\Rm}{\mu\in nDFT(\lambda)}{}$.
  
  \begin{comment}
  \begin{definition}
  	Suppose that for $i=0,1$ we have a pair $(G_{i},\cdot_{i})$, where $G_{i}$ is a group and $\cdot_{i}$ an action without inveersions of $G_{i}$ on a simplicial tree $G_{i}$.
  	A morphism from $(G_{0},\cdot_{0})$ to $(G_{1},\cdot_{1})$ consists of a pair $(\phi,f)$.
  \end{definition}
  \end{comment}
  
\section{Taut \pr-\ress }
  \begin{comment}
  Let $\Delta$ be a \gad of a group $G$ with underlying graph $(V,E)$.
  \begin{comment}
  \begin{definition}
  	By a taut structure on $\Delta$ we intend a map $\Theta$ assigning a \SEC $\TS{\Theta}{v}\in\Ess{\Sigma_{v}}$ to each $v\in 	V_{s}$.
  \end{definition}
  \begin{definition}
  	By a taut structure on $\Delta$ we intend a collection $\Theta$ of essential s.c.c. in the surfaces $\Sigma_{v}$ corresponding 	to surface type vertices of $\Delta$, such that $\Theta\rst_{\Sigma_{v}}\in\Ess{\Sigma_{v}}$ for each $v\in V_{s}$.
  \end{definition}
  \end{comment}
  The fact that a morphism factors through a resolution provides plenty of information about it. Unfortunately, not enough for the verification process of $\forall\exists$ formulas, which relies on a fine control of both completion and shortening. This is what motivates the introduction of taut \pr-resolutions; properly speaking they consist of a resolution together with some additional information
  on each of its nodes, which determines a restricted notion of factoring.
  %  These do not just only constitute a restricted class of resolutions: they contain additional information which gives a more precise control of the way morphisms factor through it.
  %  Suppose now that we are also given a \pr-structure $\qG{G}$ on $G$. Each element of $PMod(\Delta;G)$ determines some class in $Out_{\partial}(\pi_{1}(\Sigma_{v}))$ and consequently on in $Mod(\Sigma_{v})$.
  %  \newcommand{\TS}[2]{#1\cap\Sigma_{#2}}
  Suppose we are given a group $G$, a \gad $\Delta$ of $G$, an \cfam $\delta$ in $\Delta$. Recall that in this situation the image $G^{p}_{\delta}$ admits a free decomposition of the form:
  \begin{equation}
  	\label{taut factoring}(\bfrp{A_{i}}{i=1}{m})\frp(\bfrp{\pi_{1}(\Pi)}{\Pi\in Pinch^{ex}(\Delta,\delta)}{})\frp\F_{b_{1}(X^{c}_{\delta})}
  \end{equation}
  Given a homomorphism $\fun{f}{G}{H}$ we say that $f$ is taut with respect to $\delta$ if:
  \elenco{
  	\item $f$ pinches $\delta$
  	\item $\bar{f}$ does not pinch any essential simple closed curve in any of the surfaces in $Pinch^{in}(\Delta,\delta)$
  	\item The external rank of $\Delta^{p}_{\delta}$ is maximal among all possible choices of $\delta$ for which the first two properties are satisfied.
  }
  Suppose now that we are also given a \pr-structure $\qG{G}$ on $G$. For any surface type vertex $v$ the group $Mod(\qG{G},\Delta)$ restricts to some finite index subgroup of
  %  $Out^{\delta}(\Delta_{v})$,
  $Out^{\delta}(\pi_{1}(\Sigma_{v}))$, corresponding to a finite index subgroup of $Mod(\Sigma_{v})$. For each such $v$ chose a set $\mathcal{X}_{v}$ of representatives of the action of the latter on the family of all \cfams on $\Sigma_{v}$; this is finite.  Denote by $Ess(\qG{G},\Delta)$ (the result does not depend on $A$) the family of all the \cfams on $\Delta$ obtained by picking some $\delta_{v}\in\mathcal{X}_{v}$ for each $v$ of surface type.
  %  For each $v$, let $Ess(\Sigma_{v})$ be a set of representatives of the orbits of (possibly empty) \cfams in $\Sigma_{v}$ by the action of $\PPMod{\Sigma_{v}}$.
  The discussion above implies that this is finite. By a taut structure on $\Delta$ we intend an \cfam in $\Delta$ such that $\Sigma_{v}\cap\delta\in Ess(\Sigma_{v})$ for all vertices $v$.
  %    \newcommand{\TS}[2]{#1\cap\Sigma_{#2}}
  %    or equivalently, the rank of its fundamental group
  \newcommand{\TRAT}[0]{\mathcal{R}_{A}[\Theta]}
  \newcommand{\TR}[0]{\mathcal{R}[\Theta]}
  \begin{definition}
  	\label{taut resolutions} A taut \rs \pr-\res of a \rs \pr-group $\rqG{G}{A}$ is given by  $(\LRL{R}{A}{(J,r)},\delta^{\lambda})_{\lambda\in (J)}$, comprising a \ws \rs \pr-\res $\LRL{R}{A}{(J,r)}$ and for each $\lambda\in J$ a taut structure on $\delta^{\lambda}$ on $\Delta$ such that for $\lambda\in\hat{J}$:
  	\elenco{
  		\item $\eta^{R}_{\QR{\lambda}}$ is taut with respect to $\delta^{\lambda}$
  		\item $\setof{\Tm}{\mu\in nDFT(\lambda)}$ are precisely the images of those free factors in the first term of \ref{taut factoring}.
  		\item $\setof{\Tm}{\mu\in DFt(\lambda)}$ are each one a maximum rank free image of the remaining factors.
  	}
  	We will refer to $(\Theta^{\lambda})_{\lambda}$ as a taut structure on $\RRA$.
  	
  	By a \trm of $\mathscr{R}_{A}=\mathcal{R}_{A}[\Theta]$ we intend a \trm $\fcrm{(J)}{f}{\tau}$ of $\RRA$ such that for all
  	$\lambda\in (J)$
  	%  	$\lambda\in\hat{J}$
  	$f^{\lambda}\circ\sigma^{\lambda}$ is taut with respect to $\delta$ for some $\sigma^{\lambda}\in Mod^{\pi}_{A_{\lambda}}(\GRla)$,
  	where $\sigma^{\lambda}=(\tau^{\lambda})^{-1}$ in case $\lambda\in\hat{J}$. A morphism from $\rqG{G}{A}$ to $\rqG{F}{A}$ factors through $\TRAT$ if and only if it belongs
  	to some taut \trm of $\TRAT$. If we are additionally given an additional \gr structure on the underlying \pr-\res talk about a taut \gr
  	\pr-\res.
  \end{definition}
  Recall we are allowed to have one non-strict map at the top of the \pr-\res, provided the \gad at the top is trivial.
  We will say that a taut \gr \pr-\res $\TR$ is closed if $\QR{R}$ is closed as a graded \pr-\res and $\delta_{\lambda}=\nil$ at any leaf $\lambda$.
  \begin{lemma}
  	Suppose we are given a \gr \pr-limit group $\gqG{G}{A}{P}{I}$ and a family $\mathcal{G}$ of morphisms from $\rqG{G}{A}$ to some free $\rqG{K}{A}$.
  	Then there is a finite family $\mathcal{TMR}(\mathcal{F})$ of closed taut resolutions with the following properties:
  	%  	of pairs $(q,\LRL{R}{A}{I})$; there are finitely
  	\enum{i)}{
  		\item Any morphism in $\mathcal{G}$ factors solidly through one of them. \label{factorization}
  		%  		from $\rqG{G}{A}$ to $\rqG{F}{A}$
  		%  		\item $\DR{R_{\alpha}}{\lambda}$ is a graded abelian* JSJ of $\RLR{R_{\alpha}}{A}{\lambda}$ whenever the later is not free with trivial parameter and constant subgroups.
  		\item \label{constructed from G}For each $\LRL{R}{}{(J,r)}$ and any $\delta^{\lambda}$ the restriction of $\eta^{\lambda}_{\QR{R}}$ to each of the inner factors of $\Delta^{c}_{\delta}$ is
  		the limit quotient of a sequence of factorizations at $\lambda$ of morphisms in $\mathcal{G}$. \footnote{This is in fact stronger than what will be needed in the analysis of $\forall\exists$-sentences.}
  	}
  \end{lemma}
  \begin{proof}
  	To begin with, define the notion of a pseudo-taut \gr \pr-\res in the same way as a taut one, by replacement the requirement that the \pr-\res is strict with the weaker one that
  	the groups assigned to those edges of $\GRla$ adjacent to a surface type vertex group to be mapped injectively by $\eta^{\lambda}_{\QR{R}}$. We modify the definition of a \trm accordingly:
  	the map $f^{\lambda}$ is not required to be injective on boundary subgroups of surface type vertex groups of $\GRla$.
  	\begin{claim}
  		Given any $\gqG{G}{A}{P}{I}$ and $\mathcal{G}$ as above, there is a finite family $\mathcal{MR}^{qt}(\mathcal{G})$ of closed quasi-taut \gd \pr-\ress such that every morphism in $\mathcal{G}$
  		factors solidly through one of them and property (\ref{constructed from G}) holds.
  	\end{claim}
  	\begin{proof}
  		The proof is by induction on the $\leq_{rk}\times\leq_{Z}$ of $G$, using a very similar construction to that used in \ref{MR}. In view of this, we will merely underscore the differences between both constructions. By the arguments employed there we can assume that $\gqG{G}{A}{P}{I}$ is a freely indecomposable \gr \pr-limit group, with abelian* JSJ decomposition $\Delta^{JSJ}$.
  		
  		First of all, for any $e\in\Delta^{JSJ}$ one of whose endpoints is of surface type let $\mathcal{G}_{e}$ the family of all those morphism in $e$ which kill $\Delta_{e}$. Let $\mathcal{MQ}_{e}$
  		the finite family of maximal limit quotients of sequences of morphisms in $\mathcal{G}_{e}$. Applying the induction hypothesis to each of them we can easily obtain a family $\mathcal{MR}^{qt}(\mathcal{G}_{e})$ satisfying the required conditions with respect to the family $\mathcal{G}_{e}$.
  		%  		Each of them is obviously non-trivial.
  		Let $\mathcal{G}_{ndeg}$ be the family of all those $f\in\mathcal{G}$ which do not belong to any of the $\mathcal{G}_{e}$ above.
  		We can write $\mathcal{G}_{ndeg}=\bunion{\delta\in Ess}{\mathcal{G}_{\delta}}{}$,
  		where $\mathcal{G}_{\delta}$ is the family of all those $f\in\mathcal{G}$ such that $f\circ\tau$ pinches
  		$\delta\in Ess:=Ess(\qG{G},\Delta^{JSJ})$ for some $\delta\in Ess$ and the external rank associated to $\delta$ is maximal for that property.
  		If $\delta\neq\nil$, then consider the inner factors of $G_{1}\cdots G_{k}$ of $G^{p}_{\delta}$. We can assume that $A\leq G_{1}$ and $P_{i}\leq G_{j_{i}}$ for
  		any $i\in I$ and some $1\leq j_{i}\leq k$. For each $f\in\mathcal{G}$ denote by $\bar{f}$ the unique map such that $f=\bar{f}\circ p^{\delta}$
  		and for any $1\leq j\leq k$ by $\mathcal{G}_{\delta,j}$ the collection $\setof{\bar{f}\rst_{A_{j}}}{f\in\mathcal{G}}$.
  		Let $\mathcal{MR}_{j}$
  		%  		induction hypothesis we can construct a family of
  		the finite collection of \gr \pr-\ress of $\gqG{G}{A_{i}}{P}{I_{i}}$ obtained applying the induction hypothesis to the $\mathcal{G}_{\delta,j}$, where $A_{i}=A\cap G_{i}$ and $I_{i}=\setof{j}{i_{j}=i}$.
  		Given any tuple $(\mathscr{R_{j}})_{j}\in\mathcal{MR}^{qt}_{1}\times\cdots\mathcal{MR}^{qt}_{k}$ for all $1\leq j\leq k$ we can construct a pseudo-taut \pr-\res of $\gqG{G}{A}{P}{I}$ with
  		$(\gqG{G}{A}{P}{I},\Delta_{JSJ},\delta)$ at its root $r$, where $\eta^{r}_{\QR{R}}=(q_{1}\dcup q_{2}\cdots q_{k})\circ p^{\delta}$
  		and $q_{j}$ is the quotient map of $G_{j}$ onto $\mathcal{R}_{j}$. Denote by  $\mathcal{MR}^{qt}_{\delta}$ the collection of all of them.
  		If $\delta=\nil$ then we proceed exactly as in the proof of \ref{MR}. If $\gqG{G}{A}{P}{I}$ is solid, we have to take both the trivial \pr-\res where $\gqG{G}{A}{P}{I}$ is the only vertex
  		and those coming from applying the induction hypothesis to the maximal flexible $\mathcal{G}$-quotients. By the latter we mean the limit quotients of flexible sequences whose members are of the form
  		$f\circ\tau$, where $f\in\mathcal{G}$ and $\tau\in Mod(\gqG{G}{A}{P}{I})$.
  		
  		It is easy to check that the union $(\bunion{e\in E\Delta^{JSJ}}{\mathcal{MR}^{qt}_{e}}{})\cup (\bunion{e\in Ess}{\mathcal{MR}^{qt}_{\delta}}{})$ meets all of our demands.
  	\end{proof}
  	One can use the same argument as in the proof of \ref{quasi strict MR} to replace each $\mathcal{R}$ in $\mathcal{MR}^{qt}$ by finitely many taut \pr-\res capturing all the morphisms which factor through it. This time each of the nodes in the  rooted tree $(\mathfrak{T},\leq,r)$ is assigned some pseudo-taut $\mathscr{R}_{\alpha}^{J_{\alpha}}$, satisfying property (\ref{factorization}) of the statement of the lemma. Denote by $\mathcal{G}_{\alpha}$ the family of all the morphisms of $\mathcal{G}$ which factor through it.
  	
  	There is essentially one single obstruction to the termination of the process at any given node $\alpha\in\mathfrak{T}$, namely,
  	that for some $\lambda\in J_{\alpha}$ the map $\eta^{\lambda}_{\QR{R}}$ is not the limit of a sequence of factorizations at $\lambda$ of members of $\mathcal{G}_{\alpha}$; note
  	that this is necessarily the case if $\LRL{R_{\alpha}}{}{}$ fails to be strict \footnote{We can always assume the root has a trivial decomposition and a single child \pr-\res.}.
  	
  	Chose some maximal $\lambda$ witnessing this failure. Let $\Delta=\GRla$ and $\bar{G}^{\lambda}_{\alpha}$ be the set of all factorizations at $\lambda$ of morphisms from $\mathcal{G}_{\alpha}$.
  	%  	The properties of $\LRL{R_{\alpha}}{}$ imply that
  	This implies for some child $\mu$ of $\lambda$ which is the image of an inner factor $G_{\mu}$ of
  	$\Delta^{p}_{\delta^{\alpha}_{\lambda}}$ all $f\in\bar{G}^{\lambda}_{\alpha}$ factor through $q\circ \eta^{\lambda}_{\QR{R_{\alpha}}}$ for $q$ in some finite family
  	$\mathcal{MQ}_{\alpha}^{\lambda}$ of proper quotients of $\RLR{R_{\alpha}}{}{\lambda}$. At each children of $\alpha$
  	we place the taut \pr-\res obtained by replacing (the construction is described in more detail in the proof of \ref{quasi strict MR}) $\LRL{R_{\alpha}}{}{}\rst_{\lambda}$ by some member of $\mathcal{MR}^{qt}(\bar{\mathcal{G}}_{\mu})$ (composing $\eta^{\lambda}_{\QR{R_{\alpha}}}$ with the top quotient map of the latter). Just as in \ref{quasi strict MR},
  	the process must eventually stop.
  \end{proof}
  %  As a matter of fact, one could use \ref{tree preorders} directly in the proof of the theorem, by applying to partially constructed (not necessarily closed) taut \pr-\res.

% include well-structuredness into the tower definition
% allow for \ws weak resolutions
\chapter{\pr-Towers}
\label{Towers chapter}
\section{Definition}
  
  Although \pr-limit groups admit a fairly regular description in terms of resolutions, as seen in the previous section, for many purposes one must restrict oneself to the smaller class of those
  having a \pr-tower structure, which we now go on to define.
  
  \df{
  	By a \pr-tower structure we will intend a (a priory) weak \ws \pr \res $\mathcal{T}=\LRL{T}{}{(J,r)}$ where for
  	each $\lambda\in (\hat{J})$ not of free product type
  	the following conditions are satisfied:
  	\enum{a)}{
  		\item The graph $X=(V,E)$ underlying the decomposition $\GTla$ is a bipartite graph, $V=V_{1}\coprod V_{2}$, where $V_{1}\neq \nil$ and $V_{2}$ consist precisely of those vertices of rigid type.
  		\item Each $v\in V_{2}$ of abelian type has incidence number at most $1$ (incidence number $0$ implies the whole group is abelian).
  		%  		\item $\GTla$  does not consist of a single vertex of abelian t
  		\item There is a bijective correspondence between $V_{1}$ and the set $nDFT(\lambda)$ assigning to each $v\in V_{1}$
  		some $\mu_{v}\in nDFT(\lambda)$ for which $\GTla_{v}=\RLR{T}{}{\mu_{v}}$.
  		\item The map $\eta^{\lambda}_{\QR{R}}$ restricts to the identity on $\GTla_{v}$ for any $v\in V_{1}$.
  		Observe that the collection $\GTla_{v}$ for $v\in V_{1}$ generate a group isomorphic to their free product.
  	}
  	A \rs (\gr) \pr-tower is simply a \rs (\gr) \pr-\res whose underlying \pr-\res is a \pr-tower.
  	
  	If a \rs \pr-tower $\RTA$ is closed, we will refer to it as a \pr-tower structure for $\TAl$.
  	If a collection $\Lambda\subset\lvs{J}\setminus DFt(\mathcal{T})$ exists such that declaring $\TAl$ for $\lambda\in\Lambda$ as a parameters
  	makes $\RTA$ into a closed \gr resolution, we refer to it as a tower structure for $\RLR{T}{A}{r}$ relative to $\{\Tl\}_{\lambda\in\Lambda}$.
  	Conversly, we will refer to $\RLR{T}{A}{r}$ as the top of the \pr-tower.
  }
  The way towers are described in most of the literature, the fundamental group of a closed surface must belong to the very bottom of the tower, while the existence of a homomorphism with non-abelian image from this fundamental group to the free group with is expressed in terms of the surface having characteristic $\leq-2$.
  \newcommand{\tlf}[1]{tlf(#1)}
  
  Let now $\RTA$ be a \rs \pr-tower relative to the collection $\{\Tl\}_{\lambda\in\Lambda}$.
  Observe that if in the situation above $\lambda\in nDFT(\lambda)$, and $\mu\succ\lambda$, then $\Tl$ is a rigid vertex of $\DR{R}{\mu}$ and hence is contained in
  $\Tm$. Iterating the argument one can easily see that in fact $\Tl\leq \RLR{T}{}{r}$ in this case.
  Denote by $\tlf{J}$ the set of minimal nodes in $nDFT(\mathcal{T})$. The group generated by $\{\Tl\}_{\lambda\in\tlf{J}}$ is isomorphic to $\bfrp{\Tl}{\lambda\in\tlf{J}}{}$.
  
  \newcommand{\bott}[1]{bot(#1)}
  \newcommand{\FHA}[0]{free hanging non-cyclic abelian }  %tc
  \newcommand{\FHAS}[1]{FHA(#1)}
  \begin{comment}
  Observe there are several possibilities for $\lambda\in nDFT$.
  \enum{a)}{
  	\item If $\lambda\in\tlf{R}$, then either:
  	\enum{i)}{
  		\item $\lambda\in\Lambda$
  		\item $\qG{A}\leq\Tl$
  		\item $\Tl$ is free and $A\cap\Tl=\tg$
  		\item $\GTla$ contains a single vertex of surface type, $A\cap\Tl=\tg$ and $\Tl$ is isomorphic to a closed surface group
  		\item $\GTla$ contains a single vertex of abelian type, $A\cap\Tl=\tg$ and $\Tl$ is isomorphic to a non-cyclic free abelian group.
  		We willl refer to $\Tl$ in this case as an isolated abelian subgroup and to $\lambda$ as a node of isolated abelian type.
  	}
  	\item If $\lambda\nin\Lambda$ then either:
  	\enum{i)}{
  		\item $\lambda$ is of free product type
  		\item $\GTla$ contains vertex groups of rigid type, in which case be refer to $\lambda$ as a node of \emph{floor type}
  	}
  }
  \end{comment}
  
  As a convention, we denote the \pr-group on top of a closed \pr-tower with the same letter, so unless otherwise specified, by $\rqG{T}{A}$ we will always mean $\RLR{T}{A}{r}$ and so on.
  We will also sometimes improperly refer to $\rqG{T}{A}$ itself as a \pr-tower instead of as a group admitting a \pr-tower structure.
  
\section{Closures}
  
  %  Let $\LRL{T}{A}{(J,r)}$ be a \pr-tower structure relative to a family of \rs \pr-limit groups $\qG{L_{1}}_{A_{1}},\qG{L_{2}}\cdots\qG{L_{m}}$.
  Let  $\LRL{T}{A}{(J,r)}$ be a \rs \pr-tower structure and $M$ a maximal non-cyclic abelian subgroup of $T$.
  Clearly the set of those $\lambda\in nDFT$ for which some conjugate of $M$ intersects $\Rl$ non-trivially
  is a descending branch $r=\lambda_{0}\succ \lambda_{1}\succ \lambda_{2}\succ\cdots \lambda_{m}$ in $J$.
  Up to replacing $M$ by some of its conjugates, we can assume that $M_{\lambda_{j}}=M\cap \RLR{T}{}{\lambda_{j}}$ is non-trivial.
  %  	[is the first part needed?]
  Given $\lambda$, let $r_{j}=rk(M_{\lambda_{j}})$. We claim that $r_{j}>r_{j+1}$ for some $0\leq i\leq m-1$
  if and only if $M$ is conjugate to $Z_{G(A})$, where $A$ is some abelian vertex group of $\DR{T}{\lambda_{j}}$ attached to the rigid vertex group $\RLR{T}{}{\lambda_{j+1}}$ of $\DR{T}{\lambda_{j}}$. It is easy to see that there can be no growth in the case $\lambda_{j}$ is not of free-product type, so we might as well assume it is of floor type.
  All non-cyclic abelian subgroups of $\RLR{T}{}{\lambda_{j}}$ are elliptic in $\DR{T}{\lambda_{j}}$.
  
  Now, since $M\cap \RLR{T}{}{\lambda_{j+1}}\neq\tg$, if $M$ did not intersect any abelian vertex group of $\GTla$, then the rigid vertex group $\RLR{T}{}{\lambda_{j}}$ would be the only one into which one might conjugate $M_{\lambda_{j}}$, contradicting the growth in rank. On the other hand, centralizers in surface groups are cyclic, so necessarily $M$ contains a conjugate of an abelian vertex group, as desired.
  Let us note in passing that this shows that:
  \begin{lemma}\label{abelians in towers}
  	Any maximal non-cyclic abelian subgroup of $\RLR{T}{}{r}$ is either conjugate into some of the abelian vertex groups in the decomposition associated to some node of the tower
  	or can be conjugated into $L_{i}$ for some $1\leq i\leq m$.
  \end{lemma}
  Consider first the case in which $\lambda_{m}$ neither belongs to $\Lambda$ nor is of isolated abelian type (observe this includes the possibility that $\RLR{T}{}{\lambda_{m}}$ is cyclic).
  Let us give a closer look at $\RLR{T}{}{\lambda_{m}}$.
  %   if we assume that the tower is \fine, then
  In case $\lambda_{m}$ is of free-product type $M_{\lambda_{m}}$ has to be hyperbolic in the corresponding free product, hence cyclic. In the floor case the minimality of $\lambda_{m}$ implies $M_{\lambda_{m}}$ cannot intersect a conjugate of a rigid type vertex non trivially, hence ($M$ is maximal) neither a conjugate of an abelian type vertex.
  On the other hand, maximal abelian subgroups of surface groups and free groups are cyclic, so $M_{\lambda_{m}}$ has to be cyclic in every other case as well.
  %  On the other hand, centralizers hyperbolic elements of $\DR{T}{\lambda}$ or non-peripheral elements of vertex groups are cyclic.
  
  Whenever $M_{\lambda_{m}}$ is cyclic we will refer to $M$ as a pegged abelian group and to any generator of the conjugate of $M_{\lambda_{m}}$ in $M$  as a peg of $M$.
  Let $\mathcal{PM}$ denote some system of representatives of the conjugacy classes of pegged non-cyclic maximal abelian subgroups of $T$.
  % the family of abelian type vertex groups appearing along $\QR{T}$ which are maximal abelian subgroups of $\qG{T}$; observe they are a system of representatives of conjugacy classes of maximal abelian subgroups of $T$
  %  stand for the collection of all those pegged maximal abelian groups intersecting non-cyclic abelian vertex groups of $\QR{T}$ non-trivially.
  
  In case $M_{\lambda_{m}}$ is an isolated abelian vertex group, we will say that $M$ is \emph{\FHA}
  and refer to the unique conjugate of $M_{\lambda_{m}}$ contained in $M$ as the bottom of $M$, denoted by $\bott{M}$.
  Let $\mathcal{FHM}(T)$ stand for the collection of representatives of conjugacy classes of \FHA groups.
  %  for thoe \FHA groups containing an isolated abelian vertex group
  We remark that if $\Lambda=\nil$, any non-cyclic maximal abelian group is conjugate to an element of $\mathcal{PM}(T)\cup\mathcal{FHM}(T)$.
  %  	 to $M$ as a free hanging maximal abelian group (\FHMA) and
  \begin{comment}
  In general, we will say of $m\in M$ that it is \fundamental in $\qG{M}$ if it generates a cyclic summand of $\qG{M_{\lambda_{m}}}$ onto which $M_{\lambda_{m}}$ \pr-retracts, or equivalently, if it is primitive in $M$ and $\subg{\pi(m)}=\pi(Q)$.
  \end{comment}
  
  \begin{definition}
  	% change the word \emph{closed}
  	Suppose we are given a closed \rs-\pr-tower $\LRL{T}{A}{J}$. % of a \rs limit group $\rqG{G}{A}$.
  	By a \emph{closure} of $\QR{T}$ we intend a pair $(\kappa,\LRL{S}{A}{J})$ where $\LRL{S}{A}{J}$ is a \pr-tower with the same index set and dropped nodes as $\RTA$ and $\kappa$ an injective morphism
  	\begin{align*}
  		\fun{\kappa}{\RLR{T}{r}{}\frp(\bfrp{\Tl}{\lambda\in DFt(T)}{})}{\RLR{S}{r}{}\frp(\bfrp{\qG{S}^{\lambda}}{\lambda\in DFt(S)}{})}
  	\end{align*}
  	compatible with the quotient maps in the two resolutions and such that:
  	\enum{i)}{
  		\item  $\lambda$ is of free product type in $\QR{S}$ if and only if it is in $\QR{T}$
  		%	  	\item For each $\lambda\in J$ the dropping factors at $\lambda$ are the of the two resolutions coincide.
  		%    		\item For any $\lambda\in DFt(T)$ the $\kappa$ restricts to an isomorphism between $\Tl$ and a finite index subgroup of $\qG{S}^{\lambda}$, which is only a proper subgroup in case $\Tl$ is covered by an isolated abelian group.
  		\item If $\lambda\in\lvs{J}$ then $\kappa$ restricts to an isomorphism between $\Tl$ and $\qG{S}^{\lambda}$, unless the first is the image of a \FHA group, in which case the codomain of the inclusion is abelian as well and contains the image as a finite index group.
  		%    		This will depend of the proof of the solid stufff
  		\item For any other $\lambda\in J$ there is a vertex type preserving isomorphism $\phi$ between the underlying graphs of the decompositions $\Delta=\GTla$ and $\Delta^{S}=\Delta\qG{S}^{\lambda}$ such that $\kappa(t_{e})=t_{\phi(e)}$, $\kappa$ restricts to embeddings of $\Delta_{e}^{T}$ into $\Delta_{\phi(e)}^{S}$ and of $\Delta_{v}^{T}$ into $\Delta_{\phi(v)}^{S}$. Moreover, $\kappa(\Delta_{v}^{T})$ is equal to $\Delta_{v}^{S}$ if $v$ is of surface type and a finite index subgroup if it is of abelian type.
  		The image in $Q$ of $\Delta_{v}^{S}$ and $\Delta_{v}^{T}$ are the same.
  		%  		And some element of $\Delta_{v}^{s}$
  	}
  	%  By a \emph{weak closure} we intend a tower satisfying the same set of conditions, except that inclusions between abelian vertex groups are allowed to be of infinite index.
  	Incurring in a certain abuse of language we sometimes refer as the underlying group $\rqG{S}{A}$ itself as a closure of $\rqG{T}{A}$.
  \end{definition}
  
  Notice that given a system of representatives
  $\mathcal{M}$ of the conjugacy classes of non-cyclic maximal abelian subgroups of $\rqG{T}{A}$
  and $\RSA$ a closure of $\RTA$, the group $S$ is the fundamental group of a star shaped graph of groups, with $T$ in its center, and a free abelian group $\bar{M}$ amalgamated to $T$ over the finite index subgroup $M$ at each of the other vertices.
  \begin{observation}
  	In case $M$ is pegged, its peg is a peg of $\bar{M}$.
  	%  	Add a corresponding observation for the tipp of a free hanging non-cyclic abelia
  \end{observation}
  \begin{proof}
  	All we have to show is that any peg $p$ in $T$ has no proper roots in $S$ either. If some non-trivial element of $\bar{M}$ could be conjugated to a rigid vertex group of the decomposition
  	in $\GTla$ in which $p$ is hyperbolic, then the same would be true for some finite power in $M$ of said element, contradicting choice of $p$.
  	In general, suppose we are given a reduced simplicial $G$-tree $T$ in which no non-trivial element fixes an infinite axis and $H\leq G$, with minimal tree $T_{H}\subset T$, so that
  	the natural map from $H\backslash T_{H}$ to $G\backslash T$
  	is injective on edges.
  	Then any hyperbolic element $h\in H$ lacking a proper root in $H$ also lacks one in $G$.
  	Indeed, suppose this is not the case, so that there some $g\in G$ and $n\in\N$ such that $g^{n}=h$.
  	
  	We claim that $g$ preserves $T_{H}$.
  	By assumption action of $g$ preserves the orbits of edges of $A=Ax(h)$ by the action
  	of $G$. If we pick $e\in A$ this implies the existence of some $k\in H$ such that $k\cdot e=g\cdot e$. Since both $e$ and $g\cdot e$ point in the same direction along $A$, necessarily
  	$Ax(k)=A$. It follows that $k^{-1}g$ fixes $A$, which by assumption can only be the case if $g=k$.
  \end{proof}
  \newcommand{\augm}[0]{}   % amalgamation of the model to a tower
  
  We say that a morphism $f$ which factors through $\RTA$ is primitive if $f(M)=Z_{\F}(M)$ for any $M\in\mathcal{MA}$.
  We say that $f$ extends to some closure $\RSA$ of $\RTA$ if it extends to a morphism from $\rqG{S}{A}$ to $\rqG{F}{A}$ factoring through $\RSA$.
  \begin{definition}
  	\label{covering closure definition}By a covering system of closures of a \rs \pr-tower $\LRL{T}{A}{J}$ of
  	\footnote{A covering closure in Sela's terminology.} we intend a finite family
  	%  	$\{(\kappa_{i},\LRL{S^{i}}{A}{J})\}_{i=1}^{m}$
  	$\mathcal{CL}$ of closures of $\LRL{T}{A}{J}$ with the property that every morphism $\fun{f}{\rqG{T}{A}}{\rqG{F}{A}}$ factoring through
  	$\RTA$ extends to some $\fun{\bar{f}}{\RLR{S}{A}{r}}{\rqG{F}{A}}$ factoring through $\RSA$ for some $\RSA\in\mathcal{CL}$.
  \end{definition}
  %  Let $\LRL{T}{A}{(J,r)}$ be a (relative) \pr-tower structure.
  \newcommand{\pegging}[0]{pegging }   %tc
  \newcommand{\ep}[1]{#1'}
  This extends to a homomorphism
  $\fun{\bar{f}}{\rqG{S}{A}}{\rqG{H}{A}}$ if and only if for any $M\in\mathcal{FHM}$
  %  for which $f(M)$ is non trivial, contained in a maximal subgroup $Z$ of $H$,
  the restriction $f_{M}:=\hat{f}\rst_{\qG{M}}$ extends to $f_{\bar{M}}:\qG{\hat{\overline{M}}}\to\qG{H}$ (notice that these, and thus $\bar{f}$, are unique).
  The necessary and sufficient condition for this extension to factor through $\RSA$ is precisely that:
  \enum{i)}{
  	\item $f_{\bar{M}}(\bar{M})=f_{\bar{M}}(pg_{M})$ for any $M\in\mathcal{MA}$\label{clause1}
  	\item $f_{\bar{M}}(\bar{M})\leq f_{\bar{M}}(\bott{\bar{M}})$ for any $M\in\mathcal{FHM}$\label{clause2}
  }
  An obvious but important observation that can be drawn from this is that if some primitive $f$ as above extends to $\RSA$ then any $f'$
  such that the maps $\fun{f_{M}}{M}{f(M)}$ and $\fun{f'_{M}}{M}{f'(M)}$ differ by a \pr-automorphism of their images will also extend to $\RSA$.
  In particular, only primitive morphisms are needed to test whether a system of closures is covering or not.
  
  It can come in handy to modify $\QR{T}$ by turning all maximal non-cyclic abelian subgroups into pegged ones, as described below.
  Let $\hat{T}$ be the group obtained by amalgamating $\hat{M}:=M\oplus Z_{M}$ over $M$ to $T$ for each $M\in\mathcal{FHM}$, where $Z_{M}$ is isomorphic to $\Z$.
  For each $M\in\mathcal{MA}:=\mathcal{PM}\cup\mathcal{FHM}$ let $pg_{M}$ be either a peg of $M$ in $\QR{T}$ in case $M\in\mathcal{PM}$ or a fixed generator of $Z_{M}$ otherwise.
  Denote by $\pg{T}$ the collection all $q\in Q^{\mathcal{FHM}}$ with the property $\pi(Z_{M})=\pi(M)$ for any $M\in\mathcal{FHM}$. For any $\pg{T}$
  %  the extension of $\pi_{T}$ to $\hat{T}$ obtained by setting $\pi(pg_{M})=q_{M}$ for any $M\in\mathcal{FHM}$ is
  denote by $\qG{T^{q}}_{A}$ the \rs \pr-group with underlying group $\hat{T}$ in which $\pi(pg_{M})=q_{M}$ and the rest of the structure stays the same.
  There is an obvious \rs \pr-tower structure $\LRLi{T}{q}{A}{(J,r)}$, obtained by enlarging the appropriate vertex groups of $\GTla$ for $\lambda\in J$. At the bottom of the new tower, for each 	$\lambda\in\lvs{J}$ for which $\Tl=\bott{M}$ for some $M\in\mathcal{FHM}$ we find a retraction $\eta^{T^{q}}_{\QR{\lambda}}$ from
  $Z_{M}\oplus\bott{M}$ to $\qG{Z_{M}}$ instead. Here $Z_{M}$ is the unique vertex group of an obvious (degenerate) decomposition, making $pg_{M}$ a peg of $\hat{M}$.
  We will refer to any such $\qG{T^{q}}_{A}$ as a pegging of $\rqG{T}{A}$. For any morphism $f$ which factors through $\QR{T}$ there is some $q\in\pg{M}$
  such that $f$ admits a 'good' extension $\hat{f}$ to $\qG{T^{q}}_{A}$, i.e., one for which $\fun{\hat{f}}{\qG{T^{q}}_{A}}{\rqG{G}{A}}$
  factors through $\qG{T^{q}}_{A}$ and with the additional property that
  %  $Z(\hat{f}(M))=\hat{f}(\subg{pg_{M}})$ in case $M\in\mathcal{FHM}$.
  $\hat{f}(M)=\hat{f}(\subg{pg_{M}})$ in case $M\in\mathcal{FHM}$.
  %  In this situation we say in that $f$ conforms to the tuple $q$.
  
  \begin{definition}
  	\label{congruency condition definition} Let $K$ be a positive integer divided by $|\pi(M)|$ for any $M\in\mathcal{MA}$. An atomic $K$-congruence condition on $\RTA$ is given by a tuple 		$\mathcal{C}=(q,(\xi_{PM})_{PM\in\mathcal{MA}})$,
  	where $q\in\pg{T}$ and for each $M\in\mathcal{MA}$ the map $\xi_{M}$ is a homomorphism from $\hat{M}$ to $\Z/N\Z$.
  	We say that a primitive $f$ factoring through $\RTA$ satisfies $\mathcal{C}$ if there is an extension $\fun{\hat{f}}{\qG{T^{q}}_{A}}{\rqG{H}{A}}$ of $f$
  	such that $\hat{f}(M)=\hat{f}(\subg{pg_{M}})$ for all $M\in\mathcal{FHM}$ and for all $M\in\mathcal{MA}$ and $m\in M$ equality $\hat{f}(m)=f(pg_{M})^{l}$ holds for some $l\in\xi(m)\in\Z/K\Z$.
  	A $K$-congruency condition is a union of $K$-congruency conditions.
  \end{definition}
  The intersection of a $(q,K)$ and a $(q,K')$ atomic congruency conditions is a (possibly empty) $lcm(K,K')$ atomic congruency condition and that the complement of a congruency condition is again a congruency condition. Note that the condition might be empty even being satisfied by some non-primitive $f$.
  Our goal now is to show the following:
  \begin{lemma}
  	\label{cclemma}Given a finite family $\mathcal{CL}$ of closures of a \rs \pr-tower $\RTA$, a primitive morphism $\fun{f}{\rqG{T}{A}}{\rqG{H}{A}}$ factoring through $\RTA$
  	extends to one of the resolutions in $\mathcal{CL}$ precisely in case it satisfies certain congruence condition.
  \end{lemma}
  Suppose we are given a family $\mathcal{F}$ of sequences of primitive maps factoring through $\RTA$
  We say that $\mathcal{F}$ is congruence complete if and only if for any non-empty congruency condition $\mathcal{C}$
  there is $\ssq{f}{n}\in\mathcal{F}$ such that $f_{n}$ satisfying $\mathcal{C}$ for $n$ big enough.
  \begin{corollary}
  	\label{congruency and covers} Let $\RTA$ be a \rs \pr-tower and $\mathcal{F}$ a congruency complete family of sequences of morphisms from $\rqG{T}{A}$ to $\rqG{F}{A}$ factoring through $\RTA$.
  	Then any finite family $\mathcal{CL}$ of closures of $\RTA$ with the property that $f_{n}$ eventually extends to some
  	$\RSA$ in $\mathcal{CL}$ for any $\ssq{f}{n}\in\mathcal{F}$ is a covering system of closures.
  \end{corollary}
  
  \begin{proof}
  	(of \ref{cclemma})
  	Fix $\RSA\in\mathcal{CL}$ and let $f$ be primitive factoring through $\RTA$ and conforming to some $q\in Q$. Recall the discussion after definition \ref{covering closure definition}.
  	%	  Suppose we are given a primitive morphism $f$ from $\rqG{T}{A}$ to a free \rs \pr-group $\rqG{H}{A}$ factoring through $\RTA$, which conforms to $q$.
  	%  Assume that $f$ conforms to $q$
  	%  and let $\hat{f}$ be a good extension.t
  	Since $M$ has finite index in $\bar{M}$ and $pg_{M}$ is primitive in $\bar{M}$ in case $M\in\mathcal{PM}$, for $M\in\mathcal{MA}$,
  	the basic theory $\Z$-modules tells us that there are bases $\{pg_{M},m_{j}\}_{j=1}^{r}$ of $\hat{M}$ and $\{pg_{M},\bar{m}_{j}\}_{j=1}^{r}$ of $\bar{\hat{M}}$
  	and positive integers $k_{j}$ for $1\leq j\leq r$ such that $m_{j}=\bar{m}_{j}^{k_{j}}$.
  	Let $K_{M}=lcm\,\{k_{j}|\pi(M)|\,|\,1\leq j\leq r\}$. Consider $M\in\mathcal{PM}$ first.
  	For each $1\leq j\leq r$, let $e_{j}$ the unique integer such that $f(m_{j})=z^{e_{j}}$, where $z=f(pg_{M})$.
  	Clearly, whether $f_{M}$ has an extension to $\bar{M}$ satisfying (\ref{clause1}) above
  	depends only on whether each $f(m_{j})$ has a $k_{j}$-th root in $\subg{z}$, hence only on the tuple $([e_{j}])_{j=1}^{r}\in(\Z/N\Z)^{r}$.
  	The same is true fir the question of whether this extension is in fact a \pr-morphism, i.e., whether the $k_{j}$-th root of $f(m_{j})$ is mapped to $\pi(\bar{m}_{j})$ by $\pi$, since $|\pi(M)||K$.
  	
  	In the case of $M\in\mathcal{FHM}$, it is a priory not enough for us that $f(m_{j})$ has a $k_{j}$-th root,
  	in order for (\ref{clause2}) to be satisfied we need that $f(m_{j})=f(n_{j})$ for some $n_{j}\in\bott{M}$ which itself has a $k_{j}$-th root in $\bar{M}$. This follows from the fact that
  	$Z(f_{M}(M))=\bar{f}_{M}(\bar{M})=f_{M}(M)$.
  	%	  Since $n_{j}$ is not given in advanced, how to characterize the maps $f$ for which this is true using a single congruence condition might seem a priori problematic.
  	%	  Consider the subgroup $M^{[k_{j}]}\leq M$ consisting of all the
  	%	  elements of $M$ having a $k_{j}$-th root in $\bar{M}$. This has finite index in $M$, and the index $[f(M),f(M^{[k_{j}]})]$  divides $K$. Therefore, whether $f_{m}(m_{j})\in f(M^{[k_{j}]})$
  	%	  or not is controlled by a congruency condition.
  	%  Putting together the conditions coming from each of $M\in\bar{M}$ we conclude that whether a morphism $f$ factoring through $\RTA$ extends to one factoring through $\RSA$ is given by
  	%  a a $K$-congruency condition, where $K$ is the lowest common multiple of $\{K_{M}\}_{m\in\mathcal{MA}}$.
  \end{proof}

\section{Completions}
  \label{completions subsection}
  %    \newcommand{\CR}[0]{}
  % Need to do it for graded.
  %    Introduce a special notation for the embeddings $\iota$, so that one does not need to mention them explicitely every time
  %     as the normalization operation preserves well-structuredness and \pr-modular groups; though not essential, this hypothesis allows neater statements and proofs in what follows.
  %    Let $\LRL{R}{A,A}$ be a strict \ws \pr-\res.
  Let $\LRL{R}{A}{J}$ be a \gr \ws \pr-\res for which the \pr-group at any leaf is a \rs \pr-limit group.
  We will assume that all the $\GRla$ are normalized, as defined earlier (if not, the construction can be done using their normalization instead). By adding edges to groups of surface type (possibly making the tree non-minimal) we can guarantee that
  
  In this section we describe the construction of a so called completion of aa \ws \pr-\res $\mathcal{R}$: this will be a \rs (\gr) \pr-tower $\LRL{S}{A}{(J,r)}$ over the same index set, together with certain compatibility data. The standard construction for plain groups can be found in \cite{sela2}.
  First of all, for each $\lambda\in J$ we have an injective morphism $\fun{\iota_{\lambda}}{\RLR{R}{A}{\lambda}}{\RLR{S}{A}{\lambda}}$. This is an isomorphism in case $\lambda\in\lvs{J}$ or $\GRla$ consists of a single non-rigid type vertex (so the associated group is a closed surface group or a free non-cyclic abelian group.
  
  If $\lambda$ is not of free product type and $\GRla$ has at least one edge, Let $T_{\lambda}=\widetilde{\GRla}$ and $\hat{T}^{*}_{\lambda}$ the $\RLR{R}{}{\lambda}$-tree obtained by normalizing $\widetilde{\GSla}$ sliding all edges between a rigid and a surface type vertex. Then there is also an isomorphism of \gats $\Theta_{\lambda}$  $\afld{T_{\lambda}}$ (as described in the previous subsection) and some subtree of $\hat{T}^{*}_{\lambda}$ such that:
  \enum{i)}{
  	\item $\eta^{\lambda}_{\QR{S}}\circ\iota_{\lambda}=\bfrp{\iota_{\mu}}{\mu\prec\lambda}{}\circ \eta^{\lambda}_{\QR{R}}$. \label{natural transformtion}
  	\item For each \trm $(f_{\lambda},\tau_{\mu})_{\lambda\in J,\mu\in\hat{J}}$ of $\RRA$, there is a \trm $(g_{\lambda},\tau_{\mu})_{\lambda\in J,\mu\in\hat{J}}$ such that $f_{\lambda}=g_{\lambda}\circ\iota_{\lambda}$ for each $\lambda\in J$. \label{extension prooperty}
  	\item $\Theta_{\lambda}$ is equivariant with respect to $\iota_{\lambda}$, i.e., for each $g\in Sl$ and $x\in\afld{T_{\lambda}}$ we have $\Theta_{\lambda}(g\cdot x)=\iota_{\lambda}(g)\cdot\Theta_{\lambda}(x)$.
  	%    	\item For each $v\in V_{s}T$ the map $\iota_{\lambda}$ maps $Stab(v)$ onto	 $Stab(\Theta_{\lambda}(v))$.
  	\item Each abelian vertex group of $T_{\lambda}$ is mapped by $\Theta$ into an abelian group of $\hat{T}^{*}_{\lambda}$.
  	\item The star around each surface type vertex of $T_{\lambda}$ is mapped isomorphically onto the star of a surface type vertex of $\hat{T}^{*}_{\lambda}$.
  	%    	\item For each
  	%    	[add conditions on the shape and complexity of the completion]
  }
  \begin{comment}
  \begin{remark}
  	If the decomposition $\GRla$ is not in normal form for each $\lambda$, but nevertheless a \gad relative the family of non-cyclic abelian subgroups, we can replace them by their normal form and apply the construction to the resulting resolution. Normalization preserves \pr-modular groups and conjugacy classes of non-abelian groups (see \ref{}), hence the resulting resolution $QR{S}$ together with the emebeddings $\iota_{\lambda}$ still satisfies properties \ref{natural transformation} and \ref{extension property} with respect to the starting resolution.
  \end{remark}
  \end{comment}
  We construct $QR{S}\rst_{\lambda}$ and the embedding $\iota_{\lambda}$ by induction on $\leq$. For $\lambda\in\lvs{J}$ let $Sl=Rl$ and $\iota_{\lambda}$ the idenity.
  Now, suppose that $QR{S}\rst_{\lambda}$ is defined for $\mu\succ\lambda$ and that same is true for the embeddings $\iota(\mu)$.
  There are two cases to consider.
  If $\lambda$ is of free product type, then we simply let
  $\RLR{S}{}{\lambda}$ be $\bfrp{\Sm}{\lambda\succ\mu}{}$ and $\eta^{\lambda}_{\QR{S}}$ the map $\bdcup{\iota_{\mu}}{\mu\prec\lambda}{}$. In the next section we will deal with the case in which $\DR{R}{\lambda}$ is non-tivial and $\eta^{\lambda}_{\QR{R}}$ non-injective. All along the way we will use the fact that any towers having \rs \pr-limit groups at its leaves is itself a limit group.
  
  \subsection{Floor case}
    
    In order to ease the notation, we will denote $\bdcup{\iota_{\mu}}{\mu\prec\lambda}{}$ by $\kappa$, $\theta=\kappa\circ\phi$ and by $\iota$ the map $\iota_{\lambda}$ to be constructed. In the same fashion, we will rename $\eta^{\lambda}_{R}$ as $\phi$,
    the new map $\eta^{\lambda,S}_{\lambda,S}$ to be constructed as $\eta$ and $\Theta_{\lambda}$ as $\Theta$. Let also $\qG{G}:=\qG{R}^{\lambda}$ and $\qG{\hat{G}}$ the \pr-group $\qG{G}^{\lambda}$. The \gad $\GRla$ puts the group $\Rl$ into isomorphism with some $\pi_{1}(\Delta,Z)$.
    If this decomposition consists of a single isolated vertex of surface of abelian type we then can just set $\qG{S}^{\lambda}=\qG{R}^{\lambda}$.
    
    Let $X=(V,E)$ be the underlying graph of $\Delta$. Set also $\qG{H_{\mu}}=\Sm$ and $\qG{H}=\bfrp{\Sm}{\lambda\succ\mu}{}$.  % modify the notation so that the
    We now describe the construction of a \gad $\hat{\Delta}$, with underlying graph $\hat{X}=(\hat{V},\hat{E})$.
    Its set of rigid vertex groups wil be precisely $\{H_{\mu}\}_{\mu\in nDFT(\lambda)}$. We also define a \pr-structure to each vertex group a	$\hat{\Delta}_{v}$.
    %       and
    %      an injective morphism $\zeta_{v}$ of it into
    %      and a morphism from $\hat{\Delta}_{v}$ to $\qG{H}$.
    %      This determines the restriction to the subgroup of the fundamental group of $\hat{\Delta}$ generated by elliptic elements of both the \pr-structure, and the morphism $\eta$ to $H$.
    Later in the subsection we specify a maximal tree $\tilde{Z}$ in $\hat{X}$ and show that the partial \pr-structures on vertex groups extend to a \pr-structure on $\hat{G}:=\pi_{1}(\hat{\Delta},\hat{Z})$, which we take as $Sl$, and a morphism from it to $\qG{H}$. As $\GSla$ we take the decomposition as the graph of groups $\Delta$ above. We start the construction of a graph of groups $\Lambda$ with underlying graph
    \begin{align*}
    	(W,F)=(\{v_{\mu},\hat{u},\hat{w}\}_{\mu\in NDCt(\lambda),
    	u\in V_{a},	w\in V_{s}},
    	\{e_{v},\edin{e_{v}},\hat{f},\edin{\hat{f}}\}_{v\in V_{a},\alpha(f)\in V_{s}})
    \end{align*}
    %     We will allow ourselves to temporarily lift the requirement that the underlying graph of a graph of groups must be connected.
    To begin with, set $\Lambda_{\mu}=\Sm$. For any $w\in V_{a}$ let $\qG{P_{v}}=(Per(v),\pi\rst)$ and
    chose some $\qG{K_{v}}$ such that $(\Delta_{w},\pi\rst)=\qG{P_{v}}\oplus\qG{K_{v}}$.
    Since $\iota$ is injective on $P_{v}$ the morphism $\eta$ injects $\qG{P}$ into $\qG{H_{\mu_{w}}}$ for some rigid vertex $w\in V$. In correspondence, the new graph includes vertex $\hat{w}$ and a single edge $e_{w}$ from $\hat{w}$ to $v_{\mu_{w}}$. Let $O(w)$ denote the set of those edges $e$ from $w$ to a vertex distinct from the first edge in the unique path in $Z$ from $w$ to $v_{0}$.
    Let $Q_{w}=C_{H_{\mu_{w}}}(\theta(\Delta_{w}))$. As $\Lambda_{\hat{w}}$ we take the group $Q_{w}\oplus K'_{w}\oplus(\subg{\bar{t}_{e}}_{O(w)})_{ab}$, where $\qG{K'_{w}}$ is a copy of $K_{w}$ isomorphic to it via some $\alpha_{w}$. We endow $\Lambda_{\hat{w}}$ with a \pr-structure
    %      by letting $\pi_{\Lambda_{\hat{w}}}$ r
    restricting to $\pi_{H_{\mu_{w}}}$ and to $\pi_{K'_{w}}$ respectively on the first two direct summands in the expression above and sending each of the generators $\bar{t}_{e}$ to the identity.
    \begin{comment}
    Let $\qG{Q_{w}}=(C_{H_{\mu_{w}}}(\theta(\Delta_{w}))\pi_{H\mu_{w}}\rst_{Q_{w}})$. As $\Lambda_{\hat{w}}$ we take the group $Q_{w}\oplus K'_{w}\oplus(\subg{\bar{t}_{e}}_{O(w)})_{ab}$, where $K'_{w}$ is a copy of $K_{w}$ isomorphic via some $\alpha_{w}$. We endow such vertex with the \pr-structure restricting to $Q_{w}$ to the one inherited from $H_{w_{v}}$, to the pushforward via $\alpha$ of that on  and send the image in $Q$ of $\bar{t}_{e}$ to be the identity.
    \end{comment}
    As $\Lambda_{e_{w}}$ we take the leftmost summand in the previous decomposition and as $i_{e_{w}}$ simply the identity map.
    Our $\zeta_{v}$ will be the \pr-isomorphism $\theta\rst_{K_{w}}\oplus\alpha_{v}$ between $(\Delta_{w}\pi\rst_{\Delta_{w}})$ and the subgroup of $\qG{\Lambda_{w}}$ given by $\qG{P_{w}}\oplus\qG{K'_{w}}$.
    %  The edge embedding into $H_{\mu_{w}}$ is given by $\iota_{\mu_{v}}\eta$, while that into $\Lambda_{\hat{w}}$ is the natural one.
    
    To continue, we associate to each $u\in V_{s}$ a new vertex $\hat{u}$ and as $\qG{\Lambda_{\hat{u}}}$ we take some copy of $\qG{\Delta_{u}}$, isomorphic to it via some isomorphism $\zeta_{u}$. For each $e$ originating at $u$ there is $\mu_{e}\prec\lambda$ such that $\eta(\Delta_{\edin{e}})\subset H_{\mu_{e}}$. To $e$ we associate some $\hat{e}$, originating at $u$, with  $\alpha{\edin{\hat{e}}}=v_{\mu_{u}}$. We let $\Delta_{\hat{e}}=\zeta_{u}(\Delta_{e})$. As $i_{\hat{e}}$ we take the map $\theta\circ i_{e}\circ\zeta_{u}^{-1}$ to $\Delta_{e}$.
    %     surface attaching map
    The fact that $\Delta$ is a commutative transitive graph of groups does not guarantee that $\Lambda$ is one as well, as the images of two elements in non-conjugate maximal abelian subgroups might as well commute. By induction on the size of $J$, one can assume, however, that the $H_{\mu}$ are limit groups and, therefore, commutative transitive.
    
    %      we need to define what a commutative transitive group is undefineds
    Dealing with this issue involves three steps. First of all, for each $w\in V_{a}$ we chose $\delta_{v}\in H_{\mu_{w}}$ so that $Q_{w}^{\delta_{w}}=Q_{w'}^{\delta_{w'}}$ whenever $Q_{w}$ and $Q_{w'}$ are conjugate  (which implies $\mu_{w}=\mu_{w'}$). Then replace the vertex group $\Lambda_{w}$ by its formal conjugate
    $\Lambda'_{w}=\Lambda_{w}^{\delta_{w}}\subset\pi_{1}(\Lambda\rst_{v_{\mu_{w}},e_{w},\hat{w}},\{e_{w}\})$.
    By this we mean that $\Lambda'_{w}$ inherits the \pr-structure $\pi(\delta_{w})$-conjugate to that on $\Lambda_{w}$ as well. The new left summand is still a subset of $H_{\mu_{w}}$ and as $i_{e_{w}}$ we can take the inclusion again.
    %      $\Lambda'_{w}=Q_{w}^{\delta_{w}}\oplus K'_{w}^{\delta_{w}}\oplus (\subg{\bar{t}^{\delta_{w}}_{e}}_{O(w)})_{ab}$.
    %      amounts merely to a change of fundamental domain of the associated action and%
    The operation does not affect the partial \pr isomorphism type of the fundamental group. Accordingly, we can see $\zeta_{w}$ as an injective \pr-morphism of $\Delta_{w}$ into $(\Lambda'_{w})^{\delta_{w}^{-1}}$, coinciding with $\theta$ on $K_{w}$.
    \begin{comment}
    We set $\Lambda'_{e}=\Lambda_{e}^{\delta_{w}}$ and $i_{e'_{w}}=\theta^{\inn{\delta_{w}^{-1}}}=\iota_{e_{w}}^{\inn{\delta_{w}^{-1}}}$. The new vertex group comes equipped with the \pr-structure conjugate by $\pi(\delta_{w})$. The latter operation reflects merely a change of fundamental domain of the associated action which does not affect the partial \pr isomorphism type of fundamental group. The morphism $\zeta_{w}$ now can be seen as embedding of $\Delta_{w}$ into $(\Lambda_{w})^{\delta_{w}^{-1}}$.
    \end{comment}
    
    After the operation is completed, we can assume that for any two $w,w'\in V_{a}$ if two conjugates of the images of $\Lambda_{e_{w}}$ and $\Lambda_{e_{w'}}$ in $H_{\mu}$ for $\mu=\mu_{w}=\mu_{w'}$ commute, then they coincide, in which case
    we say that $w$ and $w'$ are equivalent. For each $\mu_{w}$ fold together % at $v_{\mu_{w}}$
    all the edges from $v_{\mu}$ ending in the same equivalence class $[w]$. The vertex group associated to any of the new vertices $\hat{[w]}$ is the amalgamated product of $\Lambda_{w'}$ for $w'\sim w$ over the common peripheral subgroup. At this point we replace any of those groups by their abelianization: % say what such a tree product is
    \begin{align*}
    	P_{[w]}\oplus(\bcsum{K_{w'}}{w\sim w'}{})\oplus(\bcsum{s_{e}}{e\in O(w'),\\w'\sim w}{})
    \end{align*}
    Here $s_{e}$ stands for the image of the generator $\bar{t}_{e}^{\delta_{w}}$ of $\Lambda_{w}$. We keep the notation
    $\bar{t}_{e}$ for $s_{e}^{\delta_{w}^{-1}}$.
    Notice that since peripheral subgroups are mapped injectively by this quotient, the what we are left with is still a graph of groups.
    The fact that the quotient group inherits a \pr-structure from the amalgamated free product follows from the following:
    \begin{claim}
    	Let $w,w'\in V_{a}$ such that $w\sim w'$. Then $[\pi(\Delta_{w})^{\pi(\delta_{w})},\pi(\Delta_{w'})^{\phi(\delta_{w})}]=1$.
    \end{claim}
    \begin{proof}
    	\begin{comment}
    	Observe that for any \pr-limit group $\qG{L}$ and $g,h\in L\setminus\{1\}$ such that $[g,h]=1$ the images $\pi(g)$ and $\pi(h)$ are contained in some cyclic subgroup of $Q$. Indeed, if we regard $L$ as a limit of morphisms $\fun{f_{n}}{\qG{L}}{\qG{F}}$, necessarily $f_{n}(g)$ and $f_{n}(h)$ are contained in the same subgroup of $F$ and, hence, the same holds for $\pi_{F}(f_{n}(g))=\phi_{L}(g)$ and $\pi_{F}(f_{n}(h))=\phi_{L}(h)$.
    	\end{comment}
    	Let $\phi=\iota\circ\eta$. Observe that $\pi(\Delta_{u})=(\pi\circ\theta)(\Delta_{u})\subset C_{H_{\mu_{w}}}(P_{u})$ for $u\in\{w,w'\}$. Now, we know that $H_{\mu_{w}}$ is a \pr-limit group, 	hence $CSA$.  Since $P_{w}^{\delta_{w}}$ and $P_{w'}^{\delta_{w'}}$ are contained in some maximal abelian subgroup, the same is true for $\theta(\Delta_{w}^{\delta_{w}})$ and $\theta(\Delta_{w'}^{\delta_{w'}})$. Hence
    	\begin{align*}
    		1=[\pi(\phi(\Delta_{w'})^{\delta_{w'}}),\pi(\phi(\Delta_{w})^{\delta_{w}})]=[\pi(\phi(\Delta_{w'}))^{\pi(\delta_{w'})},\pi(\phi(\Delta_{w}))^{\pi(\delta_{w})}]=[\pi(\Delta_{w'})^{\pi(\delta_{w'})},\pi(\Delta_{w})^{\pi(\delta_{w})}]
    	\end{align*}
    \end{proof}
    The previous embedding $\zeta_{w}$ corresponds now to one of $\Delta_{w}$ into $\hat{\Delta}_{\hat{[w]}}^{\delta_{w}^{-1}}$, for which we employ the same term.
    Extend the notation by letting $\zeta_{v}=\theta\rst_{\Delta_{v}}:\Delta_{v}\to\hat{\Delta}_{\mu_{v}}$ for $v\in V_{r}$.
    %      We grab the chance to make the following remark:
    \begin{remark}
    	\label{edge groups}Let $e\in E$ and $v=\alpha(e)$. Then $\zeta_{v}\rst_{\Delta_{e}}$ is equal to:
    	\enum{i)}{
    		\item $\theta\rst_{\Delta_{e}}$ in case $v\in V_{r}\cup V_{a}$
    		\item $\inn{t_{\hat{e}}^{-1}}\circ\theta\circ\inn{t_{e}}\rst_{\Delta_{e}}=\inn{\theta(t_{e})t_{\hat{e}}^{-1}}\circ\theta\rst_{\Delta_{e}}$ in case $v\in V_{s}$
    	}
    \end{remark}
    %      By $\bar{t}_{e}$ we denote by $\bar{t}_{e}=\tilde{t_{e}}\in\hat{\Delta}_{\hat{[w]}}^{\delta_{w^{-1}}}$ for an edge $e$ incident to $w\in V_{a}$, which lays in the same vertex stabilizer of $\hat{T}$ as the image of $\zeta_{w}$.
    
    We first specify the restriction of $\eta$ to each vertex group of $\hat{\Delta}$. This is the identity on every $\hat{\Delta}_{v_{\mu_{u}}}$, as required by the definition of a tower. On $\hat{\Delta}_{\hat{[w]}}$ with $w\in V_{a}$, we let $\eta$ kill any of the generators of the form $\tilde{t_{e}}$ for $w'\sim w$ and $e\in O(w)$. On the direct summand which is the image of $\qG{K_{w'}}$ in $\hat{\Delta}_{\hat{[w]}}$ we let the restriction of $\eta$ coincide with $\inn{\delta_{v}}\circ\theta\circ\zeta_{v}^{-1}\circ\inn{\delta_{v}^{-1}}$.
    The map $\eta\rst_{\Delta_{\hat{[w]}}}$ is a \pr-morphism since both $\theta$ and $\zeta_{v}$ are.
    \begin{comment}
    In the last case
    \begin{align*}
    	\pi\circ\inn{\delta_{v}}\circ\theta\circ\zeta_{v}^{-1}\circ\inn{\delta_{v}^{-1}}=\inn{\pi(\delta_{v})}\circ\pi\circ\theta\circ\zeta_{v}^{-1}\circ\inn{\delta_{v}^{-1}}=\inn{\pi(\delta_{v})}\circ\pi\circ\inn{\delta_{v}^{-1}}=\pi
    \end{align*}
    \end{comment}
    In case $v\in V_{s}$, we let $\eta\rst_{\hat{\Delta}_{\hat{v}}}=\theta\circ\zeta_{v}^{-1}$. This, again, preserves the \pr-structure.
    The identity $\eta\circ\zeta_{v}=\theta$ clearly holds for every $v\in V$.
    
    Notice that removing the set $\mathcal{O}$ of all those $\hat{e}$ associated with edges $e\in E\setminus Z$ does not disconnect $\hat{X}$. This means it is possible to chose a maximal subtree $\hat{Z}$ of $\hat{X}$ which does not intersect $\mathcal{O}$.
    For each $v\in V$ let $\Phi(v)\in\hat{V}$ be either:
    \elenco{
    	\item The unique $v_{\mu}$ such that $\eta(v)\subset G_{\mu}$
    	\item $\hat{[v]}$ in case $v\in V_{a}$
    	\item $\hat{v}$ in case $v\in V_{s}$
    }
    Given edges $e\in E\setminus Z$ and $f\in\hat{E}\setminus\hat{Z}$ denote by $t_{e}\in\pi_{1}(\Delta,Z)$ and $\hat{t}_{f}$ the Bass-Serre elements associated to $e$ and $f$ in the respective graph of groups.
    %      Let $\mathcal{N}=\hat{E}\setminus(\mathcal{O}\cup \hat{Z})$.
    
    We now have to check that both the map $\eta$ and the partial \pr-structure defined on vertex groups extend to the entire $\hat{G}$. In view of the fact that $\eta\rst_{H}=Id\rst_{H}$ and that every edge in $\hat{E}$ has an endpoint in some of the $v_{\mu}$, this is equivalent to the existence for all $e\in\hat{E}$ of $q_{e}=\eta(t_{e})\in\hat{G}$ such that the equality
    \begin{align*}
    	i_{e}=\inn{q_{e}}\circ\eta\rst_{\hat{\Delta}_{e}}:\hat{\Delta}_{e}\to H_{\mu_{e}}
    \end{align*}
    is satisfied for all $e\in\hat{X}\setminus\hat{Z}$. Note that any such edge is adjacent to some surface type vertex.
    So suppose that $e=\hat{f}$ for some $f$ originating at $w\in V_{s}$.
    %      We can assume that $\omega(e)\in\{v_{\mu}\}_{\mu\prec\lambda}$.
    %       This is trivial for $v=\hat{[w]}$ of abelian type, as then $t_{e}=1$ (hence $q_{e}=1$) and $i_{e}$ is the identity inside $G$.
    Then:
    \begin{align*}
    	i_{e}=\theta\circ i_{f}\circ\zeta_{w}^{-1}=\theta\circ\inn{t_{f}}\circ\zeta_{w}^{-1}=\inn{\theta(t_{f})}\circ\theta\circ\zeta_{w}^{-1}
    \end{align*}
    We also have $\eta\rst_{\Delta_{e}}=\theta\circ\zeta_{w}^{-1}$. So letting $q_{e}=\theta(t_{f})$ does the job.
    Note that in virtue of the way in which we chose $\hat{Z}$, in case $t_{e}=1$ we have $q_{e}=1$, while $t_{f}=1$ in this case by virtue of the choice of $\hat{Z}$.
    % And that $q_{\hat{e}}=1$ for those $e\in Z$ not adjacent to some vertex of $V_{s}$.
    In the same fashion, setting $\pi(t_{e})=\pi(\theta(t_{f}))$ extends the \pr-structure on vertex groups to a global one
    compatible with the previous collapse map.
    % , as this choice verifies the equation
    % $\pi_{H\mu_{e}}\circ i_{e}=\inn{\pi(t_{e})}\circ\pi_{\hat{\Delta}_{\alpha(e)}}\rst_{\hat{\Delta}_{e}}$ for $e\in\hat{e}$.
    
  \subsection{The morphism $\iota$}
    
    %      We keep the identification of $G$ with $\pi_{1}(\Delta,Z)$.
    %    \subsubsection{The morphism $\iota$ }
      
      Chose a reference vertex $v_{0}\in V_{r}\cup V_{a}$. If the group $A$ of constants is not-trivial, then we can assume that $A\subset\Delta_{v_{0}}$.
      The construction depends on the assignment, to each $v\in V$ of a conjugating element $\gamma_{v}\in\pi_{1}(\hat{\Delta},\hat{Z})$. With this information one defines a morphism $\iota_{v}$ from each $\Delta_{v}$ into $\hat{G}$, given by $\inn{\gamma_{v}}\circ\zeta_{v}$.
      %         where $\zeta_{v}$, regarded as a map into $\hat{G}$.
      
      The choice of $\gamma_{v}$ is completed with the assignment of an element $r_{e}\in\hat{G}$ for each $e\in E$ in such a way that the map sending $t_{e}$ to $r_{e}$ and restricting to $\zeta_{v}$ on $\Delta_{v}$ extends to a morphism $\iota$ from $\rqG{G}{A}$ to $\rqG{\hat{G}}{A}$ with the properties sought for. In particular $r_{e}=1$ for $e\in Z$.
      %        As usual for $e\in Z$ both $t_{e}$ and $r_{e}$ are trivial.
      %       Modulo identifying $\Delta_{v}$ with its restriction in $\pi_{1}(\Delta,Z)$, we can think of the $\zeta_{v}$ as restrictions of \zeta.
      The condition for an extension which is a group homomorphism to exist is that the equation
      \begin{align*}
      	\inn{r_{e}}\circ\inn{\gamma_{v}}\circ\zeta_{v}\rst_{\Delta_{e}}=\inn{\gamma_{w}}\circ\zeta_{w}\inn{t_{e}}\rst_{\Delta_{e}}
      \end{align*}
      holds for all $e\in E$ from $v$ to $w$ (recall that $t_{e}$ acts like $i_{e}$ on $\Delta_{e}$ seen as a subgroup of $G$). Equivalently:
      \begin{align*}
      	\label{embedding well defined}
      	\inn{\gamma_{v}r_{e}\gamma_{w}^{-1}}\zeta_{v}\rst_{\Delta_{e}}=\zeta_{w}\circ\inn{t_{e}}\rst_{\Delta_{e}}
      \end{align*}
      \newcommand{\bt}[0]{\bar{t}}
      
      In order to define $\gamma_{v}$, consider the unique simple path $v_{0},e_{0},v_{1},e_{1},\cdots v_{k},e_{k},v_{k+1}=v$ contained in $Z$.
      We let $\gamma_{v}=g_{k}g_{k-1}\cdots g_{0}$, where $g_{j}$ is equal to:
      \enum{i)}{
      	\item  $\bt_{e_{j}}\in\hat{\Delta}_{\hat{v_{j}}}^{\delta_{j}^{-1}}$ in case $v_{j}\in V_{a}$ and $v_{j+1}\in V_{r}$ \label{case uno}
      	\item $t_{\hat{e_{j}}}^{-1}$ in case $V_{s}\cap\{v_{j},v_{j+1}\}\neq \nil$  (possibly $1$) \label{case dos}
      	\item $1$ in case $v_{j}\in V_{r}$ and $v_{j+1}\in V_{a}$ \label{case tres}
      }
      %lololo
      
      As for $r_{e}$, it is enough to deal with the case in which $v=\alpha(e)\in V_{a}$. If $w=\omega(e)$, then we let $r_{e}$ be equal to:
      \enum{i)}{
      	\item $\gamma_{v}^{-1}\bt_{e}\theta(t_{e})\gamma_{w}$ in case $w\in V_{r}$
      	\item $\gamma_{v}^{-1}t_{\hat{e}}\gamma_{w}$ in case $w\in V_{s}$
      }
      For $e\in E$ from a vertex $v$ to a vertex $w$ we distinguish several cases:
      \newcommand{\De}[0]{\Delta_{e}}
      \newcommand{\DDe}[0]{\Delta_{\edin{e}}}
      \enum{i)}{
      	\item None of $w,v$ belongs to $V_{s}$.
      	%        	In this case let $D=\Delta_{e}$.
      	In this case we have:
      	\begin{align*}
      		\zeta_{w}\circ\inn{t_{e}}\rst_{\De}=\theta\circ\inn{t_{e}}\rst_{\De}=\\=\inn{\theta(t_{e})}\circ\theta\rst_{\De}=\inn{\gamma_{v} r_{e}\gamma_{w}^{-1}}\circ\zeta_{v}\rst_{\De}
      	\end{align*}
      	\item One of $w,v$ belongs to $V_{s}$, let's say $w$ does. Then $\zeta_{v}\rst_{D}=\theta\rst_{D}$ while $\zeta_{w}\rst_{\DDe}=\inn{\theta(t_{e}^{-1}t_{\hat{e}})}\theta\rst_{\DDe}$
      	%        	for $D'=\Delta_{\edin{e}}=D^{t_{e}}$
      	and:
      	\begin{align*}
      		\zeta_{w}\circ\inn{t_{e}}\rst_{\De}=\inn{\theta(t_{e}^{-1})t_{\hat{e}}}\circ\theta\inn{t_{e}}\rst_{\De}=\\
      		=\inn{\theta(t_{e}^{-1})t_{\hat{e}}}\circ\inn{\theta(t_{e})}\theta\rst_{\De}=\inn{t_{\hat{e}}}\circ\theta=\\
      		=\inn{\gamma_{v}r_{e}\gamma_{w}^{-1}}\circ\theta\rst_{\De}=\zeta_{w}\rst_{\De}
      	\end{align*}
      }
      We note a couple of properties of the choice of $\gamma_{v}$ and their implications for $\iota$.
      \enum{i)}{
      	\item $\gamma_{v_{0}}=1$ imples that $\iota$ is a \rs morphism,
      	%  $\eta\circ\iota=\phi$
      	\item $\pi(\gamma_{v})=1$ for any $v$ (recall that $\pi(t_{\hat{e}})=1$ in case $e\in Z$ but $\hat{e}\nin Z$), which implies the restriction of $\iota$ to any vertex group is, just as $\zeta_{v}$ is, a \pr-morphism.
      	\item $\eta(\gamma_{v})=1$ for all $v\in V$, as $\eta(t_{\hat{e_{j}}})=1$ always in case (\ref{case dos}) of the definition of $g_{j}$. This implies that $\eta\circ\iota\rst_{\Delta_{v}}=\theta\rst_{\Delta_{v}}$. Since also $\eta(\iota(t_{e}))=\theta(t_{e})$ for each $e\in E$ (as $\theta(t_{\hat{e}})=t_{e}$) we deduce that $\eta\circ\iota=\theta$.
      }
      
  \subsection{The embedding $\Psi$, injectvity of $\iota$ }
    
    The following will use the assumption that $\hat{G}$ is a limit group so, in particular, $CSA$.
    \begin{comment}
    \begin{definition}
    	A \gad of a group $G$ relative to its non-cyclic abelian subgroups is said to be in normal form if it is commutative transitive, in the dual tree $T$ any edge is adyacent to a vertex of abelian type and for every $e,f\in ET$ such that $\alpha(e)=\alpha(f)=v$, the groups $Stab(e)$ and $Stab(f)$ are conjugate in $Stab(v)$ if and only if $e=f$.
    \end{definition}
    \end{comment}
    Recall $Z$ lifts to a subtree $\lft{Z}$ of the covering $G$-tree $T=\lft{(\Delta,Z)}$ and that
    if we let $\lft{v}$ and $\lft{e}$ stand for the images of $v$ and $e$ respectively, then
    $Stab(T;\lft{v})=\Delta_{v}$ and $Stab(T;\lft{e})=\Delta_{e}$ for all $v\in V$ and $e\in E$.
    Additionally, any $e\in E\setminus Z$ can be lifted to $\lft{e}\in ET$ in such a way that if we let the tuple $(\hat{Z},\hat{Z}\cup\{\lft{e}\}_{e\in E},(g_{\lft{e}})_{e\in E\setminus E})$ is a presentation of the action of $G$ on $T$, where $g_{\lft{e}}=t_{e}$ (in particular $\alpha(\lft{e})\lft{Z}$ for such $e$).
    %        From now on we will drop the tildes in the case of vertices.
    %        reserve the notation for edges
    %        The same holds for the action of $\hat{G}$ on $\widetilde{(\hat{\Delta},\hat{Z})}$, for we will use the same notation.
    The graph of groups $\hat{\Delta}$ is not normalized, as vertices of surface type are connected to vertices of rigid type.
    However, for each $e\in E$, ending in some $v\in V_{s}$,
    %        the subgroup $\hat{\Delta}_{\hat{e}}=\theta(\Delta_{e})$, and hence is conjugate into $\theta(P_{w})$ for the peripheral subgroup $P_{w}$ of some $w\in V_{a}$.
    However, for any edge $e\in\hat{T}$ originating in a vertex of surface type there is a unique
    %        \footnote{The conjugacy classes in $H_{\mu}$ of edge groups of the form $\Delta_{e\hat{[w]}}$ are disjoint.}
    edge $f$ with $\alpha(\edin{f})$ of abelian type such that $Stab(e)\subset Stab(f)$. Sliding each such $e$ over the corresponding $f$ yields a $\hat{G}$-tree $\hat{T}^{*}$ which is normalized.
    %        and which we can think as sharing with $\hat{T}$ the same vertex set, as well as all those edges not involved in the slide.
    We denote by $e^{*}$ the edge of $\hat{T}^{*}$ resulting from the slide of the edge $e$ of $\hat{T}$. Any other edge or vertex of the tree $\hat{T}^{*}$ can be thought of as identical to the corresponding ones in $\hat{T}$.
    %lalalala
    We will define a
    %        cellular map
    map
    $\Psi$ from $VT$ to $\hat{T}^{*}$ and then prove it extends to a cellular map from $\afld{T}$ to $\hat{T}*$ for which we keep the same name
    mapping the star of any vertex $v$ injectively into the one around $\Psi(v)$, bijectively in case $v$ is a surface-type vertex.
    \begin{comment}
    \begin{lemma}
    	\label{locally injective}Let $q$ be the projection map from $T$ to $\afld{T}$. The map $\Psi$ can be written as $\Theta\circ q$, for a cellular map $\Theta$ the star of any vertex $v$ injectively into the one around $\Psi(v)$ and bijectively onto the star around $\Psi(v)$ in the case $v$ is a surface-type vertex.
    \end{lemma}
    \end{comment}
    At that point one can apply the following easy fact:  % possibly move this to the appropriate section
    \begin{fact}
    	Any locally cellular map $\fun{f}{T}{T'}$ between $G$-trees injective on stars of vertices is an isomorphic embedding.
    \end{fact}
    \begin{comment}
    \begin{proof}
    	It is enough to show $f$ is injective on $VT$. Suppose that was not the case and consider $v,w\in VT$ at minimal distance such that $f(v)=f(w)$. Since $f$ is locally injective, we must have $d(v,w)\geq 3$, but then the images of distinct vertices in the interior of the simple path from $v$ to $w$ are distinct, which means this path maps by $f$ onto a simple closed loop, contradicting the fact that $T'$ is a tree.
    \end{proof}
    \end{comment}
    From this we will prove the map $\Theta$ to be a graph embedding.
    \begin{corollary}
    	The morphism $\iota$ is injective.
    \end{corollary}
    \begin{proof}
    	If $g\in G\setminus\{1\}$ does not fix $\afld{T}$, the same is true for $\iota(g)$ and $\hat{T}*$, so $\iota(g)\neq 1$. There might be non-trivial elements $g\in G$ acting trivially on $T_{v}$. Since the restriction of $\iota$ to any stabilizer of $\afld{T}$ is injective, though, in this case $\iota(g)\neq 1$ as well.
    \end{proof}
    
    \begin{proof}[of the existence of $\Psi$]
    	We first define $\Psi$ on $\lft{V}=V$ by letting $\Psi(\lft{v})=\gamma_{v}^{-1}\cdot\lft{\Phi(v)}$ for $v\in V$. After that we extend $\Psi$ equivariantly to the whole $VT$. That is, we send any $v\in VT$ of the form $g\cdot v_{0}$ for $v_{0}\in\lft{Z}$ to $\iota(g)\cdot\Psi(v_{0})$. This is well defined since $\iota(\Delta_{v})$ stabilizes $\Psi(v)$ for $v\in \lft{V}$.
    	\begin{claim}
    		Given $e\in T$ the vertices $\Psi(\alpha(e))$ and $\Psi(\omega(e))$ are adjacent in $\hat{T}^{*}$
    		%        		lululu\overset{}{}
    	\end{claim}
    	\begin{proof}
    		Let $v=\alpha(e)$ and $w=\omega(e)$. The map $\Psi$ preserves vertex type, hence $v':=\Psi(v)$ and $w':=\Psi(w)$ are two vertices of distinct type (one of them abelian). The fact that $\tg\neq\iota(\Sb{v}{T}\cap\Sb{w}{T})\subset Stab_{\hat{G}}(v')\cap Stab_{\hat{G}}(w')$ implies $v'$ and $w'$ are adjacent. Indeed, otherwise some diameter $2$ segment of $\hat{T}^{*}$ centered at a vertex of non-abelian type would be fixed by a non-trivial element, which is impossible, since $\hat{T}^{*}$ is in normal form.
    	\end{proof}
    	For any $e\in ET$ we let $\Psi(e)$ be the edge given by the previous claim. One needs to show that
    	the map $\Psi$ folds exactly as much as the quotient map from $T$ onto $\afld{T}$, namely:
    	\begin{lemma}
    		Any two edges $e,f$ of $T$ with $\alpha(e)=\alpha(f)=v$ have the same image by $\Psi$ only if:
    		\elenco{
    			\item $v$ is of non-rigid type and $e= f$
    			\item $v$ is of abelian type and $e$ and $f$ are in the same orbit under the action of the peripheral subgroup of $Stab(v)$
    		}
    	\end{lemma}
    	\begin{proof}
    		Since the stabilizers of rigid vertices of $\hat{T}^{*}$ are $CSA$, the proof amounts to showing that for $u\in VT$ of non-rigid type:
    		\enum{i)}{
    			\item \label{stabilizers} For any $e\in ET$ the intersection $\Sb{v}{G}\cap Stab(\Psi(e),G)$ is equal to:
    			\elenco{
    				\item $Per(T;v)$ if $e$ is adjacent to $u\in VT_{r}$ and $v\in VT_{a}$
    				\item $\Sb{e}{G}$ otherwise
    			}
    			\item \label{orbits} For any edges $e,f$ of $T$ originating at a common non-rigid vertex $v$ and which are in different orbits by the action of $Stab(v)$, their images $\Psi(e)$ and $\Psi(f)$ lie in different orbits by the action of $\Sb{v}{G}$.
    		}
    		Local injectivity at rigid type vertex $v$ depends on the fact that the \gad $\Delta$ is normalized. Edges incident at $v$ in different orbits will have mutually stabilizers non-conjugate in $\Sb{v}{G}$, so the same holds for their images by $\Psi$ with respect to $\iota(Stab(v))$.
    		
    		As usual, it is enough if we check the desired properties for the lifts of vertices in $V$.
    		Start by considering $v\in V_{a}$. Then $\Psi(\lft{v})=g\cdot x$, where $x=\lft{\hat{[v]}}$ and $g=\gamma_{v}^{-1}\delta_{v}$. Of course for any two $w,w'$ adjacent to $\lft{v}$ which are of different type their images are of different type as well, so they lie in different orbits; the same is true in case
    		%$w,w'\in VT_{s}$ belong to different orbits by construction.
    		$w,w'$ are both of surface type.
    		
    		Suppose we are given an edge $e$ from $v$ to some vertex of surface type $w$. The definition of the attaching maps in $\hat{\Delta}$ reformulates as saying that $Stab(\lft{\hat{e}})=\theta(\Sb{e}{G})\subset H_{\mu_{v}}$. Now, $\Psi(\lft{e})$ is a translate of $\lft{\hat{e}}^{*}$, so the group $Stab(\Psi(\lft{e}),\hat{G})$ is a conjugate in $\hat{G}$ of $\theta(\Delta_{e})$. Since $\iota(Stab(\lft{e}))=\iota(\Delta_{e})$ is of that form as well
    		and $\hat{G}$ is $CSA$\footnote{It admits a commutative transitive \gad with $CSA$ rigid type vertex groups. } . As $\iota$ is injective on $\Delta_{v}$, this proves (\ref{stabilizers}) for those edges ending in a surface type vertex.
    		%Observe that the two arguments given yield the two desired properties for the case in which $v\in V_{s}$.
    		
    		A more careful analysis is needed to deal with adjacent rigid type vertices. To begin with, let $RE_{v}$ be the set of all $e\in E$ from $v$ to a vertex of rigid type and $RV_{v}=\{\omega(e)\}_{e\in RE_{v}}$. For each $u\in RV_{v}$ we have $\Phi(u)=v_{\mu_{v}}$,
    		lifting to some $y_{u}\in Z$. For any $e=(v,u)\in E$, with $u\in RV_{v}$ the lift $\lft{e}\in ET$ ends in the vertex $t_{e}\lft{u}$. Now, $\Psi(\lft{u})=\gamma_{u}^{-1}\cdot y_{u}$, so by equivariance we get:
    		\begin{equation}
    			\Psi(\omega(\lft{e}))=\iota(t_{e})\gamma_{u}^{-1}\cdot y_{u}  \label{image by Psi}
    		\end{equation}
    		Remember that $\Psi(\lft{v})=\gamma_{v}^{-1}\delta_{v}\cdot x$. To each $e=(v,u)\in RE_{v}$ we know there is an element $h_{e}$ which is either $1$ or of the form $\bt_{e}^{\pm}\in\Delta_{\lft{[v]}}^{\delta_{v}^{-1}}$ and:
    		\elenco{
    			\item If $e\in Z$, then $\gamma_{u}=h_{e}\gamma_{v}$. So, as $t_{e}=1$, the equation (\ref{image by Psi}) yields $\Psi(\omega(\lft{e}))=\gamma_{v}^{-1}h_{e}^{-1}\cdot y$
    			\item If $e\nin Z$, then $\iota(t_{e})=\gamma_{v}^{-1}h_{e}^{-1}\theta(t_{e})\gamma_{u}$, so in this case we get
    			\begin{align*}
    				\Psi(\omega(\lft{e}))=\gamma_{v}^{-1}h_{e}^{-1}\theta(t_{e})\cdot y=\gamma_{v}^{-1}h_{e}^{-1}\cdot y_{u}
    			\end{align*}
    			recall that $\theta(t_{e})$ belongs $H_{\mu_{v}}$.
    		}
    		%        		In particular, a Dehn twist over an edge not in $Z$ restricts to the identity on $\pi_{1}(\Delta\rst_{Z},Z)$.
    		Now, for any distinct $e,e'\in RE_{v}$ the element $h_{e}h_{e'}^{-1}\nin\iota(\Delta_{v})^{\gamma_{v}^{-1}}=\zeta_{v}(\Delta_{v})^{\delta_{v}^{-1}}\subseteq\Delta_{\hat{[w]}}^{\delta_{v}^{-1}}$, and therefore $\Psi(\lft{e})$ and $\Psi(\lft{e'})$ lie in different orbits with respect to the action of $\iota(Stab(\lft{v}))=\iota(\Delta_{v})$.
    		Remember, moreover, that  $\Delta_{v}=P_{v}\oplus K_{v}$, where $P_{v}$ is $Per(v)$, and the map $\zeta_{v}$ has been constructed in such a way that it sends $P_{v}$ to $\hat{\Delta}_{e_{[v]}}^{\delta_{w}^{-1}}\subset\hat{\Delta}_{\hat{[v]}}^{\delta_{w}^{-1}}$ and $K_{v}$ into some complementary factor $\hat{\Delta}_{\hat{[v]}}^{\delta_{w}^{-1}}$. This clearly implies the property (\ref{stabilizers}) for $\lft{v}$.
    	\end{proof}
    \end{proof}
    
  \subsection{Factoring}
    
    %%% INcompatibility between the $c$ notation !!!!!!!!!
    %% we can change things here if we want to put everything in the tree language
    We shall now prove that for any morphism $\fun{f}{\qG{G}_{A}}{\rqG{F}{A}}$ factoring through the \pr-\res $\LRL{R}{A}{(J,r)}$ there is one $\fun{\hat{f}}{\rqG{R}{A}^{\lambda}}{\rqG{G}{A}}$ factoring through $\LRL{S}{A}{(J,r)}$ such that $f=\hat{f}\circ\iota$.
    We claim the statement reduces to the following lemma:
    \begin{lemma}
    	\label{mod lifting} For each $\lambda\in\hat{J}$ (of non-free product type) and each $\sigma\in Mod(\GRla)$ there is $\hat{\sigma}\in Mod(PMod(GSla;A))$ such that
    	\begin{equation}
    		\label{mlifteq}\hat{\sigma}\circ\iota_{\lambda}=\iota_{\lambda}\circ\sigma
    	\end{equation}
    \end{lemma}
    Indeed, in that case given a \trm $(f_{\lambda},\sigma_{\mu})_{\lambda\in\lvs{J},\mu\in\hat{J}}$ of $\RRA$ let $(\hat{f}_{\lambda},\hat{\sigma}_{\mu})_{\lambda\in J,\mu\in\hat{J}}$ be a \trm of $\RSA$ such that $\hat{f}_{\lambda}=f_{\lambda}$ for $\lambda\in\lvs{J}$. Given $\lambda\in J$ it is straightforward to prove by induction that $\hat{f}_{\lambda}\circ\iota_{\lambda}=f_{\lambda}$.
    This is trivial in case $j\in\lvs{J}$. Otherwise we have:
    \begin{align*}
    	\hat{f}_{\lambda}\circ\iota_{\lambda}=\bdcup{\hat{f}_{\mu}}{\mu\prec\lambda}{}\circ \eta^{\lambda,S}_{\lambda,S}\circ\hat{\sigma}_{\lambda}\circ\iota_{\lambda}=\\
    	=\bdcup{\hat{f}_{\mu}}{\mu\prec\lambda}{}\circ \eta^{\lambda,S}_{\lambda,S}\circ\iota_{\lambda}\circ\sigma_{\lambda}=
    	\bdcup{\hat{f}_{\mu}}{\mu\prec\lambda}{}\circ\bdcup{\iota_{\mu}}{\mu\prec\lambda}{}\circ \eta^{\lambda,R}_{\lambda,R}\circ\sigma_{\lambda}=\\
    	=\bdcup{\hat{f}_{\mu}\circ\iota_{\mu}}{\mu\prec\lambda}{}\circ \eta^{\lambda,R}_{\lambda,R}\circ\sigma_{\lambda}=\bdcup{f_{\mu}}{\mu\prec\lambda}{}\circ \eta^{\lambda,R}_{\lambda,R}\circ\sigma_{\lambda}=f_{\lambda}
    \end{align*}
    Where the first equality follows from the properties satisfied any \trm of $\RSA$, the second from equality \ref{mlifteq}, the third from the properties of the family $(\iota_{\lambda})_{\lambda\in J}$, the fifth by the induction hypothesis and the last one by the equalities valid for any \trm of $RRA$.
    
    \begin{proof}
    	[of \ref{mod lifting}] First of all, notice that it is enough to prove the property for a set of generators of $PMod(\Delta)$. Indeed, suppose that extensions $\hat{\sigma}$ and $\hat{\tau}$ exist for $\sigma,\tau\in PMod(\Delta)$. Then $\iota\circ\sigma\circ\tau=\hat{\sigma}\circ\iota\circ\tau=\hat{\sigma}\circ\hat{\tau}\circ\iota$, so  $\hat{\sigma\circ\tau}=\hat{\sigma}\circ\hat{\tau}$ satisfies the condition
    	%      	\footnote{A priori, different decomposition of $\sigma\circ\tau$ as a product of modular automorphism whose lift is already defined might yield a different result, so the construction does not imply automatically the existence of an $\iota$-compatible embedding of $PMod(\Delta)$ into $PMod(\hat{\Delta})$, but is enough for our purposes.}
    	%      	e get $c_{e}=c_{\lft{e}}\in\ker{\pi}$ for $e\in E$ and $\sigma_{v}\in Aut(\Delta_{v})$ for $v\in V$ preserving the \pr-structure on $\Delta_{v}$ and the identity for $v\in V_{r}$ such that
    	%      	letting $d_{e}=[\lft{v_{0}},\omega(\lft{e})]$ and $d_{e}=[\lft{v_{0}},\lft{v}]$
    	Obviously, for $\inn{c}$ we can take $\hat{\inn{c}}=\inn{\iota(c)}$. Contrary to the definition, to make the notation lighter when denoting the generators of the \pr-modular group we will refer to edges and vertices of $V$ rather than to their lifts in the fundamental trees $T$ and $\hat{T}$.
    	In virtue of \ref{elementary automorphisms} it is then enough to prove the lemma for $\sigma$ of the following form:
    	\enum{i)}{
    		\item \label{twist case} A Dehn twist automorphism $\tau_{e,c}$ where $e\in E\setminus Z$ and $\lft{e}=(\lft{v},\lft{v'})$ points away from $\lft{v_{0}}$, where $c\in\ddn{v}$, one among $\{v',v\}$ is of rigid type and the other of abelian type.
    		\item\label{surface vga case} A natural extension $\vga{\rho}{c}{v}{E}$ of some $\rho\in Mod(\ddn{v})$ for some $v\in Z$ of surface type.
    		\item \label{abelian vga case} A vertex group automorphism $\vga{\rho}{1}{v}{E}$ extending some $\rho\in Aut(\ddn{v})$ which fixes $Per(T;v)$, for some $v\in Z$ of abelian type.
    		%     		\item \label{surface twist case} A Dehn twist $\tau_{e,c}$ over an edge $e$ from an abelian type vertex $v$ to a surface type vertex $w$.
    	}
    	Moreover, in the last two cases we can assume that $c_{e_{0}}=1$ for the unique $e_{0}\in \lft{Z}$ originating at $v$ and pointing towards $v_{0}$.
    	%      	\lft{e} points outwards from $\lft{v_{0}}$
    	Start with case (\ref{twist case}) and let $x=\abv{v}$. We will assume $v$ is of abelian type, the other one being idetical. Remember that $\DDa{v}$ can be written as $P_{\abv{v}}\oplus\bcsum{L_{u}}{u\sim v}{}\oplus\bcsum{\subg{\bar{t}_{f}}}{f\in O(v)}{}$,
    	%      	where associated with $e$ there is a generator $\bar{t}_{e}$ of a direct summand of $\DDh{x}$.
    	where $L_{u}=(K'_{v})^{\delta_{u}}$ for $\qG{K'_{v}}\subset\qG{\hat{G}}$ isomorphic via $\zeta_{u}$ to $\qG{K_{v}}$. % a complement to the peripheral subgroup of $\ddn{u}$
    	Consider the automorphism $\xi$ of $\DD{x}$ fixing the two leftmost direct summands in the expression above
    	and sending $\bar{t}_{e}$ to $\zeta_{v}(c)^{\delta_{v}}\bar{t}_{e}$ and thus $\bt_{e}$ to $\zeta_{v}(c)\bt_{e}$. Let $\bar{\xi}$ be its extension to $\hat{G}$ with trivial twisting elements.
    	Let us check that $\hat{\sigma}=\bar{\xi}$ serves as the appropriate lift.
    	We first check the behaviour of $\hat{\sigma}$ on a vertex group $\ddn{u}$ of $\Delta$. Notice that in any case $\bar{\xi}\circ\zeta_{u}=\zeta_{u}$, since $x$ has a single incident edge.
    	%      	$u$ has a single incident edge.
    	%      	 $\bar{\xi}$ is a vertex group automorphism of a vertex of $\hat{\Delta}$ with a single incident edge, while $\iota(c)=\zeta_{v}(c)^{\gamma_{v}}$.
    	%and all the elements of $\ddn{v}$ in the normal form of the image of $\zeta_{u}$ belong to conjugates of direct summands fixed by $\xi$.
    	There are two cases:
    	\elenco{
    		\item If
    		%      		$e\nin[v_{0},u]_{Z}$ (the path in $Z$ from $v_{0}$ to $u$)
    		$\lft{e}\nin[\lft{v_{0}},\lft{u}]$
    		, so that $\sigma\rst_{\ddn{u}}=Id_{\ddn{u}}$, then $\gamma_{u}$ is the product of elements of the form $\bt_{f}$ for $f\neq e$ or of the form $t_{\hat{f}}$ for $f\in E\setminus Z$, hence $\bar{\xi}(\gamma_{u})=\gamma_{u}$.
    		%      		 if $u\neq v$ this is obvious; in case $u=v$, notice that $\zeta_{v}$ sends $\ddn{v}$ to the conjugate by $c_{w}^{-1}$ of a summand of $\hDDa{v}$ fixed by $\hat{\sigma}$, so also fixed by $\hat{\sigma}$.
    		\begin{align*}
    			\bar{\xi}\circ\iota\rst_{\ddn{u}}=\bar{\xi}\circ\inn{\gamma_{u}}\zeta_{u}=\\
    			=inn{\bar{\xi}(\gamma_{u})}\circ \bar{\xi}\circ\zeta_{u}\circ\sigma\rst_{\ddn{u}}=\iota\circ\sigma\rst_{\ddn{u}}
    		\end{align*}
    		\item If $\lft{e}\in[\lft{v_{0}},\lft{u}]$, then $\sigma\rst_{\ddn{u}}=\inn{c}\rst_{\ddn{u}}$.
    		%      		Let $v_{0},e_{0},v_{1}\cdots v_{m+1}=u$ the path in $Z$ from $v_{0}$ to $u$, as in the definition of $\iota$;
    		Going back to the definition of $\gamma_{v}$, we can write $\gamma_{v}=g\bt_{e}\gamma_{v}$ where $\gamma_{v}$ and $g$ are products of Bass-Serre elements of $\hat{\Delta}$ and elements of $\ddn{w}$ for $w\neq u$ or elements of $\ddn{x}$ in a direct summand of $\ddn{u}$ fixed by $\sigma$. Hence
    		\begin{align*}
    			\bar{\xi}(\gamma_{u})=g\zeta_{v}(c)\bt_{e}\gamma_{v}=g\bt_{e}\gamma_{v}\zeta_{v}(c)^{\zeta_{v}}=\gamma_{u}\iota(c)
    		\end{align*}
    		So that:    		.
    		\begin{align*}
    			\bar{\xi}\circ\iota\rst_{\ddn{u}}
    			%      			=\bar{\xi}\circ\inn{g\bt_{e}\gamma_{v}}\zeta_{u}=
    			=\inn{\gamma_{u}\iota(c)}\circ \bar{\xi}\circ\zeta_{u}=\\
    			=\inn{\iota(c)}\inn{\gamma_{u}}\circ\zeta_{u}=\inn{\iota(c)}\circ\iota\rst_{\ddn{u}}=
    			\\=\iota\circ\inn{c}\rst_{\ddn{u}}=\iota\circ\sigma\rst_{\ddn{u}}
    		\end{align*}
    	}
    	Let us now check the equality on $t_{f}\neq 1$, for $f\in\E$ from an edge $u$ to an edge $w$.
    	If $e\neq f$, then $\sigma(t_{e})=c_{u}^{-1}t_{f}c_{w}$, where the assignment $x\mapsto c_{x}\in\{1,c\}$
    	is as in the definition of a Dehn twist.
    	For any $w\in\V$ we know from the previous case that $\bar{\xi}(\gamma_{w})=\gamma_{w}\iota(c_{w})$. On the other hand $\iota(t_{f})$ was defined as $\gamma_{u}^{-1}\epsilon_{f}\gamma_{w}$, where $\epsilon_{f}=\bt_{f}\theta(t_{f})$ in case $w\in V_{r}$, and $t_{\hat{f}}$ in case $w\in V_{s}$. Note that $\bar{\xi}(\epsilon_{f})=\epsilon_{f}$.
    	\begin{align*}
    		\bar{\xi}(\iota(t_{f}))=\bar{\xi}(\gamma_{u}^{-1}\epsilon_{f}\gamma_{w})=\\
    		=\bar{\xi}(\gamma_{u})^{-1}\bar{\xi}(\epsilon_{f})\bar{\xi}(\gamma_{w})=\\
    		%      		=\bar{\xi}(\gamma_{u})^{-1}\epsilon_{f}\bar{\xi}(\gamma_{w})=\\
    		=\iota(c_{u})^{-1}\gamma_{u}^{-1}\epsilon_{f}\gamma_{w}\iota(c_{w})=\\
    		=\iota(c_{u})^{-1}\iota(t_{f})\iota(c_{w})=\\
    		=\iota(c_{u}^{-1}t_{f}c_{w})=\iota(\sigma(t_{f}))
    	\end{align*}
    	as desired.
    	If $e=f$, then $\sigma(t_{e})=c^{-1}t_{e}$ and $u=v$. Also $c_{u}=c_{w}=1$.
    	Now $\iota(t_{e})=\gamma_{v}^{-1}\bt_{e}\theta(t_{e})\gamma_{w}$, so that:
    	\begin{align*}
    		\bar{\xi}(\iota(t_{e}))=
    		=\bar{\xi}(\gamma_{v}^{-1}\bt_{e}\theta(t_{e})\gamma_{w})=\\
    		=\gamma_{v}^{-1}\bt_{e}\theta(c)^{-1}\theta(t_{e})\gamma_{w}=\\
    		%      		=\gamma_{v}^{-1}\bt_{e}\zeta_{v}(c)^{-1}\zeta_{v}(t_{e})\gamma_{w}=\\
    		=\gamma_{v}^{-1}\bt_{e}\zeta_{v}(c^{-1}t_{e})\gamma_{w}=\\
    		=\iota(c^{-1}t_{e})=\iota(\sigma(t_{e}))
    	\end{align*}
    	\newcommand{\bhx}[0]{\bar{\hat{\xi}}_{(\hat{c}_{e})}}
    	\newcommand{\sbhx}[0]{\bar{\hat{\xi}}}
    	We now deal with case (\ref{surface vga case}). Here $\rho$ extends an automorphism of $\Delta_{v}$ for some $v\in V$ of surface type. Let $Or$ be the set of all $e\in\E$ with $\alpha(e)=v$ and $\hat{Or}=\{\hat{e}\,:\,e\in Or\}$.
    	%      	in which $\sigma=\bar{\xi}_{(c_{e})_{e}}$ extending an automorphism of some surface type vertex group $\ddn{v}$.
    	As $\hat{\sigma}$ we take $\sbhx=\bhx$, where $\xi=\zeta_{v}\circ\rho\circ\zeta_{v}^{-1}$
    	and $\hat{c}_{\hat{e}}=\iota(c_{e})^{\gamma_{v}^{-1}}$ for each $e\in\hat{Or}$. Extend the assignment $e\mapsto\hat{c}_{e}$ to $\hV$ and the whole of $\E$ as described in definition \ref{elementary automorphisms} (by lifting to the fundamental domain).
    	We start by showings the following:
    	\begin{lemma}
    		\label{twisting check} The following equalities hold for $u\neq v$:
    		\begin{align*}
    			\hat{c}_{\Phi(u)}\bhx(\gamma_{u})=\gamma_{u}\iota(c_{u})  \\
    			\bhx(\gamma_{v})=\gamma_{v}
    		\end{align*}
    	\end{lemma}
    	\begin{proof}
    		%      		Clearly $\iota(t_{e})=\gamma_{u}^{-1}t_{\hat{e}}\gamma_{w}$, and $c_{u}=c_{e}$ for $e\in Or$.
    		For $x\in\Va$ observe that $\hat{\sigma}(\bt_{e})=\bt_{e}^{c_{\abv{v}}}$ for any of the generators $\bt_{e}\in\hDDa{x}$, as $\delta_{x}$ belongs to a vertex $v_{\mu}$ for which $\hat{c}_{v_{\mu}}=c_{\abv{x}}$.
    		Let $v_{0}=undefined,u_{m}=u$ be the unique simple path in $\Z$ from $v_{0}$ to $u$. There are two possibilities:
    		\enum{i)}{
    			\item In case $v\neq u_{i}$ for all $i$, let $0\leq i_{1}<i_{2}<\cdots i_{l}\leq m$ enumerate the (possibly empty) set of those indices $i$ for which $e_{i}$ is incident to a rigid type vertex and $\hat{e_{i}}\nin\MT$; extend the notation by letting $i_{0}=0$  and $i_{l+1}=m$. Then, for $0\leq j\leq l$, the value $\hat{c}_{\Phi(u_{i})}$ is constant, say equal to $c_{j}$ for all $i\in\ocint{i_{j}}{i_{j+1}}$. If we let $f_{j}=\hat{e_{i_{j}}}$, then of course
    			\begin{align*}
    				\bhx(t_{f_{j}})=\hat{c}_{\alpha(f_{j})}^{-1}t_{f_{j}}\hat{c}_{\omega(f_{j})}=c_{j}^{-1}t_{f_{j}}c_{j+1}
    			\end{align*}
    			Now, we know that	$\gamma_{u}=a_{l}t_{f_{l}}a_{l-1}\cdots t_{f_{1}}a_{0}$, where $a_{j}$ is the product of either elements in some $\ddn{x}$ and $c_{\Phi(x)}=c_{j}$, or of the form $t_{\hat{f}}$, where $\Phi(\alpha(f))=\Phi(\omega(f))=c_{j}$. We conclude that $\bhx(a_{j})=a_{j}^{c_{j}}$ for all $j$. Since $c_{u}=1$ and $c_{v_{0}}=1$ we have:
    			\begin{align*}
    				\hat{c}_{\Phi(u)}\bhx(\gamma_{u})=c_{l}(a_{l}^{c_{l}}(c_{l}^{-1}t_{f_{l}}c_{l-1})a_{l-1}^{c_{l-1}}\cdots a_{0}^{c_{0}})=
    				\\=a_{l}t_{f_{l}}a_{l-1}\cdots t_{f_{1}}a_{0}=\gamma_{u}=\gamma_{u}\iota(c_{u})
    			\end{align*}
    			\item In case $v=u_{i^{*}}$ for some $i^{*}$, this index unique.
    			Observe also that for any $e\in\E$ originating at $v$ and such that $\hat{e}\nin\hat{Z}$, $\bhx(t_{\hat{e}})=\hat{c}_{\hat{e}}^{-1}t_{\hat{e}}$.
    			We can write $\gamma_{u}$ as $g\gamma_{v}$, where $g$ corresponds to the subpath $u_{i},e_{i}\cdots u_{m}$. Now, $\bhx(\gamma_{v})=\gamma_{v}$, as $\hat{c}_{\hat{e_{i}}}=\hat{c}_{u_{i}}=e_{i^{*}-1}=1$ for $i<i^{*}$. Similarly, $\bhx(g)=\hat{c}_{\Phi(u)}^{-1}g\hat{c}_{\hat{e_{u}}}$, where $e_{u}$ is the first edge in the path from $v$ to $u$ in $Z$. Since $\hat{c}_{\hat{e_{u}}}=\iota(c_{e_{u}})$, we get
    			\begin{align*}
    				\hat{c}_{\Phi(u)}\bhx(\gamma_{u})=g\hat{c}_{\hat{e_{u}}}\gamma_{u}=g\iota(c_{u})^{\gamma_{v}^{-1}}\gamma_{v}=g\gamma_{v}\iota(c_{u})=\gamma_{u}\iota(c_{u})
    			\end{align*}
    		}
    	\end{proof}
    	
    	We can now check the equality on $\ddn{u}$, for any $v\neq u\in\V$. Note, first of all, that $\sigma\circ\iota\rst_{\ddn{u}}=\inn{\iota(c_{v})}\circ\iota\rst_{\ddn{u}}$.
    	On the other hand, since $\zeta_{u}(\ddn{u})\subset\hDD{\Phi(u)}$, we get
    	\begin{align*}
    		\bhx\circ\iota\rst_{\ddn{u}}=\inn{\bhx(\gamma_{u})}\circ\bhx\circ\zeta_{u}=\\
    		=\inn{\bhx(\gamma_{u})}\circ\bhx\circ\zeta_{u}=\inn{\bhx(\gamma_{u})}\circ\hat{c}_{\Phi(u)}\circ\zeta_{u}=\\
    		=\inn{\hat{c}_{\Phi(u)}\bhx(\gamma_{u})}\circ\zeta_{u}=\\=\inn{\hat{c}_{\gamma_{u}}\iota(c_{u})}\circ\zeta_{u}=\\
    		=\inn{\iota(c_{u})}\circ\inn{\hat{c}_{\gamma_{u}}}\circ\zeta_{u}=\\
    		=\inn{\iota(c_{u})}\circ\iota\rst_{\ddn{u}}=\iota\circ\inn{c_{u}}\rst_{\ddn{u}}=\\
    		=\sigma\rst_{\ddn{u}}
    	\end{align*}
    	As for $v$ itself, $\bhx\circ\zeta_{v}=\zeta_{v}\circ\xi\rst_{\ddn{v}}$ by definition, hence:
    	\begin{align*}
    		\bhx\circ\iota\rst_{\ddn{u}}=\inn{\hat{c}_{\bhx(\gamma_{v})}}\circ\zeta_{u}\circ\xi\rst_{\ddn{v}}=\\
    		=\inn{\gamma_{v}}\circ\zeta_{v}\circ\xi\rst_{\ddn{v}}=\iota\circ\sigma\rst_{\ddn{v}}
    	\end{align*}
    	All is left to check is that $\bhx(\iota(t_{f})=\iota(\sigma(t_{f}))$.
    	Let us first check the case in which $f\in E\setminus Z$ originating at a vertex $u\neq v$ and ending at a vertex $w$.
    	We have $\sigma(t_{f})=(c_{u})^{-1}t_{f}c_{v}$, hence:
    	%      	 	equality \ref{twisting check} yields:
    	\begin{align*}
    		\iota(\sigma(t_{f}))=\iota(c_{u})^{-1}\iota(t_{f})\iota(c_{w})=\\
    		=(\gamma_{u}^{-1}\hat{c}_{\Phi(u)}\bhx(\gamma_{u}))^{-1}
    		(\gamma_{u}^{-1}\epsilon_{f}\gamma_{v})\gamma_{w}^{-1}\hat{c}_{\Phi(u)}\bhx(\gamma_{w})=\\
    		=\bhx(\gamma_{u})^{-1}(\hat{c}_{\Phi(u)}^{-1}\epsilon_{f}\hat{c}_{\Phi(w)})\bhx(\gamma_{w})=\\
    		=\bhx(\gamma_{u}^{-1}\epsilon_{f}\gamma_{w})=\bhx(\iota(t_{f}))
    	\end{align*}
    	Here $\epsilon_{f}$ is equal to either $\epsilon_{f}=\bt_{f}\theta(t_{f})$ or $t_{\hat{f}}$.
    	Up to taking inverses, the only case left case is that in which $f$ originates at $v$. Then $\bhx(\gamma_{v})=\gamma_{v}$ and using \ref{twisting check} again we get:
    	%      	$\sigma(t_{f})=c_{f}^{-1}t_{f}c_{w}$ and $\bhx(t_{\hat{f}})=\hat{c}_{\hat{v}}^{-1}t_{\hat{f}}$,
    	\begin{align*}
    		\iota(\sigma(t_{f}))=\iota(t_{f}c_{f})=\gamma_{v}^{-1}t_{\hat{f}}\gamma_{w}\iota(c_{w})=\\
    		=\gamma_{v}^{-1}t_{\hat{f}}\hat{c}_{\Phi(w)}\bhx(\gamma_{u})=\\
    		=\bhx(\gamma_{v}^{-1})\bhx(t_{\hat{f}})=\\
    		=\bhx(\gamma_{v}^{-1}t_{\hat{f}}\gamma_{w})=\bhx(\iota(t_{f}))
    	\end{align*}
    	
    	In case (\ref{abelian vga case}) one can take as $\hat{\sigma}$ a vertex group automorphism extending an automorphism of $\hDDa{v}$ supported in the direct summand $K_{v}^{'c_{w}}$ corresponding to $v$. This case lacks the minor subtleties of the previous two and is left for the reader.
    	
    	In all these cases, it is clear that if $\sigma$ is an \pr-modular automorphisms and the twisting element $c_{e}$ maps to the center of $Q$ then the lift $\hat{\sigma}$ constructed above is a \pr-modular automorphism of $\hat{\Delta}$. Likewise, if $\sigma$ fixes $A$ the same is true of $\hat{\sigma}$. This concludes the proof.
    	%      	where $e_{u}$ is the first vertex of the simple path in $\MT$ from $v$ to $u$ and $w_{u}=\omega(e_{u})$.
    \end{proof}
    
    \begin{remark}
    	It follows from the proof that if $\QR{R}$ is \ws (taut)
    	the same is true for its completion. In the latter case the extension of a morphism factoring in a taut fashion does as well.
    \end{remark}

\section{\pr-test sequences}
  
  \label{test sequence section} In \cite{sela2}, a certain class of sequences of homomorphisms from a group endowed with a tower structure to the free group, called test sequences, are defined in terms of a certain (long) list of combinatorial conditions. Roughly speaking, the goal of the definition is, to ensure that the image of a fixed tuple of generators by the test sequence (often itself denoted by the term 'test sequence') eventually escapes any diophantine condition not already witnessed within the tower. This will be made more precise in the last chapter.
  
  For our purposes all we need is test sequences which are morphisms of \pr-groups, what we call \pr-test sequences.
  In an attempt to make the reading easier for the reader, our presentation differs slightly from that in \cite{sela2}. The definition of a \pr-test sequence provided below is (regardless of \pr) is in fact weaker than the version presented there, as it regards the geometric properties of the sequence in isolation from which we immediately draw a couple of simple consequences. The next section is devoted to the proof of the existence of a family of test sequences with certain properties. In Sela's terms this would be phrased as the existence of test sequences, without further qualifiers.
  
  \newcommand{\refin}[0]{Res(T';A)}
  Given a pegging $\qG{T^{q}}_{A}$ of $\RTA$, it will be convenient to replace the \pr-tower structure
  inherited from $\RTA$ by what we will call a 'refinement' of $\RTA$.
  The underlying index tree $J^{*}$, node groups and graph of groups decomposition $\RLR{T}{}{\lambda}$ will be again essentially independent of $q$.
  being the for each $T$.
  Decompose the complement of the peripheral group of each abelian vertex group appearing in $\QR{T}$ into a direct sums of cyclic groups; let $\mathcal{CS}$
  the collection of all the resulting cyclic groups and $\mathcal{SG}$ be that of all the surface type vertex groups appearing along $\QR{R}$.
  The properties we need $\LRLii{T}{q}{A}$ to have are as follows. Firstly for each $\lambda\in (\hat{J})$ either:
  \enum{a)}{
  	\item $\lambda$ is of free product type, with exactly two children
  	\item $\DR{T}{\lambda}$ has exactly one vertex $v$ of non-rigid type and either $\DR{T}{\lambda}_{v}$ is conjugated to
  	a member of $\mathcal{SG}$ or else $v$ is of abelian type, $\GTla$ contains at least a rigid vertex
  	and $\GTla_{v}$ is the direct sum of its peripheral group with a cyclic group conjugated to one in $\mathcal{CS}$.
  }
  Secondly, any group in $\mathcal{SG}$ or in $\mathcal{CS}$ appears at some node of floor type in the way described above.
  
  \begin{comment}
  Technically speaking, each $\qG{T^{q}}_{A}\in \pmb{@F}$ gives rise to a different refinement, the groups at the nodes are the same in each case.
  For completeness, we sketch the construction of $\LRLii{T}{q}{}$. It is by no means canonical, but this should not bother us.
  \end{comment}
  
  We start by placing $\qG{T^{q}}_{A}$ at the root $r$ (same as before) of our tower. Assume that for some $\lambda\in J$, the group $\Tl$ already appears as $\RLR{T}{}{\lambda}$ for some $\lambda^{*}\in J^{*}$ which is a leaf of the partial resolution. We extend our partially constructed tower $\LRL{T}{}{}$ by attaching at $\lambda$ a resolution at whose leaves we find precisely all the $\RLR{\bar{T}}{}{\mu}$ for $\mu\prec\lambda, \mu\in\lvs{J}$, as described now. Assign an order $v_{1},v_{2},\cdots v_{m_{S}}$ to the non-rigid surface groups appearing in $\Delta=\GTla$. In case $v_{i}$ is of abelian type, pick up an ordered base $(m_{i}^{j})_{j=1}^{n_{i}}$ of the complement of the peripheral subgroup of $\Delta_{u_{i}}$ and let $M_{i}^{j}=\subg{m_{i}^{k}}_{k=1}^{n_{i}}$.
  
  If removing $v_{1}$ from $\abs{\Delta}$ does not disconnect the graph, after renumbering the $v_{i}$, we can assume it does not disconnect our reference spanning tree either.
  At the nodes $nDFT(\lambda)$ we take the fundamental groups of each of the resulting connected components with respect to the restriction of the spanning tree.
  We proceed by removing the vertices of non rigid type of those components in the order specified by the given numbering, working with each of the components in parallel.
  If the vertex $v_{i}$ is of abelian type, we add a chain $\lambda^{*}=\mu_{0}\succ\mu_{1}{}{\succ}\cdots\mu_{n_{1}}$, where
  $\RLR{T}{}{\mu_{l}}$ is the amalgamated product of $M_{i}^{n_{i}-l}\frp_{Per(M_{i})} H$, for $H$ is the fundamental group of the graph of groups obtained removing $v_{i}$.
  If $\lambda\in J$ is of free product type, we can dislodge free factors one by one in a series of nodes of free product type with only two factors.
  
  In addition to this, we will impose a linear order, denoted by $\lo$, on $J^{*}$.
  On the set of children of any fixed $\lambda\in J^{*}$, we define $\lo$ arbitrary, except for the constraint that the child carrying the constants must be the minimum element.
  Finally, given any pair of distinct nodes $\lambda_{1},\lambda_{2}\in J^{*}$, consider their lowest upper bound $\nu$ and suppose that $\lambda_{i}\leq\mu_{i}$ for $\mu_{i}\prec\lambda$ for $i\in\{1,2\}$. We let $\lambda_{i}\lo\lambda_{j}$ if and only if  $\mu_{i}\lo\mu_{j}$, where $\{i,j\}=1,2$. This is a standard construction: it corresponds to a depth-first search of
  the nodes of the tree.
  %  The reader not-familiar with that notion can easily check that the resulting relationship is a linear order as an exercise.
  \newcommand{\phn}[2]{f^{#1}_{#2}}
  \renewcommand{\H}[0]{H}
  \newcommand{\GG}[0]{\RLR{T}{}{\lambda}}
  %  In the situation above we say that $\ssq{f^{2}}{n}$ grows faster than $\ssq{f^{1}}{n}$.
  
  \begin{definition}
  	\label{tsdef} Let $\mathcal{T}_{A}^{(J,r)}$ be a \rs \pr-tower structure of a group $\rqG{G}{A}$ over
  	\pr-limit groups $\{\Tl\}_{\lambda\in\Lambda}$ and suppose that for each $\lambda\in\Lambda$ we are given
  	a sequence $(g^{\lambda}_{n})_{n}$ of morphisms from $\TAl$ to $\rqG{F}{A}$ with trivial limit kernel. %% Or without
  	%  	 for a sequence of morphisms $\ssq{f}{n}$
  	%  	factoring through $\RTA$
  	We will say that a sequence $(f_{n})_{n}$ of morphisms factoring through $\RTA$ is a \emph{\pr-test sequence}
  	%  	 adapted to $\qG{T^{q}}_{A}$
  	relative to the family $\{(g^{\lambda}_{n})_{n}\}_{\lambda\in\Lambda}$
  	if and only if it has trivial limit kernel,
  	$(f_{n}\rst_{\Tl})_{n}$ is equal to some subsequence $\ssq{g^{\lambda}}{n}$ up to postcomposition by an inner automorhisms
  	%  	each $f_{n}$ is primitive
  	%  	and extends to $\hat{T}$ (we keep the same name) in a way that the
  	and the following properties are satisfied.
  	
  	For any $\lambda\in J\setminus\Lambda^{*}$, where $\Lambda^{*}\subset\lvs{J^{*}}$ is the set of leaves corresponding to those
  	in $\Lambda$, we require the sequence $\phn{\lambda}{n}$ to be unbounded and metrically convergent and
  	the action of $\GG$ on a limiting tree $Y$ associated to the sequence $(\phn{\lambda}{n})_{n}$, where $\phn{\lambda}{n}:=f_{n}\rst_{(T^{*})^{\lambda}}$ to be as follows:
  	\begin{enumerate}
  		\labitem{(FG)}{free group case} If $\lambda\in\lvs{J\setminus\Lambda^{*}}$, (which implies $\GG$ is unrestricted)	$Y$ is a geometric realization of a
  		simplicial $\GG$-tree $S$, with a single orbit of trivially stabilized vertices and as many orbits of trivially stabilized edges as the rank of the free group.
  		\labitem{(FP)}{free product case} If $\GG$ is of free product type, with children $\mu_{1}$ and $\mu_{2}$, then $Y$ is equivariantly isomorphic to the geometric
  		realization of the simplicial action associated with the free product $\RLR{T}{}{\mu_{1}}\frp \RLR{T}{}{\mu_{2}}$.
  		\labitem{(Fl)}{floor type} If $\RLR{T}{}{\lambda}$ is of floor type
  		then the action of $\GG$ on $Y$ admits a decomposition $\treeact{S}{Y}{}$ as a graph of actions.
  		in which $S$ is the tree dual to $\Gamma$.
  		There are two cases, according to the type of the unique non-rigid vertex $u$ of $\DR{T}{\lambda}$:
  		\begin{enumerate}
  			\labitem{(S)}{surface case} If $u$ is of surface type, then $\Gamma=\DR{T}{\lambda}$ 	and for any surface type vertex $u\in VS$, the action of $\Sb{u}{S}$ on $Y_{u}$ is dual to an	arational measured foliation on $\Sigma_{u}$ and is not the pullback of any similar action by the map induced by a covering map $\fun{q}{\Sigma}{\Sigma'}$, while the action of the stabilizers of rigid vertices is trivial.
  			\labitem{(A)}{pegged abelian case} If $u$ is of abelian type, in which case $Y$ is isomorphic to the geometric realization of the tree dual to an HNN extension, where the vertex group is that of the
  			single rigid type vertex of $\GTla$ and the Bass-Serre a generator of a non-peripheral direct summand of $Stab(u)$.
  		\end{enumerate}
  	\end{enumerate}
  	
  	In cases \ref{free group case} and \ref{free product case} above we also impose that for any finitely generated \rs \pr-limit group $\qG{H}_{A_{\lambda}}$ containing $(\mathcal{T}_{A}^{*})^{\lambda}$ and any diverging subsequence of extensions $(h_{k_{n}})_{n}$ of $(f_{k_{n}})_{n}$ to $\rqG{H}{A}$ with trivial limit kernel and limiting tree $Z$ either $\GG$ is elliptic in $Z$ or its minimal tree $Z_{\GG}$ satisfies:
  	\property{$\dagger$}{transverse property}{
  		For any $h\in H$, if $Z\cap h\cdot Z$ is non-degenerate, then $h\in\Tl$.
  	}
  	We also require:
  	\property{$\dagger\dagger$}{growth}{
  		$\Tm$ grows faster than $\RLR{T}{}{\nu}$ under the sequence $\ssq{f}{n}$ for $\mu\lo\nu$.
  		%  		$H_{0}$ is elliptic in any limiting tree for $\ssq{h}{n}$
  	}
  	\begin{comment}
  	Suppose we are given $\lambda_{0}\lo\lambda_{1}$ in $J^{*}$ and $H_{0},H_{1}\leq H$ a finitely generated limit group, where
  	$H_{i}$ is a copy of $Gl(T^{*};\lambda_{i})$ together with a sequence $\ssq{h}{n}\subset Hom(H,\F)$ with trivial limit kernel, such that $h_{n}\rst_{H_{i}}=\inn{\gamma_{n}^{i}}\circ f_{n}\rst_{H_{i}}$ for some $\gamma_{n}^{i}$.  :
  	\property{$\dagger\dagger$}{growth}{
  		$H_{0}$ is elliptic in any limiting tree for $\ssq{h}{n}$
  	}
  	\end{comment}
  \end{definition}
  As we will see in the course of the construction, the condition that the sequence $\ssq{f}{n}$ has limit kernel is to a certain extent redundant, and can be pulled up from the bottom of the tower even if one starts with a slightly weaker (but somehow cumbersome) definition.
  %  maybe change location
  We note in passing the following:
  \begin{observation}
  	Let $\LRL{T}{A}{(J,r)}$ be a \pr-tower structure of a group $G$ and $H$ a subgroup of $G$. There is a unique minimal $\lambda\in J_{0}$ such that $H$ can be conjugated into
  	$\Tl$.
  \end{observation}
  \begin{proof}
  	It follows easily from the following two facts:
  	\enum{i)}{
  		\item Given a retraction $r$ of a group $G$ onto some $H\leq G$, two subgroups $K_{1},K_{2}\leq H$ are conjugate in $K$ if and only if
  		they are conjugate in $H$.
  		\item Given a free product $G=\bfrp{G_{i}}{i=1}{m}$, a subgroup $H\leq G$ can be conjugated into $G_{i}$ for at most on value of $i$.
  	}
  \end{proof}
  
  \begin{lemma}
  	\renewcommand{\GG}[0]{\qG{G}}
  	\renewcommand{\H}[0]{H}
  	%  	define a test sequence for a general tower as that of some refinement
  	\label{extra properties}
  	%  	Let $\LRL{T}{A}{(J,r)}$ a \fine tower structure
  	Fix $\lambda\in J^{*}$ and let $\GG$ be $\RLR{T}{}{\lambda}$, $\H$ finitely generated containing $\GG$ and $g_{n}$ an extension of $\phn{\lambda}{n}$ to some finitely generated $\H$ for each $n\in\N$ such that $\limker{g}{n}=\tg$, with limting tree $Y$.
  	
  	Assume that $\H$ is freely indecomposable relative to $\{\GG\}\cup\mathcal{A}$, where $\mathcal{A}$ is the family of all the subgroups of $\H$ elliptic in $Y$.
  	%  and the action of $\F$ on its Cayley Graph. %maybe
  	Then the action of $H$ on $Y$ admits a decomposition as a tree of actions $\treeact{S}{Y}{}$
  	in which each vertex action is either of Seifert (surface), axial or simplicial type.
  	%  with non-trivially stabilized edges.
  	If $\GG$ is elliptic in $Y$, then all the the simplicial components are non-trivially stabilized.
  	If $\GG$ is not elliptic in $Y$, denote by $L$ either the non-rigid type vertex group of $\GTla$
  	in cases \ref{surface case} and \ref{pegged abelian case} or $\GG$ in case \ref{free group case}.
  	Then there is a unique orbit of non-degenerate components and the minimal tree $Y_{L}$ is equal to one of them $Y_{v_{0}}$.
  	\elenco{
  		\item  In case \ref{free group case} we have $\GG=\H$.
  		\item  In cases \ref{surface case} and \ref{free product case} the group $Stab(v_{0})$ coincides with $L$ and all the points of intersection betweeen $Y_{L}$
  		and other non-degenerate components are non-trivially stabilized in $L$.
  		%  	\item  In case \ref{pegged abelian case} $Stab(v_{0})$ is the centralizer of $L$ and its action on $Y_{L}$ either of simplicial or axial type.
  	}
  	%  Any simplicial segment not laying in a translate of $Y_{v_{0}}$ is non-trivially stabilized (in particular any such segments in all cases except \ref{free product case} and \ref{free group case}).
  	%  Any simplicial edge not contained in a translate of $Z$ has a non-trivial edge stabilizer.
  	%  	In case \ref{pegged abelian} the
  	%  In particular, translates of $Y_{L}$ intersecting in non-degenerate segments coincide.
  \end{lemma}
  
  \begin{proof}
  	First, observe that if the group $\GG$ is itself elliptic in $Y$, then $\H$ cannot be freely indecomposable relative to the family of all point stabilizers of $Y$. If it is not elliptic, the same is true in cases \ref{surface case} and \ref{pegged abelian case} of the definition of test sequences, since then $\GG$ is freely indecomposable relative to its rigid vertex groups, which are elliptic in $Y$ in virtue of the properties of test sequences. This means that the first alternative of \ref{Rips decomposition} does not hold, and we can decompose $Y$ as a normalized tree of actions $\treeact{S}{Y}{}$ of the desired type, with no trivially stabilized simplicial edges.
  	
  	So assume now that $\GG$ is not elliptic. Consider \ref{surface case} first. Let the action on $Z$ be dual to an arational measured folitation $\mathcal{F}_{0}$ on $\Sigma_{0}$. Observe that since the action of the vertex surface group is indecomposable and non-linear $Z\subset Y_{v_{0}}$ for some component $Y_{v_{0}}$ of Seifert type.
  	Of course, this implies $L\leq Stab(Y_{v_{0}})$. Now, boundary subgroups of $L$ fix points in $Y_{v_{0}}$, hence they are contained in
  	the ones associated to boundary components of $\Sigma$. In virtue of the results of \cite{scott1978subgroups}, the embedding of
  	$L$ into $\pi_{1}(\Sigma)$ is induced by some covering map $\fun{q}{\Sigma_{0}}{\Sigma}$.
  	The defining properties of test sequences imply $q$ has to be a homeomorphism, hence $L=Stab(v_{0})$.
  	In \ref{pegged abelian case}, the minimal tree $Z=Y_{L}$ is a line which according to \ref{factoring lines} has to be fundamental in $Y$ and the conclusion follows from corollary \ref{factoring lines}.
  	In \ref{free group case} and \ref{free product case}  and corollary \ref{Rips refined} yields a decomposition $\treeact{S}{Y}{}$ of the required type which includes $Z$ as a vertex tree $Y_{v_{0}}$.
  	The full strength of (\ref{transverse property}) implies that $Stab(v_{0})=L$.
  	%  	Case \ref{cyclic tip case} can be dealt in a similar way using property \ref{weak transverse property}(only needed in case the line stabilized by $\RLR{T}{}{\lambda}$ has trivial setwise st{abilizer).
  \end{proof}
  \begin{corollary}
  	\label{nice decompositions} In the situation above, if $\GG$ is not elliptic in $Y$, there is a minimal simplicial $\H$-tree $U$ in which $\GG$ is not elliptic and a $\GG$-subtree (sub-\gat) $W$ of $Y$ such that:
  	\elenco{
  		\item If $\lambda$ is as in case \ref{free group case}, then $\GG=\H$.
  		\item If $\lambda$ is as in \ref{free product case} $U$ has trivially stabilized edges and $W$ is equivariantly isomorphic to the tree dual to the decomposition 	$\RLR{T}{}{\lambda}=(\RLR{T}{}{\mu_{1}})^{\gamma_{\lambda}}\frp (\RLR{T}{}{\mu_{2}})^{\gamma_{\lambda}}$, where $\{\mu_{1},\mu_{2}\}=Ch(\lambda)$.
  		\item In case $\lambda$ is of floor type then $W$ is equivariantly isomorphic as a \gat to the dual tree of the decomposition $\DR{T}{\lambda}$.
  		In case \ref{surface case}, $\Sb{v}{\GG}=\Sb{v}{\H}$ holds for any non-rigid vertex $v$ of $W$.
  	}
  	Any subgroup of $\H$ elliptic in $Y$ is elliptic in $U$ and stabilizes a rigid vertex of $U$ in case $\lambda$ is of floor type.
  \end{corollary}
  \begin{proof}
  	Case \ref{free group case} is already dealt with in the lemma. In case \ref{pegged abelian case} the conclusion follows directly from lemma \ref{factoring lines}.
  	In case \ref{free product case}, we can apply lemma \ref{blow and collapse} to the skeleton $S$ and the simplicial tree underlying the action of $\GG$ on $Y_{v_{0}}$.
  	Recall that the stabilizers of edges of $S$ incident to $v_{0}$ can only stabilize points of $Y_{v_{0}}$ which are non-trivially stabilized in $\GG$ and is therefore elliptic in the
  	associated simplicial tree. In case \ref{surface case} the decomposition can be obtained by collapsing all the edges of $S$ which are not incident to a translate of $v_{0}$. Each edge
  	of $S$ incident to $v_{0}$ corresponds to a point of $Y_{L}$ stabilized by a boundary subgroup of $L$ and have that same stabilizer in $\H$. This implies the union of the translates
  	by elements of $\GG$ of the star around $v_{0}$ is isomorphic to the tree dual to $\DR{T}{\lambda}$.
  	% 	 	argument in the proof of the embedding
  	A subgroup of $\H$ elliptic in $Y$ is elliptic in $S$. On the other hand any subgroup of $Stab(v_{0})$ elliptic in $Y_{v_{0}}$ is contained in the stabilizer of an incident edge group.
  	This proves the last claim for cases \ref{free group case} and \ref{surface case}. In case \ref{pegged abelian case} the property was already contained in lemma \ref{factoring lines}.
  \end{proof}
  
  \subsection{Two basic properties of test sequences}
    
    \begin{lemma}
    	\label{primitive pegs} Assume $\LRL{T}{A}{(J,r)}$ is a closed tower structure of $\rqG{G}{A}$ and that all members of a given \pr-test sequence $\ssq{f}{n}:\qG{T}^{q}_{A}\to\rqG{F}{A}$ restrict to the same embedding $\kappa$ of the \lvlgp $\rqG{H}{A}$ which carries the constants into $\rqG{F}{A}$.
    	%  	peg $p$ of a maximal abelian subgroup $M$ of $G$ either $p$ has a conjugate $p^{g}$ which belongs to the group at a leaf carrying the constants or
    	Assume that for $c\in G$ either $Z_{G}(c)=\subg{c}$ or a peg of $Z_{G(c)}$. Then either of the following happens:
    	\enum{i)}{
    		\item Some conjugate $c^{g}$ of $c$ belongs to the terminal group of $\RRA$ which carries the constants.
    		\item For $n$ big enough $f_{n}(c)$ is primitive in $\F$.
    	}
    	%  	and each $n\in\N$, let $d_{n}(p)=[Z_{\F}(f_{n}(p)):\subg{f_{n}(p)}]$. Then:
    	%  	chose a generating base $p_{i}=c_{1}^{i},c_{2}^{i},\cdots c_{rk(M_{i})}^{i}c_{i}^{j}$ of $M_{i}$,
    	%  where $p_{i}$ is a peg of $M_{i}$ and for $2\leq i\leq m$, the family $\{p_{j}^{i}\}_{2\leq i\leq j}$ is the union of some choice of complementary bases
    \end{lemma}
    \begin{proof}
    	Assume for the sake of contradiction none of the two options were true and let $\lambda\in J\setminus DFt(\mathcal{T})$ be minimal such that $c$ is conjugate to some $d\in\Tl=\qG{L}_{A_{\lambda}}$. Let $K$ the amalgamated product of $M=\subg{d,u}_{ab}$ and $L=\RLR{T}{}{\lambda}$ over $\subg{d}$. Notice that this is a $CSA$ group, since $d$ generates its own centralizer in $L$.
    	
    	Then a metrically convergent sequence $\ssq{g}{n}$ of extensions of $\ssq{f}{n}$ to $K$ exists
    	so that $g_{n}(u)^{k_{n}}=f_{n}(d)$ for some $k_{n}\nin\{-1,0,1\}$. Let $P=K/\limker{n}{g_{n}}$. The first observation is that $Z_{Q}(d)\neq\subg{z}$ as for no $n,m\in\N$ the element $d^{m}u^{-1}$ in the kernel of $g_{n}$, and hence that $L$ does not map onto $P$. Clearly $P$ is freely indecomposable relative to $H$, since free factors are centralizer closed. Consider now the action of $P$ on a limiting tree for the sequence $L$.
    	
    	Observe that $L$ is not elliptic in $Y$. Indeed, suppose for a second it is elliptic we claim that then $M$ must necessarily be elliptic as well. For $M$ can only act by translations on its minimal tree (abelianity of an action is clearly preserved by taking limits in the \egt). On the other hand $\tl{u}{Z}\leq\lim_{n}\frac{1}{k_{n}}\tl{d}{Z}$, so $M$ in this case is generated by elliptic elements.
    	As a result there are points $x$, stabilized by $M$ and $y$ whose stabilizer contains $L$ such that $Y$ is equivariantly isomorphic to a geometric realization of the tree dual to the amalgamated product of $Stab(x)$ and $Stab(y)$ over the set-wise stabilizer of $[x,y]$ (since by superstability if an element of $Stab(x)\cup Stab(y)$ fixes a non-degenerate segment of $[x,y]$ then it must fix the entire $[x,y]$). But this is impossible, since implies that either $d$ stabilizes a
    	non-degenerate tripod (if the index of the set-wise stabilizer of $[x,y]$ in $\subg{u,d}$ is greater than $2$) or that some element of $M$ acts as a reflection on some infinite line (if the index is equal to $2$).
    	
    	%    	MIStake
    	So $G$ is not elliptic in $Y$. It follows from minimality and our assumptions on $d$ that $\lambda$ falls in one of cases \ref{free product case}, \ref{free group case} or \ref{surface case}
    	and $d$ is not elliptic in $Y$.
    	%    	Since $[u,d]=1$, $u\cdot Ax(d)=Ax(d)$.
    	%  Let $P$ be either the surface vertex group of $\GTla$ in case \ref{surface case} and $L$ otherwise.
    	%   lemma \ref{extra properties} implies
    	%  for $q\in Q$ the intersection $Y_{P}\cap q\cdot Y_{P}$ can only be degenerate if
    	On the other hand, consider the decomposition of the action of $P$ on $Y$ as a tree of actions $\treeact{S}{Y}{}$ provided by corollary \ref{extra properties}.
    	In case \ref{surface case} let $S_{0}\subset S$ be the miminal tree of $L$ in $S$ and $Z=\abunion{v\in S_{0}}{Y_{v}}{}$. In the other ones let $Z=Y_{v}$ for the unique $v$ stabilized by $L$.
    	Observation \ref{fundamental subtrees} implies that $Z$ is fundamental and $Stab(Z)=L$ in both. This means that $P=L$, since $Ax(d)\subset\abunion{v\in VS_{0}}{Y_{v}}{}$ fixed by $u$ as well (in other words $ug$ belongs to the limit kernel of $\ssq{g}{n}$ for some $g\in G$ ).
    \end{proof}
    
    As a warm up before the proof of Merzlyakov's theorem, let us give one of the simplest and most natural examples of the idea behind test sequences.
    \begin{comment}
    principle behind the construction of (\pr-)test sequences:
    if the image of some tuple of elements of the tower by the maps in the test sequence satisfies some diophantine (\pr-constrained in our case) diophantine condition, up to the addition of some roots
    the tuple itself must satisfy the condition within the tower.
    \end{comment}
    \begin{lemma}
    	Let $\LRL{T}{A}{(J,r)}$ be a closed tower structure
    	where $A$ is full and malnormal in $\F$,
    	$u,u_{2}\in \hat{T}$ non-trivial elements of $\hat{T}$ and for some \pr-test sequence $\ssq{f}{n}$
    	relative to a fixed embedding $\kappa$ of the \lgp group carrying the constants into $\F$
    	and $n$ big enough $f_{n}(u_{1})=f_{n}(u_{2})^{t_{n}}$ for some $t_{n}\in\F$.
    	Let $\kappa$ be image of
    	$u_{1}$ and $u_{2}$ are conjugate in $\hat{T}$.
    \end{lemma}
    \begin{proof}
    	Assume without loss of generality that $u_{j}$ belongs to $\RLR{T}{}{\lambda_{j}}$ for the minimal $\lambda_{j}\in J^{*}$ such that $u_{j}$ can be conjugated to $\RLR{T}{}{\lambda}$
    	(this is well defined in virtue of \ref{nice decompositions}).
    	%    	If $\lambda_{1}$ and $\lambda_{2}$ are comparable, say $\lambda_{1}\leq\lambda_{2}$ consider the quotient of
    	Consider first the case in which $\lambda_{1}\neq\lambda_{2}$ and
    	%    	 for no $j\in\{1,2\}$ does $\RLR{T}{}{\lambda_{j}}$ carry the constants. Let
    	%    	$\lambda_{1}\neq\lambda_{2}$
    	assume without loss of generality that $\lambda_{1}\lo\lambda_{2}$. For simplicity let $H_{i}=\RLR{T}{}{\lambda_{i}}$ for $i=1,2$ up to replacing $\ssq{f}{n}$ by a subsequence, a geometrically convergent sequence $\ssq{g}{n}$ of homomorphisms from $H_{1}\frp H_{2}$ to $\F$ exist with the properties:
    	\enum{i)}{
    		%    		\item $g_{n}\rst_{H_{1}}=f_{n}\rst_{H_{1}}$.
    		\item $g_{n}\rst_{H_{1}}=f_{n}\rst_{H_{1}}$
    		\item $g_{n}\rst_{H_{2}}=\inn{t_{n}}\circ f_{n}\rst_{H_{2}}$ for some $t_{n}\in\F$ satisfying the property $f_{n}(u_{1})=f_{n}(u_{2})^{t_{n}}$ for which $\sl{*^{1}_{n}}{f_{n}(x)^{t_{n}s}}{}$ is minimal, where $x$ is a fixed tuple of generators of $H_{2}$.
    	}
    	Let $\fun{q}{H_{1}\frp H_{2}}{Q}$ be the limit quotient of $\ssq{g}{n}$ and $Y$ the associated limiting tree. We can think of $H_{1},H_{2}$ as a subgroups of $Q$, which necessarily satisfies equality $u_{1}=u_{2}^{t}=u_{1}^{ts}$. Property \eqref{growth} implies $\RLR{T}{}{\lambda_{1}}\leq G$ (hence $u_{1}$) fixes some $*_{1}\in Y$ which is the limit of the sequence $(*^{1}_{n})_{n}$. We claim that $H_{2}$ is elliptic as well. Indeed, otherwise $u_{2}$, not being conjugate to a non-rigid vertex group of $\DR{T}{\lambda}$, would have to act hyperbolically on $Y$, something clearly incompatible with being conjugate to $u_{1}$.
    	In view of this, let $*_{2}\in Y$ be fixed by $H_{2}$.
    	%    	Hence the sequence $(*^{2}_{n})_{n}$ converges to some
    	Since the action of $H$ on $Y$ is non-trivial, necessarily $*_{1}\neq *_{2}$. This implies that the action of $\subg{H_{1},H_{2}^{ts}}$ on its minimal tree in $Y$ must be isomorphic to some geometric realization of the tree dual to an amalgamated free product of the form $H_{1}\frp_{E}H_{2}^{ts}$, where $E$ is the abelian segment stabilizer of $I=[*_{1},*_{2}]$, since superstability of the action implies for any $g\in H_{1}$ either $g\in E$ or $g\cdot I\cap I=\{*_{1}\}$. But then \ref{general shortening} implies there is $c\in E$ and $n$ big enough $\sl{*^{1}_{n}}{g_{n}(x)^{g_{n}(c)}}{}<\sl{*^{1}_{n}}{g_{n}(x)}{}$, and the element $t'_{n}=t_{n}g_{n}(c)$ contradicts the minimality of $t_{n}$.
    	
    	If $\lambda_{1}=\lambda=\lambda_{2}$, then if $\RLR{T}{}{\lambda}$ carries the constants the claim is clear.
    	If it does not, let $H=H_{1}$ again, up to replacing $\ssq{f}{n}$ by a subsequence, a geometric convergent  $\ssq{g}{n}$ of homomorphisms form $H\frp\subg{t}$ to $\F$ exist such that:
    	\enum{i)}{
    		\item $g_{n}\rst_{H_{1}}=f_{n}\rst_{H_{1}}$
    		\item $g_{n}(u_{2})^{g_{n}(t)}=g_{n}(u_{1})$ and $t_{n}:=g_{n}(t)$ minimizes $\pl{*^{1}_{n}}{t_{n}}{}$ among those with that property.
    	}
    	Let $\fun{q}{H\frp\subg{t}}{Q}$ be the limit quotient of $\ssq{g}{n}$ and $Y$ the associated limiting tree. We can identify $H$ with $q(H)$, so that $u_{1}=u_{2}^{q(t)}$; this equation implies that $Q$ is freely indecomposable relative to $H$. If $H$ is elliptic in $Y$, then it fixes the limit $*_{1}$ of the sequence $(*^{1}_{n})$ and an argument entirely analogous to the one used above shows that $Y$ is isomorphic to some geometric realization of an HNN extension with vertex group $H$ and Bass-Serre element $t$. 	An application of \ref{general shortening} contradicts the minimality condition for $t_{n}$ in this case. If $H$ is not elliptic in $Y$ then neither are $u_{1}$ and $u_{2}$ and the intersection $Y_{H}\cap q(t)^{-1}\cdot Y_{H}$ must contain the whole axis of $u_{1}$. If $\lambda$ falls in cases \ref{free group case},\ref{surface case} or \ref{free product case} this implies right away that $q(t)\in u_{1}$ and we are done. If $\lambda$ is as in case \ref{pegged abelian case}, we know that there is a decomposition of $Q$ obtained from $\DR{T}{\lambda}$ by enlarging the two vertex groups (the non-rigid one stays abelian), dual to some tree $S$. If one among $u_{1},u_{2}$ fixes an abelian-type vertex of $S$, then the other must fix that same vertex as well and by the $CSA$ property $[q(t),u_{2}]=1$, implying that $u_{1}=u_{2}$.
    	
    	If one of them is hyperbolic in $S$, then the other one must be as well. By inspecting their normal form one can easily check that the only possible case in which they are not conjugate in $H$ is
    	that in which, up to conjugation in $Q$ (by different elements) we have $u_{1}=r_{1}a_{1}$, $u_{2}=r_{2}a_{2}$, where $a_{i}\in A$ and the two elements $r_{1},r_{2}\in R$ are not mutually
    	conjugate in $R$ but $r_{1}=r_{2}^{b}$ for some $b\in B\setminus H$ (boolean '$\setminus$'). This can be ruled out by induction hypothesis, so we are done.
    	
    	\begin{comment}
    	Without loss of generality we can assume that $g$ fixes some edge $e\in Ax(u_{1})\cap Ax(u_{2})$.
    	Let $v$ be the origin of $e$ of rigid type and $f_{1},f_{2}$ the other edges adjacent to $v$ in $Ax(u_{1})$ and $Ax(u_{2})$ respectively. Clearly $Stab(f_{1})^{g}=Stab(g_{2})$
    	by induction hypothesis this implies either $Stab(f_{1})$ is conjugated into the
    	\end{comment}
    	
    	\begin{comment}
    	a way that B X H equals the edge group the latter. If one among u1 , u2 fixes an abelian-type vertex of
    	S, then the other must fix that same vertex as well and by the CSA property rqptq, u2s “ 1, implying
    	
    	By \ref{elliptic expansion}
    	\end{comment}
    	
    \end{proof}

\section{Existence of \pr-test sequences}
  \label{TS construction section}
  Fix now a \rs \pr-tower structure $\LRL{T}{A}{(J,r)}$ relative to the the groups at the nodes of $\Lambda\neq (J)$ and suppose that for each
  $\lambda\in\Lambda$ a sequence $\ssq{g^{\lambda}}{n}$ of morphisms from $\TAl$ to $\qG{F}_{A_{\lambda}}$ with trivial limit kernel is given.
  %Fix some refinement $\LRL{T'}{A}$ of $\RTA$, as described in the previous section.
  For the remaining of the section, let $\RTA$ be a \rs \pr-tower which has already been refined, as described at the beginning of section \ref{test sequence section}
  The goal of this section is to prove the following result:
  \begin{proposition}
  	\label{test sequence existence} For any $q\in\pg{T}$ there is a family $\mathcal{F}$ of test sequences
  	of $\mathcal{T}_{A}$ relative to $(\ssq{g^{\lambda}}{n})_{\lambda\in\Lambda}$ such that:
  	\elenco{
  		\item $\mathcal{F}$ is closed under taking subsequences
  		%    		\item For each $\qG{T^{q}}_{A}\in \pmb{@F}$ the subfamily of all those memebers of $\mathcal{F}$ extending to $\qG{T^{q}}_{A}$ is closed under taking diagonals.
  		\item $\mathcal{F}$ is closed under diagonal subsequences
  		\item $\mathcal{F}$ is congruency complete
  	}
  \end{proposition}

  \subsection{Proof of proposition \ref{test sequence existence}}
    
    The proof is by induction on the size of $J$; the induction step involve showing the existence of an adequate family of \pr-test sequences of $\LRL{T}{A}{(J,r)}$ provided such a family exists for
    $\RTA\rst_{\lambda}$ for any $\lambda\in Ch(r)$.
    We will refer to a family of divergent sequences satisfying the first two points above as a CSD family.
    It is not hard to see that a CSD family $\mathcal{F}$ corresponds to a descending chain
    \begin{align*}
    	Hom(G,\F)\supseteq[\mathcal{F}]_{0}\supset[\mathcal{F}]_{1}\cdots
    	%	\supset[\mathcal{F}]_{2}\supset\cdots$
    \end{align*}
    with trivial intersection, from which  $\mathcal{F}$ is recoverable as the set of all sequences $\ssq{f}{n}$ such that $f_{n}\in[\mathcal{F}]_{n}$.
    %	\label{unboundedness remark}
    Clearly all sequences in such an $\mathcal{F}$ have trivial limit kernel if and only if for any $g\in G\setminus\{1\}$ there is some
    $n\in\N$ such that $f_{m}(g)\neq 1$ for all $\ssq{f}{n}\in\mathcal{F}$ and $m\geq n$.
    If all sequences in the family share a common limiting tree we say that the family is convergent. Note that in this case there is some $g\in G$ which grows uniformly fast, by which we mean that for any $L>0$ some
    $n$ exists such that $\tl{g}{\lambda^{f_{m}}}>L$ for any $m\geq n$ and $\ssq{f}{m}\in\mathcal{F}$.
    
    \begin{observation}
    	\label{asymmetrical growth} Suppose we are given groups $G_{1},G_{2}\cdots m$ and for $1\leq i\leq m$ a CSD family $\mathcal{F}_{i}$ of convergent sequences in $Hom(G_{i},\F)$.
    	Then there is a CSD family $\mathcal{F}$ of sequences of homomorphisms from $G_{1}\frp G_{2}\cdots G_{m}$ to $\F$	such that:
    	\elenco{
    		\item  $(f_{n}\rst_{G_{i}})\in\mathcal{F}_{i}$ for $1\leq i\leq m$
    		\item Given any tuple $(\ssq{f^{j}}{n})_{1\leq j\leq m}\in\mathcal{F}_{1}\times\mathcal{F}_{2}\cdots\mathcal{F}_{m}$
    		there is some $\ssq{f}{n}$ in $\mathcal{F}$ such that $(f_{n}\rst_{G_{i}})_{n}$ is a subsequence of $\ssq{f^{i}}{n}$ for $1\leq i\leq m$.
    		\item $\ssq{f}{n}$ makes $G_{i}$ grow faster than $G_{j}$ for $i<j$
    	}
    \end{observation}
    \begin{proof}
    	Simply chose a finite set of generators $\mathcal{X}_{i}$ of $G_{i}$ for $1\leq i\leq m-1$ and take as
    	$\mathcal{F}$ the collection of all sequences $\ssq{g}{n}$ of maps from $G_{1}\frp G_{2}\cdots G_{m}$ such that $(g_{n}\rst_{G_{j}})_{n}\in\mathcal{F}_{i}$ and that
    	$\tl{xy}{\lambda^{f^{j}_{n}}}\leq n\cdot\tl{g_{j+1}}{\lambda^{f^{j+1}_{n}}}$ for all $n>0$, $1\leq j\leq m-1$, $x,y\in\mathcal{X}_{j}$ and some $g_{j+1}\in G_{j+1}$.
    	%	with respect to  $\mathcal{F}_{j}$.
    \end{proof}
    
    \subsubsection{Free groups, free products and small cancellation sequences}
      Fix a base $\mathcal{B}$ of a free group $\F$. We say that $x\in\F$ is cyclically reduced if the translation length of its normal form is minimal among those of all of its conjugates.
      Alternatively, if the the length of its normal form coincides with the translation length of its action on the Cayley graph of $\F$, i.e. if the orbit of $1\in Cayl(\F)$ by powers of $x$ is contained in the axis of $x$. Any element $y\in\F$ can be expressed in a unique way as a product without cancellation of the form $u\cdot x\cdot u^{-1}$ so that $x=y^{u}$ is cyclically reduced (see \cite{lyndon2015combinatorial} for more details. Given elements $x,y$ (possibly equal) of a free group $\F$, a common piece of $x$ and $y$ is a word in $\mathcal{B}^{\pm}$ appearing as a subword of the normal form of cyclically reduced conjugates $x^{g}$ and $y^{h}$ of $x$ and $y$. If $x=y$, we additionally require that $gh^{-1}\nin Z(x)$.
      A tuple $x=(x_{j})_{j=1}^{n}$ of cyclically reduced elements of a free group $\F$ is said to have the small cancellation property $C'(k)$ for $k\in \N$ if
      any common piece of $x_{i}$ and $x_{j}^{\pm}$ for $1\leq i,j\leq n$ has length smaller than $\min\{|x_{i}|,|x_{j}|\}$.
      %     If some of the $x_{j}$ are not cyclically reduced, we say that
      %    $x$ has property $C'(k)$ if and only if the tuple $(y_{i})_{i=1}^{n}$, where $y_{n}$ is any cyclically reduced conjugate of $x_{n}$.
      %    $\frac{1}{n}\min_{\leq j\leq n}l_{j}$, where $l_{j}$ is the length of any cyclically reduced conjugate of $x_{j}$.
      %  Check the name of the small cancellation property.
      \begin{lemma}\label{small cancellation embeddings}
      	Let $(F^{1},\pi^{1})$ and $(F^{2},\pi^{2})$ be free \pr-groups, where $\pi^{i}$ and chose a base $(x_{1},\cdots x_{m})$ for $F^{1}$.
      	Then for each $m\in\N$ there is a morphism $\phi^{m}$ from $(F^{1},\pi^{1})$ to $(F^{2},\pi^{2})$ such that if we let $y^{m}_{j}=\bar{\phi^{m}(x_{j})}$ and $z_{i,m}$ be a cyclically reduced conjugate of 	$y_{i,m}$ then:
      	\elenco{
      		\item The tuple $(z^{m}_{1},z^{m}_{2},\cdots z^{m}_{n})$ has property $C'(m)$
      		\item $|y_{i,m}|>m$ for all $1\leq i\leq n$
      		\item $\frac{|y_{i,m}|}{|y_{j,m}|}<1+\frac{1}{m}$ for all $i,j\in\{1,2\cdots m\}$
      		\item $|z_{i,m}|\geq\frac{(m-1)}{m}|y_{i,m}|$.
      	}
      	\begin{proof}
      		Let $\mathcal{R}\subset\F$ map bijectively onto $Q$ via $\pi$ and take $L=\max_{r\in\mathcal{R}}|r|$. Let $n$ be the rank of $n$.
      		%  	a finite set of representatives of the preimages of points in $Q$.
      		It is a well known fact that for any $m$ a tuple $\bar{y}'_{m}=(y'_{1,m}\cdots y'_{n,m})$ cyclically reduced elements of $F^{2}$ of length $\geq m$
      		satisfying the three conditions above exists. The only thing left is to ensure the compatibility of those tuples with the \pr-structure on both sides.
      		
      		In order to do so, for each $1\leq i\leq n$, let $r_{i,m}\in\mathcal{R}$ be such that $\pi^{1}(x_{i})=\pi^{2}(y_{i,m})$, where $y_{i,m}=y'_{i,m}r_{i,m}$.
      		Let $u$ be such that $y'_{i,m}=w_{i,m}\cdot u$, $r_{i,m}=u^{-1}s_{i,m}$ and $y_{i,m}=w_{i,m}\cdot s_{i,m}$. Here the notation $u\cdot v$ indicates the fact that
      		the product is reduced, i.e., no cancellation between a sub-word of the normal form for $u$ and that for $v$ takes place.
      		As in the statement, let $z_{i,m}$ be a cyclically reduced conjugate of $y_{i,m}$.
      		Since $\bar{y}_{m}$ satisfies the small cancellation property, if $y_{i,m}$ starts with the inverse of some final subword $w$ of $z_{i,m}$, then $|w|\leq|s_{i,m}|+\frac{1}{m}|y_{i,m}|$. In short, $|z_{i,j}|\geq \frac{m-2}{m}|y_{i,j}|-L$. Furthermore, any piece of length $C$ contained in some $z_{i,m}$ has to contain another one
      		$y_{i,m}$ of length $\geq C-L$. Hence if we let $P(\bar{u})$ stand for the length of the biggest piece of $\bar{z}_{m}$ not intersecting any of the $u_{i,m}$ then:
      		%  	the quotient between the length of the smallest piece of $\bar{z,m}$ and that of the shortest element of $\bar{z,m}$ is smaller than
      		\begin{align*}
      			\frac{P(\bar{z,m})}{\min_{1\leq i\leq n}|z_{i,m}|}\leq \frac{P(\bar{y,m})+L}{\frac{m-1}{m}\min_{1\leq i\leq n}|z_{i,m}|-L}\leq
      			\frac{\min_{1\leq i\leq n}|z_{i,m}|}{\min_{1\leq i\leq n}|z_{i,m}|-L}(\frac{1}{m}+\frac{L}{m})
      		\end{align*}
      		Clearly this term goes to $0$ as $m$ goes to infinity, so we are done.
      		%  	We can take $\phi^{l}$ sending the base of $F^{1}$ to $z_{i}^{m_{l}}$ for sufficiently big $m_{l}$.
      	\end{proof}
      \end{lemma}
      
      \begin{observation}
      	\label{injectivity Nielsen}If a tuple $(x_{1},x_{2}\cdots x_{n})$ of cyclically reduced elements of a free group satisfies $C'(3)$, it is automatically a base of the subgroup of the free group it generates, as in particular it satisfies Nielsen property (see \cite{lyndon2015combinatorial} for a definition).
      \end{observation}
      %  (alternatively, the image of the tuple of generators)
      We will call a sequence of morphisms $f_{m}$ satisfying the first two properties listed in the statement a small cancellation sequence with respect to the base $(x_{1},\cdots x_{m})$. If it also satisfies the third we will say it is balanced. Note that given a sequence of balanced small cancellation sequences $(\ssq{f^{m}}{n})_{m\in\N}$ with respect to a certain base $x$ of $F_{1}$, the sequence
      $(f_{n}^{n})_{n\in\N}$ is also a balanced small cancellation sequence.
      \begin{lemma}
      	\label{small cancellation and limits} Let $H$ be a free group, $x$ a base of $H$ and $\ssq{f}{n}$ a balanced small cancellation sequence of morphisms from $H$ to some other free group $F$ with respect to $x$. Suppose we are given a finitely generated extension $L$ of $H$ and a sequence $\ssq{h}{n}$ of maps extending $f_{n}$ with trivial limit kernel and $Y$ a limiting tree for $\ssq{h}{n}$. Suppose that $H$ is not elliptic in $Y$ and let $Y_{H}$ be the minimal tree of $H$. Then $h\cdot Y_{H}\cap Y_{H}$ can only be a non-degenerate segment in case $h\in H$ and the action of $H$ on $Y_{H}$ is a geometric realization of the action of $H$ on $Cayl_{x}(S)$ which assigns the same length to all the edges of $Cayl(S,x)$.
      \end{lemma}
      \begin{proof}
      	It all comes down to the following easy observation:
      	\begin{observation}
      		\label{characterization of small cancellation}
      		Let $(x_{1},\cdots x_{n})$ be a tuple of cyclically reduced elements of the free group $\F$. Let $l_{i}$ be the translation length of $x_{i}$ for its action on the Cayley graph of $\F$.
      		The tuple has property $C'(n)$ if for any $h\in\F$ whenever $Ax(x_{i})$ and $h\cdot Ax(x_{j})$ overlap on a segment of length greater than $\frac{1}{n}\mn{1\leq i\leq m}\tl{x_{i}}{Cayl(\F)}$ necessarily $i=j$ and $h\in Z(x_{i})$.
      	\end{observation}
      	%    	We will start with statement \ref{first part}.
      	Assume that $H$ does not fix a point in $Y$. Let $\nu_{n}=\mn{1\leq i\leq |x|}f_{n}(x_{i})$. Small cancellation implies that $1$ cannot be at distance greater than $\frac{1}{n}\nu_{n}$ from
      	the axis of any $f(x_{i})$ for $1\leq i\leq m$. On the other hand, for any $h\in H$, we have $\frac{n-4}{n}|h|\nu_{n}\leq|f_{n}(h)|\leq\frac{m+1}{m}|h|\nu_{n}$, the left inequality following from
      	that given a reduced word $x_{i_{1}}^{\epsilon_{1}}x_{i_{2}}^{\epsilon_{2}}\cdots x_{i_{r}}^{\epsilon_{r}}$ at most $\frac{4}{n}$ of the reduced word for $f_{n}(x_{i_{j}}^{\epsilon_{j}})$ cancels out in the word
      	$f_{n}(x_{i_{1}}^{\epsilon_{1}})f_{n}(x_{i_{2}}^{\epsilon_{2}})\cdots f_{n}(x_{i_{r}}^{\epsilon_{r}})$. This implies that the functions $\frac{1}{\nu_{n}}\pl{1}{f_{n}(-)}{Cayl(\F)}$ converge point-wise to $\pl{1}{-}{Cayl_{x}(H)}$ and the action of $H$ on $Y_{H}$ is exactly as desired.
      	Suppose now that for some $g\in L$ the intersection $Y_{H}\cap g\cdot Y_{H}$ is a non-degenerate segment.
      	
      	Since $Y_{H}$ is covered by translates of the axis of the $x_{i}$, this implies $Ax(x_{i})$ and	$h\cdot Ax(x_{j})$ overlap in a non-degenerate segment for $1\leq i,j\leq m$. This in turn implies that for some $\alpha>0$
      	and any $n$ in some infinite subsequence $|Ax(f_{n}(x_{i}))\cap f_{n}(h)\cdot Ax(f_{n}(x_{j}))|>\alpha\tl{f_{n}(x_{i})}{}$. The small cancellation property together with the observation above implies that $i=j$ and $f_{n}(h)=f_{n}(x_{i})^{r_{n}}$ for some $r_{n}\in\Z$ and any such $n$. Now, $r_{n}=\frac{\tl{f_{n}(h)}{}}{\tl{f_{n}(h)}{}}$ converges to $r_{\infty}=\frac{\tl{x_{j}}{Y}}{\tl{x_{i}}{Y}}$, so it must be constantly equal to $r_{\infty}\in\Z$ after a certain point, which implies that $h\in\subg{x_{i}}$.
      \end{proof}
      \begin{comment}
      \begin{corollary}
      	\label{discriminating retractions} In the situation above, $L= \frp H$ and $h_{n}\rst_{\F}=Id_{F}$ then the limit kernel of $h_{n}$ is trivial.
      \end{corollary}
      \end{comment}
      
      \begin{corollary}
      	\label{discriminating retractions} Suppose we are given a group $L= K\frp H$ and a $\ssq{g}{n}$ a sequence of homomorphisms from $L$ to $\F$
      	such that $g_{n}\rst_{H}=f_{k_{n}}$ for some balanced small cancellation sequence $f_{n}$ and $\frac{\pl{1}{g_{n}(y_{i}y_{j})}{}}{k_{j}}$ tends to $0$ for
      	some finite tuple $y$ of generators of $L$ and $1\leq i,j\leq |y|$. Then $\limker{g_{n}}{n}=1$.
      	and the	limiting tree is a geometric realization of the tree dual to a graph of groups decomposition with a single vertex group $K$ and a one-edge trivially stabilized loop
      	for each of the given generators of $H$.
      \end{corollary}
      \begin{proof}
      	Let $\fun{q}{L}{Q}$ be the limit quotient of the sequence $\ssq{f}{n}$.
      	%  	we can think of $K$ and $H$ as subgroups of
      	Consider the action of $Q=L/\limker{h}{n}$ on any limiting tree $Y$ for $\ssq{h}{n}$. The last property in the definition of small cancellation sequences implies that
      	$q(K)$ fixes a point in the minimal tree $Y_{H}$\footnote{The ratio between the distance between the axis of the small cancellation elements in $\F$ from the identity and the rescaling factor tends to $0$.} . We know that $Y_{H}\cap g\cdot Y_{H}$ can only be non-degenerate for $g\in q(H)$ which then, if non-trivial, cannot fix any point of $Y$. It follows that the limit kernel of $\ssq{g}{n}$ is trivial and $Y$ of the desired form.
      	\begin{comment}
      	We know also 	$\{l\cdot Y_{H}\}_{l\in L}$ is a transverse family. Let $S$ be the skeleton of the associated decomposition of $Y$. Of course, $K$ fixes a point of $Y$, it fixes a point $v_{0}$ of $S$ as well, necessarily distinct the one, say $v_{1}$, fixed by $H$. Now, any $e$ incident to $v_{0}$ is trivially stabilized, since the set-wise stabilizer of $Y_{H}$ is $H$ itself. This automatically implies that $Q\cong q(K)\frp q(H)\cong L$.
      	\end{comment}
      \end{proof}
      \begin{comment}
      \begin{corollary}
      	\label{discriminating retractions} Suppose we are given a group $L= K\frp H$ and a $\ssq{g}{n}$ a sequence of homomorphisms from $L$ to $\F$
      	such that $g_{n}\rst_{H}=f_{k_{n}}$ for some balanced small cancellation sequence $f_{n}$ and $\frac{\tl{g_{n}(y_{i}y_{j})}{}}{k_{j}}$ tends to $0$ for
      	some finite tuple $y$ of generators of $L$ and $1\leq i,j\leq |y|$. Then $\limker{g_{n}}{n}=1$.
      	%    		limiting tree is a geometric realization of the tree dual to a graph of groups decomposition with a single vertex group $K$ and a one-edge trivially stabilized loop
      	%    	for each of the given generators of $H$.
      \end{corollary}
      \begin{proof}
      	Let $\fun{q}{L}{Q}$ be the limit quotient of the sequence $\ssq{f}{n}$.
      	%  	we can think of $K$ and $H$ as subgroups of
      	Consider the action of $Q=L/\limker{h}{n}$ on any limiting tree $Y$ for $\ssq{h}{n}$. We know	$\{l\cdot Y_{H}\}_{l\in L}$ is a transverse family. Let $S$ be the skeleton of the associated decomposition of $Y$. Of course, $K$ fixes a point of $Y$, it fixes a point $v_{0}$ of $S$ as well, necessarily distinct the one, say $v_{1}$, fixed by $H$. Now, any $e$ incident to $v_{0}$ is trivially stabilized, since the set-wise stabilizer of $Y_{H}$ is $H$ itself. This automatically implies that $Q\cong q(K)\frp q(H)\cong L$.
      \end{proof}
      \end{comment}
      
      %    The previous discussion provides the base case of our proof (almost) and the family
      \begin{corollary}
      	Any \pr-tower without constants and a single node of type \ref{free group case} admits a CSD family of \pr-test sequences.
      \end{corollary}
      \begin{corollary}
      	Suppose that the root $r$ is as in case \ref{free product case} and  that for each $\mu\in Ch(r)$ a congruency complete CSD family of \pr-test sequences
      	of $\RTA\rst_{\mu}$ exists. Then  a familiy $\mathcal{TS}$ with such properties exists for $\RTA$ as well.
      \end{corollary}
      \begin{proof}
      	Let $Ch(r)=\{\mu_{1},\mu_{2}\}$ and assume that $\mu_{1}$ carries the constants. Let $\mathcal{TS}_{0}$ be the family of sequences of morphisms from
      	$T$ to $\F$ obtained by applying \ref{asymmetrical growth} to the families $\mathcal{TS}^{\mu_{1}}$ and $\mathcal{TS}^{\mu_{2}}$.
      	%    	implies the existence of a CSD family
      	Let $\qG{Z}$ be cyclic, with $\pi(Z)=\tg$ and $\ssq{g}{n}$ a \pr-test sequence. Take $L=T\frp Z$,
      	%    	 this can be also decomposed as $K\frp Z$, where $K=\RLR{T}{\mu_{1}}\frp \RLR{T}{\mu_{2}}^{t}$.
      	Let $\iota$ be the embedding of $T$ into $L$ given by $\iota\rst_{\RLR{T}{\mu_{1}}{}}=Id$ and $\RLR{T}{\mu_{2}}=\inn{t}$.
      	Fix a tuple of generators $x$ of $T$. As $\mathcal{TS}$ we take the set of all sequences of the form $(g_{n}\circ\iota)\subset Mor(\rqG{L}{A},\rqG{F}{A})$,
      	where:
      	\elenco{
      		\item $(g_{n}\rst_{T})_{n}\in\mathcal{TS}_{0}$
      		\item $g_{n}\rst_{Z}=h_{k_{n}}$ for some \pr-test sequence $\ssq{h}{n}$ of $\qG{Z}$ and $k_{n}\in\N$ is such that
      		$\tl{g_{n}(x_{i}x_{j})}{}\leq k_{n}$ for any $1\leq i,j\leq |x|$.
      	}
      	It is easy to check that $\mathcal{TS}$ is CSD sequence of morphisms factoring through $\RTA$. The action of $T$ on any limiting tree for $\ssq{f}{n}\in\mathcal{TS}$ is simply
      	the minimal tree of the pullback by $\iota$ of the action of $L$ on the limiting tree for $\ssq{g}{n}$. The result we are looking for is immediate, given what the we know about the latter from the previous lemma.
      	It is also true, although not needed for this, that $\limker{g_{n}}{n}$.
      \end{proof}
      \newcommand{\De}[1]{\Delta_{#1}}
      \newcommand{\pM}[0]{\mathcal{A}}
    \subsubsection{Approximating actions by twisting}
      \begin{lemma}
      	\label{Dehn twist approximation}
      	Suppose we are given a group $G$, $\mathcal{F}$ a finite subset of $G$,
      	$\BS{Z_{0}}{Z}{t}$ a presentation of an action of $G$ without inversions on a simplicial tree $S$ with abelian edge stabilizers. Then there is a finite subset $\mathcal{H}\subset\bunion{v\in VS}{Stab(v)}{}$ such that given any $\fun{f}{G}{H}$ and any action $\lambda$ of $H$ by isometries on a real tree $(X,d)$, if the following properties are satisfied:
      	\enum{i)}{
      		\item $\lambda$ is $D$ acylindrical for some $D>0$
      		\item  $\mathcal{H}\cap ker\,f=\nil$
      		\item  $f(Stab(e))$ acts hyperbolically on $X$ for all $e\in S$
      		\item \label{sita} Given $e,f$ incident at $v\in S$, if $Stab(e)\cap Stab(f)^{g}\neq\tg$ then either
      		\footnote{ Equivalently, at the level of graph of groups:
      		\elenco{
      			\item No two conjugates of the groups $\Delta_{[e]}$ and $\Delta_{[f]}$ intersect non-trivially for distinct and non mutually inverse $[e],[f]$ originating at $[v]$.
      			\item The two edge groups associated to a one-edge loop coincide and the corresponding Bass-Serre generator commutes with them.
      			\item $\Delta_{[e]}$ is malnormal in $\Delta_{[v]}$ for any edge $[e]\in G\backslash S$
      		}
      		If $G$ is $CSA$ the last bullet is equivalent to $\Delta_{[e]}$ being a maximal abelian subgroup of $\Delta_{[v]}$. }:
      		\elenco{
      			\item $e=f$
      			\item $e=t\edin{f}$ for some Bass-Serre element $t$ commuting with $Stab(e)$
      		}
      	}
      	\begin{comment}
      	or any $v_{0}\in S$
      	if each $e\in E$ is assigned some $c_{e}=c_{\edin{e}}^{-1}$ in such a way that
      	%    	$c_{e}^{t_{e}}=c_{\edin{e}}^{-1}$
      	\elenco{
      		\item $c_{t_{e}\cdot e}=c_{e}^{t_{e}^{-1}}$
      		\item $c_{e}$ acts hyperbolically on $Y$ for each $e\in E$   \label{hypothesis1}
      		%    		\item For any $v\in V$, the action of $Stab(v)$ on $Y$ is faithful. \label{hypothesis2}
      		\item If $\alpha(e)=\alpha(f)=v$, and $g\in\mathcal{G}_{v}$ the elements $c_{e}$ and $c_{f}^{g}$ commute only if $e=f$ and $g=1$.   \label{hypothesis3}
      	}
      	\end{comment}
      	Then for each $x\in X$, $v_{0}\in S$ and any geometric realization $(|S|,d')$ of $S$ there is a sequence $\ssq{\tau}{n}\subset\pM$ of products of Dehn twists over the edges of $e$
      	%    	 any equivariant assignment of a length $\alpha_{e}\in(0,1)$ to each of the edges of $S$
      	such that the sequence $\frac{1}{n}\pl{x}{g}{\lambda\circ f\circ\tau_{n}}$ converges to $\pl{v_{0}}{g}{(\rho,d')}$ for any $g\in\mathcal{F}$,
      	%    	in which $e$ has length $\alpha_{e}$
      	where $\rho$ is the obvious action of $G$ on $|S|$.
      \end{lemma}
      In particular, if we take $H$ to be free and as $\lambda$ the action of $H$ on its Cayley graph we deduce:
      \begin{corollary}
      	\label{discrimination corollary} Suppose we are given a finitely generated \pr-group $\qG{G}$ a simplicial $G$-tree $S$ with abelian edge stabilizers
      	satisfying condition (\ref{sita}) or of the lemma above and a sequence $\ssq{f}{n}$ of morphisms from $G$ to $\F$ such that for any $v\in S$ the sequence
      	$(f_{n}\rst_{Stab(v)})_{n}$ has trivial limit kernel. Then there is a sequence $\ssq{\tau}{n}\subset Mod^{\pi}(\qG{G},S)$
      	%    	of products of Dehn twists over edges of $S$
      	such that $(f_{n}\circ\tau_{n})_{n}$ has trivial limit kernel.
      \end{corollary}
      \begin{comment}
      \begin{proof}
      	Indeed, given any $n\in\N$, let $\mathcal{F}_{n}$ be the ball of radius $n$ in some fixed Cayley graph of $G$ minus the identity and $\mathcal{H}_{n}\subset G\setminus\{1\}$ the finite set associated to it by the lemma, applied to the action $\lambda$ of $\F$ on a geometric realization of its Cayley graph. There is some $k_{n}$ such that for $m\geq k_{n}$ the map $f_{m}$ does not kill those elements of $\mathcal{F}_{n}$ elliptic in	$S$ nor any $h\in\mathcal{H}_{n}$.
      	The lemma implies, in particular, that for any such $m$ some product $\tau_{m}\in Mod^{\pi}(\qG{G},S)$ of Dehn twists over edges of $S$ exists such that any element
      	of $\mathcal{F}$ hyperbolic with respect to $S$ is hyperbolic for the action $\lambda\circ f_{m}\circ\tau$, so that its image by $f_{m}\circ\tau_{m}$ is non-trivial.
      \end{proof}
      \end{comment}
      \begin{example}
      	The corollary is not true if one simply removes condition (\ref{sita}). Consider the amalgamated product of a free group $\F(a,b,c)$ with the free abelian group with base $\{c,d\}$ over their common subgroup $\subg{c}$. Let $F'=\F(x,y)$ and let $\ssq{f}{n}$ be the sequence of maps given by:
      	\begin{align*}
      		f_{n}(a)=y, \,\,\, f_{n}(b)=y^{x}, \,\,\, f_{n}(c)=x^{n+3},\,\,\, f_{n}(d)=x^{n+2}
      	\end{align*}
      	It is easy to check that the restriction of $f_{n}$ to both vertex groups has trivial limit kernel, as the triple $(f_{n}(a),f_{n}(b),f_{n}(c))$ is Nielsen reduced. Now, any $g_{n}$ obtained form $f_{n}$ by precomposing by a Dehn twist of the decomposition above	is equal to $f_{n}$ up to post-composing by an automorphism of $\F$ (some power of the one fixing $x$ and sending $y$ to $y^{x}$), so in particular $ker\,f_{n}=ker\,g_{n}$. But $\limker{f_{n}}{n}\neq\tg$, as the element $b^{-1}a^{cd^{-1}}\neq 1$ is sent by the identity by all the $f_{n}$.
      \end{example}
      The hypothesis can be slightly weakened, though:
      \begin{corollary}
      	\label{second discrimination corollary} Suppose we are given a finitely generated \pr-group $\qG{G}$ and $S$ a a simplicial $G$-tree with non-trivial abelian stabilzers such that $VS=V_{0}\dcup V_{1}$ and
      	\elenco{
      		\item Each vertex of $V_{0}$ is adjacent only to vertices of $V_{1}$
      		\item There are no non-trivial elements fixing two distinct edges adjacent to a vertex of $V_{1}$.
      		\item The stabilizer of any edge of $S$ coincides with its centralizer in the stabilizer of any of its endpoints.
      	}
      	Assume we are also given a sequence $\ssq{f}{n}$ of homomorphisms from $G$ to $\F$ such that its restriction to the stabilizer of any vertex of $S$ has trivial limit kernel.
      	Then there is a sequence $\ssq{\tau}{n}\subset Mod^{\pi}(\qG{G},S)$ such that $\limker{f_{n}\circ\tau_{n}}{n}$ is trivial.
      \end{corollary}
      \begin{proof}
      	First of all notice that the conditions listed above imply that $S$ is $2$-acylindrical and $Stab(e)$ is self-centralizerd in $G$ for any edge $e$ in $S$ \footnote{That elliptic elements commuting with $Stab(e)$ are in $Stab(e)$ is easy to check. The axis of an hyperbolic element commuting with $Stab(e)$ has to be fixed by any non-trivial element of $Stab(e)$, against acylindricity. } .
      	Our proof is by induction on the number of edges of the associated decomposition $\Delta$. Pick any edge $e\in S$ and let $S^{*}$ be the tree resulting from the collapse of any edge outside the orbit of $e$. The induction hypothesis implies the existence of a sequence $\ssq{\sigma}{n}\subset\mathcal{A}$ such that $\limker{f_{n}\circ\sigma_{n}\rst_{Stab(u)}}{n}=\tg$
      	for any vertex $u\in S^{*}$.
      	
      	The desired result will follow as soon as we manage to prove that $S^{*}$ satisfies condition (\ref{sita}). Since $Stab(e)$ is self-centralized and
      	$Stab(u)$ is $CSA$, the only thing left to show is that $Stab(e)\cap Stab(f)=\tg$ for $e,f$ originating at $u$ in distinct orbits, this is clear, since the lifts $\tilde{e}$, $\tilde{f}$
      	in $S$ of $e$ and $f$ respectively must span a segment of diameter at least $3$.
      	%    	There are at most two orbits of edges originating at $u$. We claim that the st
      \end{proof}
      %    Recall that the \egt is metrizable. Let $d^{GT}$ be some metric on it.
      \begin{proof}
      	\newcommand{\expe}[1]{exp_{#1,n}}
      	(of lemma \ref{Dehn twist approximation}) Let $f$ be any homomorphism from $G$ to $H$. Chose a \BSP $\BS{Z_{0}}{Z}{t}$ of the action of $G$ on $S$	for which $v_{0}\in Z_{0}$. For any edge $e\in Z$ mapping onto a loop in $G\backslash S$, chose a Bass-Serre element $t_{e}$ commuting with
      	$Stab(e)$, in accordance to (\ref{sita}). Assign to every edge $e\in S$ some $c_{e}\in Stab(e)$ in such a way that:
      	%    	$c_{e}^{t_{e}}=c_{\edin{e}}^{-1}$
      	\elenco{
      		\item $c_{\edin{e}}^{-1}=c_{e}$
      		\item $c_{g^{-1}\cdot e}=c_{e}^{g}$
      		\item $c_{e}$ is hyperbolical with respect to $\lambda^{f}$ for each $e\in E$; let $tl_{e}$ be the translation length of $f(c_{e})$ with respect to $\lambda$
      	}
      	For each $e\in S$ let $\alpha_{e}$ be the length of the segment corresponding to $e$ in $(|S|,d')$ and $tl_{e}$ the translation length of $f(c_{e})$ with respect to $\lambda$. For each $n\in\N$ chose an integer $exp_{e}(n)$ closest to $\flr{\frac{n\alpha}{tl_{e}}}$, with the property that the Dehn twist over $e$ by $c_{e}^{exp_{e}(n)}\in\pM$. Clearly there is a constant $C$ independent of $e$ and $n$
      	such that $exp_{e}(n)$ and $\frac{n\alpha}{tl_{e}}$ differer at mostby $C$ regardless of $n$.
      	As $\tau_{n}$ we take the composition of the Dehn twists $\tau_{e,c_{e}}^{exp_{e}(n)}$, for $e$ ranging among all edges of $Z$ pointing away from $v_{0}$.
      	Represent any $g\in\mathcal{G}$ as a product of the form $\looprep{g}{e}{m}$, associated to a sequence $v_{0},e_{0},v_{1},e_{1}\cdots e_{m-1},v_{m}=v_{0}\subset Z$
      	%  , where $v_{i}\in V$ and $e_{i}\in E$,
      	projecting to a closed path in  $G\backslash T$, where $\alpha(e_{i})\in V$ for each $i$.
      	%    	Let $\mathcal{G}_{v}$ be the collections of all the $g_{i}$ appearing in such an expression with $v_{i}=v$.
      	Let $R$ be some constant bigger than:
      	\enum{i)}{
      		\item Any of the $tl_{e}$.
      		\item $\pl{x}{g_{i}}{\lambda}$ for $0\leq i\leq m$
      		\item The distance between $x$ and $Ax(c_{e})$ for any $e\in S$.
      		\item   $(C+1)(\serprop{e\in E}{\tl{c_{e}}{\lambda}})$
      		%    		\item $K$ŋcI am going to completely change the construction. It proceeds by proving the result on discrimination first. I need to modify the \twist lemma so as to not talk about actions
      	}
      	Then $\tau_{n}(g)$ can be written as
      	$g_{0}s_{0}t_{e_{0}}g_{1}s_{1}t_{e_{1}}g_{2}\cdots s_{m-1}g_{m}$
      	where $s_{i}=c_{e_{i}}^{-\expe{e_{i}}}$. We want to estimate the distance between the points $x$ and $y=\lambda(\tau(g))\cdot x$.
      	Let:
      	$I_{0}=[x,g_{0}\cdot x]$,
      	$J_{j}=\bar{g}_{j}\cdot[x,s_{j}\cdot x]$ for $0\leq j\leq m$,
      	$I_{j}=\bar{h}_{j}\cdot[x,t_{e_{j-1}}g_{j}\cdot x]$ for $1\leq j\leq m$, where $\bar{g}_{j}=g_{0}s_{0}t_{e_{0}}\cdots g_{j}$
      	and $\bar{h}_{j}=g_{0}s_{0}t_{e_{0}}g_{1}\cdots s_{j-1}$ (the reader is reminded that although we have decided to leave the subscripts out, all of this is dependent of $n$).	Notice that for $0\leq j\leq m-1$ the right endpoint of the oriented segment $I_{j}$ is the left endpoint of $J_{j}$, while the right endpoint of $J_{j}$ is the left endpoint of $I_{j+1}$ and the right endpoint of $I_{m}$ is $y$.
      	%    	Add  n  to the notation all over the place.
      	Now, clearly $|I_{j}|,|K|\leq 2R$ for $0\leq j\leq m$. On the other hand $|J_{j}-n|\leq R$ and the subsegment $J'_{j}\subset J_{j}$ consisting of those points at distance $\geq R$ from $\partial J_{j}$ is contained in $\bar{h}_{j}Ax(s_{j})=\bar{h}_{j}Ax(c_{e_{j}})$, since $x$ is at distance $\leq R$ from $Ax(c_{e_{j}})$.
      	For each $e\in S$, let $l_{e}(g)$ be the number translates of $e$ or $\edin{e}$ contained in the path $[v_{0},g\cdot v_{0}]$.
      	We deduce that:
      	\begin{align*}
      		\frac{1}{n}\pl{x}{\tau_{n}(g)}{\lambda}\leq
      		\frac{1}{n}(\serprop{e\in E}{\alpha_{e}|\{0\leq j\leq m\,|\,[e_{j}]=e\}|}+(2m+1)R)
      		=l_{\alpha}(g)+\frac{2m+1}{n}
      	\end{align*}
      	\newcommand{\djj}[0]{c_{\edin{e_{j-1}}}}
      	\newcommand{\cjj}[0]{c_{e_{j-1}}}
      	\newcommand{\cj}[0]{c_{e_{j}}}
      	\renewcommand{\dj}[0]{c_{\edin{e_{j}}}}
      	\newcommand{\tjj}[0]{t_{e_{j-1}}}
      	\newcommand{\gj}[0]{g_{j}}
      	\newcommand{\sj}[0]{s_{j}}
      	\newcommand{\ctj}[0]{c_{t_{e_{j-1}}}^{-1}}
      	% check this
      	\begin{comment}
      	Suppose for a second that the subsegment $J''_{j}$ of $J_{j}$ consisting of those points at distance $\geq 3R+D+1$ from $\partial J_{j}$ intersects some $I_{k}$ or $J_{l}$ for $j\neq l$. Since $|I_{j}|,|I_{j+1}|\leq R$, for $n$ big enough $J_{j}$ shares a subsegment of length $\geq 2R+D+1$ with either $J_{j-1}$ or $J_{j+1}$, let us say with the first of them. Then $J'_{j-1}$ and $J'_{j}$ share a subsegment of length $\geq R+D+1$ and hence the same holds for $(t_{e_{j-}}g_{j}s_{j})^{-1}\cdot Ax(c_{e_j-1})=Ax(c_{e_{j-1}}^{t_{e_{j-1}}g_{j}s_{j}})$ and $Ax(c_{e_{j}})$.
      	\end{comment}
      	Assume for a moment that the subsegment $J_{j}''$ of $J_{j}$ consisting of those points at distance greater or equal than $ 3R+D+1$ from $\partial J_{j}$ intersects some $I_{k}$ or $J_{l}$ for $j\neq l$. Since $|I_{j}|,|I_{j+1}|\leq R$, for $n$ big enough some subsegment of $J_{j}$ of length greater or equal than $2R+D+1$ is contained in $J_{j-1}$ or $J_{j+1}$ with the opposite orientation, let us say $J_{j-1}$.
      	This implies $J'_{j-1}$ and $J'_{j}$ share a subsegment of length $\geq R+D+1$ and hence the same is true for $(\tjj\gj\sj)^{-1}\cdot Ax(\cjj)=Ax(\cjj^{\tjj\gj\sj})$ and $Ax(\cj)$ (in the sense of $\lambda^{f}$).
      	%  	Now, $c_{e_{j-1}}^{t_{e_{j-1}}g_{j}s_{j}}=(c_{\edin{e_{j-1}}}^{g_{j}s_{j}})^{-1}\in\Delta_{v_{j}}$.
      	Now, $\cjj^{\tjj\gj\sj}=(\ctj)^{\gj\sj})\in Stab(v_{j})$.
      	The following fact is well known and can be easily proven by inspection:
      	\begin{fact}
      		Given an action of a group $G$ on a real tree, $L>0$ and two hyperbolical elements $h,g\in G$ such that $|Ax(g)\cap Ax(h)|\geq\tl{g}{}+\tl{h}{}+L$, the commutator $[h,g]$ fixes a segment of length at least $L$ of 	$Ax(g)\cap Ax(h)$.
      	\end{fact}
      	It follows from this, together with the $D$-acylindricity of $Y$, that $f(([\ctj)^{\gj\sj},\cj])=1$, which implies that
      	$f([(\ctj)^{\gj},\cj])=1$. We claim that $[(\ctj)^{\gj},\cj]$ is non-trivial in $G$. Indeed, assumption \ref{sita} together with the fact that the expression $\looprep{g}{e}{m}$ is reduced imply
      	that the only way it can be trivial is if $[\gj,\cj]=1$ and $e_{j-1}=e_{j}$. This implies the existence of two inversely oriented subsegments of $Ax(c_{j})$ both of which are oriented in the same direction in which
      	$\cj$ acts by translations, which is obviously an absurdity.
      	%  	Observe that assumption \ref{sita} and the fact  implies the commutator $[c_{\edin{e_{j-1}}}^{g_{j}},c_{j}]$ is non-trivial in $G$.
      	We deduce from this that there is a finite collection of non-trivial elements of $G$ such that if $f$ does not kill any of them for each $n$ big enough
      	the segment $[x,y]$ will contain the disjoint union of the segments $J''_{j}$, which in turn implies $\frac{1}{n}\pl{x}{\tau_{n}(g)}{Y}\geq\nicefrac{nl_{\alpha}(g)-m(6R+3D+D))}{n}$.
      	Both this and the lower bound above tend to $l_{\alpha}$ with $n$, giving us the desired result.
      	%  and the theorem follows from the fact the \pegt and the length function topology coincide (see \ref{}).
      \end{proof}
      Let $\Sigma$ be a big compact surface with boundary.  We say that $\delta$ cuts the surface into small pieces if and only if
      any $\Sigma_{0}\in [\Sigma]_{\delta}$ is either a a sphere minus three discs (a.k.a. 'pair of pants') or a projective plane with two disks removed.
      \begin{remark}
      	\label{remark on small pieces} Observe that such a \cfam has the property that given any homotopy class $\alpha$ of essential s.c.c. either $\alpha\in\delta$ or else any of the cyclic subgroups associated with $\alpha$ are hyperbolic in $T_{\delta}$.
      \end{remark}
      We say that two \cfams $\delta$ and $\epsilon$ fill $\Sigma$ if for suitable representatives of $\delta,\epsilon$, each connected component of $\Sigma\setminus(\delta\cup\epsilon)$ is either simply connected or a boundary-parallel annulus. In particular, for any $\alpha\in\gamma$ there is some $\beta\in\delta$ such that $i(\alpha,\beta)\neq 0$ and viceversa.
      \begin{remark}
      	This is equivalent to the fact that any non-trivial element of $\pi_{1}(\Sigma,*)$ must be hyperbolic in either $T_{\gamma}$ or $T_{\gamma}$.
      \end{remark}
      %    for some basepoint $*\in\Sigma$
      Given a \cfam $\gamma$ on $\Sigma$ let $\Delta_{\gamma}$ be a graph of groups decomposition
      corresponding to the tree $T_{\gamma}$, as defined in \ref{surfaces section}, whose vertex groups are the image in $\pi_{1}(\Sigma,*)$ of the fundamental groups of the subsurfaces in
      $[\Sigma]_{\delta}$ (this image is well defined up to conjugation).
      %    Let also $\Sigma^{*}_{\gamma}\in [\Sigma]$ be the unique piece containing $*$.
      %    As we know $\pi_{1}(\Sigma^{*}_{\gamma},*)$ can be identified canonically with a subgroup $S_{\gamma}^{*}$ of $S$ fixing a unique
      %    vertex $u_{\gamma}\in VT_{\gamma}$.
      \begin{fact}
      	Given a \cfam $\gamma$ in a compact bounded surface $\Sigma$ there is another \cfam $\delta$ such that $\gamma$ and $\delta$ fill $\Sigma$.
      \end{fact}
      The curve graph of a big surface $\Sigma$, which we denote by $\mathcal{C}^{1}(\Sigma)$ is the graph whose vertices are homotopy classes of essential simple closed curves in $\Sigma$, where an edge is added between
      $\alpha$ and $\beta$ in case $i(\alpha,\beta)$ is minimal. Clearly a \cfam is a subset $\mathcal{C}^{1}(\Sigma)$ of diameter $\leq 1$. On the other hand, any two curves at distance $\geq 3$ in $\mathcal{C}^{1}(\Sigma)$ fill $\Sigma$, as there is no essential simple closed curve disjoint with both. The statement above follows therefore from the very well-known fact that $\mathcal{C}^{1}(\Sigma)$ has infinite diameter.
      A proof of the following can be extracted from \cite[4.17]{wilton2009solutions}
      \begin{lemma}
      	%    	Suppose we are given $G$ and  \label{special system of curves}
      	\label{small pieces injected}Let $S$ be the fundamental group of a big bounded compact surface $\Sigma$ .
      	Suppose we are given a homomorphism $f$ from $S$ to some limit group $G$ which injects the groups corresponding to boundary components of $S$ and has non-abelian image. Then there is a \cfam $\delta$ on $\Sigma$ cutting $\Sigma$ into small pieces such that for each $\alpha\in\delta$ such that the fundamental group of any $\Sigma_{0}\in [\Sigma]_{\delta}$ is mapped injectively by $f$ into $G$.
      \end{lemma}
      
    \subsubsection{Surface type vertex groups}
      
      %    In order not to clutter the (inductive) construction step for surface type floors, let us first proof the first partial result,
      As a warm up, let us start by proving the following weaker and well-understood result:
      \begin{lemma}
      	\label{surface discrimination} Suppose we are given a \pr-tower structure $\LRL{T}{A}{(J,r)}$ and $\lambda\in J$ of floor type, where $\GTla$ contains a unique vertex of non-rigid type,
      	which is of surface type, as well as a sequence of morphisms $\fun{\ssq{f}{n}}{\bfrp{\TAl}{\mu\prec\lambda}{}}{\rqG{F}{A}}$ with trivial limit kernel.
      	Let $\delta$ be a \cfam on the corresponding surface $\Sigma$ which cuts it into small pieces.
      	%    	Then there is a sequence $\ssq{\sigma}{n}\subset PMod(GRla;A_{\lambda})$ such that the sequence $(f_{n}\circ \eta^{R,A,\lambda}_{R,A,\lambda}\circ\sigma_{n})_{n}$ from $G$ to $\F$ has trivial limit kernel as well.
      \end{lemma}
      \begin{proof}
      	Let $\delta$ be the \cfam provided by lemma \ref{small pieces injected}, applied to the unique surface group of $\Delta=\GTla$, isomorphic to $\pi_{1}(\Sigma,*)$, and $\eta:=\eta^{\lambda}_{\QR{R}}$.
      	%    	or maybe not; then there
      	Let $S$ be the tree dual to $G:=\GTla$
      	and $S_{\delta}$ the one obtained by blowing up those vertices in the unique orbit $G\cdot v_{0}$ of $S$ to a copy of the action on
      	$T_{\delta}$. Let $S'_{\delta}$ be the tree obtained from $S_{\delta}$ by folding together $(e,g\cdot e)$ for any $e$ originating at a rigid type vertex $v$ and any $g$ in the centralizer of $Stab(e)$ in $Stab(v)$. The resulting decomposition is identical to that associated to $S_{\delta}$, except that at each of the vertices where we saw $H:=\pi_{1}(\Sigma_{0})$ for some $\Sigma_{0}\in [\Sigma]_{\delta}$
      	we now see the amalgamated product of $H$ with the centralizers of groups corresponding to each of the one or two boundary components of $\Sigma_{0}$. It is clear now that the tree $S'_{\delta}$ satisfies
      	the list of conditions of corollary \ref{second discrimination corollary} (with the set of rigid vertices as $V_{0}$).
      	
      	We need to check that the restriction of $\ssq{f}{n}$ to the stabilizer of a vertex $u$ of $S'_{\delta}$ not already in $S_{\delta}$ has trivial limit kernel. Observe that in that case $Stab(u)=A\frp B$, with $A$ and $B$ are the centralizers of incident edges.
      	\begin{claim}
      		Let $G=A\frp B$ for free abelian groups $A$ and $B$ and let $\ssq{f}{n}$ be a sequence of homomorphisms from $G$ to a free group
      		such that $(f_{n}\rst_{A})_{n}$ and $(f_{n}\rst_{B})_{n}$ both have trivial limit kernel and $f_{n}(G)$ is non-abelian for any $n$.
      		Then $\ssq{f}{n}$ has trivial limit kernel.
      	\end{claim}
      	\begin{proof}
      		This follows easily from the well-known fact that any pair of two non-commuting elements of the free groups form a basis of the subgroup of the free group they generate; in particular here
      		$f_{n}(G)\cong f_{n}(A)\frp f_{n}(B)$ for any $n$. For any $g\in G$ with normal form (let's say) $a_{1}b_{1}a_{2}b_{2}\cdots b_{n}$, where $a_{i},b_{i}\neq 1$,
      		for $n$ big enough both $f_{n}(a_{i})$ and $f_{n}(b_{i})$ are non-trivial for all $i$, which together with the fact above clearly implies $f_{n}(g)\neq 1$.
      	\end{proof}
      \end{proof}
      \renewcommand{\pM}[0]{R}  %tc
      \newcommand{\close}[2]{\sim_{#1}^{#2}}
      %    An elementary way of proving this is through the following lemma, together with proposition \ref{Rips decomposition}.
      \begin{lemma}	\label{generic test sequences}
      	Let $G$ be a group admitting an action on a \gat $S$,
      	such that $G\backslash S$ is bipartite, with edges joining a single vertex $[v_{0}]$ of surface type with (finitely many) vertices of rigid type.
      	Let $\Sigma$ be the surface associated to $v_{0}$, so that $Q_{0}:=Stab(v_{0})\cong\pi_{1}(\Sigma)$.
      	Suppose we are given some finite index subgroup $\pM\leq Mod(\Delta)\subset Aut(G)$.
      	and \cfam $\delta$ on $\Sigma$ cutting it into small pieces.
      	Let $W$ be the tree obtain by blowing up $v_{0}\in S$ to the $Stab(v_{0})$-tree $T_{\delta}$ dual to the \cfam $\delta$\footnote{In this particular case there is a unique way of performing such a blow-up. } and then collapsing all edges inherited from $S$ and let $(|W|,d,\rho)$ a geometric realization of it. Then there is a sequence $\ssq{\theta}{n}\in\pM$ and a diverging sequence $\ssq{\mu}{n}$ of positive constants such that
      	if we let $W_{0}$ be the minimal tree of $W$ the sequence $(|W_{0}|,\frac{1}{\mu_{n}}d,\rho\circ\theta_{n})$ converges
      	in the equivariant Gromov topology
      	to some action on a real tree $(Y,d_{lim},\rho_{lim})$
      	%  	the sequence of functions $\frac{1}{\mu_{n}}\pl{x}{\--}{\rho\circ\theta_{n}}$ converges to the length function $\pl{_{lim}}{\--}{\rho_{lim}}$ for some pointed action $(Y,d_{lim},x_{lim},\rho_{lim})$
      	of $G$ on a real tree
      	%	   on a real tre
      	which can be decomposed as a tree of actions $\treeact{S}{Y}{}$, where
      	\elenco{
      		\item $Y_{u}$ is a point for $u\in G\cdot v_{0}$
      		\item $Y_{v_{0}}$ is dual to an arational measured foliation on $\Sigma$
      		and for any covering map $\fun{q}{\pi_{1}(\Sigma,*)}{\pi_{1}(\Sigma',*)}$ the action of $Stab(v_{0})$ on its minimal tree does not extend via $q_v{*}$
      		to any action of $\pi_{1}(\Sigma',*)$ on a real tree.
      	}
      \end{lemma}
      \begin{proof}
      	Let $\gamma$ be a \cfam such that $\gamma$ and $\delta$ fill $\Sigma$.
      	%    	Take $\gamma$ and $\delta$ \cfams on $\Sigma$ two mutually intersecting \cfams, and assume moreover that $\delta$ cuts $\Sigma$ into small pieces.
      	Take any collection $\{\phi_{n}\}_{n\in\N}$ of non-trivial Deck transformations of $\Sigma$ such that (up to an homeomorphism isotopic to the identity)
      	any proper covering map with domain $\Sigma$ admits $\phi_{n}$ as a Deck transformation for some $n$.
      	Let us first show how to reduce the lemma to the following claim:
      	\newcommand{\tld}[1]{\tl{#1}{(\rho,d)}}
      	\newcommand{\lnf}[2]{|ln(\frac{#1}{#2})|}
      	\newcommand{\lnft}[2]{\lnf{\tld{\theta_{n}(#1)}}{\tld{\theta(#2)}}}
      	\begin{claim}
      		%	Fix some presentation $\BS{U_{0}}{U}{t}$ of the action of $G$ on $W$, where $v_{0}\in U_{0}$.
      		There are sequences $\ssq{\theta}{n}\subset\pM$ restricting to an automorphisms of $Q_{0}$ and $\ssq{\mu}{n}\subset\R^{+}$ such that:
      		\enum{i)}{
      			\item For any $g\in G$ which does not stabilize a rigid type vertex of $W$ and $n$ big enough $\theta_{n}(g)$ is hyperbolic in $W$.
      			\item $\lnft{g}{h}$ is convergent for $h,g$ none of which stabilizes a rigid type vertex of $W$.
      			\item $\tld{g}$ diverges for some $g_{0}\in Q_{0}$ \label{divergence}
      			\item For any two elements $h,g\in G$ which do not stabilize a rigid type vertex of $W$
      			$\lmt{n}{\lnft{g}{h}}<\infty$  \label{arationality}
      			\item For any $n$ there is some $a_{n}\in Q_{0}$ such that $\lmt{m}{\lnft{a_{n}}{a'_{n}}}>0$, where $a'_{n}$ stands for $(\phi_{n})_{*}(a_{n})$ by $a'_{n}$ \footnote{Note that $(\phi_{n})_{*}(a_{n})$ 	being defined up to conjugacy is enough for the expression above to be well defined.}. \label{asymmetry}
      		}
      	\end{claim}
      	%    	NOOOOOOOOOOOOOOOOOOOOOOOOOOOOOOOOOOOOOOOOOOOOOOOO
      	Indeed,
      	%  up to replacing $\ssq{\theta}{n}$ with a subsequence, we can assume
      	this clearly implies some divergent (by \ref{divergence}) sequence $\ssq{\mu}{n}$ of positive constants exists such that
      	the sequence $\frac{1}{\mu_{n}}\tl{--}{(d,\rho\circ\theta_{n})}$ converges pointwise to $\tl{--}{(d_{lim},\rho_{lim})}$ for some
      	isometric action $\rho_{lim}$ of $G$ on real tree $Y_{lim}$. As we know, this implies convergence in the equivariant Gromov topology.
      	%      	\footnote{As we could have extracted a subsequences, this is not essential.}.
      	%  	 and hence, as we know, in the equivariant Gromov topology.
      	%    	$ ((|W_{0}|,\frac{1}{\mu_{n}}d,\rho\circ\theta_{n}))$
      	%    	the sequence $((|W_{0}|,\frac{1}{\mu_{n}}d,\rho\circ\theta_{n}))$ converges to some minimal
      	%  	action $(Y,d_{lim},\rho_{lim})$ in the Gromov topology.
      	%    	Observe that for some fixed set of generators the maximum translation length  with respect to $(d,\rho\circ\theta_{n})$ of some of them must converge to infinity.
      	Observe the following:
      	\enum{a)}{
      		%    		\item Clearly $Stab(u)$ is elliptic in $Y$ for any rigid type vertex $u\in S$, since $\rho_{n}(Stab(u))$ is conjugate to $Stab(u)$, which is obviously elliptic in $W$.
      		\item For any $u$ of rigid type any $g\in Stab(u)\setminus\{1\}$ is elliptic with respect to each of the actions $\rho\circ\theta_{n}$ and therefore the point it fixes in the limit action is unique.
      		\item For such $u$ the group $Stab(u)$ itself is elliptic.
      		let $x_{u}$ be the unique point it fixes.
      		\item Observe that $G$ is freely indecomposable with respect to the family of all its subgroups  which are elliptic in $Y$. \footnote{Consider the action of $G$ on any tree $U$
      		dual a free decomposition that would contradict this. Stabilizers of rigid vertices of $T$ must act elliptically on $U$, as their image is elliptic in $Y$ and all stabilizers of vertices in $T$ inject in $Q$. Since surface groups are freely indecomposable relative to the family of their boundary subgroups, this means that $Stab(u)$ fixes a unique vertex of $U$ for any $u\in VT$ so that the action of $G$ on $U$ is trivial (free factors are malnormal). }
      		\item The tree $W$ is $2$-acylindrical.
      		%    		since no non-trivial element can fix the star around some vertex of surface type.
      	}
      	Let $Z$ stand for either singleton containing the point fixed by $Stab(v_{0})$ (unique since the group is not abelian) or else the minimal tree of $Stab(v_{0})$.
      	For any segment $v_{0},u,g\cdot v_{0}$ of length $2$ the sets $Z$ and $g\cdot Z$ have to intersect in the unique point fixed by $Stab(u)$.
      	If $Z=\{x\}=Fix(Stab(v_{0}))$ it follows easily by induction that $Stab(v_{1})=Stab(v_{0})^{g^{-1}}$ fixes $x$ for any $v_{1}=g\cdot v_{0}$ of surface type,
      	so that in fact every $g\in G$ does, contradicting the non-triviality of the action. Therefore  $Stab(v_{0})$ is not elliptic in $Y$.
      	In virtue of the results of section \ref{limit actions section}  we know that the action of $Y$ decomposes as a normalized tree of actions $\treeact{S}{Y}{}$,
      	where each non-degenerate $Y_{v}$ is either of Seifert, axial or simplicial type.
      	
      	Assume for a moment that the action of $H:=Stab(v_{0})$ on $Y$ were not elliptic in $S$.
      	%    	and consider the	miinimal tree $S_{Stab(v_{0}})$ of $Stab(v_{0})$ in $S$.
      	%    	this must be properly quoted
      	Since the peripheral subgroups of $Stab(v_{0})$ are elliptic in $Y$, they are elliptic in the minimal tree $S_{Stab(v_{0})}$
      	and hence there is a non-empty \cfam $\delta$ and an equivariant map from the tree $T_{\delta}$ dual to $\delta$
      	%    	admits an equivariant map onto
      	to $S_{H}$.	In particular the conjugacy class corresponding to any of the curves in $\delta$ is mapped into an edge stabilizer of $S$
      	and it must hence be elliptic in $Y$, contradicting hypothesis (\ref{arationality}) listed above, which implies $\tl{g}{(d_{lim},\rho_{lim})}\neq 0$ for any $g\in G$. We conclude that $Stab(v_{0})$ is contained in $Stab(u)$ for some component $Y_{w}$, which cannot be other than of Seifert type, the other ones being either lines or segments.
      	
      	Assume that $Stab(w)\cong\pi_{1}(\Sigma')$ with subgroups of $Stab(w)$ elliptic in $S$ in correspondence with peripheral subgroups of $\Sigma'$.
      	Since boundary subgroups of $Stab(w)$ are sent to peripheral groups of $Stab(w)$, a result due to Scott \cite{scott1978subgroups}
      	that the inclusion of $Stab(v_{0})$ into $Stab(w)$ is induced by a finite covering map $\fun{q}{\Sigma}{\Sigma'}$.
      	That this is in fact a homeomorphism can be proven directly\footnote{As $Stab(w)$ contains $Stab(v_{0})$ as a subgroup of finite index, it must be elliptic in the tree $S$ as well. Since $v_{0}$ is the unique vertex fixed by $Stab(v_{0})$ in fact $Stab(v_{0})=Stab(w)$}, but it also follows from the last property in the statement of the lemma, which is all is now left to show.
      	If the action of $Stab(v_{0})$ on its minimal tree is the pullback of an action of $\pi_{1}(\Sigma',*)$ on a real tree for any covering map $\fun{q}{\Sigma}{\Sigma'}$,
      	pick a deck transformation $\phi_{n_{0}}$ associated to the covering map $q$.
      	%		 and let $(\phi_{n})_{*}$ be the automorphism of $Stab(u)$ it induces.
      	The fact that for any loop $\alpha$ in $\Sigma$ we have $q\circ\phi_{n_{0}}\circ\alpha=q\circ\alpha$ implies that $a$ and $(\phi_{n})_{*}(a)$ are conjugate in $Stab(w)$
      	for any element $a\in Stab(v)$. In particular, this implies that $\tld{a_{n}}=\tld{(\phi_{n})_{*}(a_{n})}$, contradicting property (\ref{asymmetry}) listed above and the continuity of
      	translation length with respect to the Gromov topology.
      	
      	\begin{proof}(of the claim)
      		For each $\alpha\in\delta$ chose $w_{\alpha}\in(0,1)$ such that $w_{\alpha}\neq w_{\beta}$ for any distinct $\alpha\beta\in\delta$ and let $d^{w}$ be the distance on
      		$|S_{\delta}|$ obtained by assigning length $w_{\alpha}$ to each of the edges in the orbit corresponding to $\alpha$, $(|T_{\gamma}|,d^{1})$ the standard geometric realization of $T_{\delta}$ and  	$\epsilon=2\mn{\alpha,\beta\in\gamma\\\alpha\neq\beta}|ln\frac{w_{\alpha}}{w_{\beta}}|$.
      		
      		Given an automorphism $\sigma$ of $S$ we let $\tl{g}{\lambda^{w}_{\gamma}}$ stand for the translation length of $\sigma(g)$ with respect to the action of $g$ on $(T_{\delta},d^{w})$.
      		Denote by $Ell(\lambda)$ the set of elements which are elliptic for an action $\lambda$ and by $Hyp(\lambda)$ the set of those which are hyperbolic.
      		For each $R$ let $B_{R}$ stand for the ball of radius $R$ around the identity in the Cayley graph of $S$ with respect to some fixed finite tuple of generators.
      		In the course of the construction, we generate a sequence $\ssq{A}{n}$ of finite subsets of $S$, where $B_{n}\subset A_{n}$.
      		Assume we have successfully constructed $\theta_{n}$ with the property that $A_{n}\cap Ell(\lambda_{\delta}\circ\theta_{n})=\tg$.
      		
      		Let $\bar{\theta_{n}}$ be some homeomorphism of $\Sigma$ (fixing the basepoint $*$) which induces $\theta_{n}$. Consider the \cfam $\delta'=\bar{\theta_{n}}^{-1}$.
      		Since $\delta'$ splits $\Sigma$ into small pieces and $\phi_{n}$ is not homotopic to the identity, it might behave in two different ways with respect to the deck transformation $\phi_{n}$. The first 	possibility is that $i(\phi_{n}(\alpha),\delta')\neq 0$ for some $\alpha\in\delta'$. Given a non-trivial element $a_{n}$ in the conjugacy class associated to $\alpha$, this is the same as to say that $a_{n}\in Ell(\lambda_{\delta}\circ\theta_{n})$ but $\phi_{n}(a_{n})\in Hyp(\lambda_{\delta}\circ\theta_{n})$.
      		The other one is that $\phi_{n}(\delta')=\delta'$. In that case some $\alpha\in\delta'$ is sent by $\phi$ to some other $\beta\in\delta'$.
      		%    	represented by some element $b_{n}$ of the fundamental group.
      		Chose as $a_{n}$ in this case an element with translation length $2$ in $T_{\gamma'}$ (the original simplicial tree) and whose axis consists entirely of translates of $e_{\alpha}$ of $\edin{e_{\alpha}}$.
      		Clearly the same property holds for $\phi_{n}(a_{n})$ only with respect to $e_{\beta}$ instead of $e_{\alpha}$.
      		In virtue of lemma \ref{Dehn twist approximation} there is a power $\tau^{\delta}_{n}\in\pM$ of Dehn twists over edges of $T_{\delta}$
      		satisfying the following inequalities for any $h,g\in (Hyp(\lambda_{\delta}\circ\theta_{n})\cap B_{n+1})\cup A_{n}$:
      		\begin{align*}
      			|ln(\frac{(\tl{g}{\lambda_{\gamma}\circ\tau^{\delta}_{n}\circ\theta_{n},d^{1})}
      			\tl{h}{(\lambda_{\delta}\circ\theta_{n},d^{w})}}
      			{\tl{g}{(\lambda_{\delta}\circ\theta_{n},d^{w})}
      			\tl{h}{(\lambda_{\gamma}\circ\tau^{\delta}_{n}\circ\theta_{n},d^{1})}}
      			)|<\frac{\epsilon}{2^{n+3}} \\
      			\tl{g}{\lambda_{\gamma}\circ\tau^{\delta}_{n}\circ\theta_{n}}>n
      		\end{align*}
      		For $\eta_{n}=\tau^{\delta}_{n}\circ\theta_{n}$ this implies in particular that
      		$Hyp(\lambda_{\delta'})\cap A_{n}=Hyp(\lambda_{\delta}\circ\theta_{n})\cap A_{n}\subset Hyp(\lambda_{\gamma}\circ\eta_{n})$.
      		On the other hand if $g\in Ell(\lambda_{\delta}\circ\theta_{n})=Ell(\lambda_{\delta}\circ\eta_{n})$ then $g\in Hyp(\lambda_{\gamma}\circ\eta_{n})$, since $\delta'=\bar{\eta_{n}}^{-1}(\delta)$ and 	$\bar{\eta_{n}}^{-1}(\gamma)$ fill $\Sigma$.
      		% $\eta_{n}(g)$ and $\theta_{n}(g)$ musts be conjugate to each other.
      		If $\phi_{n}$ is as in the first of the two possibilities discussed above, we can also take:
      		\begin{align*}
      			|ln(\frac{\tl{(\phi)_{*}(a_{n})}{(\lambda_{\gamma}\circ\eta_{n},d^{1})}}{\tl{a_{n}}{(\lambda_{\gamma}\circ\eta_{n},d^{1})}})|>\epsilon
      		\end{align*}
      		\newcommand{\lama}[0]{(\lambda_{\delta}\circ\tau^{\gamma}_{n}\circ\eta_{n},d^{w})}
      		\newcommand{\nlama}[0]{(\lambda_{\gamma}\circ\eta_{n},d^{1})}
      		If the second one is the case, we can assume that
      		\begin{align*}
      			|ln(\frac{\tl{a_{n}}{\nlama}}{\tl{\phi_{n})_{*}(a_{n})}{\nlama}})|>\epsilon
      			%    		\frac{\epsilon}{2^{n+3}}
      		\end{align*}
      		for $g\in\{a_{n},(\phi_{n})_{*}(a_{n})\}$ as well. Let $A_{n+1}=A_{n}\cup B_{n}\cup\{a_{n},(\phi_{n})_{*}(a_{n})\}$.
      		%    	$\eta_{n+1}=min\{\eta_{n},\frac{1}{2}|ln(\frac{w_{\alpha}}{w_{\beta}})|\}$ in the first case and $\eta_{n}$ otherwise.
      		Another application of lemma \ref{Dehn twist approximation} yields the existence of a product $\tau^{\gamma}_{n}$ of Dehn twists over edges of $T_{\gamma}$ such that
      		if we let $\theta_{n+1}:=\tau^{\gamma}_{n}\circ\eta_{n}$ then $Ell(\lambda_{\delta}\circ\theta_{n+1})\cap A_{n+1}=\{1\}$ and:
      		\begin{align*}
      			|ln(\frac{\tl{g}{\lama}}{\tl{h}{\lama}})|<\frac{1}{2^{n+3}}
      		\end{align*}
      		for any $h,g\in A_{n+1}$. At this point we can let $\theta_{n+1}=\tau^{\gamma}_{n}\circ\tau^{\delta}_{n}\circ\theta_{n}$ and iterate the process.
      		Given any $h,g\in B_{n_{0}}$ for $n\geq n_{0}$ obviously $h,g\in Hyp(\lambda_{\delta}\circ\theta_{n})$ and:
      		\renewcommand{\lama}[0]{(\lambda_{\delta}\circ\theta_{n},d^{w})}
      		\newcommand{\lamao}[0]{(\lambda_{\delta}\circ\theta_{n+1},d^{w})}
      		\begin{align*}
      			\serprop{n\geq n_{0}}{|ln(\frac{\tl{g}{\lama}\tl{h}{\lamao}}{\tl{h}{\lama}\tl{g}{\lamao})|}}
      			\leq\serprop{n\geq n_{0}}{\frac{\eta_{n}}{2^{n+2}}}\leq\frac{\epsilon}{2}
      		\end{align*}
      		Using a standard diagonal argument, one can obtain a subsequence $(\theta_{n_{k}})_{k}$ and any $h,g$ not contained in a boundary subgroup of $S$
      		the sequence of values
      		\renewcommand{\lama}[0]{(\lambda_{\delta}\circ\theta_{n_{k}},d^{w})}
      		\renewcommand{\lamao}[0]{(\lambda_{\delta}\circ\theta_{n_{k+1}},d^{w})}
      		$(|ln(\frac{\tl{g}{\lama}}{\tl{h}{\lama}})|)_{n}$ converges to some limit, which is greater than $\frac{\epsilon}{2}$ in case $(g,h)=(a_{n},(\phi_{n})_{*}(a_{n}))$.
      		There is a small catch: we claimed to be able to take as $(|S_{\delta}|,d^{w})$ any geometric realization of $S_{\delta}$
      		however this restriction can be easily removed by postcomposing each $\theta_{n}$ with a sufficiently big power of Dehn twists over $\delta$ and a big powers of Dehn twists over $\gamma$.
      	\end{proof}
      	This concludes the proof of \ref{generic test sequences}.
      \end{proof}
      Another useful lemma is the following:
      \begin{lemma}
      	\label{free product and growth}Suppose we are given a (\rs) \pr-tower $\LRL{T}{}{(J,r)}$, a test sequence $\ssq{f}{n}$ of $\QR{T}$, a limit group $L$ and
      	a convergent sequence $\ssq{g}{n}$ of homomorphisms from $G:=L\frp T$ to $\F$ such that:
      	\elenco{
      		\item	 $g_{n}\rst_{T}=f_{n}$
      		\item  $\limker{g_{n}}{n}\cap G=\tg$
      		\item $T^{\lambda}$ grows faster than $G$ under $\ssq{f}{n}$ for any $\lambda\in J$
      	}
      	Then $\limker{g_{n}}{n}=\tg$.
      \end{lemma}
      \begin{proof}
      	This follows from the arguments used in the proof of Merzlyakov's theorem, so we won't exercise excesive rigor in the proof.
      	Let $N:=L/\limker{g_{n}}{n}$. Let $Y$ be the limiting tree for the sequence $\ssq{g}{n}$ and $Y_{T}$ the minimal tree of $T$ in $Y$.
      	If no non-trivial element of $G$ fixes a point of $Y_{T}$ with non-trivial stabilizer in $T$,
      	then it is easy to show that $N$ is isomorphic to the free product of $G$ and $T$ (see the proof of \ref{adding a node}).
      	Otherwise $N$ is isomorphic to the amalgamated product $G$ and some group $B_{\lambda}$ obtained from $T$ by enlarging $\Tl$ for some $\lambda\prec r$.
      	By induction we can write $B_{\lambda}$ as the free product of $G\cap B_{\lambda}$ and $\Tl$. This easily implies that $N$ is the free product of $T$ and $G$.
      \end{proof}
      \begin{corollary}
      	\label{surface test sequences} Suppose we are given a \pr-tower structure $\LRL{T}{A}{(J,r)}$, where $r$ is of floor type and $\DR{T}{r}$ contains a unique vertex of non-rigid type, 	which is of surface and that for any $\mu\prec r$ we are given a family $\mathcal{TS}^{\mu}$ of (relative) \pr-test sequences $\RTA\rst_{\mu}$, closed under taking subsequences and diagonal 	subsequences and such that $\TAm$ is elliptic in any limit tree of a sequence extending both $\ssq{f^{1}}{n}$ and $(f_{n}^{\mu_{2}}\rst_{\RLR{T}{}{\nu}})$ for $\mu_{1}\lo\mu_{2}$, 	$\ssq{f^{j}}{n}\in\mathcal{TS}^{\mu_{j}}$ and $\nu\leq\mu$.
      	Then, given any tuple $(\ssq{f^{\mu}}{n})_{\mu}\in\times_{\mu\prec r}\mathcal{TS}^{\mu}$ there is a \pr-test sequence $\ssq{f}{n}$ of $f_{n}$ extending
      	$\bdcup{\mu\prec r}{f^{\mu}_{n}}{}$ and such that the collection $\mathcal{TS}$ of all of them is closed under taking subsequences and diagonal sequences.
      	Given any fixed increasing sequence $\ssq{\kappa}{n}$ of positive constants we can require that
      	$\frac{\tl{f_{n}(g)}{r}}{\kappa_{n}}$ tends to infinity with $n$ for any $\ssq{f}{n}\in\mathcal{TS}$ and some $g\in G$.
      	%    		 $\ssq{f^{\mu}}{n}$ and $\ssq{f^{r}}{n}$ . Then each
      	%as well as a sequence of morphisms $\fun{\ssq{f}{n}}{\bfrp{\RLR{T}{A}{r}}{\mu\prec r}{}}{\rqG{F}{A}}$ with trivial limit kernel.
      \end{corollary}
      \begin{proof}
      	Let $\delta$ be the \cfam provided by lemma \ref{small pieces injected}, applied to the unique surface group of $\Delta=\DR{T}{r}$, isomorphic to $\pi_{1}(\Sigma)$, and $\eta:=\eta^{r}_{\QR{R}}$. Take $\gamma$ a \cfam on $\Sigma$ which together with $\delta$ fills $\Sigma$. Let $\ssq{\theta}{n}\subset Mod^{\pi}(\qG{G},S)$ and $(Y_{lim},d_{lim},*_{lim},\rho_{lim})$ as provided
      	by the last lemma.
      	%    	where $(|S^{c}_{\delta}|,d,\rho)$ is the standard geometric realization of the tree $S^{c}_{\delta}$.
      	%    	We would like to apply corollary \ref{Dehn twist approximation} to approximate convergent sequence. We face, however, two potential obstacles:
      	Let $\mathcal{TS_{0}}$ be the family of sequences obtained by applying \ref{asymmetrical growth} to $(\mathcal{TS}^{\mu})_{\mu\in Ch(r)}$,
      	in a way that makes $\Tm$ grow faster than $\RLR{T}{}{r}$ in case $\mu\prec r$.
      	%  	collection of all the sequences of the form $(f_{n}\circ\eta)_{n}$, where $\ssq{f}{n}=\bdcup{\mu\prec r}{f^{\mu}_{n}}{}$
      	%  	for a tuple $((\ssq{f^{\mu}}{n})_{n})_{\mu\prec r}\in\times_{\mu\prec r}\mathcal{TS}^{\mu}$.
      	In virtue of \ref{free product and growth}, we know that $\limker{g_{n}}{n}$ for all $\ssq{g}{n}\in\mathcal{TS}_{0}$.
      	Let $\mathcal{UTS}_{0}$ be the collection of all the $f\in Hom(G,\F)$ appearing as $f_{n}$ for some $\ssq{f}{n}\in\mathcal{TS}_{0}$.
      	For any $f\in \mathcal{UTS}_{0}$
      	denote by $N_{f}$ the maximum value of $n$ such that $f=\inn{c}\circ f_{n}$ for for some $c\in\F$ and $\ssq{f}{n}\in\mathcal{F}$.
      	Let $A_{0}\subset A_{1}\subset A_{2}\cdots A_{n}\cdots$ some exhaustion of $G$. For each $n$ let $B_{n}\subset G\setminus\{1\}$ be the set of $S$-elliptic elements obtained by applying
      	lemma \ref{Dehn twist approximation} to $\sigma_{n}(A_{n})$
      	%    	Up to replacing the sequence $\ssq{\theta}{n}$ by a subsequence, we can assume that
      	and $k_{n}$ a positive integer such that $ker\,f_{m}\cap B_{n}=\nil$ for any $m\geq k_{n}$ and $\ssq{f}{n}\in\mathcal{F}$.
      	Pick any basepoint $*'$ in the Cayley graph of $\F$. Given $f\in\mathcal{UTS}_{0}$
      	let $u_{f}:=max\,\setof{n\in\N}{k_{n}\leq N_{f}}$, so that $ker\,f\cap B_{u_{f}}=\nil$.
      	In virtue of lemma \ref{Dehn twist approximation} we know that there is some
      	%    	product
      	$\tau_{f}\in Mod^{\pi}(\qG{G},S)$
      	%    	of Dehn twists over edges of $S$ (which we can take in any given finite index subgroup of $Aut(G)$)
      	and $\nu_{f}\in\R^{+}$ such that
      	\begin{align*}
      		|\frac{1}{\nu_{f}\mu_{u_{f}}}\pl{*'}{\sigma_{u_{f}}(g)}{r\circ f\circ\tau_{f}}-\frac{1}{\mu_{u_{f}}}\pl{*}{\sigma_{u_{f}}(g)}{\rho}|\leq
      		|\pl{*_{lim}}{g}{\rho_{lim}}-\frac{1}{\mu_{u_{f}}}\pl{*}{\sigma_{u_{f}}(g)}{\rho}|
      	\end{align*}
      	for every $g\in B_{u_{f}}$, where $r$ is the action of $\F$ on its Cayley graph.
      	We claim that the family $\mathcal{TS}$ consisting of all the sequences of the form $(g_{n}\circ\tau_{g_{n}})_{n}$, where $\ssq{g}{n}\in\mathcal{TS}_{0}$ satisfies the desired properties.
      	First of all, the family clearly inherits being closed under subsequences and taking diagonal sequences from $\mathcal{TS}_{0}$.
      	Now, given some $\ssq{f}{n}\in\mathcal{TS}_{0}$ and $g\in B_{m}$
      	%    	observe the sequence of values $(u_{f_{n}})_{n}$ tends to infinity.t This means that
      	for any $n\geq k_{m}$ we have
      	\begin{align*}
      		|\frac{1}{\nu_{f_{n}}\mu_{u_{f_{n}}}}\pl{*'}{g}{r\circ f_{n}\circ\tau_{f_{n}}\circ\sigma_{u_{f_{n}}}}-
      		\pl{*_{lim}}{g}{\rho_{lim}}|\leq
      		%	    	\frac{1}{\mu_{u_{f_{n}}}}\pl{*}{\sigma_{u_{f_{n}}}(g)}{\rho_{lim}}|\leq
      		2|\pl{*_{lim}}{g}{\rho_{lim}}-\frac{1}{\mu_{u_{f_{n}}}}\pl{*}{\sigma_{u_{f_{n}}}(g)}{\rho_{lim}}|
      	\end{align*}
      	since $u_{f_{n}}$ tends to infinity with $n$, the right hand side tends to $0$, which implies the sequence $(f_{n}\circ\tau_{f_{n}}\circ\sigma_{f_{n}})_{n}$
      	has trivial limit kernel and admits $Y_{lim}$ as a limiting tree, which implies it satisfies the test sequence condition at $r$.
      	\begin{comment}
      	For the last claim, chose some $g$ which cannot be conjugated into the stabilizer of a rigid type vertex group of $S$. For $n$ bigger than some $C$
      	$\sigma_{n}(g)$ is hyperbolic in $S_{\delta}^{c}$. So as soon as $u_{f}\geq C$ we can clearly chose $\tau_{f}$ so that
      	$\tl{f(\sigma_{u_{f}}(g))}{\lambda}\geq u_{f}\kappa_{N_{f}}$. Given $\ssq{f}{n}\in\mathcal{TS}_{0}$,
      	eventually) $\frac{\tl{f_{n}\circ\tau_{f_{n}}\circ\sigma_{u_{f_{n}}}(g)}{\lambda}}{\kappa_{n}}\geq n$  (of course $N_{f_{n}}\geq n$).
      	Observe that postcomposing both $\tau_{f}$ and $\sigma_{n}$ can be taken to fix the rigid type stabilizer containing $A$, so that the resulting sequence factors through $\RTA$.
      	%    	 with an inner automorphism does not alter any of the properties above (provided we change the basepoint accordingly)
      	\end{comment}
      \end{proof}
      Note that in the theorem above, if each $\mathcal{TS}^{\mu}$ is congruency complete, the same is true for $\mathcal{TS}$.
      
  \subsection{The induction step: abelian vertex groups}
    
    \begin{lemma}
    	Let $G$ be a group of the form $L\frp_{M}(M\oplus\subg{t})$, where $M$ is a non-trivial free abelian group and suppose we are given a sequence $\ssq{f}{n}\neq Hom(G,\F)$ with trivial limit kernel.
    	Let $x$ be some finite set of generators of $G$. Suppose that for each $n$ we are given an extension $\fun{g_{n}}{G}{\F}$ of $f_{n}$
    	such that $\tl{g_{n}(t)}{}\geq n\cdot\tl{g_{n}(x_{i}x_{j})}{}$ for any $1\leq i<j\leq|x|$. Then then any limiting tree for the sequence $\ssq{g}{n}$ is isomorphic to a geometric realization of the simplicial tree $S$ dual to an HNN extension with vertex group $L$ and Bass-Serre element $t$.
    \end{lemma}
    \begin{proof}
    	This can be proven directly without any advanced technology using the same techniques of the proof of \ref{Dehn twist approximation}. But it also follows straightforwardly from the fact that $L$ is elliptic and $t$ (its image in the limit quotient of the sequence $g_{n}$ ) hyperbolic in any limiting tree $Y$, together with all what we know about the action on $Y$. In particular, the point-wise stabilizer of $Ax(t)$ coincides with the centralizer of $t$ in $G$.
    \end{proof}
    
    \begin{corollary}
    	Let $\LRL{T}{A}{(J,r)}$ a refined \pr-tower where $r$ falls in case \ref{pegged abelian case}. If there is a congruency complete CSD family $\mathcal{TS}_{\lambda}$ of \pr-test sequences
    	of $\RTA\rst_{\lambda}$, where $\lambda$ is the unique child of $r$, then there is a congruency complete CSD family of \pr-test sequences of
    	$\RTA$.
    \end{corollary}
    \begin{proof}
    	Let $M\oplus\subg{t}$ be the abelian type vertex group of $\DR{R}{r}$, where $M\leq\Tl=:L$ and $x$ some finite tuple of generators of $L$. As $\mathcal{TS}$ we take the family of sequences
    	$\ssq{g}{n}$ of morphisms factoring through $\RTA$ such that $g_{n}\rst_{L}\in\mathcal{TS}$ and $\tl{g_{n}(t)}{}\geq n\cdot\tl{f_{n}(x_{i}x_{j})}{}$ for $1\leq i<j\leq |x|$. Any such sequence is a \pr-test sequence, by the previous lemma, and the resulting family is clearly congruency complete.
    \end{proof}

\section{Characterization of \rs \pr-limit groups}
  \label{discrimination section}
  %  The most immediate consequence of the existence of test sequences is the following.
  
  % \subsection{Elementary \pr-limit groups }
    \begin{comment}
    The most basic examples of limit groups, other than free groups, are finitely generated abelian groups and closed surface groups. The family of all \pr-structures on such groups making any of those into \pr-limit groups admits a simple description, which we present here, together with a couple of additional facts.
    \end{comment}
    The most basic examples of limit groups, other than free groups, are finitely generated 	abelian groups and closed surface groups. It turns out that the family of all \pr-structures on such a group making it into \pr-limit groups admits a simple description.
    %    The results here require some notions and results from later chapters, but we include them here for the sake of illustration.
    \newcommand{\rcl}[1]{\sqrt{#1}}
    \begin{comment}
    \begin{lemma} \label{free plimits}
    	Let $\pmb{@F}$ be a \pr-group which is free, where $H$ is a non-abelian free group. Then $\qG{H}$ is a \pr-limit group.
    \end{lemma}
    \end{comment}
    %  We need a little discussion earlier on normal forms in the free group.
    \begin{comment}  % 2016/01/27/10:52:17_:
    \begin{lemma} \label{free plimits}
    	Let $\pmb{@F}_{A}$ be a \rs \pr-group, where $H$ is a non-abelian free group. Then $\qG{F_{l}}_{A}$ is a \rs \pr-limit group
    \end{lemma}
    \end{comment}
    The following is a consequence of observation \ref{injectivity Nielsen} below.
    \begin{lemma} \label{free plimits}
    	Let $\qG{H},\qG{H}'$ be  \pr-groups and suppose that $H,H'$ are non-abelian free groups. Then there is an injective morphism from $\qG{H}$ to $\qG{H}'$.
    	Let $\rqG{L}{A}=\qG{A}\frp\qG{H}$, where $H$ is free. Then $\rqG{L}{A}$ is a \rs \pr-limit group.
    \end{lemma}
    
    \begin{corollary}
    	\label{reference model invariance}The notion of a \pr-limit group depends only on $Q$ and not on the particular epimorphism $\fun{\pi}{\F}{Q}$.
    \end{corollary}
    
    %    Hence lemma \ref{free plimits} immediately follows from the much stronger lemma \ref{small cancellation embeddings}.
    
    %    Given a commutative transitive group $G$ and $N$ an abelian subgroup of $G$, by the root closure of $N$ in $G$, denoted by $\rcl{N}$, we intend the set of all elements of $G$ with powers in $N$. By commutative transitivity, this is a subgroup of $Z_{G}(N)$. Observe that if the latter is a free abelian group and $N$ finitely generated, then $\rcl{N}$ is a direct summand of $Z_{G}(N)$.
    \begin{lemma} \label{free abelian plimits}
    	Let $\qG{M}$ a \pr-group whose underlying group $M$ is free abelian. Then $\qG{M}$ is a \rs \pr-group if and only if $\pi(M)\leq C$ for some cyclic $C\leq \F$.
    \end{lemma}
    \begin{proof}
    	The implication from left to right is clear. For the opposite direction, take $x\in M$ such that $\pi(x)$ generates $\pi(M)$ in $Q$. We can assume
    	$Z=\subg{x}$ is root closed in $M$. Write $M$ as $Z\oplus N$ and let $y_{1},\cdots y_{m}$ a base for $N$. Pick $K\in\N$ dividing the order of $\pi(C)$ in $Q$.
    	Fix some tuple $\bar{r}=(r_{i})_{i=1}^{m}\in\{0\cdots K\}$ of distinct integers such that $\pi(y_{i})=x^{r_{i}}$. For each $n\in \N$ consider the \pr-retraction $f_{n}$ of $\qG{M}$ onto $\qG{Z}$ sending $y_{i}$ to $x^{e_{i,n}}$, where $e_{i,n}=K(n+i)!+r_{i}$.
    	It is easy to check that given $0\leq i_{1}< i_{2}\cdots< i_{r}\leq m$ and non-zero integers $k_{1},k_{2}\cdots k_{r}$
    	the last term of the sum $k_{1}e_{i_{1},n}+k_{2}e_{i_{2},n}\cdots k_{r}e_{i_{r},n}$ grows in absolute value faster with $n$ than all others and that therefore the total number is non-trivial for $n$ big enough.
    	This easily implies taht the limit kernel of the $f_{n}$ is trivial. Since $\qG{Z}$ is itself a \pr-limit group we are done.
    	%  	As given any sequence $f_{n}$ of restricted morphisms with limit kernel from $M$ to $\F$ there are finitely many
    \end{proof}
    
    \begin{lemma} \label{closed surface plimits}
    	Let $\qG{S}$ be a \pr-group whose underlying group $S$ is isomorphic the fundamental group of a  closed surface $\Sigma$ of Euler characteristic $\geq-2$. Then $\qG{S}$ is a \pr-limit group if and only if for some maximal \cfam $\delta$ in $\Sigma$ each of the root closure of the conjugacy classes in $S$ in correspondence with the curves is killed by the homomorphism $\pi$.
    \end{lemma}
    \begin{proof}
    	For the 'only if' all is needed is the fact that $\pi_{S}=\pi\circ h$ for some homomorphism $h$ from $S$ to $\F$. In virtue of \ref{pinching to a free group}, $h$ can be written as $h'\circ p$, where $p=q\circ pinch_{\delta}$ for an \cfam $\delta$ in $\Sigma$ and $q$ kills the fundamental group of any projective plane in $Pinch(\Sigma,\delta)$.
    	
    	Notice that the image of $p$ is the fundamental group of an finite graph, hence free.
    	For the opposite direction notice that viceversa, the fact that some $\delta$ satisfies the property of the statement implies $\pi_{S}$ factors through the quotient $p$ above.
    	
    	Now, by virtue of the discussion in \ref{pinching subsection} the image of $p$ is a free group $F$ and the assumption that $\xi(\Sigma)\leq-2$ implies it is non-abelian. Push the \pr-structure on $S$ forward onto $F$.
    	%    	We know that there is some morphism of $\F$ in
    	%    	 system the corresponding pinching map (see \ref{pinching section} for details) is a morphism with non-abelian free image.
    	The result follows from \ref{reference model invariance} together with the existence of \pr-test sequences (the \pr-tower in this case consists of a unique vertex group of surface type).
    \end{proof}

  \begin{lemma}
  	Let $\rqG{G}{A}$ be a \rs \pr-group, $S$ a \gad $\Delta$ of $\rqG{G}{A}$. Then any morphism $f$ from $\rqG{G}{A}$ to a \rs \pr-limit group $\rqG{H}{A}$ which is formally strict with respect to $S$.
  	is strict (in particular $\rqG{G}{A}$ is also a \rs \pr-limit group ).
  \end{lemma}
  \begin{proof}
  	Let $\ssq{h}{n}$ be a sequence of morphisms from $\rqG{H}{A}$ to $\rqG{F}{A}$ with trivial limit kernel. To begin with, let $G'$ be the retract of $G$ obtained by replacing any abelian type vertex group of $\Delta$ by its extended peripheral subgroup.
  	%  	Using lemma \ref{small surfaces injeceted} and corollary \ref{discrimination corollary} one can show, exactly as in the construction of \pr-test sequences
  	Corollary \ref{second discrimination corollary}, together with the fact that $f$ is formally strict with respect to $S$ implies the existence of a sequence of automorphisms
  	$\ssq{\tau}{n}\supset Mod^{\pi}(\qG{G},S)$ such that the limit kernel of  $(h_{n}\circ f\circ\tau_{n})$ has trivial intersection with $G'$ (up to subsequence extraction we can assume that this sequence is convergent).
  	
  	Now, up to the fact that the image of an abelian vertex group of $\Delta$ by $h_{n}$ might be larger than that of its extended peripheral subgroup, the result is now a consequence of regarding $\rqG{G}{A}$ as a \pr-tower over $\rqG{G'}{A}$.
  	A quick look at the proof of case \ref{pegged abelian case} of the construction of \pr-test suffices to check that this is not really an obstruction.
  	%  	 the result now follows from the abelian floor case of the construction of test sequences.
  	%  	 enlarging the abelian vertex groups of
  \end{proof}
  \begin{corollary}
  	A \rs \pr-group $\rqG{G}{A}$ is am $A$-\rs \pr-limit group if and only if it admits a strict closed weak \rs \pr-\res.
  \end{corollary}
  The latter can be also seen as a direct consequence of proposition \ref{test sequence existence} together with the embedding constructed in section \ref{completions subsection}.

\chapter{Merzlyakov-type theorems}
\label{Merzlyakov chapter}
\section{\pr-towers and formal solutions}
  \newcommand{\cosetdet}[4]{\bwedge{i=0}{#2_{#1}\in #3_{#1}}{#4}}
  In this section we provide a generalization of the following theorem, originally due to Merzlyakov:
  %    In what follows we will designate tuples of variables by non-barred small letters $x,y,z$
  \begin{theorem}
  	%    	\renewcommand{\F}[0]{\mathbb{F}}
  	% 	Add a section on the meaning of $\subg{x}$ in the family of morphisms.
  	Let $\F$ be a free group and $\Sigma(x,y,a)=1$ a finite system of equations in the tuple $x$ over a tuple $a$ of constants from $\F$. Assume that the group generated by $a$ in $\F$
  	is a free factor of $\F$ and suppose that:
  	\begin{align*}
  		\F\modelof\forall x\exists y\,\,\,\syeq{\Sigma}{x,a}
  	\end{align*}
  	Then there is a word $w(x,a)$ such that the equality $\Sigma(x,w(x,a),a)=1$ holds in $\subg{x}\frp\subg{a}$.
  	%  	So that for each group $G$ one has $\F\modelof\forall x\; \Sigma(x,w(x,a),a)=1$.
  \end{theorem}
  We would like to remark that the result does not hold in case $\subg{a}$ is not a free factor. For a trivial example of how this can fail if $A$ take
  $\F=\F(\alpha,\beta)$, $a=(\alpha,\alpha^{\beta})$, $x$ and $y$ single variables and our system of equations simply:
  \begin{equation*}
  	\begin{cases}
  		x=1 \\
  		a_{1}^{y}=a_{2}
  	\end{cases}
  \end{equation*}
  The corresponding $\forall\exists$ sentence is clearly true, but there is no way that some word $w(x,a)$ as above can exist (in this case really $w(a)$) since $a_{1}$ and $a_{2}$ are not conjugate in $\subg{a_{1},a_{2}}$.
  One refers to the word $w$ as a \emph{formal solution} for the sentence above. Observe that the existence of such a word is equivalent to that of a retraction from the group
  $G_{\Sigma}$ to the subgroup generated by $x$ and $a$.
  %  The condition on $A$ cannot be weakened. Indeed, let $G_{\sigma}$ be the group
  
  Merzlyakov theorem can be generalized to finite disjunctions of systems of equations. If one lets systems of inequations enter the hypothesis, though, the picture becomes much more complex. Indeed, given a sentence with constants of the form:
  \begin{align*}
  	\forall x\exists y\,\syeq{\Sigma}{x,y,a}\wedge\syneq{\Pi}{x,y,a}
  \end{align*}
  which is valid in $\F$, in general there is no way of finding an algebraic expression $w(x,a)$ such that $\syeq{\Sigma}{x,w(x,a),a}\wedge\syneq{\Pi}{x,w(x,a),a}$ holds for any possible value of $x$.
  In case that $\F(x)\frp\subg{a}$ is non abelian one can guarantee however that $\syneq{\Pi}{x_{0},w(x_{0},a),a}$ for at least some value $x_{0}$ of $x$, in oder words, that the system $\syeq{\Sigma}{x,a}$ does not imply $\syeq{\Pi}{x,w(x,a),a}$ in $\F$ and hence the set of those values of $x$ for which the condition fails is a proper variety of $\F^{|x|}$.
  
  The natural way to proceed is trying to construct an ad-hoc formal solution that works in the subvariety of the points where first formal solution fails, except for those points in a yet smaller proper variety. Equational noetherianity would then imply that the process eventually stops, yielding a piecewise algebraic function that witnesses the validity of the sentence.
  
  The first, but inessential, complication lies in the fact that in this case there is no single formal solution that can cover all values of $x$ in the specified variety.
  The main obstacle is however, that these more specialized formal solutions cannot be found in terms of $x$ and $a$ alone. They necessarily involve some new tuple $t$ of variables, which are linked to $x$ and $a$ by certain relationships. And in that case, knowing the failure of the formal solution of a particular $(x_{0},t_{0})$ might fail to give any new information on $x_{0}$ in isolation, breaking the previous argument.
  \begin{comment}
  That is, for any value of $x$ and of $t$, setting $y=w(x,t,a)$ ensures the satisfaction of the system of equalities. And tehe system of inequalities for at least one value $(x_{0},t_{0})$. It can be the case, however, that for all values of $x$ some value of
  which prevents the previous argument to work.
  \end{comment}
  %    That is for each value of $x$ and any value of $t$ satisfying these constraints, writing the
  %    This implies the argument above fails, for although the set of values of $(x,t)$ for which the formal solution satisfies the
  Roughly speaking, the way out of this is to show that the relationships linking $x$ are $t$ are tight enough for a certain notion of complexity to properly decrease from one step to the next of the procedure, when adequately refined.
  
  Before being able to tackle such difficulties, a subject that unfortunately exceeds the scope of this dissertation, one needs to adapt the generalization of Merzlyakov theorem from \cite{sela2}), which handles $\forall\exists$ formulas with algebraic conditions on universal variables to our particular context.
  %   We will use an expression like $\forall x\in q\,\Psi(x)$ as abbreviations for $\forall x\,(x\in q\sra\Psi(x))$ and $\exists x\in q\,\Psi(x)$ as one
  %	for \exists x\,x\in q\wedge\Psi(x)
  %    Given a \pr-group with an $x$-marking, the notation $q^{x}$ will denote the tuple $\pi(x)\in Q^{|x|}$, etc.
  \begin{comment}
  Our aim is to analyze $\forall\exists$ sentences of the following form:
  \begin{align*}
  	\pVAFor{\Delta}{\Sigma}{\Pi}{p}{q}
  \end{align*}
  \end{comment}
  %    It is not possible to analyze all solutions to the \CEQ $\ceq{\Delta}{x,a}{}$ at once.
  \begin{comment}
  Since in the common theory of \pr-groups the predicate $P_{q}$ is equivalent to $\bvee{q'\neq q}{P_{q}}{}$, any $\LQ$-formula is equivalent to one of the form:
  \begin{align*}
  	\forall x^{1}
  \end{align*}
  \end{comment}
  
  The class of sentences dealt with in the next theorem are almost general $\forall\exists$ sentences with constants, as
  any $\forall\exists$ sentence can be written as a finite conjunction of formulas of that type.
  \begin{theorem}
  	\label{clean Merzlyakov}
  	Suppose that $\qG{F}$ satisfies the sentence with constants:
  	\begin{align*}
  		\pVAFor{\Delta}{\Sigma}{\Pi}{q^{x}}{q^{y}}{r}
  	\end{align*}
  	Then there is a tuple of variables $s$, some $m\geq 1$ and for any $1\leq l\leq m$:
  	\enum{i)}{
  		\item A \EQ $\syeq{\Lambda_{l}}{x,t,a}$
  		\item Some $q^{s}(l)\in Q^{|s|}$.
  		\item A $|y|$-tuple $w^{l}(x,s,a)$ of words in $x,s,a$
  		\item Some $k_{l}\in\{1,\cdots r\}$
  	}
  	so that the following list of conditions are satisfied:
  	\enum{i)}{
  		\item \label{equations satisfied}The following holds in any group containing $A$:
  		\begin{align*}
  			\forall x\,\forall 	s\,(\syeq{\Delta}{x,a}\wedge\syeq{\Lambda_{l}}{x,s,a})\sra\syeq{\Sigma_{k_{l}}}{x,w_{l}(x,s,a),a}))
  		\end{align*}
  		%    		Is valid in any group.\\
  		%    		\item For each word $u(x,y,a)$ in $\Sigma$, $u(x,w^{l}(x,t,a),a)$ is the trivial element in the quotient of the group
  		%    		$A\frp\F(x,t,y)$ by the normal closure of the relationships $\syeq{\Delta}{x,a},\syeq{\Lambda_{l}}{x,t,a}$,
  		\item  \label{retraction in Q} The following equality holds in $Q$ for $1\leq l\leq m$:
  		\begin{align*}
  			w^{l}(q^{x},q^{s}(l),q^{a})=q^{y}(k_{l})
  		\end{align*}
  		\item  \label{extension}The following holds in $\qG{F}$:
  		\begin{align*}
  			\forall x\,\exists s\,\,(\bvee{i=1}{s\in q^{s}(l)\wedge\syeq{\Lambda_{l}}{x,s,a}}{m})
  		\end{align*}
  		\item \label{proper subvariety}For all $1\leq l\leq s$ there is $(x_{0},s_{0})\in\F^{|s|+|x_{0}|}$:
  		\begin{align*}
  			x_{0}\in q^{x}\wedge s_{0}\in q^{s}(l)\wedge \syeq{\Delta}{x_{0},a}\wedge\syeq{\Lambda_{l}}{x_{0},s,a}\wedge\syneq{\Pi}{x_{0},w^{l}(x_{0},s_{0}),a}
  		\end{align*}
  		
  	}
  	
  	If $\syeq{\Delta}{x,a}$ is empty and $\subg{\pi(a),q^{x}}=Q$, then all conditions except possibly the last one can be attained in the case in which $s$ the empty tuple. If in addition to this the free product of $\F(x)*\subg{a}$ is non-abelian, then the last condition can be achieved for $s=\nil$ as well.
  \end{theorem}
  This result illustrates neatly what the obstruction to a full generalization of Merzlyakov's theorem for sentences with algebraically constrained universal variables is, but as it stands it is not precise enough to be of use for the ultimate purpose of analyzing all $\forall\exists$ sentences with parameters. In order to achieve this goal one needs to keep track of a greater amount on the information provided by the proof. This will be phrased in the language of \pr-groups morphisms and resolutions introduced in chapter \ref{MR chapter}.
  \begin{comment}
  Formal solutions in the sense of the theorem above can be seen as retractions of  extensions of $\qG{G}$, in the following sense.
  Suppose that a morphism $\kappa$ is given mapping $\qG{G}_{A}(x)$ to a \pr-limit group $\qG{H}_{A}(x,t,y)$, so that both systems
  $\syeq{\Sigma}{x,y,a}$ and $\syneq{\Pi}{x,y,a}$ hold in $\qG{H}$ as statements on the tuples $x,y,a$ of elements of $H$ and let $r=\pi(y)$. Assume that there is a \pr-retraction from
  $\qG{H}$ onto the \pr-subgroup generated by $x,t,a$ in $H$. Then, writing the image of the
  \end{comment}
  
  Let $\qG{G}=\qG{G}^{\Delta}$. We know that there is a finite family $\mathcal{MR}^{ws}$ of strict \ws \rs \pr-\res and for each $\RRA\in\mathcal{MR}^{ws}$ some surjective morphism $q^{R}$ from $\rqG{G}{A}$ to its top group $\rqG{L}{A}$, so that each $f\in Mor(\rqG{L}{A},\rqG{F}{A})$ can be written as $\bar{f}\circ q^{R}$, for some $\RRA\in\mathcal{MR}^{ws}$, and some $\bar{f}$ which factors through $\RRA$.
  To each $\LRL{R}{A}{(J,r)}\in\mathcal{MR}^{ws}$ we can associate a completion
  $\LRL{T}{A}{(J,r)}=Compl(\QR{R})$, as constructed in \ref{completions subsection} and a morphism which embeds $\RLR{R}{A}{r}$ into $\RLR{T}{A}{r}$.
  
  Let $\qG{H}$ be a free \pr-group such that $\subg{\pi(\rqG{G}{A}),\pi(H)}=Q$ and $t$ be a tuple of variables, such that for each of the completions $\rqG{T}{A}$ obtained above the group $T\frp H$ admits a generating $t$-marking. We can assume that $x\subset t$.
  
  We claim that the theorem above follows from the following proposition:
  %  Now, we can see $\qG{T}$ as endowed with a marking $x,t$ compatible with this map.
  \begin{proposition}
  	\label{the meat of Merzlyakov} Let $\RTA$ be a closed \rs \pr-tower structure for a group $\qG{T}_{A}(t)$, where $A$ is a free factor of $\F$.
  	%  	 and $\mathcal{TS}$ a congruency complete family of \pr-test sequences of $\RTA$ closed under subsequences and diagonal sequences.
  	%  	Suppose that for any $\ssq{f}{n}\in\mathcal{TS}$ the tuple $t_{n}=f_{n}(t)$
  	Suppose that any solution $t_{0}$ of $\qG{T}_{A}(t)$ satisfies:
  	\begin{align*}
  		\exists y\,\bvee{j=1}{(y\in q(j)\wedge\syeq{\Sigma_{j}}{t_{0},y,a}\wedge\syneq{\Pi_{j}}{t_{0},y,a})}{s}
  	\end{align*}
  	for a system of equalities $\syeq{\Sigma_{j}}{t,a}$ and a system of inequalities $\syeq{\Pi_{j}}{x,a}$ and tuples $q(j)\in Q^{|y|}$. Assume moreover that for all $1\leq j\leq s$ the tuple $q(j)$ is contained in $\pi(T)\leq Q$.
  	
  	Then there is a covering system of closures $\mathcal{CL}$ of $\RTA$ and for each
  	$\RSA\in\mathcal{CL}$ some retraction:
  	\begin{align*}
  		\fun{r^{S}}{\rqG{S}{A}\frp\F(y)}{\rqG{S}{A}\frp\qG{H}}
  	\end{align*}
  	and some $1\leq j_{S}\leq s$ such that $r^{S}(y)=q(j_{S})$ and $\syeq{\Sigma_{j_{S}}}{x,r^{S}(y),a}$.
  	
  	If in addition $\qG{T}$ is non-abelian then $r_{j}$  can be chosen to preserve the validity of
  	the system $\syneq{\Pi_{j_{S}}}{x,y,a}$.
  \end{proposition}
  %  Switch to the terminology 'solution' all over this thing.
  Indeed, we can always take a free \pr-group $\qG{H}$ mapping onto $Q$ and replace the completion of each of the resolutions in $\mathcal{MR}^{ws}$
  by its product with $\qG{H}$, which can be given a \pr-tower structure extending that of $T$ in an obvious way. Apply the proposition to them and take a tuple $s$ of variables such that each of the
  resulting closures admits a (generating) $s$ marking.
  Any solution of $\qG{G}_{A}(x)$ extends to one which factors through at least one of the completions above, and hence to one factoring through one of the
  closures \footnote{It is important that it factors through the \pr-tower, not any morphism from the \pr-tower group $\rqG{T}{A}$ to $\rqG{F}{A}$ extends to a closure.},
  a property which translates into property (\ref{extension}) of the theorem; each $\Delta_{l}$ appearing in the statement corresponds to the presentation of
  one of those closures and each $q(l)$ to the image in $Q$ of the associated tuple of generators by the corresponding $\pi$-map. If we let $w_{l}(s,a)$ be the expression of $r^{S}(y)$ in terms of $s,a$, items (\ref{equations satisfied}) and (\ref{retraction in Q}) follow immediately.
  
  For the last item, observe that since each of the closures $\rqG{S}{A}$ are themselves \rs \pr-limit groups, if the system of inequalities $\syneq{\Pi_{l}}{x,w(x,t,a),a}$ holds in $\rqG{S}{A}$,
  then necessarily $\syneq{\Pi_{l}}{x_{0},w(x_{0},t_{0},a),a}$ for some particular solution $(x_{0},t_{0})$ of $\rqG{S}{A}$, since the latter is a \rs \pr-limit group.
  %  DEFINE WHAT $A$ SOLUTION\,OF\,A\,GROUP\,IS\,
  %  refer to the group underlying a closure as the closure
  
  \subsection{Formal Makanin-Razborov diagrams: proving propostion \ref{the meat of Merzlyakov}}
    %    For the entirety of this subsections we assume all the hypothesis of \ref{the meat of Merzlyakov}
    
    For each $1\leq l\leq r$ and $q\in\pg{T}$ let $\qG{(K_{l}^{q})}_{A}$ be the quotient of $\rqG{T}{A}\frp\F(y)$ by the relation given by the system $\syeq{\Sigma_{l}}{t,y,a}$, where $\pi(y)=q(l)$.
    
    We can assume for any $l$ this assignment induces an extension of the map $\pi_{\hat{T}}$, since otherwise there are not any tuples $(t_{0},y_{0})$ in $\F$ satisfying both $\syeq{\Sigma_{l}}{t_{0},y_{0},a}$
    and $y\in q(l)$, and we can remove the corresponding term from the disjunction. Likewise, we can assume that $T$ injects in the quotient obtained this way.
    
    By a $(q,\Sigma_{l})$-formal \pr-group we intend a free factor containing $\hat{T}$ of a \pr-quotient of $\qG{(K^{q}_{l})}_{A}$
    onto which $T$ maps injectively and by.
    
    Fix some family $\mathcal{TS}$ of test sequences of $\QR{R}_{A}$. By a formal sequence of morphisms from $\qG{L}_{A}(z)$ a sequence $\ssq{g}{n}$ of morphisms from $\rqG{L}{A}$ to $\qG{F}$ such that $(g_{n}\rst_{T})_{n}\in\mathcal{TS}$.
    We refer to it as merely 'formal' if when not wanting to specify the pair $(q,\Sigma_{l})$.
    
    Given a \pr-\res $\LRL{R}{}{(J,r)}$ of a marked \rs \pr-group $\qG{G}_{A}(x)$ and a system of inequalities $\syneq{\Sigma}{x,a}$, with $a\subset A$, we say that $\QR{R}$ preserves
    the validity of $\syneq{\Sigma}{x,a}$ if the map from $\rqG{G}{A}$ to $\bfrp{\Rl}{\lambda\in\lvs{J}}{}$ induced by the quotient maps along the resolution does.
  %   \subsection{Proving proposition \ref{formal MR}}
    %    for which $\qG{T^{q}_{l}}_{A}$ is well-defined
    \begin{lemma}
    	\label{full free subgroups} For any \pr-limit group $\qG{L}$ there is some free subgroup $H\leq L$ such that
    	\begin{align*}
    		\pi(H)=\pi(G)
    	\end{align*}
    \end{lemma}
    \begin{proof}
    	We will use the following well-known fact:
    	\begin{claim}
    		Given an isometric action of a group $G$ on a real tree $Y$ and elements $g_{0},g_{1},\cdots g_{m}$ of $G$ acting hyperbolically on $Y$ and such that
    		$Ax(g_{i})\cap Ax(g_{i})$ has bounded diameter for $0\leq i<j<m$, there is a positive integers $N$ such that
    		for any choice of $k_{i}\geq N$ for $0\leq i\leq m$ the set $\{g_{i}^{k_{i}}\}_{0\leq i\leq m}$ generates a free group, of which it is a base.
    	\end{claim}
    	\begin{proof}(sketch)
    		The boundary of $Y$, denoted by $\partial Y$ is a Hausdorff space whose underlying set is the family of all infinite geodesic rays in $Y$, quotiented by the equivalence relationship
    		that identifies two such rays with infinite intersection. The action of $G$ on $\partial Y$ induces an action of $G$ on $\delta Y$. Any bi-infinite geodesic $\alpha$ is associated with
    		a pair of endpoints $\partial\alpha\subset\delta Y$. Given any hyperbolic element $h$ of $G$ the pair $\partial Ax(h)$ is of the form
    		$\{e^{+},e^{-}\}$, with the property that given any neighbourhoods $U^{+}$ and $U^{-}$ of $e^{+}$ and $e^{-}$ respectively there is some positive $N>0$ such that
    		$h^{k}\cdot(\delta Y\setminus U^{-})\subset U^{+}$ for any $k\geq N$ and $h^{k}\cdot(\delta Y\setminus U^{+})\subset U^{-}$ for any $k\leq-N$.
    		
    		Now, the fact that $Ax(h_{i})\cap Ax(h_{j})$ is bounded for $i\neq j$ implies that the intersection of $\partial Ax(h_{i})$ and $\partial Ax(h_{j})$ is empty.
    		
    		For each $0\leq i\leq m$ one can chose neighbourhoods $U^{+}_{i},U^{-}_{i}$ of the endpoints of $Ax(h_{i})$, as above, in a way that $\{U_{i}^{\pm}\}_{0\leq i\leq m}$
    		are mutually disjoint. If we let $A_{i}=U_{i}^{+}\cup U_{i}^{-}$
    		%  and $B_{i}=\partial Y\setminus A_{i}$,
    		then there is $N>0$ such that for any choice of $k_{j}\geq N$,
    		and $l\in\Z\setminus\{0\}$ we have $h_{i}^{k_{i}l}\cdot A_{j}\subset A_{i}$ for any $j\neq i, 0\leq j\leq m$. An application of the ping-pong lemma yields the desired conclusion.
    	\end{proof}
    	
    	Let $Y$ be the limiting tree for some sequence of morphisms from $\qG{L}$ to $\qG{F}$ with trivial limit kernel and let $\mathcal{X}$ be any finite set of generators of
    	$\qG{L}$. Perform the following two operations until none of them is possible:
    	\elenco{
    		\item Suppose $x\in\mathcal{X}$ is elliptic and let $h$ be hyperbolic. If
    		$xh$ is elliptic then $Fix(x)$ must intrsect the axis of $h$ (indeed, consider
    		the segment $I=[p,xh^{k}\cdot p]$ for some $p$ fixed by $xh$). Now, $Fix(x)$ is linear; since our action is irreducible, this implies the existence of some element $h$ hyperbolic in $Y$ such that $Fix(x)\cap Ax(y)=\nil$, so that $xh$ is hyperbolic. Up to replacing $h$ with a proper power we can assume $\pi(h)=1$. Replace $x$ by $xh$ for some such $h$.
    		\item If two elements $x,y\in\mathcal{X}$ commute, replace $\{x,y\}$ by some $z\in Z_{L}(x)$ such that $\pi(z)$  generates  $\pi(Z_{L}(x))$.
    	}
    	We end up with some set $\mathcal{X}\subset L$ of hyperbolic elements such that $\subg{\pi(\mathcal{X})}=\pi(L)$ and $[x,y]\neq 1$ for
    	any distinct $x,y\in\mathcal{X}$, which implies that $Ax(x)\cap Ax(y)$ has bounded diameter (see \ref{factoring lines}). On the other hand any element of $\qG{L}$ has
    	arbitrarily big powers with the same image in $Q$, so the sought for result follows from the fact above.
    \end{proof}
    \newcommand{\FMR}[2]{\mathcal{FMR}^{#1}_{#2}}
    \begin{proposition}
    	\label{formal MR} Suppose that $\mathcal{TS}$ is congruency complete and closed under subsequences and diagonal subsequences.
    	
    	Then for each $1\leq j\leq s$ and $q\in\pg{T}$ a set $\FMR{q}{l}$ of finitely many closed graded \pr-\res of
    	$\qG{(K_{l}^{q})}_{A}$ with respect to the parameter group $\hat{T}$ exists such that for any formal sequence $\ssq{g}{n}$ of morphisms from $\qG{(K_{l}^{q})}_{A}$ to $\rqG{F}{A}$
    	with the following property.
    	
    	Given any formal sequence $\ssq{g}{n}$, if $t_{n}=g_{n}(t)$ and $y_{n}=g_{n}(t)$ satisfy
    	\begin{align*}
    		y_{n}\in q(l)\wedge\syeq{\Sigma_{l}}{t_{n},y_{n},a}\wedge\syneq{\Pi_{l}}{t_{n},y_{n},a}
    	\end{align*}
    	%      	some $\QR{R}\in\mathcal{FMR}_{l}$ exists through which all the members of some infinite subsequence $(g_{n_{k}})_{k}$ factor.
    	for $n$ big enough, then all the members of some infinite subsequence of $\ssq{g}{n}$ factor through some \pr-\res in $\FMR{q}{l}$.
    	
    	Moreover, for any $\LRL{T}{}{(J,r)}\in\FMR{q}{l}$, for any $\lambda\in \lft{J}$ either
    	$\Rl$ is free or a closure of $\LRLii{T}{q}{A}$ and the morphism $\fun{\psi}{\qG{K_{l}^{q}}_{A}}{\bfrp{\Rl}{\lambda\in\lvs{J}}{}}$ induced by
    	the projections along $\QR{R}$ preserves the validity of $\syneq{\Pi_{l}}{t,y,a}$.
    \end{proposition}
    \begin{proof}(of \ref{the meat of Merzlyakov} from proposition \ref{formal MR})
    	%    	Here we have the additional assumption that $T$ covers $Q$.
    	Let $\LRL{R}{J}{}\in\FMR{q}{l}$ and $\lambda_{0}$ the leaf of $J$ carrying $\hat{T}$.
    	
    	Since $\pi(K_{l}^{q})=\pi(T)$ and $\Rl$ for is free for $\lambda\in\lvs{J}\setminus\{\lambda_{0}\}$,
    	there must be some \pr-retraction $\bar{r}:\qG{S}\frp(\bfrp{\Rl}{\lambda\in\lvs{J}\\\lambda\neq\lambda_{0}}{})\to\qG{S}$  sending
    	$\Rl$ into $T$ for $\lambda\in\lvs{J}$ into
    	If in addiiton $T$ is non-abelian, then in virtue of lemma \ref{full free subgroups} there is some $H\leq\hat{T}$ such that $\pi(H)=\pi(\hat{T})$;
    	it follows from
    	lemma \ref{free plimits}
    	%    	the theorem on the existence ot \pr-test sequences, by taking $\qG{H}$ as a reference \pr-free group and group of constants,
    	that there is some \pr-retraction from $\RLR{R}{A}{\lambda_{0}}\frp(\bfrp{\Rl}{\lambda\in\lvs{J}\\\lambda\neq\lambda_{0}}{})$ to $\RLR{R}{A}{\lambda_{0}}$ preserving the validity of $\syneq{\Psi}{t,\psi(y),a}$ and sending $\bfrp{\Rl}{\substack{\lambda\in\lvs{J}\setminus\{\lambda_{0}\}}}{}$ to $\qG{T}$.
    	
    	%    By \label{free plimit} there is some \pr-retraction from $\qG{H}\frp\bfrp{\lambda\in\lvs{J}\setminus\{\lambda_{0}\}}{\Rl}{}$
    	Now, $\rqG{U}{A}:=\RLR{R}{A}{\lambda_{0}}$ has the structure of a closure of $\LRLii{T}{q}{A}$, which is the result of adding a peg to the original tower \pr-\res $\RTA$.
    	Given $M\in\mathcal{FHM}(T)$, contained in $\tilde{M}\in\mathcal{FHM}(S)$, we can write $\tilde{M}$ as $\bar{M}\oplus Z'_{M}$, where $\bar{M}$ is the root closure of $Z \oplus M$.
    	
    	Of course $M$ has finite index in $\tilde{M}$, so using lemma \ref{blow and collapse} one can easily show that $U$ admits a star shaped graph of groups decomposition with a closure of $\rqG{S}{A}$ in its center
    	and for any $M\in\mathcal{FHM}(T)$, a vertex group $\tilde{M}$ amalgamated to it over $\bar{M}$.
    	
    	We know $\rqG{U}{A}$ is a limit of a sequence of morphisms which eventually send $p_{M}$ to a primitive element of $\F$, so in particular $\pi(\tilde{M})=\pi(M)=\subg{\pi(pg_{M})}$, so there is a
    	%    the construction of test sequences for amalgamated abelian subgroups we can find a family of subgroups
    	%    completeeeeeeeeeeeeeeeeeeeeeeeeeeeeeeeeeeee
    	a \pr-retraction from $\qG{U}$ to $\qG{S}$; existence of test sequences in case \ref{pegged abelian case} implies we can assume that the latter preserves the validity of the system $\syneq{\Psi}{t,\psi(y),a}$.
    	
    	Composing all the maps above we obtain
    	a homorphism $r^{\QR{R}}$ from $\rqG{T}{A}\frp\subg{y}$ to $\qG{S}$ restricting to the identity on $\qG{T}$. The fact that it can be written as $r'\circ h$, where $r'$ is a morphism from $(\qG{K}_{l}^{q})_{A}$ to $\rqG{S}{A}$ implies that $\pi(r^{\QR{R}}(y))=q(l)$, and $\syeq{\Sigma_{l}}{x,r^{S}(y),a}$.
    	
    	%    	Now, $\syneq{\Pi_{l}}{t,y,a}$ holds in $\pi(K_{l}^{q})$ and is preserved by $r'$, hence
    	The discussion above implies that the system $\syneq{\Pi_{l}}{t,r^{\QR{R}}(y),a}$ is also valid.
    	
    	To conclude the proof, notice that the fact that the collection of all the closures $\RSA$ associated to the resolutions in $\abunion{1\leq l\leq s\\q\in\pg{T}}{\FMR{q}{l}}{}$
    	constitute a covering system follows from the fact that $\mathcal{TS}$ is congruency complete and each
    	$\ssq{f}{n}\in\mathcal{TS}$ contains an infinite subsequence factoring through  $\QR{R}\in\FMR{q}{l}$ for some $1\leq l\leq s$ and $q\in\pg{T}$ which automatically implies it extends to the closure associated with $\QR{R}$.
    \end{proof}
    \newcommand{\pcl}[1]{$#1$-tight }    %tc
    \newcommand{\flp}[0]{(q,\Sigma_{l})}
    \newcommand{\apo}[0]{\sqsubset}
    \newcommand{\epo}[0]{\sqsubseteq}
    \newcommand{\napo}[0]{\nsqsubset}
    \newcommand{\nepo}[0]{\nsqsubseteq}
    \newcommand{\alsuc}[0]{\Subset}
    
    %    The rest of this subsection is dedicated to the proof of \ref{formal MR}.
    %    \begin{proof}
  	\subsection{Proof of proposition \ref{formal MR}}
    	\begin{definition}
    		%    	make this consistent with the definition of closures
    		Let $\rqG{L}{A}$ be a formal \rs \pr-limit group, $J':=nDFt(J^{*})$ and
    		$J_{0}\subset J'$. Let $\bar{J_{0}}$ stand for $J_{0}\cup\abunion{\lambda\in J_{0}}{Ch(\lambda)}{}$ and let $n(J_{0})$ stand for the uppermost node of $J'\setminus J_{0}$ with respect to the linear order $\lo$.
    		\newcommand{\pos}[0]{\mathcal{PO}}
    		\newcommand{\pair}[0]{$(J_{\lambda_{0}},\apo)$}
    		%    	Possibly change the ^{*} into ^{q}     I\,M\,P\,O\,R\,T\,A\,N\,T
    		
    		By a \pcl{J_{0}} structure for $\RRA$ we intend a collection:
    		\begin{align*}
    			((H_{\lambda})_{\lambda\in J_{0}\cup\{n(J_{0})\}},
    			(S^{\kappa})_{\lambda\in J_{0}},
    			(\gamma_{\lambda})_{\lambda\in \bar{J_{0}}\setminus J_{0}}
    			,(v_{\lambda})_{\lambda\in\bar{J_{0}}\setminus\{r\}})
    		\end{align*}
    		where $H_{\lambda}\leq L$, $\gamma_{\lambda}\in L$ and $S^{\lambda}$ is a simplicial $H_{\lambda}$-tree such that the following properties are satisfied, where
    		we extend the notation by letting $\gamma_{\lambda}=1$ for
    		$\lambda\in J_{0}$ and $K_{\lambda}:=(\RLR{T}{}{\lambda})^{\gamma_{\lambda}}$ for $\lambda\in\bar{J_{0}}$:
    		\begin{enumerate}[i)]
    			\item \label{technical property} For any $\mu\in\bar{J_{0}}$,
    			$K_{\mu}$ is contained in $H_{c(\mu)}$, where $c(\mu)\leq\lambda$ is the minimal node of $J_{0}$ into which $\RLR{T}{}{\mu}$ can be conjugated.
    			\item \label{vertex property} For $\lambda\in\bar{J_{0}}\setminus\{r\}$ the group $K_{\lambda}$ stabilizes a vertex $v_{\lambda}$ of $U^{c(\lambda)}$ of rigid type in case $c(\mu)$ is of floor type.
    			\item $H_{r}=L$ and $H_{\lambda}=\Sb{v_{\lambda}}{H_{p(\lambda)}}$ for $\lambda\in J_{0}\cup\{n(J_{0})\}\setminus\{r\}$
    			%    		\item $v_{\mu}=v_{\nu}$ in case $c(\mu)=c(\nu)$, $v_{\mu}\in H_{c(\mu)}\cdot v_{\nu}$ and $\nu\in J_{0}$.
    			\item \label{separation property} For distinct $\mu_{1},\mu_{2}\in J_{0}\cap Ch(\lambda)$ the vertices $v_{\mu_{1}}$ and $v_{\mu_{2}}$ are in distinct orbits by the action of $H_{\lambda}$.
    			%    		\item  $K_{\lambda}:=\RLR{T}{}{\lambda}^{\gamma_{\lambda}}$ for all $\lambda\in J_{0}\cup\{n(J_{0})\}$.
    			\item If $\lambda\in J_{0}$ falls in case \ref{free product case}, the simplicial $H_{\lambda}$-tree $S^{\lambda}$ has trivially stabilized edges. The minimal tree $U^{\lambda}$ of
    			$K_{\lambda}$ in $S^{\lambda}$ is equivariantly isomorphic to the tree dual to the free product decomposition $\qG{K}_{\lambda}=(\RLR{T}{}{\mu_{1}})^{\gamma_{\lambda}}\frp (\RLR{T}{}{\mu_{2}})^{\gamma_{\lambda}}$ where $\{\mu_{1},\mu_{2}\}=Ch(\lambda)$ and its translates cover $S^{\lambda}$.
    			\item For any $\lambda\in J_{0}$ is of floor type:
    			\enum{i)}{
    				\item  The tree dual to $\DR{T}{\lambda}{}$, on which $K_{\lambda}$ acts via conjugation by $\gamma_{\lambda}^{-1}$, is equivariantly isomorphic to some subtree $U^{\lambda}$ of $S^{\lambda}$
    				and the translates of its image cover $S^{\lambda}$.
    				\item  	The stabilizer in $H_{\lambda}$ of the non-rigid type vertices of $U^{\lambda}$ is contained in $K_{\lambda}$ in case \ref{surface case}. In case \ref{pegged abelian case}
    				given a vertex $v\in U^{\lambda}$ if we let $P$ be the peripheral subgroup of $\Sb{v}{H_{\lambda}}$ then $\Sb{v}{K_{\lambda}}+P$ has finite index in
    				$\Sb{v}{H_{\lambda}}$.
    			}
    			\item If $\lambda\in J_{0}$ falls in case \ref{free group case}, then $K_{\lambda}=H_{\lambda}$
    		\end{enumerate}
    	\end{definition}
    	\begin{observation}
    		\label{first observation} For $\lambda\in\bar{J_{0}}\setminus J_{0}$ the minimality of $c(\lambda)$ implies that $v_{\lambda}$ is never a translate of $v_{\mu}$ for any $\mu\in Ch(c(\lambda))\cap J_{0}$.
    		In particular $n(J_{0})\lo c(\lambda)$ for any $\lambda\in\bar{J_{0}}\setminus J_{0}$, since otherwise $Ch(c(\lambda))\subset J_{0}$. This implies that $n(J_{0})\succ c(n(J_{0}))$.
    	\end{observation}
    	\begin{comment}
    	Also that if $\rqG{L}{A}$ admits a \pcl{J_{0}} structure as above we can always assume that in the definition above:
    	\elenco{
    		\item For $\mu_{1},\mu_{2}\in\bar{J_{0}}$ if $c_{\mu_{1}}=c_{\mu_{2}}$ and $v_{\mu_{1}}$ and $v_{\mu_{2}}$ are in the same orbit then $v_{\mu_{1}}=v_{\mu_{2}}$.
    		\end{comment}
    		\begin{observation}
    			For a formal \pr-limit group $\rqG{L}{A}$ admitting $\pcl{J'}$ structure means precisely admitting the structure of a closure of $\LRLii{T}{q}{A}$.
    		\end{observation}
    		\begin{proof}
    			All that is left to prove is the existence of a suitable retraction in case \ref{pegged abelian case}.
    			%    			For this it is enough that the centralizer $\bar{M}$ in $\rqG{L}{A}$ of any maximal non-cyclic abelian group $M$ of $\hat{T}$ has the same image in $Q$ as it has in $M$. This follows right away from the fact that $L$ is the limit quotient of some formal sequence, which in virtue of
    			This follows right away from lemma \ref{primitive pegs} will eventually send each peg of $\hat{T}$ to
    			a primitive element.
    		\end{proof}
    		Arguing in exactly the same way as it was done in the construction of Makanin-Razborov diagrams, proposition \ref{formal MR} can be easily reduced to the following:
    		\begin{lemma}
    			\label{adding a node} Let $\qG{L}_{A}(z)$ a formal limit group admitting a \pcl{J_{0}} structure and $\syneq{\Psi}{z,a}$ a system of inequalities true in $\qG{L}_{A}(x)$, where $A=\subg{a}$. Suppose that $J_{0}\neq J'$.
    			%    	 $\lambda_{0}$ is not minimal
    			%    	and let $\lambda_{1}$ be the maximum element of $J'\setminus J_{0}$ with respect to $\lo$
    			Then there is a finite collection $\mathcal{FMR}^{q}_{\lambda,L}$ of closed $T$-graded resolutions of $\rqG{L}{A}$ such that from any formal sequence $\ssq{g}{n}$ of $\rqG{L}{A}$ preserving the validity of $\syneq{\Psi}{z,a}$ one can extract a subsequence factoring through  one of $\LRL{R}{A}{(I,r)}\in \mathcal{FMR}^{q}_{\lambda_{0},L}$. Moreover, for any $\LRL{R}{A}{J}\in\mathcal{FMR}^{q}_{\lambda,L}$:
    			\elenco{
    				\item  For $\LRL{R}{A}{J}\in\mathcal{FMR}_{\lambda,L}$, the \pr-limit group at a leaf of $J$ carrying $T$ admits a \pr-\res which is \pcl{J_{0}\cup\{n(J_{0})\}}.
    				\item  $\QR{R}$ preserves the validity of the system $\syneq{\Psi}{z,a}$
    			}
    		\end{lemma}
    		
    		\subsubsection{Proof of \ref{adding a node}: uncovering a node through shortening}
      		
      		Take a \pcl{J_{0}} structure $(\ssq{H}{\lambda},\ssq{\gamma}{\mu},v_{\mu},S^{\nu})_{\lambda,\mu,\nu}$ for $\rqG{L}{A}$. Notice that the way the lineaar order $\lo$ was defined, $n(J_{0})$ is either $r$ or the the child of a node in $J_{0}$. The result can be proved by yet another application of Rips' and Sela's shortening argument.
      		\begin{comment}
      		used in the construction of Makanin-Razborov diagrams to this context, using the limiting 	action
      		on a real tree associated to the sequence of restrictions to $H_{\lambda_{1}}$ of a formal sequence $\ssq{g}{n}$.
      		\end{comment}
      		
      		Let $\mathcal{FMA}$ be the family of all conjugates in $L$ of non-cyclic maximal abelian subgroups of $\hat{T}$ and $\mathcal{BR}:=\setof{\mu\in\bar{J_{0}}}{v_{\mu}=v_{n(J_{0})}}$.
      		Given any $M\in\mathcal{FMA}$, let $\bar{M}=Z_{L}(M)$.
      		
      		Given $\lambda\in \bar{J_{0}}$ denote by $\bar{K}_{\lambda}$
      		be the subgroup generated by the union of $K_{\lambda}$ and $\bar{M}\cap Stab(v_{\lambda})$ for those $M\in\mathcal{FMA}$ such that $M\cap K_{\lambda}\neq\tg$.
      		%      	%      This requires some way guaranteeing the satisfaction by the subgroup $H_{\lambda_{1}}$ appearing at any given stage of the process of a certain indecomposibility condition.
      		%      	Given $\lambda\in\bar{J_{0}}\setminus\{r\}$ let $\bar{K}_{\lambda}$ be the group generated by:
      		%      	\begin{align*}
      			%      	 K_{\lambda}\bigcup \setof{\bar{M}\cap H}{M\in\mathcal{FMA}\,,\,M\cap K_{\lambda}\neq\tg}
      		%      	\end{align*}
      		%   Denote by $\bar{\mathcal{A}}$ the family $\{\bar{K}_{n(J_{0})}\}\cup\{K_{\mu}\}_{\substack{c(\mu)=n(J_{0})\\\mu\nin Ch(\lambda)}}$ and let $\bar{B}=K_{n(J_{0})}$.
      		The following lemma will be useful. It consists essentially in an iterative application of lemma \ref{blow and collapse}.
      		\begin{lemma}
      			\label{tower lifting trees} Let $W$ be a simplicial $H$-tree $\mathcal{O}$ a minimal invariant family of vertices of $W$ such that
      			each member of the family $\mathcal{A}:=\{\bar{K}_{n(J_{0})}\}\cup\{\bar{K}_{\mu}\,|\,v_{\mu}=v_{n(J_{0})}\}$ fixes a vertex of $\mathcal{O}$.
      			
      			Then $W$ equivariantly embeds	into a simplicial $\hat{H}$-tree $W'$ with the following properties:
      			\enum{a)}{
      				\item The pegged tower group $\hat{T}$ fixes a vertex $v_{T}$ in the image of $\mathcal{O}$ such that $Stab(v_{T})$ admits a \pcl{J_{0}} structure. \label{prop-1}
      				\item The translates of the image of $W$ in $W'$ cover $W$. \label{property0}
      				\item For any edge $e$ in $W$ the stabilizer of $e$ (in $H$) is the same as the stabilizer of $\phi(e)$ in $H_{\lambda_{j}}$. \label{property1}
      				\item If two edges $e,e'\in W$ are in distinct orbits, then the same is true for $\phi(e),\phi(e')$. \label{property2}
      				\item Suppose we are given an equivariant family $\mathcal{O}$ of vertices of $W$ such that each member of $\mathcal{A}$ fixes some $u\in\mathcal{O}$. \label{property3}
      				Then for any vertex $v\in W\setminus\mathcal{O}$ we have:
      				\begin{align*}
      					\Sb{v}{H}=\Sb{\phi(v)}{H_{\lambda_{j}}}
      				\end{align*}
      				\item  Given any two vertices $v,v'\nin\mathcal{O}$ in distinct orbits their images $\phi(v),\phi(v')$ are also in distinct and do not belong to the orbit of any
      				vertex in $\mathcal{O}$. \label{property4}
      			}
      			%      		Moreover, if $W$ is a \gat relative to $\bar{\mathcal{A}}$, then $W$ can be seen as a \gat relative to $\hat{T}$ so that the map $\phi$ preservs the vertex type.
      		\end{lemma}
      		
      		\begin{proof}
      			Let us say that $n(J_{0})=\lambda_{0}\prec\lambda_{1}\prec\cdots\lambda_{m}=r$. Let $W$ be the tree dual to the decompositon $\Delta$.
      			For $1\leq j\leq m$ we construct a $H_{\lambda_{j}}$-simplicial tree $W_{j}$ and an equivariant embedding $\phi_{j}$ of $W$ into $W_{j}$. Let
      			$\mathcal{O}_{j}$  the union of translates of $\phi_{j}(\mathcal{O})$.
      			Additionally, we require the following properties to hold:
      			\elenco{
      				\item  Each subgroup in the family
      				%      			The $H_{\lambda_{j}}$-tree $W_{j}$ is relative to the family
      				\begin{align*}
      					\mathcal{A}_{j}:=\{\bar{K}_{\lambda_{j}}\}\cup\{Stab(w)\,|\,w\in VS^{\lambda_{j}}\setminus(H^{\lambda_{j}}\cdot v_{\lambda_{j-1}}),\,\}
      					%      				\cup\{\bar{K}_{\lambda}\}_{v_{\lambda}=v_{\lambda_{j}}}
      					%      				\cup\{\bar{K}_{\lambda}\}_{v_{\lambda}=v_{\lambda_{j}}}
      					%      				\mathcal{A}_{j}:=\setof{\bar{K}_{\mu}}{\mu=}
      				\end{align*}
      				fixes some vertex in the orbit $\mathcal{O}_{j}$.
      				\item The triple $(W,\phi_{j},W_{j})$ satisfies properties (\ref{property1})-(\ref{property4}) of the statement of the lemma. .
      			}
      			We proceed by induction. Suppose that $W_{l}$ has already been constructed for $0<l<j$. Our goal is to construct $W_{j}$ as an application of \ref{blow and collapse}	to the trees $S^{\lambda_{j}}$ and 	$W_{j-1}$ with $\bar{K}_{\lambda_{j}}$ in the role of $G$.
      			%      		Since $\lambda_{j}\in J_{0}$, in fact it is easy to describe $\bar{K}_{\lambda_{j}}$.
      			Note that there are three possibilities for how some $M\in\mathcal{FMA}$  which intersects $K_{\lambda_{j}}$ non-trivially and its centralizer $\bar{M}$ can be placed with respect to $S^{\lambda_{j}}$.
      			\elenco{
      				\item The first one is that $\bar{M}\cap H_{\lambda_{j}}$ is hyperbolic in $U^{\lambda_{j}}$. Acylindricity of the tree $U^{\lambda_{j}}$ do not allow any non-trivial element elliptic in $U^{\lambda_{j}}$
      				to commute with another one which is hyperbolic. This means the group $M\cap K_{\lambda_{j}}$ is hyperbolic in $U^{\lambda_{j}}$ as well and
      				\begin{align*}
      					\bar{M}\cap H_{\lambda_{j}}=M\cap K_{\lambda_{j}}=\subg{pg_{M}}
      				\end{align*}
      				%      			prove the implication
      				\item Secondly, it is possible for $\bar{M}\cap H_{\lambda_{j}}$ to fix a rigid type vertex	of $U^{\lambda_{j}}$.
      				\item Finally, it might be that $\lambda_{j}$ is of type \ref{pegged abelian case} and $M\cap K_{\lambda_{j}}$ is the stabilizer of an abelian type vertex of $U^{\lambda_{j}}$.
      				Of course, this can be the case for at most one $M$.
      			}
      			Thus $\bar{K}_{\lambda_{j}}$ is generated by $K_{\lambda_{j}}$ together with the union $\mathcal{E}$ of $\bar{M}\cap H_{\lambda_{j}}$ for all the $M$ as in the last two bullets.
      			%      		 of subgroups of $H_{\lambda_{j}}$ elliptic in $U^{\lambda_{j}}$.
      			Let $\bar{U}^{\lambda_{j}}$ be the tree spanned by the translates of $U^{\lambda_{j}}$ by $\bar{K}_{\lambda_{j}}$ (it clearly contains the minimal tree of $\bar{K}_{\lambda_{j}}$ in $S^{\lambda_{j}}$).
      			
      			A first consequence of the property above is that the map
      			$\bar{q}:\bar{K}_{\lambda_{j}}\backslash \bar{U}^{\lambda_{j}}\to H_{\lambda_{j}}\backslash S^{\lambda_{j}}$
      			is injective on edges. The analogous quotient $q:K_{\lambda_{j}}\backslash U^{\lambda_{j}}\to H_{\lambda_{j}}\backslash S^{\lambda_{j}}$ certainly is, so all we need to show is that the quotient
      			$q':K_{\lambda_{j}}\backslash U^{\lambda_{j}}\to \bar{K}_{\lambda_{j}}\backslash U^{\lambda_{j}}$ is surjective. This follows from lemma \ref{elliptic expansion} as $\bar{K}_{\lambda_{j}}$
      			is generated by $K_{\lambda_{j}}$ and elements of $H_{\lambda_{j}}$ fixing vertices of $U^{\lambda_{j}}$. (Since $q'$ is injective on edges, we also know from \ref{elliptic expansion} is $q'$ an isomorphism.)
      			
      			Lemma \ref{elliptic expansion} implies as well that for any $v\in U^{\lambda_{j}}$ the stabilizer of $v$ in $\bar{K}_{\lambda_{j}}$ is generated by $\Sb{v}{K_{\lambda_{j}}}\cup (\mathcal{E}\cap \Sb{v}{H_{\lambda_{j}}})$. Observe any $M\in\mathcal{FMA}$ can intersect $\Sb{v}{K_{\lambda_{j}}}$ non-trivially for at most one rigid type vertex $v$ of $U^{\lambda_{j}}$.
      			Also that if $w$ is a vertex of abelian type of $U^{\lambda_{j}}$ and $\Sb{w}{H_{\lambda_{j}}}$ intersects $\Sb{v}{H_{\lambda_{j}}}$ non-trivially for some $v\in U^{\lambda_{j}}$ of rigid type,
      			then $\Sb{w}{H_{\lambda_{j}}}=\bar{M}\cap H_{\lambda_{j}}$ for some $M\in\mathcal{FMA}$ which intersects $\Sb{v}{K_{\lambda_{j}}}$ non-trivially.
      			It follows that $\Sb{v_{\lambda_{j-1}}}{\bar{K}_{\lambda_{j}}}=\bar{K}_{\lambda_{j-1}}$.
      			
      			We claim that $\Sb{v_{\mu}}{\bar{K}_{\lambda_{j}}}$ is conjugate into $\bar{K}_{\mu}$ in $H_{\lambda_{j}}$ for any $\mu\in Ch(\lambda_{j})\setminus J_{0}$ as well.
      			%      		 there is $h\in H_{\lambda_{j}}$ such that $K_{\mu}=\Sb{v}{K_{\lambda_{j}}}^{h}$. We
      			. Indeed, take $\nu\leq\lambda_{j}$, $\nu\in J_{0}\cup n(J_{0})$ such that $c(\mu)=\nu$.
      			\begin{comment}
      			For any $\nu\leq\kappa\leq\lambda_{j}$ the group $K_{\mu}$
      			stabilizes a rigid vertex $w_{\kappa}$ of $U^{\kappa}$. We can prove by induction that $\bar{K}_{\mu}$ stabilizes $w_{\kappa}$ for any such $\kappa$.
      			%      			 is conjugate into a rigid vertex $w_{\kappa}$ of $U^{\kappa}$.
      			\end{comment}
      			If that was not the case, then for some $\nu\leq\kappa\leq\lambda_{j}$
      			there should be conjugate $M,M'\in\mathcal{FMA}$ such that $K_{\mu}\subset M$ is non-trivial and
      			$M'$ intersects the stabilizer of some abelian type vertex of $S^{\kappa}$ non-trivially. However, this cannot be the case, since
      			it implies the existence of $M_{0},M'_{0}\in\mathcal{MA}(\hat{T})$
      			such that $pg_{M_{0}}\in \RLR{T}{}{\mu}$ and $pg_{M}\in \RLR{T}{}{\lambda_{j}}$. Now, $\mu\lo\lambda_{j}$
      			and $\mu$ and $\lambda_{j}$ are not $\leq$ comparable, so for any $n$ big enough $\gamma_{n}\in\F$ exists
      			such that the images $f_{n}(pg_{M_{0}})$ and $f_{n}(pg_{M'_{0}})^{\gamma_{n}}$ by some test sequence $\ssq{f}{n}$ commute.
      			Since for $n$ big enough these two are primitive elements of different length, this is impossible.
      			
      			Consider any translate $v=h\cdot v_{\lambda_{j-1}}\in S^{\lambda_{j}}_{\bar{K}_{\lambda_{j}}}$ for $h\in H_{\lambda_{j}}$. For any edge $e\in \bar{U}^{\lambda_{j}}$ incident at $v$
      			clearly $\Sb{e}{H_{\lambda_{j}}}\leq \bar{K}_{\lambda_{j}}$.
      			%      		 In fact, unless $\lambda_{j}$ falls in case \ref{pegged abelian case} $\Sb{e}{H_{\lambda_{j}}}\leq K_{\lambda_{j}}$.
      			On the other hand, as a consequence of the last discussion $\Sb{v}{\bar{K}_{\lambda_{j}}}^{h}$ is a conjugate in $H_{\lambda_{j}}$ of either $\bar{K}_{\lambda_{j-1}}$ or
      			$\bar{K}_{\mu}$ for some $\mu\in Ch(\lambda_{l})\setminus J_{0}$, where $0< k<j-1$. Now, we know $\bar{K}_{\lambda_{j-1}}$ to fix a vertex of $\mathcal{O}_{j-1}$, while any such $\bar{K}_{\mu}$ is fixes one in $\mathcal{O}_{l}$ for some $0\leq l<j$ (by the initial assumption in case $l=0$ and by the induction hypothesis otherwise), hence one in $\mathcal{O}_{j-1}$ as well.
      			
      			Trees $S^{\lambda_{j}}$ and $W_{j-1}$ thus fall under the hypothesis of lemma \ref{blow and collapse}, with $\bar{K}_{\lambda_{j}}$ in the role of $G$. The lemma produces a tree, that we take as
      			our $W_{j}$, and an equivariant embedding $\psi_{j-1}$ of $W_{j-1}$ into $W_{j}$. We defne $\phi_{j}$ as $\psi_{j-1}\circ\phi_{j-1}$.
      			
      			The construction (in view of the previous paragraph) guarantees the ellipticity of the family $\{\bar{K}_{\lambda_{j}}\}\cup\{Stab(w)\,|\,w\in VS^{\lambda_{j}}\setminus (H_{\lambda_{j}}\cdot v_{\lambda_{j-1}})\}$. Since translates of $\phi_{j-1}(W)$ cover $W_{j-1}$ and translates of $\psi_{j-1}(W_{j-1})$ cover $W_{j}$ clearly translates of $\phi_{j}(W)$ cover $W_{j}$.
      			Likewise, lemma \ref{blow and collapse} guarantees that $\psi_{j-1}(VW_{j-1}\setminus\mathcal{O}_{j-1})\subset VW_{j}\setminus\mathcal{O}_{j}$, and $\phi_{j-1}(VW\setminus\mathcal{O})\subset VW_{j-1}\setminus\mathcal{O}_{j-1}$ by induction hypothesis, so clearly $\phi_{j}(VW\setminus\mathcal{O})\subset VW_{j}\setminus\mathcal{O}_{j}$. Transitivity can be easily checked for the remaining properties in order to complete the proof.
      			%      			 Likewise, for $v\in VW\setminus\mathcal{O}$ we have $\phi_{j-1}(v)\nin\mathcal{O}_{j}$
      			%      		and $Stab(v)=Stab(\phi_{j-1}(v))=Stab(\psi_{j-1}(\phi_{j-1}(v)))$, the first equality deriving from the induction hypothesis and the second from the lemma.
      			
      			%      			Let $R_{\lambda_{j-1}}$ the full stabilizer of the vertex of $W_{\lambda_{j}}$ fixed by $K_{\lambda_{j-1}}$ and $P_{\lambda_{j}}$ the tree spanned by all translates of $U^{\lambda}$ by elements of $R_{\lambda}$.
      			
      			Pick some $v^{T}_{0}\in\mathcal{O}$ fixed by $K_{n(J_{0})}$. Let $v^{T}_{j}:=\phi_{j}(v^{T}_{j})$ and for $1\leq j\leq m$ let $R_{\lambda_{j}}$ the full stabilizer of $v^{T}_{j}$ and by $P_{\lambda_{j}}$ the tree spanned by all translates of $U^{\lambda}$ by elements of $R_{\lambda}$. We claim that the stabilizer of $v_{T}=v^{T}_{j}$ admits a \pcl{J_{0}} structure.
      			Let us see how to extend this data to a \pcl{J_{0}} structure
      			\begin{align*}
      				((R_{\lambda})_{\lambda\in J_{0}\cup\{n(J_{0})\}},
      				(P^{\lambda})_{\lambda\in J_{0}},
      				(\gamma'_{\lambda})_{\lambda\in \bar{J_{0}}\setminus J_{0}}
      				,(v'_{\lambda})_{\lambda\in\bar{J_{0}}\setminus\{r\}})
      			\end{align*}
      			%      			 $((R_{\lambda})_{\lambda},(P^{\lambda})_{\lambda},(\gamma'_{\lambda})_{\lambda},(v'_{\lambda})_{\lambda})$,
      			where $v'_{\lambda}=v_{\lambda}$ in case $\lambda\in J_{0}$.
      			For $1\leq j\leq m$, consider the two graph maps induced by inclusion:
      			$\fun{q_{1}}{K_{\lambda_{j}}\backslash U^{\lambda_{j}}}{R_{\lambda_{j}}\backslash P_{\lambda_{j}}}$ and
      			$\fun{q_{2}}{R_{\lambda_{j}}\backslash P_{\lambda_{j}}}{H_{\lambda_{j}}\backslash S^{\lambda_{j}}}$.
      			
      			Lemma \ref{blow and collapse} tells us that $q_{2}$ is injective on edges and that the only vertex in its image admitting several preimages is $v_{\lambda_{j-1}}$.  $n(J_{0})\leq\lambda_{j-1}\prec\lambda_{j}$. On the other hand $q_{2}\circ q_{1}$ is bijective on edges by definition. It follows that $q_{1}$ must be surjective as well,
      			so clearly the tree $P_{\lambda}$ satisfies the requirements of the definition.
      			%      So the action of $H'_{\lambda}$ on $S^{'\lambda=}:=U^{\lambda}_{H'_{\lambda}}$ satisfies the conditions in the definition of a \pcl{J_{0}} structure.
      			
      			Consider now any $\mu\in Ch(\lambda_{j})\setminus\{\lambda_{j-1}\}$. If $\mu\in J_{0}$ there is no need to modify the given data for $\nu\leq\mu$, as $H_{\mu}\subset R_{\lambda}$.
      			Suppose now that $\mu\nin J_{0}$, let $u_{\mu}$ be the vertex stabilized by $\RLR{T}{}{\mu}$ and $c'(\mu)$ the minimum $\kappa\in J_{0}$ such that
      			$K'_{\mu}:=(\RLR{T}{}{\mu})^{\gamma'_{\lambda}}\subset R_{\lambda_{j}}$ for some $\gamma'_{\lambda}\in R_{\lambda_{j}}$. If $c(\mu)=\lambda_{j}$, then we can take $v'_{\mu}=u_{\mu}$.
      			One can show that $\gamma_{\mu}(\gamma'_{\mu})^{-1}\in H_{\kappa}$ \footnote{By repeatedly using the general fact that given $\lambda\in J_{0}$ an element $h\in H_{\lambda}$ conjugates some non-trivial element in the stabilizer of a (rigid) vertex of $S^{\lambda}$ to another one, then $h$ itself must belong to the stabilizer.}, so that $K'_{\mu}$ must fix a rigid vertex $w$
      			of $S^{\mu}$. If $w\nin P_{\mu}$, since $K'_{\mu}$ is contained in $R_{\mu}$, it must fix the vertex in $P_{\mu}$ closest to $w$. All vertices in the boundary of $P_{\mu}$ in $S^{\mu}$ are rigid, so an adequate $v'_{\mu}$ exists.
      		\end{proof}
      		\begin{corollary}
      			\label{modular extension} With the notation of the lemma, if $W$ is endowed with a geometric abelian marking relative to $\{\bar{K}_{n(J_{0})}\}\cup\{\bar{K}_{\mu}\,|\,v_{\mu}=v_{n(J_{0})}\}$,
      			then the $L$-tree $W'$ can be regarded as a \gat in such a way  $\phi$ respects vertex type. Every $\sigma\in Mod^{\pi}_{K_{n(J_{0})}}(W)$ extends to some
      			$\tau'\in Mod^{\pi}_{\hat{T}}(W')$.
      		\end{corollary}
      		\begin{proof}
      			In this case we can take as $\mathcal{O}$ any minimal invariant family of rigid vertices such that each member of $\mathcal{A}$ stabilizes one of them.  Any vertex $v\in W$ not belonging to $\mathcal{O}':=L\cdot\mathcal{O}$ in $W'$ is a translate of
      			$\phi(w)$ for some vertex $w\in W\setminus\mathcal{O}$, whose orbit is uniquely determined. The action of $Stab(w)$ on the star around $v$ is isomorphic to the action of
      			$Stab(w)$ on the star around $w$, so one can assign to $v$ the type of $w$ (obviously well-defined).
      			
      			Notice that $\mathcal{O}'$ consists of a unique orbit, as any $v\in\mathcal{O}$ has to become an attaching point at some stage in the iterative construction.
      			%      		There is
      			%      		a subtree $Z_{0}$ of $W$, containing a vertex fixed by $K_{n(J_{0})}$, such that
      			%      		\BS{\pi()}{S}{t} that its union $Z$ with some edges originating at $Z_{0}$
      			The properties of the embedding $\phi$ imply\footnote{Take a vertex stabilized by $\hat{T}$ in $\phi(\mathcal{O})$. It is easy to find a Bass-Serre presentation $\BS{Z_{0}}{Z}{t}$ of $W'$ with $w\in Z_{0},Z\subset\phi(W)$, by adding edges to the subtree in a greedy fashion until a transversal is found.} that a \gad $\Delta'$ of $L$ associated to $W$ can be obtained (in a way compatible with the inclusion map) from a \gad $\Delta$ of $L$ by means of:
      			\enum{i)}{
      				\item Enlarging $\Delta_{u}$ for some rigid vertex $u$ for which $K_{n(J_{0})}\leq\Delta_{u}$.
      				\item Removing some of the other rigid vertices and attaching all edges incident to them to $u$ instead. This eequires adding Bass-Serre elements for some of those edges, which have to drop from the maximal tree.
      			}
      			For any of the removed vertices $w$, the group $\Delta_{u}$ is contained in $(\Delta'_{u})^{s_{u}}$, where $s_{u}$ is one of the newly introduced Bass-Serre elements.
      			Given any $\tau\in PMod(\Delta)$ in the set of generators given in section \ref{trees and automorphisms section} which fixing $\Delta_{u}$ (either a Dehn twist or an extension of a vertex group automorphisim) 	one can easily extend $\tau$ to some $\tau'\in Mod^{\pi}(W')$ fixing $\Delta'_{u}$ and sending any of the $s_{u}$ above to $s_{u}c_{u}$, whenever  $\tau$ restricts to conjugation to $c_{u}$ on $\Delta_{u}$.
      		\end{proof}
      		\begin{corollary}
      			\label{shedding free products} Let $\rqG{L}{A}$ be a \pcl{J_{0}} formal \pr-limit group.
      			Suppose that there is a non-trivial free decomposition of $H_{\lambda_{n(J_{0})}}$ relative to
      			$\{K_{n(J_{0})}\}\cup\setof{K_{\mu}}{v_{\mu}=v_{n(J_{0})}}$. Then there is a free decomposition of $L=L_{0}\frp L$
      			such that $\hat{T}\leq L_{0}$ and $\qG{L_{0}}$  admits a \pcl{J_{0}} structure as well.
      		\end{corollary}
      		\begin{proof}
      			Notice that any free splitting of $H_{n(J_{0})}$ relative to $\{K_{n(J_{0})}\}\cup\{K_{\mu}\}_{\substack{c(\mu)=n(J_{0})}}$ is also relative to
      			$\{\bar{K}_{n(J_{0})}\}\cup\{\bar{K}_{\mu}\}_{v_{\mu}=v_{n(J_{0})}}$, (if an abelian group intersects a free factor non-trivially then it is contained in it) so the previous lemma can be applied to the tree $W$ dual to the said decomposition.
      		\end{proof}
      		
      		Using the corollary one can reduce the theorem to the case in which $H:=H_{n(J_{0})}$ is freely indecomposable relative to $\{K_{\mu}\,|\,v_{\mu}=v_{n(J_{0})}\}\cup\{K_{n(J_{0})}\}$.
      		Suppose it is not; then we know there is a free decomposition of $L$ of the form $L=L_{0}\frp L_{1}$, where $\hat{T}\leq L_{0}$.
      		Chose a new tuple of variables $u=(u_{0},u_{1})$ such that each $L_{i}$ is endowed with a surjective $u_{i}$-marking.
      		
      		By collecting the elements of each $L_{i}$ appearing in the normal form of the words in $\Theta$ one can find \NEQs $\syneq{\Theta_{0}}{u^{0},a}$ and $\syneq{\Theta_{1}}{u^{1}}$ such 	that given any \rs homomorphism $q_{0}$ with domain $L_{0}$
      		and $q_{1}$ with domain $L_{1}$, if each $q_{j}$ preserves the validity of $\Theta_{j}\neq 1$ then the homomorphism from $L$ to $(\bfrp{f_{j}(L_{j})}{j=0}{1})$ induced by them preserves $\syneq{\Psi}{z,a}$.
      		\begin{comment}
      		By the same argument used in the construction of the regular Makanin-Razborov diagrams, existence of the sought resolution is clearly implied
      		by the induction hypothesis on $\qG{L_{1}}$ and $\syeq{\Theta}{u^{0},a}$ and the following lemma:
      		\end{comment}
      		It is a routine task to derive the following result given \ref{morphism shortening} and the argument following result, just as in \ref{MR}.
      		%      	 Its proof uses the same techniques of the current proof, exempt from the additional complications that having to deal with formal \pr-limit groups entails.
      		%      	 We will indicate how to modify the Makanin-Razborov procedure at the end of the subsection.
      		\begin{lemma}
      			\label{inequality MR} Given a \pr-limit group $\qG{L}(z)$ and a system of inequalities $\syneq{\Psi}{z,a}$
      			there is a finite family of closed resolutions $\mathcal{MR}_{\Psi}$ preserving the validity of the system such that any morphism from $\qG{L}$ to $\qG{F}$
      			preserving the validity of $\syneq{\Psi}{z,a}$ factors through at least one of $\QR{R}\in\mathcal{MR}_{\Psi}$
      		\end{lemma}
      		Let $\mathcal{MR}_{1}$ be the set of \pr-\ress which results from applying lemma \ref{inequality MR} to the pair $(\qG{L_{1}}(u^{1}),\syneq{\Theta_{1}}{u^{1}})$. And $\mathcal{MR}_{0}$ that obtained by applying the induction
      		hypothesis to $(\qG{L_{0}},\syneq{\Theta_{0}}{u^{0},a})$. Let $\mathcal{MR}^{nd}$ be the result of combining each resolution in the first set with each one in the second in the obvious way into a resolution of $\qG{L}$. Each $\QR{R}\in\mathcal{MR}^{nd}$ preserves the validity of $\syneq{\Psi}{z,a}$. And each formal sequence whose members eventually preserve
      		the validity of $\syneq{\Theta_{0}}{u^{0},a}$ and $\syneq{\Theta_{1}}{u^{1}}$ contains an infinite subsequence factoring through one of them.
      		Now, take any element $v\in\Theta_{0}\cup\Theta_{1}$ and consider the family $\mathcal{FS}_{v}$ of all formal sequences of $L$ which kill $v$ and preserve $\syneq{\Psi}{z,a}$. If $\mathcal{FS}_{v}=\nil$, then just ignore $v$. If not, since this family is clearly closed under diagonal subsequences, every $\ssq{g}{n}\in\mathcal{FS}_{v}$ eventually factors (in a strict sense)
      		through one of a collection in a finite collection $\mathcal{Q}_{v}$ of maximal limit quotients of sequences in $\mathcal{FS}_{v}$. For each $\rqG{Q}{A}\in\mathcal{Q}_{v}$ let $\mathcal{MR}^{Q}$
      		the finite collection of \pr-\ress given by applying the induction hypothesis to the pair $(\rqG{Q}{A},\syneq{\Psi}{z,a})$. We can complete each one of them into a \pr-\res of
      		$\rqG{L}{A}$ in an obvious way, by adding the quotient form $\rqG{L}{A}$ onto $\rqG{Q}{A}$ on top. Let $\mathcal{MR}^{d}$ be the collection of all the resolutions so obtained, by letting $v$ range in $\Theta_{0}\cup\Theta_{1}$. The family $\mathcal{MR}^{nd}\cup\mathcal{MR}^{d}$ satisfies all the required conditions.
      		
      		\newcommand{\pM}[0]{\mathfrak{M}}
      		Let $\Delta_{JSJ}$ be a JSJ decomposition of $H_{n(J_{0})}$ relative to the family $\{\bar{K}_{\mu}\}_{v_{\mu}=v_{n(J_{0})}}\cup\{\bar{K}_{n(J_{0})}\}$
      		and let $\pM$ be the group consisting of all the extensions $\tilde{\tau}$ fixing $\hat{T}$, where $\tau$ ranges among all elements of $Mod^{\pi}_{K_{n(J_{0})}}(\Delta_{JSJ})$.
      		Fix $w$ a tuple of generators for $H$. Consider the action of $\F$ on $Cayl(\F)$. Given a morphism $\fun{f}{\qG{L}_{A}(z)}{\rqG{F}{A}}$ and $\lambda\in J$, let $*_{B}^{f\rst_{H}}\in Cayl(\F)$ the basepoint chosen as described in \ref{basepoints}. We say that $f$ is short among those preserving the validity of $\syneq{\Psi}{z,a}$ if for any $\sigma\in\pM$ for which $f\circ\sigma$ preserves the validity of $\syneq{\Psi}{z,a}$ the inequality $\sl{*_{B}^{f\rst_{H}}}{w}{}\leq\sl{*_{B}^{f\rst_{H}}}{\sigma(w)}{}$ holds. By a short formal sequence $\ssq{g}{n}$ for $L$ we intend a formal sequence such that $g_{n}$ satisfies this condition for any $n$.
      		%      	defined affasdfasdfasdfassdaSFASDDADSFASDFSDFSDFASDFADFADFASDFASFASDFASF
      		%      	We now claim that the proof reduces to the following lemma:
      		\begin{lemma}
      			\label{Merzlyakov shortening analysis} Let $\qG{L}$ be a formal \pr-limit group which is freely indecomposable relative to $\mathcal{K}=\{K_{\mu}\}_{v_{\mu}=v_{n(J_{0})}}\cup\{K_{n(J_{0})}\}$ and admits a \pcl{J_{0}} structure and let $\ssq{g}{n}$ be a short formal sequence consisting entirely of short morphisms in the sense above from $\rqG{L}{A}$ to $\rqG{F}{A}$. Then one of the following takes place:
      			\enum{a)}{
      				\item  $\limker{g_{n}}{n}\neq\tg$
      				\item \label{abelian line}There is a decomposition of $H_{n(J_{0})}$ of the form $H_{n(J_{0})}=H'_{n(J_{0})}\frp_{N\cap H'_{n(J_{0})}}N$, where $K_{n(J_{0})}\leq H'$,
      				$N=Z_{H_{n(J_{0})}}(M)$ for some $M\in\mathcal{FMA}$	and $N=N\cap H'\oplus E$ for some $E$.
      				\item $\rqG{L}{A}$ admits a \pcl{J_{0}\cup n(J_{0})} structure.
      			}
      		\end{lemma}
      		\begin{proof}
      			%      		[of \ref{Merzlyakov shortening analysis}]
      			Let $J_{1}=J_{0}\cup\{n(J_{0})\}$. Suppose that $\limker{g_{n}}{n}=1$. Let $\rho$ be the action of $H$ on the limiting tree for the sequence $(g_{n}\rst_{H})_{n}$.
      			%      				In virtue of corollary \ref{morphism shortening} we know that $B$ cannot fix a point of $Y$.
      			The sequence of points $*^{g_{n}}_{B}$ necessarily converges to a point $*$ fixed by $B$ or in the minimal tree of $B$ in the limiting tree $Y$ associated to the restriction of the sequence
      			$\ssq{g}{n}$ to $H_{n(J_{0})}$.
      			
      			Suppose first that $K_{n(J_{0})}$ acts elliptically on $Y$.
      			
      			Lemma \ref{morphism shortening} implies that some $\bar{K}_{\mu}\in\mathcal{A}$ is not elliptic in $Y$. In case $\mu\neq n(J_{0})$ the properties of test sequences however imply that
      			$K_{\mu}$ does fix some $x_{\mu}\in Y$. We conclude there is some $M\in\mathcal{FMA}$ intersecting $K_{\mu}$ non-trivially such that $N:=Z_{M}(H_{n(J_{0})})$ does
      			not fix $x_{\mu}$. This implies it does not fix a point $y\neq x_{\mu}$ either.
      			If not, then notice that in this case
      			$M$ must also fix $y$, so for some $g\in N$ either the union $[x_{\mu},y]\cup g\cdot[x_{\mu},y]\cup g^{2}\cdot[x_{\mu},y]$ is a non-degenerate non-trivially stabilized tripod or $g$ acts as a symmetry on a non-degenerate segment $[x_{\mu},y]\cup g\cdot[x_{\mu},y]$ fixed, both of which cannot occur.
      			%      			This implies that $N\cap H_{n(J_{0})}$ is root closed in $N$.
      			
      			So $N$ acts hyperbolically on $Y$, while $M$ and $K_{n(J_{0})}$ act elliptically, and by \ref{factoring lines} there is a decomposition of $H_{n(J_{0})}$ as an amalgamated product of the required form.
      			%      		form $H'\frp_{N\cap Stab(x_{\mu}})\frp N$.
      			%      		Notice that for any $\mu\in\bar{J_{0}}$ with $v_{\mu}=v_{n(J_{0})}$
      			%      		or maybe an easily quotable lemma
      			
      			We are left with the case in which $K_{n(J_{0})}$ does not act elliptically on $Y$.
      			If the node $\lambda$ carries the constants then $\RLR{T}{}{\lambda}=\qG{A}$, so this is impossible.
      			In all other cases corollary \ref{nice decompositions} implies some $H_{n(J_{0})}$-tree $S^{n(J_{0})}$ exists with the desired properties.
      			
      			All is left is to check the existence of appropriate $v'_{\lambda}$ and $\gamma'_{\lambda}$ for $\lambda\in\bar{J_{1}}\setminus J_{1}$. For $\lambda\in J_{1}$ we don't effect any
      			change. Consider first some $\lambda\in\bar{J_{0}}\setminus\{n(J_{0})\}$ such that $c(\lambda)$ is the parent of $n(J_{0})$.
      			By assumption either $v_{\lambda}$ is a translate of $v_{n(J_{0})}$ or $v_{\lambda}$ and $v_{n(J_{0})}$ are not in the same orbit. In the latter case the parent $\mu$ of $n(J_{0})$ is remains the lowest node of $J_{1}$ into which $\RLR{T}{}{\lambda}$ can be conjugated.
      			In the former one we can replace $\gamma_{\lambda}$ by $\gamma'_{\lambda}$ in such a way that the resulting $K'_{\lambda}$ is contained in $Stab(v_{n(J_{0})})$. Now, $K'_{\lambda}$ is elliptic in the real tree $Y$, so  \ref{nice decompositions} implies that it stabilizes a vertex (of rigid type when it applies) of $S^{\lambda}$ which we can take as our new $v'_{\lambda}$.
      			The choice of $v_{\lambda}$ for $\lambda\in\bar{J_{0}}\setminus\bar{J_{1}}=Ch(n(J_{0}))$ is clear.
      			All we need to prove in order to verify (\ref{separation property}) is that $v_{n(J_{0})}$ is not a translate of $v_{\nu}$ for any $\nu\in Ch(\mu)\cap J_{0}$, which follows from observation \ref{first observation}.
      		\end{proof}
      		
      		Equipped with this result, let us finish the proof of \ref{formal MR}.  Assume first that $H_{n(J_{0})}$ admits a decomposition as in \ref{abelian line} of the previous lemma.
      		\begin{claim}
      			Such a decomposition $H'\frp_{N\cap H'}N$ inducess one of the form $L=L'\frp_{\bar{M}\cap L'}\bar{M}$, where $H'\leq M$, $\bar{M}=Z_{L}(M)$ and $L'$ admits a \pcl{J_{0}} structure.
      		\end{claim}
      		\begin{proof}
      			The proof is a small variation from that of \ref{tower lifting trees}.
      			%      		We start with the $H_{n(J_{0})}$-tree $W$ dual to the given decompostition and for each $1\leq j\leq m$ we iteratively construct an $H_{\lambda_{j}}$ tree $W_{\lambda_{j}}$.
      			We start with the given decomposition, $\Delta_{0}$, and for $1\leq j\leq m$ we construct one $\Delta_{j}$ of the form $H_{\lambda_{j}}=H'_{\lambda_{j}}\frp_{\bar{M}\cap H'_{\lambda_{j}}}(\bar{M}\cap H_{\lambda_{j}})$. These are the different possibilities for $0\leq j\leq m$:
      			\elenco{
      				\item $\lambda_{j}$ falls in case \ref{pegged abelian case} so that there is a decomposition
      				\begin{align*}
      					H_{\lambda_{j}}=H_{\lambda_{j-1}}\frp_{\bar{M}\cap H_{\lambda_{j-1}}}((\bar{M}\cap H_{\lambda_{l}})\oplus F)
      				\end{align*}
      				%      			$\bar{M}_{l}=\bar{M}\cap H_{\lambda_{l}}$
      				%      			$\bar{M}_{l}=\bar{M}\cap H_{\lambda_{l}}$
      				this means we can write $H_{\lambda_{j}}$ as an amalgamated product $H_{\lambda_{j-1}}\frp_{H_{\lambda_{j-1}}\cap\bar{M}}((H_{\lambda_{j-1}}\cap\bar{M})\oplus E\oplus F)$ and we can take
      				$H'_{\lambda_{j}}:=H'_{\lambda_{j-1}}\oplus F$. Clearly the action of $H'_{\lambda_{j}}$ on the corresponding subtree of
      				\item If $\lambda$ falls in any other case, we can use lemma \ref{blow and collapse} to lift the amalgamated product from $H_{\lambda_{j-1}}$ to $H'_{\lambda_{j}}$.
      				%      			and to prove that the action
      				%      			of the vertex group $H'_{\lambda_{j}}$ of the resulting decomposition acts on the tree spanned by $H'_{\lambda_{j}}\cdot U_{\lambda_{j}}$ the right way.
      			}
      			At this point one can check that $H'_{\lambda_{j}}$ appear in place of $H_{\lambda_{j}}$ in some \pcl{J_{0}} structure for $L'=H'_{\lambda_{m}}$ in an entirely analogous way as in the proof of
      			\ref{tower lifting trees}.
      			%      		The same argument as the one used in \ref{tower lifting trees} how to find a \pcl{J_{0}} structure with $H'_{\lambda_{j}}$ in place of $H_{\lambda_{j}}$.
      		\end{proof}
      		In this case the usual argument provides some \pr-retraction $\fun{r}{\qG{L}}{\qG{L'}}$ preserving the validity of the sytem $\syneq{\Psi}{z,a}$. For any formal sequence $\ssq{g}{n}$ of $\qG{L}$, in virtue of \ref{primitive pegs} eventually $g_{n}(N)=g_{n}(N\cap L')$, and $g_{n}$ can be written in the form $g'_{n}\circ r\circ\sigma$, where $\sigma$ is a \pr-modular automorphism of the decomposition $L'\frp_{N\cap L'}N$ (taking $N$ as an abelian type vertex group) fixing $L'$.
      		We are done by extending each member of the family of \pr-\ress provided by the induction hypothesis for the pair to $(\qG{L'},\syneq{\Psi}{z,a})$ in the obvious way.
      		
      		In the other case, observe that since the family of short formal sequences of $\qG{L}$ is closed under taking diagonal sequences, the collection of limit quotients by
      		short formal sequences contains finitely many maximal elements $\mathcal{MFQ}$.
      		
      		By \ref{Merzlyakov shortening analysis} we might as well assume that they are proper.
      		Given any formal sequence $\ssq{g}{n}$ preserving the validity of $\syneq{\Psi}{z,a}$, there if for all $n$ we take  $\tau_{n}\in\pM$ such that $g'_{n}=g_{n}\circ\tau_{n}$ is short then the members of any
      		convergent subsequence of $\ssq{g'}{n}$ must factor through $\phi\in\mathcal{MFQ}$ after a certain point.
      		
      		The sought family of resolutions are obtained in the usual way by combining the quotient $\fun{p}{\qG{L}}{\qG{L'}}\in\mathcal{MFQ}$ with any resolution resulting from applying the induction hypothesis to $(\qG{Q},\syneq{\Psi}{z,a})$.
      		%      	As the \gad at the root we can take simply the lift in the sense of \ref{tower lifting trees} of a JSJ decomposition of $H_{n(J_{0})}$ relative to $K_{n(J_{0})}$.
      	%      	\end{proof}
      	
	\section{The positive theory}
  	\begin{comment}
  	\newcommand{\aeae}[0]{\exists x^{0}\forall x^{1}\,\exists y^{1}\cdots\forall x^{m}\,\exists y^{m}}
  	\newcommand{\baeae}[0]{\exists x^{0}\forall x^{1}\in p^{1}\,\exists y^{1}\in q^{1}\cdots\forall x^{m}\in p^{m}\,\exists y^{m}\in q^{m}}
  	\newcommand{\abaeae}[0]{\exists v^{1}\forall u^{1}\in r^{1}\,\exists v^{1}\in s^{1}\cdots\forall u^{m}\in r^{m}\,\exists v^{m}\in s^{m}}
  	\end{comment}
  	\newcommand{\aeae}[0]{\forall x^{1}\,\exists y^{1}\cdots\forall x^{m}\,\exists y^{m}}
  	\newcommand{\baeae}[0]{\forall x^{1}\in p^{1}\,\exists y^{1}\in q^{1}\cdots\forall x^{m}\in p^{m}\,\exists y^{m}\in q^{m}}
  	\newcommand{\abaeae}[0]{\forall u^{1}\in r^{1}\,\exists v^{1}\in s^{1}\cdots\forall u^{m}\in r^{m}\,\exists v^{m}\in s^{m}}
  	\newcommand{\disu}[1]{
  	\bvee{i=1}{((x\in p(i)\wedge y\in q(i)\wedge\syeq{#1^{i}}{x,y,a})}{k}}
  	\newcommand{\atom}[1]{(x,y)\in q\wedge\syeq{#1}{x,y,a}}
  	As before, we add constants for the members of some finite tuple $a\leq\F$ generating some free factor $A$ of $\F$.
  	If we let $(\BTQ)_{A}$ stand for the theory of all $A$-\rs \pr-groups, it is clear
  	how any positive $\LQ_{A}$-sentence is $(\BTQ)_{A}$-equivalent to one of the following form:
  	\begin{align*}
  		\aeae \,\,\disu{\Sigma}
  	\end{align*}
  	Here $p(i)\in Q^{|x|}$ and $q(i)\in Q^{|y|}$,  $x=(x^{1},x^{2},\cdots x^{m})$ and $y=(y^{1},y^{2},\cdots y^{m})$.
  	
  	As a matter of fact, it is more convenient to work with a different class of formulas, which we will refer to as simple constrained positive (SCP) formulas, namely those of the form:
  	\begin{align*}
  		\baeae \,\,\syeq{\Sigma}{x,y,a}
  	\end{align*}
  	For some system $\syeq{\Sigma}{x,y,a}$ of equations with parameters in $a$. Here the following abbreviations have been used, for a tuple $x$ of variables and $q\in Q^{|x|}$:
  	\begin{align*}
  		(\forall x\in q\,\,\phi)\equiv \forall x\,\,(x\in q\sra\phi) \\
  		(\exists x\in q\,\,\phi)\equiv \exists x\,\,(x\in q\wedge\phi)
  	\end{align*}
  	Given an SCP sentence $\phi$ as above, by a formal solution of $\phi$ we intend a tuple $(w^{1},w^{2}\cdots w^{m})$
  	with the following properties:
  	\begin{comment}
  	\elenco{
  		\item $w^{j}$ is a $|y_{j}|$-tuple of $|y_{j}|$ words in $x^{1},y^{1}\cdots x^{j},a$ for $1\leq j\leq m$
  		\item $w^{j}_{l}(p^{1},q^{1}\cdots p^{j-1},\pi(a))=q^{j}$ for $1\leq j\leq m$
  		\item For any word $u(x,y,a)$ in the system $\syeq{\Sigma}{x,y,a}$ the expression
  		\begin{align*}
  			u(x,w^{1}(x_{1},a),w^{2}(x_{1},y_{1},x_{1},a),\cdots w^{m}(x_{1},y_{1},x_{2}\cdots x^{m},a))
  		\end{align*}
  		is trivial as an element of the group $A\frp\F(x,y)$.
  	}
  	\end{comment}
  	\elenco{
  		\item $w^{j}$ is a $|y_{j}|$-tuple of words in $x^{1},\cdots x^{j},a$ for $1\leq j\leq m$
  		\item $w^{j}_{l}(p^{1},p^{2},\cdots p^{j},\pi(a))=q^{j}_{l}$ for $1\leq j\leq m$ and $1\leq l\leq |y|$
  		\item For any word $u(x,y,a)$ in the system $\syeq{\Sigma}{x,y,a}$ the term:
  		\begin{align*}
  			u(x,w^{1}(x^{1},a),w^{2}(x^{1},x^{2},a),\cdots w^{m}(x^{1},x^{2}\cdots x^{m},a),a)
  		\end{align*}
  		represents the trivial element in the free group $A\frp\F(x)$.
  	}
  	Clearly, if an SCP formula admits a formal solution, then it is valid in any $A$-\rs \pr-group. So, in particular it is valid in $\rqG{F}{A}$. Later we will prove a strong converse of this result with implications for general positive formulas, under the assumption that $A$ is a free factor.
  	\begin{comment}
  	We shall show a stronger converse, namely:
  	\begin{proposition}
  		\label{formal positive} Let $\psi$ be a positive sentence of the form:
  		\begin{align*}
  			\baeae \syeq{\Sigma}{x,y,a}
  		\end{align*}
  		Suppose that $\phi$ is valid in $\rqG{F}{A}$. Then there are finitely many things
  	\end{proposition}
  	\end{comment}
  	
  	In order to start, we need to go back to the language of \pr-groups. Associated with any
  	SCP formula as above there is a \rs \pr-group $\qG{G_{\phi}}_{A}$, whose underlying group is the quotient of
  	$A\frp\F(x,y)$ by the normal subgroup generated by the words in $\syeq{\Sigma}{x,y,a}$ and $\pi=\pi_{G_{\phi}}$ maps each tuple $x^{j}$ (which we now see as elements of $G_{\phi}$) to the tuple $p^{j}$ and $y^{j}$ to $q^{j}$.
  	For any $1\leq j\leq m$ we let $G_{\phi}^{j}$ be the subgroup of $G_{\phi}$ generated by the tuple $(a,x_{1},y_{1},\cdots x_{j})$
  	and $H_{\phi}^{j}$ be the one generated by $(a,x_{1},y_{1},\cdots x_{j},y_{j})$.
  	We say that a $\phi$ is free on the on the universal variables if for any $0\leq l\leq m-1$ the subgroup
  	$\subg{G_{\phi}^{l},x_{l+1},\cdots x_{m}}$ is isomorphic to the free product of $G_{\phi}^{l}$ and $\F(x_{l+1},\cdots x_{m})$. In particular,
  	$\subg{A,x}\cong A\frp\F(x)$.
  	Viceversa, of course, any finitely presented \rs \pr-group $\rqG{G}{A}$ endowed with a surjective $(x,y)$-marking can be seen as a $G_{\phi}$, by taking as $\syeq{\Sigma}{x,y,a}$ the collection of relators in some presentation, expressed in terms of $x$ and $y$.
  	\newcommand{\GP}[0]{G_{\phi}}
  	\newcommand{\HP}[0]{H_{\phi}}
  	\newcommand{\PGP}[0]{\qG{G_{\phi}}}
  	\newcommand{\PHP}[0]{\qG{H_{\phi}}}
  	\begin{observation}
  		A SCP $\phi$ free on the universal variables admits a formal solution if and only if there is a \pr-retraction
  		$\fun{f}{G_{\phi}}{A\frp\F(x)}$ such that $f(H_{\phi}^{j})\subset A\frp\F(x_{1},\cdots x_{j})$ for any $1\leq j\leq m$.
  	\end{observation}
  	\begin{corollary}
  		\label{quotient positive reduction}Suppose that we are given SCP formulas:
  		\begin{align*}
  			\phi\equiv\baeae \syeq{\Sigma}{x,y,a} 	\\
  			\psi\equiv\abaeae \syeq{\Pi}{u,v,a}
  		\end{align*}
  		which are free in the universal variables and that some surjective \rs morphism $f$ from $\qG{G_{\phi}}_{A}$ to $\qG{G_{\psi}}_{A}$ exists such that
  		$f(H_{\phi}^{j})\subset H_{\psi}^{j}$ for all $1\leq j\leq m$. Then the existence of a formal solution for $\psi$ implies that of one for $\phi$.
  	\end{corollary}
  	\begin{comment}
  	The following is a trivial observation:
  	Define the partial order $\prec_{var}$ among SCP formulas as follows.
  	We say that
  	\begin{align*}
  		\baeae\syeq{\Sigma}{x,y,a}\equiv\phi\prec\psi
  	\end{align*}
  	if
  	$\psi$ can be
  	The following is clear:
  	\begin{observation}
  		If $\phi\prec_{var}\psi$, then $\prove\phi\sra\psi$.
  	\end{observation}
  	Notice also:
  	\begin{observation}
  		The order $\prec_{var}$ is well-founded.
  	\end{observation}
  	\end{comment}
  	\begin{definition}
  		Given an SCP formula:
  		\begin{align*}
  			\phi\equiv\baeae \syeq{\Sigma}{x,y,a}
  		\end{align*}
  		by a $\phi$-formal sequence we intend a metrically convergent sequence of morphisms
  		\begin{align*}
  			\fun{\ssq{f}{n}}{\qG{G_{\phi}}_{A}}{\rqG{F}{A}}
  		\end{align*}
  		such that:
  		\enum{i)}{
  			\item The restriction of $\ssq{f}{n}$ to $\subg{x^{j}}$ is a \pr-test sequence for the trivial \pr-tower with  $(\subg{x^{j}},\pi_{G_{\phi}}\rst_{\subg{x^{j}}})$ at its single node (i.e. a small cancellation sequence).
  			\item The group $\subg{x^{j}}$ grows faster than $H_{\psi}^{l}$ for $1\leq l<j\leq m$.
  		}
  	\end{definition}
  	
  	\begin{proposition}
  		\label{testing formulas}Suppose we are given an SCP formula
  		\begin{align*}
  			\phi\equiv \baeae\,\, \syeq{\Sigma}{x,y,a}
  		\end{align*}
  		where $A$ is a free factor of $\F$ and $\subg{A,p^{1},p^{2},\cdots p^{m}}=Q$.
  		Then $\phi$ admits a formal solution if and only if some $\phi$-formal sequence exists.
  	\end{proposition}
  	\begin{proof}
  		The only if direction is clear.
  		%  		 we just need to precompose with a $\phi$-formal solution extends any sequence
  		%  		of morphisms from $(\F(x),\pi\rst_{\F(x)})$ with the desired properties (this clearly exists, by \ref{asymmetricla growth}).
  		We prove the opposite direction by induction on the pair $(m,G_{\phi})$ with respect to the partial order
  		$\leq\times\leq_{rk}\times\leq_{Z}$. Case $m=1$ is a particular case of the main result of the previous subsection.
  		
  		To begin with, let us discard the case in which the $\phi$-formal sequence $\ssq{f}{n}$ has non-trivial limit kernel. In this situation let $\rqG{H}{A}$ be its limit quotient. Since $H$ is finitely presented
  		$\rqG{H}{A}$ equals $\qG{G_{\phi'}}_{A}$ for some other SCP formula $\phi'$ and some subsequence
  		of $\ssq{f}{n}$ pushes forward to a $\phi'$-formal sequence. By induction $\phi'$ admits a formal solutions and therefore so does $\phi$.
  		
  		We can also assume that $G_{\phi}$ is freely indecomposable relative to $G_{\phi}^{m}$. Indeed, otherwise we can write $G_{\phi}=L_{1}\frp L_{2}$, where $L_{2}\neq\tg$ and $G_{\phi}^{m}\leq L_{1}$.
  		Take any morphism $h$ from $\qG{L}$ to $\F$ and let $\iota$ be some morphism from $\qG{F}$ to $\subg{A,x^{1},x^{2},\cdots x^{m}}$. This exists due to the fact that the latter group maps onto $Q$.
  		The map $\iota\circ h$ extends to a \pr-retraction from $\qG{G_{\phi}}$ onto $\qG{L_{1}}$.
  		
  		So assume now that $\ssq{f}{n}$ has trivial limit kernel and $G_{\phi}$ is freely indecomposable relative to $G_{\phi}^{j}$. And consider the limiting tree $Y$ for the sequence $\ssq{f}{n}$.
  		
  		Let $\lambda$ be the action of $\F$ on its Cayley graph.
  		For any homomorphism $\fun{f}{G_{\phi}}{\F}$ not killing $\subg{a,x^{1}}$ (there is no harm in assuming none of the $f_{n}$ does),
  		\newcommand{\bpf}[0]{\frp^{\lambda_{f}}_{G_{\phi}^{m}}}
  		let	$\bpf$ be the basepoint associated to $G_{\phi}^{m}$ in the minimal tree of $G_{\phi}$ in $X$ with respect to the action $\lambda_{f}$, chosen as in section \ref{basepoints section} and
  		$\mu(f)=(\sl{\\bpf}{y^{m}}{\lambda_{f}})_{j=1}^{m}$. We can assume that
  		$\mu(f_{n})\leq\mu(f')$ for any $\fun{f'}{G_{\phi}}{\F}$ which coincides with $f_{n}$ on $G_{\phi}^{m}$ (the condition for being a $\phi$-formal sequence depends only on $f_{n}\rst_{G_{\phi}^{m}}$).
  		
  		\renewcommand{\baeae}[0]{\forall x^{1}\in p^{1}\,\exists y^{1}\in q^{1}\cdots\forall x^{m}\in p^{m}\,\exists y^{m}\in q^{m}}
  		We claim that this implies $\subg{x^{m}}$ is not elliptic in
  		$Y$. Indeed, if it was, then $G_{\phi}^{m}$ would fix the limit of the sequence $(\bpf)_{n}$ and there would thus be an
  		automorphism $\sigma\in Mod^{\pi}_{G_{\phi}^{j}}(G)$ such that for $n$ large enough $\mu(f_{n}\circ\sigma)\leq\mu(f_{n})$,
  		contradicting the choice of $\ssq{f}{n}$.
  		
  		By \ref{extra properties} we know then that the minimal tree $Y_{0}$ of $\subg{x^{m}}$ in $Y$ is fundamental, while $H_{\phi}^{m-1}$ fixes some point $z$ in $Y$.
  		The free indecomposibility hypothesis implies that $z$ can be chosen inside $Y_{0}$. Just as in the proof of Merzlyakov theorem, using 	observation \ref{free skeleton} one can see that then any intersection of $Y_{0}$ and one of its translates must belong to the orbit of $z$ and $G_{\phi}=G'\frp\subg{x^{m}}$, where $G'=Stab(z)$.
  		Let $u$ be a finite tuple of variables marking elements of $G'$ such that
  		$G'$ is generated by $(x^{1},y^{1},\cdots y^{m-1},u)$. We can see $\rqG{G'}{A}$ as  $\qG{G_{\psi}}$ for some
  		SCP formula:
  		\begin{align*}
  			\psi\equiv\forall x^{1}\in p^{1}\,\,\exists y^{1}\in q^{1}\cdots x^{m-1}\,\,\exists y^{m-1}\in q^{m-1}\exists u\in r\,\,\syeq{\Pi}{x^{1},\cdots x^{m-1},y^{1},\cdots y^{m-1},u,a}
  		\end{align*}
  		with the variables in $\psi$ in agreement with the marking of $G_{\phi}$. Now, clearly $(f_{n}\rst_{G'})_{n}$ is a $\psi$-formal test sequence, so by the induction hypothesis $\phi'$ admits a formal solution.
  		$\psi$ admits a formal solution. This implies $\phi$ does, by \ref{quotient positive reduction} (after adding innermost universally quantified variables to $\psi$).
  	\end{proof}
  	
  	\renewcommand{\disu}[1]{\bvee{i=1}{(x\in p(i)\wedge y\in q(i)\wedge\syeq{#1^{i}}{x,y,a})}{k}}
  	\renewcommand{\baeae}[0]{\forall x^{1}\in p^{1}\,\exists y^{1}\in q^{1}\cdots\forall x^{m}\in p^{m}\,\exists y^{m}\in q^{m}}
  	Let us now go back to our general positive formula:
  	\begin{align*}
  		\phi\equiv\aeae \disu{\Sigma}
  	\end{align*}
  	For any
  	%  	$(p,q)=(p^{1}\cdots p^{1},q^{1},\cdots q^{m})\in Q^{|(x,y)|}$
  	$p=(p^{1}\cdots p^{m})\in Q^{|x|}$
  	denote by $\mathcal{D}_{p}$ the collection of all the SCP formulas with free universal variables of the form:
  	\begin{align*}
  		\baeae\,\,\syeq{\Sigma^{i}}{x,y,a}
  	\end{align*}
  	for some $1\leq i\leq k$ such that $p(i)=p$.
  	%  	We can construct formal $\phi$ easily
  	%  	The construction of \pr-test sequences, together with an application of the pidgeon hole principle yields:
  	\begin{lemma}
  		Suppose that $\rqG{F}{A}\modelof\phi$. Then for each $p\in Q^{|x|}$ there is some sentence $\psi\in\mathcal{D}_{p}$ admitting a formal $\psi$-sequence
  		(in particular, $\mathcal{D}_{p}$ is not empty).
  	\end{lemma}
  	\begin{proof}
  		Indeed, given any such $p$, the construction of test sequences implies the existence of a $\phi$-formal sequence $\ssq{f}{n}$ such that $\pi(f_{n}(x))=p$ for all $n\in\N$. The sought result follows from the validity of $\phi$, together with a straightforward application of the pidgeon hole principle.
  	\end{proof}
  	
  	Now, let $\phi^{surj}$ the axiom stating that each of the predicates $P_{q}$ is non-empty. Modulo the theory $T^{\pi}_{gp}\cup\phi^{surj}$, any SCP sentence $\phi$ is equivalent to the one obtained adding any dummy constrained quantifier expression $\forall x'\in p'$ just before the atomic part of $\psi$, where the variables of $x'$ are disjoint from those appearing in $\psi$. In particular, we can always assume that the condition $\subg{\pi(x,a)}=Q$ is satisfied.
  	Using the previous lemma, the last observation and proposition \ref{testing formulas} we can deduce:
  	\begin{corollary}
  		If in the situation above $\rqG{F}{A}\modelof\phi$, then there is some finite collection $(\psi_{i})_{i=1}^{r}$ of SCP formulas, each of which admits a formal solutions, such that
  		$\BTQ\prove(\bwedge{i=1}{\psi_{i}}{r}\sra\phi)$.
  	\end{corollary}
  	From this we one can readily deduce:
  	\begin{theorem}
  		\label{equivalent} Let $A$ be a free factor of non-abelian free groups $F_{1}$ and $F_{2}$ and for $i=1,2$ let $\pi_{i}$ be a homomorphism from $F_{i}$ to the finite group $Q$.	Then
  		\begin{align*}
  			Th^{+}_{A}(F_{1},\pi_{1})=Th^{+}_{A}(F_{2},\pi_{2})
  		\end{align*}
  		%  		In particular, $Th^{+}(F,\pi)$ depends only on $Q=Im(\pi)$.
  		%  		In particular, if the free factor $A$ is non-abelian and maps onto $Q$, then $(A,\pi\rst_{A})$ is a positive elementary sub-model of $(F_{i},\pi_{i})$.
  	\end{theorem}
  	\begin{proof}
  		Take any positive $A$-sentence $\phi$ such that  $(F_{i},\pi_{i})\modelof \phi$. By the previous corollary, $\phi$ is implied
  		in $\BTQ$ by some conjunction $\bwedge{i=1}{\psi_{i}}{r}$ of SCP formulas admitting a formal solution. This  implies the validity of each $\psi_{i}$, and therefore that of $\phi$, in $(F_{3-i},\pi_{3-i})$.
  		%  		Obviously, the positive theory does distinguish between abelian and non-abelian free groups. Notice that in the proof of \ref{testing formulas} non-abelianity of the reference model $\F$ is required to construct test sequence.
  	\end{proof}

    \bibliographystyle{alpha}
    \bibliography{thesis}{}

\newcommand{\etalchar}[1]{$^{#1}$}
\begin{thebibliography}{DGH05}

\bibitem[Bes02]{bestvina2002r}
Mladen Bestvina.
\newblock R-trees in topology, geometry, and group theory.
\newblock {\em Handbook of geometric topology}, pages 55--91, 2002.

\bibitem[BF95]{bestvina1995stable}
Mladen Bestvina and Mark Feighn.
\newblock Stable actions of groups on real trees.
\newblock {\em Inventiones mathematicae}, 121(1):287--321, 1995.

\bibitem[BF09]{bestvina2009notes}
Mladen Bestvina and Mark Feighn.
\newblock Notes on sela’s work: Limit groups and makanin-razborov diagrams.
\newblock {\em Geometric and cohomological methods in group theory}, 358:1--29,
  2009.

\bibitem[CG05]{champetier2005limit}
Christophe Champetier and Vincent Guirardel.
\newblock Limit groups as limits of free groups.
\newblock {\em Israel Journal of Mathematics}, 146(1):1--75, 2005.

\bibitem[Chi01]{chiswell2001introduction}
Ian Chiswell.
\newblock {\em Introduction to $\Lambda$-trees}.
\newblock World Scientific, 2001.

\bibitem[DD89]{dicks1989groups}
Warren Dicks and Martin~John Dunwoody.
\newblock {\em Groups acting on graphs}, volume~17.
\newblock Cambridge University Press, 1989.

\bibitem[DGH05]{diekert2005existential}
Volker Diekert, Claudio Gutierrez, and Christian Hagenah.
\newblock The existential theory of equations with rational constraints in free
  groups is pspace-complete.
\newblock {\em Information and Computation}, 202(2):105--140, 2005.

\bibitem[Dun85]{dunwoody1985accessibility}
Martin~J Dunwoody.
\newblock The accessibility of finitely presented groups.
\newblock {\em Inventiones mathematicae}, 81(3):449--457, 1985.

\bibitem[Dun98]{dunwoody1998folding}
MJ~Dunwoody.
\newblock Folding sequences.
\newblock {\em Geometry and Topology monographs}, 1:143--162, 1998.

\bibitem[FLP12]{fathi2012thurston}
Albert Fathi, Fran{\c{c}}ois Laudenbach, and Valentin Po{\'e}naru.
\newblock {\em Thurston's Work on Surfaces (MN-48)}, volume~48.
\newblock Princeton University Press, 2012.

\bibitem[FM12]{farb2012primer}
Benson Farb and Dan Margalit.
\newblock A primer on mapping class groups, volume 49 of princeton mathematical
  series, 2012.

\bibitem[Fuj02]{fujiwara2002}
Koji Fujiwara.
\newblock On the outer automorphism group of a hyperbolic group.
\newblock {\em Israel Journal of Mathematics}, 131(1):277--284, 2002.

\bibitem[GL09]{guirardel2009jsj}
Vincent Guirardel and Gilbert Levitt.
\newblock Jsj decompositions: definitions, existence, uniqueness. i: The jsj
  deformation space.
\newblock {\em arXiv preprint arXiv:0911.3173}, 2009.

\bibitem[GL10]{guirardel2010jsj}
Vincent Guirardel and Gilbert Levitt.
\newblock Jsj decompositions: definitions, existence, uniqueness. ii.
  compatibility and acylindricity.
\newblock {\em arXiv preprint arXiv:1002.4564}, 2010.

\bibitem[Gub86]{guba1986equivalence}
Victor~Sergeevich Guba.
\newblock Equivalence of infinite systems of equations in free groups and
  semigroups to finite subsystems.
\newblock {\em Mathematical Notes}, 40(3):688--690, 1986.

\bibitem[Gui08]{guirardel2008actions}
Vincent Guirardel.
\newblock Actions of finitely generated groups on $\mathbb{R}$-trees.
\newblock In {\em Annales de l'institut Fourier}, volume~58, pages 159--211,
  2008.

\bibitem[KM06]{kharlampovich2006elementary}
Olga Kharlampovich and Alexei Myasnikov.
\newblock Elementary theory of free non-abelian groups.
\newblock {\em Journal of Algebra}, 302(2):451--552, 2006.

\bibitem[Lev05]{levitt2005automorphisms}
Gilbert Levitt.
\newblock Automorphisms of hyperbolic groups and graphs of groups.
\newblock {\em Geometriae Dedicata}, 114(1):49--70, 2005.

\bibitem[LP97]{levitt1997geometric}
Gilbert Levitt and Fr{\'e}d{\'e}ric Paulin.
\newblock Geometric group actions on trees.
\newblock {\em American Journal of Mathematics}, pages 83--102, 1997.

\bibitem[LS15]{lyndon2015combinatorial}
Roger~C Lyndon and Paul~E Schupp.
\newblock {\em Combinatorial group theory}.
\newblock Springer, 2015.

\bibitem[Mak83]{makanin1983equations}
Gennady~S Makanin.
\newblock Equations in a free group.
\newblock {\em Izvestiya: Mathematics}, 21(3):483--546, 1983.

\bibitem[Mer66]{merzlyakov1966positive}
Yuri~I Merzlyakov.
\newblock Positive formulae on free groups.
\newblock {\em Algebra i Logika}, 5(4):25--42, 1966.

\bibitem[MKS07]{magnus2007combinatorial}
W~Magnus, A~Karras, and D~Solitar.
\newblock Combinatorial group theory: Presentations of groups in terms of
  generators and relations (interscience, new york, 1966).
\newblock {\em Received: June}, 27, 2007.

\bibitem[MS84]{morgan1984valuations}
John~W Morgan and Peter~B Shalen.
\newblock Valuations, trees, and degenerations of hyperbolic structures, i.
\newblock {\em Annals of Mathematics}, 120(3):401--476, 1984.

\bibitem[Pau88]{paulin1988topologie}
Fr{\'e}d{\'e}ric Paulin.
\newblock Topologie de gromov {\'e}quivariante, structures hyperboliques et
  arbres r{\'e}els.
\newblock {\em Inventiones mathematicae}, 94(1):53--80, 1988.

\bibitem[Pau89]{paulin1989gromov}
Fr{\'e}d{\'e}ric Paulin.
\newblock The gromov topology on r-trees.
\newblock {\em Topology and its Applications}, 32(3):197--221, 1989.

\bibitem[Per08]{perin2008elementary}
Chlo{\'e} Perin.
\newblock {\em Elementary embeddings into a torsion-free hyperbolic group}.
\newblock PhD thesis, PhD thesis, Universit{\'e} de Caen, 2008.

\bibitem[Per11]{perinElementary}
Chloé Perin.
\newblock Elementary embeddings in torsion-free hyperbolic groups.
\newblock {\em Annales scientifiques de l'École Normale Supérieure},
  44(4):631--681, 2011.

\bibitem[PS]{perinforking}
C~Perin and R~Sklinos.
\newblock Forking and jsj decompositions in the free group, to appear in j.
\newblock {\em Eur. Math. Soc.(JEMS)}.

\bibitem[Raz85]{razborov1985systems}
Alexander~A Razborov.
\newblock On systems of equations in a free group.
\newblock {\em Mathematics of the USSR-Izvestiya}, 25(1):115, 1985.

\bibitem[Rez97]{reznikov1997quadratic}
Alexander Reznikov.
\newblock Quadratic equations in groups from the global geometry viewpoint.
\newblock {\em Topology}, 36(4):849--865, 1997.

\bibitem[RS94]{rips1994structure}
Eliyahu Rips and Zlil Sela.
\newblock Structure and rigidity in hyperbolic groups i.
\newblock {\em Geometric \& Functional Analysis GAFA}, 4(3):337--371, 1994.

\bibitem[RS97]{rips1997cyclic}
Eliyahu Rips and Zlil Sela.
\newblock Cyclic splittings of finitely presented groups and the canonical jsj
  decomposition.
\newblock {\em Annals of Mathematics}, pages 53--109, 1997.

\bibitem[RW]{reinfeldtmakanin}
Cornelius Reinfeldt and Richard Weidmann.
\newblock Makanin-razborov diagrams for hyperbolic groups, preprint, 201 0.

\bibitem[SB77]{serre1977arbres}
Jean-Pierre Serre and Hyman Bass.
\newblock Arbres, amalgames, sl2= asr{\'e}risque 46.
\newblock 1977.

\bibitem[Sco78]{scott1978subgroups}
Peter Scott.
\newblock Subgroups of surface groups are almost geometric.
\newblock {\em Journal of the London Mathematical Society}, 2(3):555--565,
  1978.

\bibitem[Sel97]{Sela1997acylindrical}
Zlil Sela.
\newblock Acylindrical accessibility for groups.
\newblock {\em Inventiones mathematicae}, 129(3):527--565, 1997.

\bibitem[Sel01]{Sela1}
Zlil Sela.
\newblock Diophantine geometry over groups i: Makanin-razborov diagrams.
\newblock {\em Publications Mathematiques de l'IHES}, 93:31--105, 2001.

\bibitem[Sel03]{sela2}
Zlil Sela.
\newblock Diophanting geometry over groups ii: Completions, closures and formal
  solutions.
\newblock {\em Israel Journal of Mathematics}, 134(1):173--254, 2003.

\bibitem[Sel06]{Sela6}
Zlil Sela.
\newblock Diophantine geometry over groups vi: The elementary theory of a free
  group.
\newblock {\em Geometric \& Functional Analysis GAFA}, 16(3):707--730, 2006.

\bibitem[Sho67]{shoenfieldbook}
Joseph~R. Shoenfield.
\newblock {\em Mathematical logic}.
\newblock Addison Wesley, 1967.

\bibitem[SL77]{schupp1977combinatorial}
RC~Lyndon~PE Schupp and Roger~C Lyndon.
\newblock Combinatorial group theory.
\newblock {\em Ergebnisse der Mathematik und ihrer Grenzgebiete Bd}, 89, 1977.

\bibitem[Sta65]{stallings1965topological}
John~R Stallings.
\newblock A topological proof of grushko's theorem on free products.
\newblock {\em Mathematische Zeitschrift}, 90(1):1--8, 1965.

\bibitem[Sta91]{stallings1991foldings}
John~R Stallings.
\newblock Foldings of g-trees.
\newblock In {\em Arboreal Group Theory}, pages 355--368. Springer, 1991.

\bibitem[T{\etalchar{+}}88]{thurston1988geometry}
William~P Thurston et~al.
\newblock On the geometry and dynamics of diffeomorphisms of surfaces.
\newblock {\em Bulletin (new series) of the american mathematical society},
  19(2):417--431, 1988.

\bibitem[TZ12]{tent2012course}
Katrin Tent and Martin Ziegler.
\newblock {\em A course in model theory}, volume~40.
\newblock Cambridge University Press, 2012.

\bibitem[Wei02]{weidmann2002nielsen}
Richard Weidmann.
\newblock The nielsen method for groups acting on trees.
\newblock {\em Proceedings of the London Mathematical Society}, 85(1):93--118,
  2002.

\bibitem[Wil09]{wilton2009solutions}
Henry Wilton.
\newblock Solutions to bestvina \& feighn's exercises on limit groups.
\newblock {\em Geometric and cohomological methods in group theory, London
  Mathematical Society Lecture Note Series}, 358:30--62, 2009.

\end{thebibliography}
    \pagestyle{empty}
    \afterpage{\null\newpage}

    \clearpage
    %    \newpage\null\thispagestyle{empty}\newpage
    %    \pagestyle{empty}
    
    \pagestyle{empty}

\centerline{\bf{Curriculum Vitae}}
\vspace{0.7cm}

%	\centerline{\noindent Javier de la Nuez Gonz\'alez}

\begin{flushleft}
	\begin{tabular}{p{.2\textwidth}  p{.7\textwidth}}{}
		Name: &Javier   \\
		First surname: & de la Nuez \\
		Second surname: & Gonz\'alez\\
		Date of birth:  & 28/07/1988 \\
		Place of birth:  & Madrid, Spain
	\end{tabular}
\end{flushleft}
\vspace{0.5cm}
\noindent\hrulefill

\vspace{0.5cm}

\noindent
\begin{tabular}{p{.08\textwidth} p{.02\textwidth} p{.08\textwidth} p{.7\textwidth}}
	Sep 11 &--& & PhD studies at the Instit\"ut f\"ur mathematische Logik und Grundlagenforschung, Universit\"at M\"unster, advisor Prof. Katrin Tent  \vspace*{0.2cm}\\
	Sep 09 &--& Jun 10 & Erasmus stay at the Universit\`{a} degli Studi di Pisa, Pisa (Italy) \vspace*{0.2cm}\\
	Sep 06 &--& Jun 11 & Diploma studies in mathematics (Licenciatura en Ciencias Matem\'aticas) at Universidad Complutense, Madrid (Spain) \vspace*{0.2cm}\\
	& -- & Jun 06 & Highschool studies at Liceo Scientifico Italiano \textit{Enrico Fermi}, Madrid
\end{tabular}

\end{document}